\documentclass[11pt]{amsbook}
\usepackage{graphicx}
\usepackage{rotating}

 \usepackage{caption}
 \usepackage{subcaption}

\newtheorem{definition}{Definition}[section]

\newtheorem{remark}{Remark}[section]

\def\rat#1{\mathcal R_{#1}}

\def\Forall{\quad \hbox{ for all }}
\def\bal#1\eal{\begin{aligned} #1 \end{aligned}}
\def\beq#1\eeq{\begin{equation} #1 \end{equation}}
\def\cI{\widetilde {\mathcal I}}
\usepackage[colorlinks,bookmarksopen,bookmarksnumbered,citecolor=red,urlcolor=red]{hyperref}
\usepackage{amsmath, amssymb, amsfonts, amsbsy, latexsym}
\usepackage{multirow}
\usepackage{epsfig}
\usepackage{graphicx}
\def\argmin{\mathop{\rm argmin}}

\def\calA{{\mathcal A}}
\def\calAt{{\mathbb A}}

\def\calIt{{\mathbb I}}

\def\wcalSt{ \widetilde{\mathbb S}}
\def\wcalMt{ \widetilde{\mathbb M}}
\def\wcalAt{ \widetilde {\mathbb A}}
\def\wcalIt{\widetilde{\mathbb I}}

\newcommand{\RR}{\mathbb{R}}

\def\argmin{\mathop{\rm argmin}}

\def\tiluh{{\widetilde u}_h}
\def\tilfh{{\widetilde f}_h}
\def\tilwh{{\widetilde w}_h}
\def\tilvh{{\widetilde v}_h}
\def\tilF{{\widetilde F}}

\def\calSt{{\mathbb S}}
\def\calSt{\widetilde{\mathbb S}}
\def\calIt{{\mathbb I}}

\newcommand{\R}{\mathbb{R}}

\newcommand{\calI}{{\cal I}}

\def\argmin{\mathop{\rm argmin}}

\def\calIt{{\mathbb I}}

\def\calI{\mathcal I}

\def\f{g}
\def\t{t}
\def\tt{\xi}
\def\q{b}

\title[Numerical Methods for Fractional Laplacian Using Rational Approximation]
{The Best Uniform Rational Approximation: \\
Applications to Solving Equations Involving 
Fractional powers of Elliptic Operators
}
\author{S.~Harizanov, R.D.~Lazarov, S.~Margenov, P.~Marinov} 
\date{started July 13, 2019, today is \today}

\address{Institute of Information and Communication Technologies,
Bulgarian Academy of Sciences, Sofia, Bulgaria}
\email{sharizanov@parallel.bas.bg}

\address{Department of Mathematics, Texas A\&M University, College Station, TX 77843} 
\email{lazarov@math.tamu.edu}

\address{Institute of Information and Communication Technologies,
Bulgarian Academy of Sciences, Sofia, Bulgaria}%Svetozar Margenov}
\email{margenov@parallel.bas.bg}

\address{Institute of Information and Communication Technologies,
Bulgarian Academy of Sciences, Sofia, Bulgaria}%Svetozar Margenov}
\email{marinov@parallel.bas.bg}

\begin{document}

\begin{abstract}

In this paper we consider one particular mathematical problem of this large area of
fractional powers of self-adjoined elliptic operators, defined either by Dunford-Taylor-like integrals or 
by the representation through the spectrum of the elliptic operator. 
%The elliptic operator generated by the Laplacian with homogeneous Dirichlet data is 
%know in the mathematical literature as spectral fractional Laplacian.

Due to the mathematical modeling of various non-local phenomena using such operators recently a 
number of numerical methods for solving equations involving operators of fractional order were introduced, studied, and tested. 
Here we consider the discrete counterpart of such problems obtained from finite difference or 
finite element approximations of the corresponding elliptic problems. In short, these are  linear equations for  $ \tiluh \in \R^N$
of the type $\wcalAt^{\alpha} \tiluh = \tilfh$ or $\tiluh = \wcalAt^{-\alpha} \tilfh$, 
where $\alpha \in (0,1)$, and $\wcalAt  \in \R^{N \times N}$ is an SPD matrix.

Among the existing methods is 
a method based on the best uniform rational approximation (BURA) introduced and analyzed in \cite{HLMMV18,harizanov2019analysis}.
%This approach  is based on the representation of the solution through the spectrum of the elliptic operator. 
In fact, the method of Bonito and Pasciak, \cite{bonito2017numerical,BP15,BP17}, which uses 
exponentially convergent sinc-quadratures for the Dunford-Taylor integrals, results in a rational approximation of the 
corresponding kernel. Thus theoretically, the BURA approach should be as good as the method of Bonito-Pasciak.

In the simplest case, to implement the BURA method one needs to generate the best uniform rational 
approximation of $t^{-\alpha}$ on the spectrum of $\wcalAt$.  In order to make the method feasible, 
instead we seek the BURA on the interval $[\lambda_1, \lambda_N]$,
where $\lambda_1 \le \dots \le \lambda_N$ are the eigenvalues of $\wcalAt$. This is further simplified by rescaling the 
system so the solution is sought in the form $  \tiluh = (\lambda_1^{-1} \wcalAt)^{-\alpha}  \lambda_1^{\alpha} \tilfh$, 
so we need the find the  BURA of $t^{-\alpha}$ on $[1, \lambda_N/\lambda_1]$. 
If we introduce a parameter $0 \le \delta \le \lambda_1/\lambda_N$, an estimate of the ratio $ \lambda_1/\lambda_N$
from below, we can avoid the necessity to know the spectrum of $\wcalAt$.

In this report we provide all necessary information regarding the best uniform rational approximation (BURA) 
$r_{k,\alpha}(t) := P_k(t)/Q_k(t)$ of  $t^{\alpha}$ on $[\delta, 1]$
% with rational functions  in $ $ % \in \mathcal R_k$ 
for various $\alpha$, $\delta$, and $k$. The results are presented in 160 tables containing the 
coefficients of $P_k(t)$ and $Q_k(t)  $, the zeros 
and the poles of $r_{k,\alpha}(t)$, the extremal point of the error  $t^\alpha - r_{k,\alpha}(t)$,  
the representation of $r_{k,\alpha}(t)$ in terms of partial fractions, etc.
In short,  we provide all necessary data to compute efficiently the approximation 
$\tilwh = r_{k,\alpha}((\lambda_1^{-1} \wcalAt)^{-1}) \lambda_1^{-\alpha} \tilfh$ of the exact solution 
$  \tiluh = (\lambda_1^{-1} \wcalAt)^{-\alpha}  \lambda_1^{\alpha} \tilfh$. 
Moreover,  we provide links to the files with the data that characterize $r_{k,\alpha}(t)$ which 
are available with enough significant digits so one can use them in his/her own computations. 
The presented numerical results  use Remez algorithm for computing BURA, see, e.g. \cite{Driscoll2014,PGMASA1987}.
It is well known that this method is highly sensitive to the computer arithmetics precision, so 
precomputing (or off-line computations) of BURA seems to be an important and necessary step, 
\cite{Driscoll2014,PGMASA1987}.
Here we report the results in  performing this step and we provide the data associated with handling this problem.

Further, similar data (the poles, the extreme points, the the partial fraction 
representation) is generated and presented for the BURA of the function $\f(q,\delta,\alpha;\t)=\t^{\alpha}/(1+q\,\t^{\alpha})$ 
(of three parameters $ q,\delta,\alpha$,  such that $0\le q$, $0 \le \delta < 1$, $0 < \alpha <1$ 
and a variable $ \t \in [\delta,1]$). The information is organized in the way as the case BURA of $t^\alpha$.
In the end on two examples of model problems we go through all necessary steps of extracting the 
necessary information in order to solve approximately the problems.

\end{abstract}

\maketitle

%\pagebreak

\tableofcontents

%
%\begin{center}
% {\bf Acknowledgments}
% \end{center}
%
%\noindent
%This work has been partly support by the Grant No BG05M2OP001-1.001-0003, financed by the 
%Science and Education for Smart Growth Operational Program (2014-2020) and co-financed by 
%the EU through the European structural and Investment funds.
%
%%\pagebreak

%\part{Theoretical background}\label{part:theory}
\chapter{Theoretical background}\label{chapt:theory}

\section{Introduction}\label{se:intro}

\subsection{
%\paragraph{\it 
General comments}
"Fractional" order differential operators appear naturally in many areas in mathematics and physics, e.g.
trace theory of functions in Sobolev classes (Sobolev imbedding, elliptic type),
the theory of special classes analytic functions 
(Dzhrbashyan, \cite{Djrbashian:1993}, Riemann-Liouville fractional derivative, \cite{KilbasSrivastavaTrujillo:2006,Podlubny:1999}),
modeling various phenomena, e.g. particle movement in heterogeneous media, \cite{MetzlerKlafter:2004} and/or  
modeling dynamics in fractal media (Tarasov, \cite{tarasov2011fractional}),
modeling materials with memory (e.g. viscoelasticity, Bagley-Torvik equation, \cite{torvik1984})
heavily tailed Levy flights of particles, \cite{HatanoHatano:1998}, peridynamics (deformable media with fractures), 
image reconstruction using non-local operators, \cite{gilboa2008nonlocal}, etc.
The most important property of these operators is that they are non-local.

There are two main definitions of "fractional Laplacian" (and more general steady-state sub-diffusion problem) used 
modeling of various non-local diffusion-like problems. For a reveling and thorough discussion on 
this topic we refer to \cite{Karniadakis2018fractional}. Here we shall consider the case of so-called 
spectral fractional Laplacian.

In this work we consider numerical methods and algorithms for solving equations involving
fractional powers of self-adjoined elliptic operators, defined either by Dunford-Taylor-like integrals or 
by the representation through the eigenvalues and the eigenfunctions of the elliptic operator.
%Note, that the elliptic operator generated by the Laplacian with homogeneous Dirichlet data is 
%know in the mathematical literature as spectral fractional Laplacian.
The numerical methods are done in two basic steps: \\
(1) approximation of the corresponding elliptic operator by
finite elements in a finite dimension space $V_h$ (of dimension $N$) or similar 
approximation by finite differences on a rectangular mesh leading to a matrix 
acting in the Euclidean space  $\R^N$; this results in a semi-discrete scheme;\\
(2) approximation of the fractional powers of the discretized elliptic operator 
%(by finite elements or finite differences) 
using the best uniform rational approximation (BURA) of a certain function on $[0,1]$,
which results in a fully discrete scheme.

\subsection{
%\paragraph{\it 
Semi-discrete and fully discrete approximations in the elliptic case}
First, in Subsection \ref{ss:problem}  we consider systems of equations generated by fractional powers of elliptic operators 
of the type  $   \calA^\alpha u = f $. The fractional powers of $\calA$ are defined by Dunford-Taylor integrals, which can be transformed when $\alpha \in (0, 1)$ to the Balakrishnan integral \eqref{bal}. 
Further we also consider other problems like $ \calA^\alpha u + \q u =f $ and initial value problem
$\frac{\partial u (t)}{\partial t} + \calA^\alpha u(t) =f(t) $, $t >0$, with $  u(0)=v $.

The approximation is done in two steps.
In the first step we discretize the elliptic operator $ \calA $.  
In the case of finite element  approximation we get a 
symmetric and positive definite operator $\calAt: V_h \to V_h$, that  results in 
an operator equation $\calAt^{\alpha} u_h = f_h$ for $f_h \in V_h$ given and $u_h \in V_h$ unknown. 
In the case of finite difference approximation %we get the following algebraic problem: 
we get a symmetric and positive
definite matrix $\wcalAt \in \R^{N\times N}$ and a vector $\tilfh \in \R^N$, so that the approximate solution  
$\tiluh \in \R^N$ %is obtained from solving the algebraic problem 
satisfies $\wcalAt^{\alpha} \tiluh = \tilfh$. The fractional powers of the operator $\calAt^\alpha$ and the 
matrix $\wcalAt^\alpha$ are defined by \eqref{bal} or \eqref{eq:spectral-0} (with finite summation), correspondingly.
These equations generate the so-called semi-discrete solutions
$$
 u_h = \calAt^{-\alpha} f_h  \ \mbox{or/and} \ \   \tiluh = \wcalAt^{-\alpha} \tilfh.
$$

In the second step we essentially approximate the Balakrishnan integral \eqref{bal}.
This is done by introducing the rational function $r_{k,\alpha}(\t)$, which is the 
best uniform rational approximation (BURA) of $\t^{\alpha}$ on $[0,1]$ and apply it to 
produce the fully discrete approximations, see \eqref{wh}, 
$$
w_h 
=\lambda_1^{-\alpha} r_{\alpha,k} (\lambda_1 \calAt^{-1}) f_h
\quad \mbox{and}  \quad 
\tilwh 
=\lambda_1^{-\alpha} r_{\alpha,k} (\lambda_1 \wcalAt^{-1}) \tilfh.
$$
The paper is devoted to the characterization of the rational functions $r_{k,\alpha}(\t)$
and providing  their extremal point, poles, partial fraction representation, etc.

%\subsection
\paragraph{\it Semi-discrete and fully discrete approximations in the parabolic sub-diffusion case}
Further, in Sections \ref{ss:subdiffusion-react} and \ref{ss:transient-subdiff} we apply the same strategy to the
sub-diffusion-reaction problem \eqref{eq:sub-reaction} and time-stepping procedure for the transient 
sub-diffusion problem \eqref{eq:parabolic}.  Both results in solving the following type 
semi-discrete problem \eqref{FEM-sub-reaction}:
$$
u_h = (\calAt^\alpha + \q \calIt)^{-1} f_h. 
$$
In this case, we consider two possible fully discrete schemes. The first one is based on the BURA of the 
function $ t^\alpha/(1+q t^\alpha)$ on the interval $[\delta, 1]$. The second fully discrete scheme is based on
a rational approximation of the same function, but NOT the best one, see Subsection \ref{ss:ura-res} 
and Definition \ref{pdef0ura} and call further URA-method.

%\begin{enumerate}
%
%\item Abstract
%\item Short introduction with origins of the corresponding problems;
%\item Elliptic problems leading to such equations;
%\item Parabolic problems producing relevant problems;
%\item The issue of BURA and URA as solution methods;
%\item What are the relevant approximation problems we need to consider
%\item Our findings
%\item Data gathering for some particular problems
%\end{enumerate}

\section{Examples of problems involving fractional powers of elliptic operators}\label{s:ell-operators}

\subsection{Spectral fractional powers of elliptic operators} %Laplacian} % setting}
\label{ss:problem}
In this paper we consider the following second order elliptic equation with homogeneous Dirichlet data:
\beq
\bal
 - \nabla \cdot( a(x) \nabla v(x)) + \q(x) u&= f(x),\qquad \hbox{ for }
x\in \Omega,\\
v(x)&=0, \qquad \hbox{ for } x\in \partial \Omega.
\eal
\label{strong}
\eeq
Here $\Omega $ is a bounded domain in $\RR^d$,  $d\ge 1$. %  and 
We assume that $0<a_0 \le a(x) $ and $\q(x) \ge 0$ for $x\in \Omega$.

The fractional powers of the elliptic operator associated with the problem \eqref{strong} are defined in 
terms of the weak form of
\eqref{strong}, namely, $v(x)$ is the unique function in  $V= H^1_0(\Omega)$
satisfying
\beq
a(v,\theta)= (f,\theta)\qquad \Forall \theta\in V.
\label{weak}
\eeq
Here
$$
a(w,\theta):=\int_\Omega  \left (a(x) \nabla w(x) \cdot \nabla  \theta(x) + \q w \theta \right )\, dx  \quad \hbox{ and } \quad
(w,\theta):=\int_\Omega w(x) \theta(x)\, dx.
$$
For $f\in L^2(\Omega):= X$, \eqref{weak} defines a solution operator $T f :=v$.
Following \cite{kato}, we define an unbounded operator $\calA$ on $X$  as follows.
The operator $\calA$  with domain
 $$ D(\calA) = \{ Tf\,:\, f\in X\}$$
 is defined by  
 \beq\label{def-A}
 \calA v=g \ \ \forall  v\in D(\calA), \ \ \mbox{ where } \ \  g\in X \ \ \mbox{ with } \ \  Tg=v.  
 \eeq
 The operator $\calAt$ is well defined as $T$ is injective.

Thus, the focus of our work in this paper is numerical approximation and algorithm development for the
equation: 
\beq\label{frac-eq}
 \calA^\alpha u = f \quad  \mbox{with a solution} \quad u =\calA^{-\alpha} f.
\eeq
Here  $\calA^{-\alpha}=T^\alpha$ for $\alpha>0$
is defined by Dunford-Taylor integrals which can be transformed  when
$\alpha\in (0,1)$, to the  Balakrishnan integral, e.g. \cite{balakrishnan}: for $f\in X$,
\beq
u=\calA^{-\alpha} f =\frac {\sin(\pi \alpha)} \pi 
\int_0^\infty \mu^{-\alpha} (\mu \calI +\calA)^{-1}f\, d\mu.
\label{bal}
\eeq
This definition is sometimes referred to as the spectral definition of fractional powers.   
One can also use an equivalent definition through the 
expansion with respect to the eigenfunctions $\psi_j$ and the eigenvalues $\lambda_j$ of $\calA$, 
e.g. \cite{Karniadakis2018fractional,Acosta-Borthagaray2017}:
\beq\label{eq:spectral-0}
\calA^{\alpha} u= \sum_{j=1}^\infty \lambda_j^{\alpha} (u,\psi_j) \psi_j \ \ \mbox{so that }
u = \sum_{j=1}^\infty \lambda_j^{-\alpha} (f,\psi_j) \psi_j 
\eeq
Since the bilinear form $a(\cdot, \cdot)$ is symmetric on $ V \times V$ and $\calA$ is a unbounded 
operator we can show that $\lambda_j $ are real and positive and $\lim_{j \to \infty} \lambda_j =\infty$.

%We note that there are also problems on bounded domains involving
%fractional powers, for example, those related to L\'evy diffusion 
%\cite{Acosta-Borthagaray2017,ros-oton2015}.   These problems
%involve the restriction of non-local operators defined on $\RR^d$
%applied to bounded domain functions extended, e.g., by 0, outside of $\Omega$.  
%However, in this paper, we focus on the spectral definition \eqref{bal} and the corresponding
%approximations by the finite element or finite difference methods.

\begin{remark}
An operator $\calA$ is positivity preserving if $\calA f\ge 0$ when $f\ge 0$.
We note that by the maximum principle, $(\mu \calI+ \calA)^{-1}$ is a positivity
preserving operator for $\mu\ge 0$
and the formula \eqref{bal} shows that $\calA^{-\alpha}$ is  also.  In many
applications, it is important that the discrete approximations share
this property.
\end{remark}

%From the solution representation \eqref{bal} of the sub-diffusion problem 
%we can easily see that if $f \ge 0$ then $ u \ge 0$,
%i.e. the solution preserves the positivity in the same way as in the standard diffusion problem.
%{\it Positivity preservation} is a feature that often is needed to be preserved by the approximate solution as well.
%

\subsection{Sub-diffusion-reaction problems}\label{ss:subdiffusion-react}
Another possible model of sub-diffusion reaction is give by the operator equation: find $u \in V$ s.t. 
\beq\label{eq:sub-reaction}
\calA^\alpha u + \q u =f \quad \mbox{where} \quad \q=const \ge 0.
%(Qu,\theta) =  \int_\Omega q u \theta \, dx \  \  \forall \theta \in V.
\eeq

\subsection{Transient sub-diffusion-reaction problems}\label{ss:transient-subdiff}
For time dependent problems one can consider: find $ u(t) \in V$ for $t \in (0,t_{max}]$  such that 
\beq\label{eq:parabolic}
\frac{\partial u (t)}{\partial t} + \calA^\alpha u(t) =f(t) \quad \mbox{and }  u(0)=v,
\eeq
with $v$ given initial data, $t_{max} >0$ is a given number, and the finite dimensional operator 
(matrix) $\calAt$ defined by \eqref{def-A}.

\section
[Semi-discrete approximations of equations involving fractional powers of ...]
{Semi-discrete approximations of equations involving fractional powers of elliptic operators}
\label{s:semi-discrete}

We study approximations to $u= \calA^{-\alpha} f$ defined in terms  of finite difference or finite element
approximation of the operator  $T$.   We shall use the following convention regarding 
the approximate solutions by the  finite element, \cite{ern-guermond},  and finite difference, \cite{samarskii2001theory}, methods.   

The finite element solutions are  functions  in $V_h$,
an $N$-dimensional space of continuous piece-wise linear functions over a partition  ${\mathcal T}_h$ of
the domain. Such functions will be denoted by $u_h$, $ v_h$, etc.   
Also we shall denote by $\calAt$, $\calIt$, etc operators acting on the elements $u_h, \theta_h$, etc 
in the finite dimensional space of functions $V_h$.  When a nodal basis of the finite 
element space is introduced, then the vectors coefficients in this basis
are denoted  $\tiluh$, $\tilvh$, etc. Under this convention  operator equations in $V_h$ such as 
$\calAt u_h=f_h$ will be written as a system of linear algebraic equations $ \wcalAt  \tiluh=\tilfh$  in $\RR^N$.

In the finite difference case, discrete solutions are vectors in $\R^N$
and are also denoted $\tiluh$, $\tilvh$, etc. Then the corresponding 
counterparts of operators action on these vectors are denoted by $\wcalAt, \wcalIt,$ etc.

%\subsubsection{
\subsection{The finite difference approximation}\label{ss:FD}
%\label{ss:FD}
%%%%%%%%%%%%%%%%%%%%%%%%%%%%%
%{\it  In the finite difference the approximation}  
In this case the approximation $\tiluh \in \R^N$  of $u$ is given by
\beq
\wcalAt^{\alpha} \tiluh = \cI_h f := \tilfh, \ \ \mbox{or equivalently} \ \ \ 
\tiluh = %\calAt^{-\alpha} \cI_h f:=
\wcalAt^{-\alpha}  \tilfh, 
% \quad \tiluh \in \R^N
\label{fda}
\eeq
where $\wcalAt$ is an $N\times N$ symmetric and positive definite matrix
coming from  a finite difference approximation to the differential
operator appearing in \eqref{strong}, $\tiluh$ is the vector in $\R^N$ of the approximate solution at the 
interior $N$ grid points,  and $\cI_h f:= \tilfh \in \R^N$ denotes the vector of the values of
the data $f$ at the grid points.  
%Examples of such matrices are given in Subsection \ref{ss:FD-examples}.

\paragraph{\it Example 1} We first consider the one-dimensional equation \eqref{strong} with
variable coefficient, namely, we study the following boundary value problem
$ - ( a(x) u^{\prime})^{\prime} =f(x),$ $ u(0)=0, \ u(1)=0, \ $ for $ 0<x<1$,
where $a(x)$ is uniformly positive function on $[0,1]$. On a uniform mesh $x_i =ih$, $i=0, \dots, N+1$,
$ h=1/(N+1)$,  we consider the three-point  approximation of the second derivative
\begin{equation*}
\begin{split}
 (a(x_i) u'(x_i))' &
%\approx-\frac{2f(x_{i-1})}{(x_i-x_{i-1})(x_{i+1}-x_{i-1})}+\frac{2f(x_{i})}{(x_i-x_{i-1})(x_{i+1}-x_{i})} -
%\frac{2u(x_{i+1})}{(x_{i+1}-x_{i})(x_{i+1}-x_{i-1})}\\
\approx %\frac{1}{\widetilde h_i} \left ( 
 \frac{1}{h}\left (a_{i+\frac12}\frac{u(x_{i+1}) - u(x_i)}{h} - a_{i-\frac12}\frac{u(x_i) - u(x_{i-1})}{h} \right )
%:=\calAt_{i,i-1}u(x_{i-1})+\calAt_{i,i}u(x_{i})+\calAt_{i,i+1}u(x_{i+1}),
\end{split}
\end{equation*}
Here $a_{i-\frac12}=a(x_i - h/2)$ or  $a_{i - \frac12}=\frac{1}{h}\int_{x_{i-1}}^{x_i} a(x) dx $. Note that the former is 
the standard finite difference approximation obtained from 
the balanced method (see, e.g. \cite[pp. 155--157]{samarskii2001theory}). 
%, while the latter is a result of finite element method with mass lumping, see Subsection \ref{ss:FEM}.
%    This is discussed later...
%>From the representation of the solution by \eqref{bal} with $\calA$ replaced by $ \wcal%At$ and 
%$f$ replaced by $\tilfh$ we see that if $f \ge 0$ then $\tilfh \ge 0$, i.e. 
%the finite difference method preserves the positivity.

Then the finite difference approximation of the differential equation $- (a(x) u'(x))' =f(x)$ is given by the matrix equation \eqref{fda}
with 
%$$%\calAt =\calAt_{i,i-1}, \calAt_{i,i}, \calAt_{i,i+1} ) 
\begin{equation}\label{FD-matrix-1D-k} 
\wcalAt= \frac{1}{h^2} \left[
\begin{array}{ccccc} a_{\frac12}+ a_{\frac32}& -  a_{\frac32}  &&&\\
 -  a_{\frac32} &  a_{\frac32} +  a_{\frac52} & -  a_{\frac52}&&\\
\cdots &\cdots &\cdots &\cdots &\cdots\\ 
& -  a_{i-\frac12}& a_{i-\frac12} +  a_{i+\frac12} & a_{i+\frac12} &\\\cdots &\cdots &\cdots &\cdots  &\cdots\\
&&& - a_{N-\frac12} &  a_{N-\frac12} + a_{N +\frac12}
\end{array}\right], \ \ 
%  \cI_h f =\tilfh = \left[
%\begin{array}{c}
% f(x_1) \\
% f(x_2) \\
% \dots  \\
% f(x_i) \\
%\dots \\
% f(x_N)
%\end{array} \right ] .
\end{equation}
%$$
The eigenvalues $\lambda_i$ of the matrix $\calAt $ are all real and positive and satisfy %the following inequalities 
$$ 
4 \pi^2 \min_x a(x) \le \lambda_i \le  4 \max_x a(x)/h^2,  \ \  i=1, \dots, N.
$$

%%%%%%%%%%%%%%%%%%%%%%%%%%%%%%%%%%%%%%%%%%%%%%
%\subsubsection{
\subsection{The finite element approximation}\label{ss:FEM}
%\label{ss:FEM}
The approximation in the finite element case is defined in terms of a
conforming finite dimensional space $V_h\subset V$ of piece-wise linear 
functions over a quasi-uniform partition ${\mathcal T}_h$ of $\Omega$ 
into triangles or tetrahedrons.
Note that the  construction \eqref{bal}  of negative fractional powers
carries over to the finite dimensional case,  replacing $V$ and $X$ by
$V_h$ with $a(\cdot,\cdot)$ and $(\cdot,\cdot)$ unchanged.  

The discrete operator $\calAt$ is defined to be the inverse of
$T_h:V_h\rightarrow V_h$ with
$T_h g_h:=v_h$  where $v_h\in V_h$ is the unique solution to
\beq
a(v_h,\theta_h) = (g_h,\theta_h),\Forall \theta_h\in V_h.
\label{T_h}
\eeq
The finite element approximation $ u_h \in V_h$ of $u$ is then given by
\beq
\calAt^{\alpha} u_h = \pi_h f, \ \ \mbox{or equivalently} \ \ \ 
u_h = \calAt^{-\alpha} \pi_h f:=\calAt^{-\alpha}f_h,
\label{fea}
\eeq
where $\pi_h$ denotes the $L^2(\Omega) $ projection into $V_h$.  In this
case, $N$ denotes the dimension of the space $V_h$ and equals the number
of (interior) degrees of freedom.  The operator $\calAt$ in the finite element
case is a map of $V_h$ into $V_h$ so that 
$\calAt v_h:=g_h$,  where $g_h\in V_h$ is the unique solution to
\beq
 (g_h,\theta_h)=a(v_h,\theta_h),\Forall \theta_h\in V_h.
\label{calAt-f}
\eeq
Let $\{\phi_j\}$ denote the standard ``nodal" basis of $V_h$.
In terms of this basis  %$\calAt$ corresponds to the matrix
\beq\label{FEM-matrices}
\calAt \mbox{ corresponds to the matrix } \wcalAt = \wcalMt^{-1} \wcalSt, \ \  
\mbox{where} \ \ \wcalSt_{i,j} =a(\phi_i, \phi_j), \ \ \ \wcalMt_{i,j} =(\phi_i, \phi_j). 
\eeq
In the terminology of the finite element method, $\wcalMt$ and $ \wcalSt$ are the mass 
(consistent mass) and  stiffness matrices, respectively.

Obviously, if $\theta =\calAt \eta$ and $\widetilde \theta,\widetilde \eta \in \R^N$ are the coefficient vectors
corresponding to $\theta,\eta\in V_h$, then $\widetilde \theta = \wcalAt
\widetilde \eta$.  Now, for the coefficient vector $\tilfh$ corresponding to $f_h=\pi_hf$ we have 
 $\tilfh = \wcalMt^{-1} \tilF$, where $\tilF$ is the vector with entries
$$
{\tilF}_j=(f,\phi_j), \qquad \hbox{ for }j=1,2,\ldots,N.
$$
%Obviously the solution of the discrete problem 
%\eqref{fea} is expressed through \eqref{bal} by
%replacing $\calA$ by $\calAt$ and $f$ by $\pi_h f$. 
%Therefore, if $\calAt$ preserves the positivity then $\calAt^\alpha$ also preserves the positivity.
Then using vector notation so that $\tiluh$ is the coefficient vector representing the solution $u_h$ 
through the nodal basis, we can write the finite element approximation of \eqref{strong} 
in the form of an algebraic system 
\beq \label{classic-FEM}
\wcalAt \tiluh = \wcalMt^{-1} \tilF \ \ \mbox{which implies} \ \ \calSt \tiluh = \tilF.
\eeq
Consequently, the finite element approximation of the sub-diffusion problem \eqref{fea} becomes
\beq\label{mat-FEM}
\wcalMt \wcalAt^\alpha \tiluh = \tilF  \quad \hbox{or}\quad \tiluh=
\wcalAt^{-\alpha} \wcalMt^{-1} \tilF.
\eeq
\vspace{2mm}

%%%%%%%%%%%%%%%%%%%%%%%%%%%%%%%%%%%%%%%%%%%
%\subsubsection{
%\paragraph{\it \underline
\subsection{The lumped mass finite element approximation}\label{ss:lumped-mass}
%\label{ss:lumpedFEM}
%%%%%%%%%%%%%%%%%%%%%%%%%%%%%%%%%%%%%%%%%%%%

We shall also introduce the finite element method with ``mass lumping" for two reasons.
First, it leads to positivity preserving fully discrete methods.  % (see, Section \ref{ss:fem-lumped}).
Second, it is well known that on uniform meshes lumped mass schemes for linear elements are equivalent
to the simplest finite difference approximations. This will  
be used to study the convergence of the finite difference method 
for solving the problem \eqref{frac-eq}, an outstanding issue in this area.
%To produce positivity preserving fully discrete methods  we shall apply ``mass lumping"   

We introduce the lumped mass (discrete) inner product $ (\cdot,\cdot)_h$ on $V_h$ in the
%Then the corresponding lumped mass matrix is defined in the 
following way (see, e.g. \cite[pp.~239--242]{Thomee2006})  for $d$-simplexes in $\R^d$:
\beq\label{mass-lumping}
(z,v)_h = \frac{1}{d+1} \sum_{\tau \in {\mathcal T}_h } \sum_{i=1}^{d+1} |\tau|  z(P_i) v(P_i)
\ \ \mbox{and  } \ \  {\wcalMt}_h =\{ (\phi_i,\phi_k)_h\}_{i,k}^N.
\eeq
%that introduces the lumped ``mass" matrix $ \ {\wcalMt}_h =\{ (\phi_i,\phi_k)_h\}_{i,k}^N$.
Here  $P_1, \dots, P_{d+1}$ are the vertexes of the $d$-simplex $\tau$ and $|\tau|$ is its 
$d$-dimensional measure.  The matrix $\wcalMt_h$ is called lumped mass matrix.
Simply, the ``lumped mass" inner product 
is defined by replacing the integrals determining the finite element mass matrix by local 
quadrature approximation, specifically, the quadrature defined by
summing values at the vertices of a triangle weighted by the area of the triangle. 
%The resulting lumped  mass matrix is diagonal with positive diagonal entries.  

In this case, we define $\calAt$ by 
$\calAt v_h:=g_h$  where $g_h\in V_h$ is the unique solution to
\beq
 (g_h,\theta_h)_h=a(v_h,\theta_h),     \Forall \theta_h\in V_h
\label{calAt-fm}
\eeq 
so that 
\beq\label{A-lumped}
\calAt   \mbox{  corresponds to the matrix } \wcalAt = {\wcalMt}_h^{-1} \wcalSt,  \quad {\wcalMt}_h =\{
(\phi_i,\phi_k)_h\}_{i,k}^N. 
\eeq
Here  $ {\wcalMt}_h$ is the lumped mass matrix which is diagonal with
positive entries.
We also replace $\pi_h$ by $\calI_h$ so that the lumped mass semi-discrete
approximation is given by
\beq u_h = \calAt^{-\alpha} \calI_h f := \calAt^{-\alpha} f_h\quad \hbox{or} \quad   \widetilde u_h = 
\wcalAt^{-\alpha} \widetilde F.
\label{lumped-semi}
\eeq
Here $\widetilde F$ is the coefficient vector in the representation of the function $\calI_hf$
with respect to the nodal basis in $V_h$.
We shall call $\tiluh$ in \eqref{fda} and $ u_h$ in \eqref{fea} and
\eqref{lumped-semi}   {\it semi-discrete   approximations} of $u$.

\subsection{Discretization of sub-diffusion-reaction problem \eqref{eq:sub-reaction}}\label{ss:sbreac-dscrete}

If we use \eqref{fda}, a finite difference approximation $\wcalAt$ of the operator $\calA$, then the 
corresponding discrete problem is: find $\tiluh \in \R^N$ such that
$$
\left ( \wcalAt^\alpha  +  \q \wcalIt  \right ) \tiluh =   \tilfh.
$$

In a similar  way one can introduce the corresponding finite element 
discretizations of the problem \eqref{eq:sub-reaction}:
\beq\label{FEM-sub-reaction}
\left ( \calAt^\alpha  +  \q \calIt  \right ) u_h =   f_h,
\eeq
where $\calAt$ is defined by \eqref{calAt-f} and $\calIt: V_h \to V_h$ is the identity operator in $V_h$.  
Using consistent mass matrix evaluation the operator $  \calIt$ has matrix representation  $\wcalMt$
defined by \eqref{FEM-matrices}, while using lumped mass evaluation, the operator $  \calIt$ is represented
by the lumped mass matrix $\wcalMt_h$ defined by \eqref{mass-lumping}.

By using these matrix representations % of the operator $\calAt$,
this equation can be written as  
a system of linear algebraic equations in $\R^N$.
Taking into account the matrix representation of the operator $\calAt$ we get the following systems,
corresponding to the consistent mass and lumped mass finite elements approximations of the $L^2$-inner product in $V_h$:
$$
\wcalMt \left ( \wcalAt^\alpha  + \q \wcalIt \right ) \tiluh = \tilF \quad \mbox{or} \quad 
 \wcalMt_h \left ( \wcalAt_h^\alpha  + \q \wcalIt \right ) \tiluh = \tilfh.
$$

\subsection{Time-dependent problems}\label{ss:parabolic}
Similarly, the discretization of the time-dependent problem \eqref{eq:parabolic} with 
implicit Euler approximation in time, for $t_n=n \tau$, $n=1,2, \dots, M$, $\tau = t_{max}/M$ and 
$u_h^n \in V_h$ an approximation of $u(t_n)$,
will lead to the operator equation
%following algebraic problem %related to the following time stepping procedure
\beq\label{eq:time-stepping}
\left (\frac{1}{\tau} \calIt + \calAt^\alpha \right ) u_h^{n} = \frac{1}{\tau} u_h^{n-1} + f_h^n, \ \ n=1, \dots, M,
\eeq
where $ f_h^n$ is the $L^2$-projection of $f(t_n)$ on $V_h$. Denoting by $v_h^n = \frac{1}{\tau} u_h^{n-1} + f_h^n$,
we have the following representation of the solution $ u_h^{n} $:
\beq\label{eq:sol-time-n}
u_h^{n} = \left (\frac{1}{\tau} \calIt + \calAt^\alpha \right )^{-1} v_h^n.
\eeq

Having in mind the matrix representations 
 \eqref{FEM-matrices} (for the consistent mass FEM), \eqref{mass-lumping} (for the lumped-mass FEM), and \eqref{A-lumped} 
of the operator $\calAt: V_h \to V_h$ 
we get the following systems of algebraic equations:
\[
\left (\frac{1}{\tau} \wcalMt + \wcalSt^\alpha \right ) \widetilde u_h^{n} = \frac{1}{\tau} \wcalMt \widetilde u_h^{n-1} + \widetilde F^n
\quad
\mbox{and }
\quad
\left (\frac{1}{\tau} \wcalMt_h + \wcalSt^\alpha \right ) \widetilde u_h^{n} = \frac{1}{\tau} \wcalMt_h \widetilde u_h^{n-1} + \widetilde F^n,
\] 
for the standard and lumped mass FEM, correspondingly.

%\begin{remark}

\section{Brief review of methods for solving equations involving fractional powers of elliptic operators} \label{sec:other}

We note that computing the solutions of \eqref{fea}, \eqref{mat-FEM}, \eqref{eq:sol-time-n} involves inverting fractional powers of elliptic operators or their shifts. This is a computationally intensive task and the aim of this paper is to provide
a methodology that results in fast and efficient methods.  Further in the paper these are called fully discrete schemes.
%\end{remark}

Due to the serious interest of the computational mathematics and physics communities in modeling and simulations involving
fractional powers of elliptic operators, a number of approaches and algorithms has been developed, studied and tested on 
  various problems,  \cite{aceto2018efficient,nochetto2015pde,nochetto2016pde,HOFREITHER2019}. We survey some of these approaches by splitting them into four basic categories. These are methods based on:

\begin{enumerate}

\item An extension of the problem from $ \Omega \subset \R^d$  to a problem in 
$\Omega \times (0,\infty) \subset \R^{d+1}$, see, e.g. \cite{caffarelli2007extension}. Nochetto and co-authors 
in \cite{nochetto2015pde,nochetto2016pde}  developed efficient computational method based on  finite element 
discretization of the extended problem and subsequent use of multi-grid technique. The main deficiency of the 
method is that instead of problem in $\R^d$ one needs to work in a domain in one dimension higher which adds to the 
complexity of the developed algorithms.

\item Reformulation of the problem as a pseudo-parabolic on the cylinder $(0,1) \times \Omega$
by adding a time variable $t \in (0,1)$. Such methods were proposed, developed, and tested by Vabishchevich in
\cite{vabishchevich2015numerical,vabishchevich2018numerical}. As shown in the numerical experiments in 
\cite{HOFREITHER2019}, while using uniform time stepping, this method is  very slow. However, the improvement made in \cite{DuanLazarovPasciak,CIEGIS2019} makes the method quite competitive.

\item Approximation of the Dunford-Taylor integral representation of the solution of equations involving fractional powers of 
elliptic operators, proposed in the pioneering paper of Bonito and Pasciak \cite{BP15}. Further the idea was extended and 
augmented in various directions in \cite{BP17,bonito2019sinc,bonito2019numerical}. These methods use exponentially 
convergent sinc quadratures and are the most reliable and accurate in the existing literature.

\item Best uniform rational approximation of the function $t^\alpha$ on $[0,1]$, proposed in \cite{HMMV2016,HLMMV18},
and further developed in \cite{harizanov2017positive,harizanov2019cmwa,harizanov2019analysis}
and called BURA methods. 

\end{enumerate}

As shown recently in \cite{HOFREITHER2019},
these methods, though entirely different, are interrelated and all seem to involve certain rational approximation of
the fractional powers of the underlying elliptic operator. As such, from mathematical point of view,
those based on the best uniform rational approximation should be the best. However, one should realize that BURA
methods involve application of the Remez method of finding the best uniform rational approximation,
\cite{PGMASA1987,Driscoll2014}, a numerical algorithm for solving certain min-max problem that is highly 
nonlinear and sensitive to the precision of the computer arithmetic. For example, in \cite{varga1992some}
the errors of the best uniform rational approximation of $t^\alpha$ for six values of $\alpha  \in [0,1]$ 
are reported for degree $k\le 30$ by using computer arithmetic with $200$ significant digits.

%\newpage

%%ssssssss%%%%%%%%%%%%%%%%%%%%%%%%%%%%
\section{Fully discrete schemes}\label{se:fully-discrete}
%%%%%%%%%%%%%%%%%%%%%%%%%%%%%%%%%%%

\subsection{Explicit representation of the solution of $ u_h = \calAt^{-\alpha}  f_h $}
\label{ss:explicit}

Now consider the spectral properties of the  operator $\calAt$:
$$
\calAt \psi_j =\lambda_j \psi_j \ \ \mbox{or in matrix form} \ \    \wcalAt  \psi_j = \lambda_j \psi_j, \ \ j=1, \dots, N.
$$
Note that if $ \wcalAt $ is defined from the finite difference approximation, then this results in a standard matrix eigenvalue problem.
In case of finite element approximation \eqref{FEM-matrices} or \eqref{A-lumped}  this results in 
corresponding  generalized eigenvalue problems
\beq\label{eq:spectral}
\wcalSt \psi_j =  \lambda_j \wcalMt   \psi_j \ \ \mbox{or} \ \ \wcalSt \psi_j =  \lambda_j \wcalMt_h   \psi_j, \ \ j=1, \dots, N.
\eeq
 Using the eigenvalues and the eigenfunctions we have explicit representation of the
 solution of the operator equation:
 \beq\label{eq:sol-spectral}
 u_h = \sum_{j=1}^N \lambda_j^{-\alpha} (f_h,\psi_j) \psi_j.
 \eeq
 This representation can be used as a direct method for solving the equation $ u_h = \calAt^{-\alpha}  f_h $. 
 Moreover, using FFT-like technique, in the cases when possible (rectangular domain and
 constant coefficients), this could be an efficient numerical method.
 However, this will limit substantially the applicability of such approach.
 
 Similarly, the solution of the problem resulting in time-stepping method \eqref{eq:sol-time-n} can be expressed through the
 eigenvalues and the eigenfunctions we have
\beq\label{eq:sol-spectral-t}
 u_h^n = \sum_{j=1}^N \left (\lambda_j^{\alpha} + \frac{1}{\tau} \right )^{-1} (v_h^n,\psi_j) \psi_j 
 = \sum_{j=1}^N  \lambda_j^{-\alpha} \left (1 + \frac{1}{\tau} \lambda_j^{-\alpha} \right )^{-1}  (v_h^n,\psi_j) \psi_j .
\eeq
 
 \subsection{The idea of the fully discrete schemes}\label{ss:fully-discrete-i}
 To explain our approach we consider $P_k(t)$, the polynomial of degree $k$ that approximates $t^{\alpha}$ on the 
 interval $[1/\lambda_N, 1/\lambda_1] $ in the maximum norm. Then the vector $w_h$ defined by 
 $$
 w_h= \sum_{j=1}^N P_k(\lambda_j^{-1}) (f_h,\psi_j) \psi_j \ \ \mbox{is an approximation of} 
 \ \ u_h = \sum_{j=1}^N \lambda_j^{-\alpha} (f_h,\psi_j) \psi_j. 
 $$ 
 Moreover, we can express the error of this approximation through the approximation properties of $P_k(t)$ \cite{Lund17}:
 $$
 \| u_h - w_h \| \le \max_{t = \lambda_1, \dots, \lambda_N} |t^{-\alpha} - P_k(t^{-1}) | \| f_h\| =
 \max_{t = 1/\lambda_N, \dots, 1/\lambda_1} |t^{\alpha} - P_k(t) | \| f_h\|.
  $$
 Once we have the polynomial $P_k(t)$ then we can find its roots $\xi_i$, $i=1, \dots, k$, so that we have  
 $$P_k(t^{-1})=c_0\prod_{i=1}^k (t^{-1} - \xi_i^{-1}) $$ %\ \ \mbox{and consequently} \ \ 
 and consequently
 $$
 w_h= c_0 \prod_{i=1}^k (\calAt^{-1} - \xi_i^{-1} \calIt) f_h := c_0 \prod_{i=1}^k  w_h^{(i)}.
 $$
 Thus, to find $w_h$ we need  to find $ w_h^{(i)} =  (\calAt^{-1} - \xi_i^{-1} \calIt) f_h $,  $i=1, \dots, k$
  which results in  solving $k$ systems $\calAt w_h^{(i)}=( \calIt - \xi_i^{-1} \calAt) f_h$. %, $i=1, \dots, k$. 
 
%\red {
This idea will produce a computable approximate solution, but will not lead 
to an efficient method since the required polynomial degree $k$ will depend on 
the spectral condition number $\kappa (\calAt) = \lambda_N/\lambda_1$.
Namely, the best uniform polynomial approximation of $t^\alpha$, 
$t\in [1/\lambda_N,1/\lambda_1]$, is given by the scaled and shifted 
Chebyshev polynomial $\widetilde{P}_k(t)$. Then, the error estimate
$$
\max_{t = \lambda_1, \dots, \lambda_N} |t^{-\alpha} - \widetilde{P}_k(t) | < 
2\left (
\frac{\sqrt{\kappa^\alpha} +1} {\sqrt{\kappa^\alpha}-1}
\right )^{{k}}
$$
holds true, and therefore a polynomial of degree 
$
k \approx 2 \kappa^{\frac{\alpha}{2}} \log\frac{2}{\epsilon}
%\left (\sqrt{\left (\frac{\lambda_N}{\lambda_1}\right )^\alpha}\right )
$
will be needed to guarantee the the relative error less than $\epsilon$ of $\| u_h - w_h \|$.
%For {\it Example 1}  
Note, that the matrix $ \calAt$ defined by \eqref{FD-matrix-1D-k} has a condition number 
$\kappa (\calAt) = O\left (\max a(x)/\min a(x) h^{-2}\right )$ and the degree of the 
polynomial of best uniform approximation will grow as $1/h^{\alpha/2}$ as $h \to 0$. 
%}

Our aim is to produce a method that involves much smaller $k$,  so we need to solve fewer systems 
of the type  $( \calAt - \xi_i \calIt) v= f$ with $\xi_i \le 0$.

%We need to make a small review of the existing results, like Clemens work, \cite{HOFREITHER2019}.

 %%%%%%%%%%%%%sssssssssss
 \subsection{Fundamental properties of BURA of $\t^\alpha$ on $(0,1]$}
 \label{ss:BURA}
 %%%%%%%%%%%%%%%
 Instead of polynomial approximation, we shall seek a rational approximation of the function $t^{-\alpha}$. In order to make the 
 computations uniform and to use the known results for the approximation theory, we first rewrite the solution 
 of the \eqref{fda} in the form
 \beq\label{eq:sol-function}
 u_h = \lambda_1^{-\alpha} (\lambda_1 \calAt^{-1})^\alpha f_h.
 \eeq
 The scaling by $\lambda_1$ maps the eigenvalues of $ \lambda_1 \calAt^{-1}$ to the interval 
 $(\lambda_1/\lambda_N,1] :=(\delta, 1]  \subset (0,1]$. Here $\delta=\lambda_1/\lambda_N$ is a small parameter.
 Often below we shall take even $\delta =0$.

Similarly, the solution \eqref{eq:sol-spectral-t} of the time dependent problem after the scaling 
by $\lambda_1$ maps the eigenvalues of $\lambda_1 \calAt^{-1}$ to the interval 
$(\lambda_1/\lambda_N, 1] :=[\delta,1] \subset (0,1]$.
Then % and the solution has a representation
$$
u_h^n  = \sum_{j=1}^N  \lambda_j^{-\alpha}(1 + \q \lambda_j^{-\alpha})^{-1} (v_h^n,\psi_j) \psi_j 
%= \lambda_j^{-\alpha} +q )^{-1} (v_h^n,\psi_j) \psi_j 
= \Big (\q \lambda_1^{-\alpha} (\lambda_1^{-1} \calAt)^{-\alpha} + \calIt \Big )^{-1} (\lambda_1^{-1} \calAt)^{-\alpha} \lambda_1^{-\alpha} v_h^n.
 $$
Introducing the parameters $q=\q \lambda_1^{-1}$ and defining  the function
$
\f(\t):= \f(q,\delta, \alpha;\t ) =\frac{\t^{\alpha}}{ 1 + q \t^{\alpha}  } \ \ \mbox{for} \ \ t \in (\lambda_1/\lambda_N, 1]
$
we get
\beq\label{eq:sol-function-q}
u_h^n  = \f(\lambda_1 \calAt^{-1}) \lambda_1^{-\alpha} v_h^n 
= \sum_{j=1}^N  \f(\lambda_j) \lambda_1^{-\alpha}  (v_h^n,\psi_j) \psi_j .
\eeq

Now we consider BURA along the diagonal of
the Walsh table and take $\rat k$ to be the set of rational functions 
% of the form  $P_k(t)/Q_k(t)$ with $P_k(t)$ and $Q_k(t)$ polynomials of degree $k$: 
%and $Q_k(0)=1$.   We shall also use the following notation  for a class of rational functions:  
$$
\mathcal R_k= \bigl\{ r(\t): r(\t)=P_k(\t)/Q_k(\t),   
\ P_k \in {\mathcal P}_k, \mbox{ and  } \ Q_k \in {\mathcal P}_k, \ \mbox{monic}  \bigr\}
$$
with  ${\mathcal P}_k$ set of algebraic polynomials of degree $k$. % on the interval $ [\delta,1] $.

The best {\it discrete} uniform rational
approximation (discrete BURA) of $t^\alpha$ is the rational function $R_{\alpha,k}\in \rat k$ 
satisfying
\begin{equation}\label{bura-discrete}
 R_{\alpha,k}(t) := \argmin_{s(t)\in \rat k}\,  
 \max_{t \in \{\frac{\lambda_1}{\lambda_N}, \frac{\lambda_2}{\lambda_N}, \dots, 1\}} |  s(t) - t^{\alpha} |.
\end{equation}

Unfortunately, such approximation depends of the knowledge of the eigenvalues, something we would like to avoid.
Now we shall show how to avoid in our computation such dependence for both solutions defined by \eqref{eq:sol-function} 
and \eqref{eq:sol-function-q}.

%Defining  the small parameter $\delta=\lambda_1/\lambda_N$ w
To find a computable approximation to 
\eqref{eq:sol-function} we introduce the following 
best uniform rational approximation (BURA) $r_{\delta,\alpha,k}(t) $ of $t^\alpha$ on $[\delta,1]$
$$
r_{\delta,\alpha,k}(t) := \argmin_{s(t)\in \rat k}\,   \sup_{t \in [\delta,1]} |  s(t) - t^{\alpha} |.
$$

Obviously we have 
$$
 \max_{t \in \{\frac{\lambda_1}{\lambda_N}, \frac{\lambda_2}{\lambda_N}, \dots, 1\}} |  s(t) - t^{\alpha} |
 \le 
 \sup_{t \in [\delta,1]} |  s(t) - t^{\alpha} |.
  $$

Often, for practical considerations, we would like to get rid of the parameter $\delta=\lambda_1/\lambda_N$ by 
using the best uniform rational approximation $ r_{\alpha,k}(t)$ of $t^\alpha$ on the whole interval $[0,1]$, namely
\begin{equation}\label{bura}
 r_{\alpha,k}(t) := \argmin_{s(t)\in \rat k}\,  
 \max_{t \in [0,1]} |  s(t) - t^{\alpha} | = \argmin_{s(t)\in \rat k}\, \| s(t) - t^{\alpha} \|_{L^\infty(0,1)}.
\end{equation}

The problem \eqref{bura}
has been studied extensively in the past, e.g. \cite{saff1992asymptotic,Stahl93,varga1992some}.  
Denoting the error by  
\beq\label{BURA-error}
E_{\alpha,k}:=\|r_{\alpha,k} (t) - t^{\alpha} \|_{L^\infty[0,1]},
\eeq
and applying Theorem 1 of \cite{Stahl93} we conclude that there is a constant 
$C_\alpha>0$,  independent of $k$, such that 
\begin{equation}\label{error-bound}
E_{\alpha,k}  \le C_\alpha e^{-2 \pi \sqrt{k \alpha}}.
\end{equation}
Thus, the BURA error converges exponentially to zero as $k$ becomes large. 

 It is known that the best rational approximation
$r_{\alpha,k}(t)=P_k(t)/Q_k(t)$  of $t^\alpha$  for $\alpha \in (0,1)$ is
non-degenerate, i.e., the polynomials $P_k(t)$ and $Q_k(t)$ are of full
degree $k$.   Denote  the roots of $P_k(t)$ and $Q_k(t)$ by
$\zeta_1, \dots, \zeta_k$ and  $d_1, \dots, d_k$, respectively.    
It is shown in
\cite{saff1992asymptotic,stahl2003}  that the roots are real, interlace and satisfy
\begin{equation}\label{interlacing}
0 > \zeta_1 > d_1 > \zeta_2 > d_2 > \cdots > \zeta_k > d_k.
\end{equation}
We then have
\beq
r_{\alpha,k} (t)=b\prod_{i=1}^k\frac{t-\zeta_i}{t-d_i}\label{prodr}
\eeq
where, by \eqref{interlacing} and the fact that
$r_{\alpha,k}(t)$ is a best approximation to a non-negative function, $b>0$ and
$P_k(t)>0$ and $Q_k(t)>0$ for $t\ge 0$. 

Knowing the poles $d_i$, $i=1, \dots, k$ we can give an equivalent representation of \eqref{prodr}
as a sum of partial fractions, namely
\beq
  r_{\alpha,k}(\t)=c_0+\sum_{i=1}^k \frac{c_i}{\t- d_i}
  \label{trg}
  \eeq
  where $c_0>0$ and $c_i<0$ for $i=1,\ldots,k$.
\def\trg{\widetilde r_{\alpha,k}}

Now changing the variable $\tt =1/\t$ in $r_{\alpha,k}(t)$ we get a rational function 
%We consider the rational function 
$\trg(\tt)$ defined by
\beq\label{eq:error-tilde}
\trg(\tt):=r_{\alpha,k}(1/\t)=\frac {\widetilde P(\tt)}{\widetilde Q(\tt)}.
\eeq
Here $\widetilde P_k(\tt)=\t^k P_k(\t^{-1})$ and $\widetilde Q_k(\tt)= \t^k Q_k(\t^{-1})$  and hence their coefficients 
are defined by reversing
the order of the coefficients in $P_k$ and $Q_k$ appearing in $r_{\alpha,k}(t)$.  
In
addition, \eqref{interlacing} implies that we have the following properties for  the roots of $\widetilde P_k$ and $\widetilde Q_k$, 
$\widetilde d_i=1/d_i$ and $\widetilde \zeta_i = 1/\zeta_i$,  respectively.
\begin{equation}\label{interlacing1}
  0 >  \widetilde d_k > \widetilde \zeta_k > \widetilde d_{k-1} > \widetilde
  \zeta_{k-1}\cdots >   \widetilde d_1>\widetilde \zeta_1.
\end{equation}  

\begin{remark}\label{lemma:BURA}  For $\alpha\in (0,1)$, 
%\beq
\beq\label{r-compute}
 \trg(\tt)=\widetilde c_0+\sum_{i=1}^k {\widetilde c_i}/{(\tt-\widetilde d_i)}
  % \label{trg-tilde}
   \eeq
 where 
  $$
  \widetilde c_0 = c_0 - \sum_{i=1}^k c_i \widetilde d_i=r_{\alpha,k}(0)=E_{\alpha,k}>0 \ \ \mbox{with} \  
  \widetilde c_i = - c_i { d}_i^{-2}  = - c_i {\widetilde d}_i^{2} >0,  \    i=1,\ldots,k.
  $$
\end{remark}

Indeed, % having in mind that $\widetilde d_i=1/d_i$, we get
$$
\widetilde r_{\alpha,k}(\tt) = r_{\alpha,k}(1/\tt) = c_0+\sum_{i=1}^k \frac{\ c_i}{1/\tt- d_i} 
= c_0 -\sum_{i=1}^k c_i d_i^{-1}  - \sum_{i=1}^k \frac{c_i d_i^{-2}}{\tt -  d_i^{-1}}
$$
and   having in mind that $\widetilde d_i=1/d_i$, we get \eqref{r-compute}.

%%%%%%%%%%%%%sssssssssss
 \subsection{Fully discrete schemes based on BURA}\label{ss:fully-discrete}
 %%%%%%%%%%%%%%%
 
 Now we introduce the {\it \underline{fully discrete}} approximations:  
$ w_h  \in V_h$ of the finite element approximation $u_h \in V_h$,  defined by \eqref{eq:sol-function},  and 
$\tilwh \in \R^N$ of the finite difference approximation $\tiluh  \in \R^N$ by
\beq
w_h 
=\lambda_1^{-\alpha} r_{\alpha,k} (\lambda_1 \calAt^{-1}) f_h
\quad \mbox{and}  \quad 
\tilwh 
=\lambda_1^{-\alpha} r_{\alpha,k} (\lambda_1 \wcalAt^{-1}) \tilfh.
\label{wh}
\eeq
Here $\calAt$ and $f_h$ are as in \eqref{fea} or \eqref{lumped-semi}
and $\wcalAt$ and $\tilfh$ are as in \eqref{fda}.

%corresponding to the finite element and finite difference methods.
In the paper \cite{harizanov2019analysis}, we studied the error of these fully discrete solutions.
For the finite element case we obtain the error estimate 
\beq
\|u_h - w_h\| \le \lambda_1^{-\alpha} E_{\alpha,k} \|f_h\|
\label{pbest}
\eeq
with $\|\cdot\|$ denoting the norm in $L^2(\Omega)$.
In the finite difference case,
we have
\beq
\|\tiluh - \tilwh  \|_{\ell_2} \le \lambda_1^{-\alpha} E_{\alpha,k} \|\tilfh \|_{\ell_2},
\label{pbest-fd}
\eeq
where the norm $\| \cdot \|_{\ell_2} $ %= (\tilwh^T \tilwh)^\frac12$  
denotes the  Euclidean norm in $\RR^N$.

%\input sect1ra-1.tex

%\sloppy

%\def\t{\lambda}
%\def\f{g}

\subsection{BURA approximation of ${\t^{\alpha}\over{1+q\,\t^{\alpha}}}  \mbox{ on }   [\delta,1] $}
\label{ss:bura-ura-res}
Now we consider the solution \eqref{eq:sol-function-q} and introduce 
the function $\f(q,\delta,\alpha;\t)$, 
of the variable $\t$ on $[\delta, 1]$,  $0 \le \delta < 1$
%$ \f_{q,\delta,\alpha}(\t):=
% ={\t^{\alpha}\over{1+q\,\t^{\alpha}}}  $  
and two parameters 
$q \in [0, \infty)$ and $0 < \alpha <1$: 
%  defined on the interval $[\delta, 1]$,  $0 \le \delta < 1$: 
$$
\f_{q,\delta,\alpha}(\t):=\f(q,\delta,\alpha;\t)={\t^{\alpha}\over{1+q\,\t^{\alpha}}} \  \mbox{ on }  \ \t \in [\delta,1]. 
$$
%Note that for $q=0$ we get the 
%standard fractional powers of an elliptic operator. 
Note that for $q=1/\tau$ we get the corresponding problem from time-discretization 
of sub-diffusion equation \eqref{eq:time-stepping}.
The role of this function is clear from the representation of the solution by \eqref{eq:sol-function-q}.
As before, our goal is to approximate this function using 
%
%  Here we employ   
the {\it best uniform rational approximation} (BURA). To find BURA of 
  $\f(q,\delta,\alpha;\t)$ we employ  Remez algorithm, cf. \cite{PGMASA1987}. 
  %={\t^{\alpha}\over{1+q\,\t^{\alpha}}}$ on $\t\in [\delta,1]$.
 
\medskip

\begin{definition}\label{pdef0bura}
The best uniform rational approximation ${r}_{q,\delta,\alpha,k}(\t) \in \mathcal R_k$
 of $ \f(q,\delta,\alpha; \t)$ on $[\delta,1]$, 
% and its approximation error $\widetilde{E}_{q,\delta,\alpha,k}$ 
called further $(q,\delta,\alpha,k)$-BURA,  is the rational function 
%are defined as follows:
\beq\label{BURA}
%\begin{aligned}
{r}_{q,\delta,\alpha,k}(\t):= % {r}(\f;\t) & =
\argmin_{ s \in \mathcal R_k}
 \left \| \f(q,\delta,\alpha;\t)- s(\t)  \right \|_{L^\infty[\delta,1]}.
% \end{aligned}
 \eeq
Then the error-function is denoted by
 \beq\label{BURA-er-eps}
% \begin{aligned}
%\widetilde {E}_{q,\delta,\alpha,k}  & = \left \| f(q,\delta,\alpha)- \widetilde{r}(q,\delta,\alpha,k)  \right \|\,\\
{\varepsilon}(q,\delta,\alpha,k;\t)  = {r}_{q,\delta,\alpha,k}(\t)-  \f(q,\delta,\alpha;\t), \\ 
%E_{q,\delta,\alpha,k} &:= \max_{\t\in [\delta,1]} \left| \varepsilon(q,\delta,\alpha,k;\t) \right|=
% \bigl\| f(q,\delta,\alpha)- r(q,\delta,\alpha,k) \bigr\|_{L^\infty[\delta,1)}
%\end{aligned}
\eeq
and its $L^\infty$-norm is denoted by
\beq\label{BURA-error-def}
 \begin{aligned}
{E}_{q,\delta,\alpha,k}  & =  %\sup_{t \in [\delta,1]} | | =
\sup_{\t \in [\delta,1]} | \f(q,\delta,\alpha;\t)- {r}_{q,\delta,\alpha,k}(\t) | 
 = \left \| {\varepsilon}(q,\delta,\alpha,k;\t) \right \|_{L^\infty[\delta,1]}. 
% \left \| f(q,\delta,\alpha;\t)- \widetilde{r}(q,\delta,\alpha,k;\t)  \right \|_{L^\infty[\delta,1]}. 
%{\varepsilon}(q,\delta,\alpha,k;\t) & = r(q,\delta,\alpha,k;\t)- f(q,\delta,\alpha;\t), \\ 
%E_{q,\delta,\alpha,k} &:= \max_{\t\in [\delta,1]} \left| \varepsilon(q,\delta,\alpha,k;\t) \right|=
% \bigl\| f(q,\delta,\alpha)- r(q,\delta,\alpha,k) \bigr\|_{L^\infty[\delta,1)}.
\end{aligned}
\eeq
%\red{Pencho, some comments about the outcome of this ???}
\end{definition}

We observe that the zeros and poles of ${r}_{q,\delta,\alpha,k}$ are again real, nonnegative, and interlacing for all considered choices of the four parameters $q,\delta,\alpha,k$. Furthermore, there seem to be $2k+2$ (the theoretically maximal possible number) points where ${\varepsilon}(q,\delta,\alpha,k;\t)$ achieves its extremal value $\pm {E}_{q,\delta,\alpha,k}$. However, we are not aware of a theoretical proof for any of the above observations in the general setting $q>0$ and/or $\delta>0$, which is not covered by Section~\ref{ss:BURA}.   
%  As shown in \cite[Lemma 2.2]{harizanov2019cmwa}, there are at lest $2k+2$ points where ${\varepsilon}(q,\delta,\alpha,k;\t)$
%  achieves its extremal value $\pm {E}_{q,\delta,\alpha,k}$. 
% {\red {Slavi, check if this is correct for $\delta \ne 0$. If not, please, make a comment about the experimental evidence !!!}}
 
% \red{We need to introduce the poles and the partial fraction representation !!!}
 
 Obviously, for $\tt =1/ \t$, $0 < \t <1$,
 \beq\label{eq:tilde}
  \widetilde \varepsilon(q,\delta,\alpha,k;\tt):= \varepsilon(q,\delta,\alpha,k;\t), \ \ \tt \in [1, 1/\delta]
  \eeq
 we have
 $$
 {E}_{q,\delta,\alpha,k}  
%\sup_{\t \in [\delta,1]} | f(q,\delta,\alpha;\t)- \widetilde{r}(q,\delta,\alpha,k;\t) | 
 = \left \| {\varepsilon}(q,\delta,\alpha,k;\t) \right \|_{L^\infty(\delta,1)} 
 = \left \| {\widetilde \varepsilon}(q,\delta,\alpha,k;\tt) \right \|_{L^\infty(1, 1/\delta)}.
% \le \left \| {\widetilde \varepsilon}(q,\delta,\alpha,k;\t) \right \|_{L^\infty(1, \infty)}.
 $$

\medskip

\subsection{URA %and URA 
approximation of $\f(q,\delta,\alpha;\t)={\t^{\alpha}\over{1+q\,\t^{\alpha}}}  \mbox{ on }   \t \in [\delta,1] $}
\label{ss:ura-res}
Now we shall introduce a rational approximations of function $\f(\t)$ that has simpler appearance and 
that is based on the BURA of $t^{\alpha}= \f(0,\delta,\alpha;\t)$.

\begin{definition}\label{pdef0ura}
The  function %(0-URA) 
\beq\label{0-URA}
%\bar{r}(q,\delta,\alpha,k;\t)  
\bar{r}_{q,\delta,\alpha,k}(\t)  := {{r}_{0,\delta,\alpha,k}(\t)\over {1+q\,{r}_{0,\delta,\alpha,k}(\t)}} \in \mathcal R_k
\eeq
is an approximation of $ \f(q,\delta,\alpha; \t)$ on $[\delta,1]$. 
% \[
%\begin{aligned}
%\bar{r}(q,\delta,\alpha,k;\t) & := {{r}(0,\delta,\alpha,k;\t)\over {1+q\,{r}(0,\delta,\alpha,k;\t)}}.
%\bar{E}_{q,\delta,\alpha,k} & = \left \| f(q,\delta,\alpha)- \bar{r}(q,\delta,\alpha,k)  \right \|_{L^\infty[\delta,1)} 
% = \max_{\t\in [\delta,1]} \left| \bar{\varepsilon}(q,\delta,\alpha,k;\t) \right|
%\end{aligned}
%\] 
Then the error-function is defined as
\[
\bar{\varepsilon}(q,\delta,\alpha,k;\t) =  \f(q,\delta,\alpha;\t)- \bar{r}_{q,\delta,\alpha,k}( \t) 
\]
 and 
 %its approximation error $\bar{E}_{q,\delta,\alpha,k}$ 
%(called further $(q,\delta,\alpha,k)$-URA or 0-$(q,\delta,\alpha,k)$-URA), are defined as:
\[
\begin{aligned}
%\bar{r}(q,\delta,\alpha,k;\t) & := {{r}(0,\delta,\alpha,k;\t)\over {1+q\,{r}(0,\delta,\alpha,k;\t)}},\\
\bar{E}_{q,\delta,\alpha,k} & = \left \| \f(q,\delta,\alpha;\t)- \bar{r}_{q,\delta,\alpha,k}(\t)  \right \|_{L^\infty[\delta,1)} 
 = \max_{\t\in [\delta,1]} \left| \bar{\varepsilon}(q,\delta,\alpha,k;\t) \right| .
\end{aligned}
\] 
%where ${r}(0,\delta,\alpha,k)$ is $(0,\delta,\alpha,k)$-BURA for $\f(0,\delta,\alpha,k).$
The rational function  $  \bar{r}_{q,\delta,\alpha,k}(\t)$ will be called $(q,\delta,\alpha,k)$-0-URA  
approximation of $ \f(q,\delta,\alpha;\t)$.
\end{definition}
Further, we present this rational function as a sum of partial fractions
\beq
  {\bar r}_{q,\delta,\alpha,k}(\t)=\bar c_0+\sum_{i=1}^k \frac{\bar c_i}{\t-\bar d_i}
  \label{trg-bar}
  \eeq
  where $\bar c_i>0$ and $\bar d_i$ are the poles, $i=0,1,\ldots,k$.
  
%  \red{Is this representation for $\t \in (0,1)$ or in $(1, \infty)$ ?????}

As shown in \cite[Theorem 2.4]{harizanov2019cmwa} the approximation error $ {E}_{q,\delta,\alpha,k} $ and 
$ \bar{E}_{q,\delta,\alpha,k}$ are related by 
$$
 \bar{E}_{q,\delta,\alpha,k}/( 1+q)^2  < {E}_{q,\delta,\alpha,k} <   \bar{E}_{q,\delta,\alpha,k}.
$$
%\red{Slavi, please, check this claim, I am not sure ?????}
The importance of this approximation is that by using corresponding
BURA of $t^\alpha$ we reduce the number of the parameters involved.  

\begin{remark}\label{prem0ura}
Now consider  
$q > q_0 >0$  and take 
%with  $q_0 >0$, $q_1 >0 $ and 
$\bar{r}_{q_0,\delta,\alpha,k}(\t)$ as 
$(q_0,\delta,\alpha,k;\t)$-0-URA approximation of  $\f(q_0,\delta,\alpha;\t)$. 
Then
\[ %\beq\label{0-URA}
\bar{r}_{q,\delta,\alpha,k}(\t)  % = \bar{r}(q_0 + q_1,\delta,\alpha,k;\t)
={\bar{r}_{q_0,\delta,\alpha,k}(\t)\over {1+(q - q_0)\bar{r}_{q_0,\delta,\alpha,k}(\t)}} 
={{r}_{0,\delta,\alpha,k}(\t)\over {1+q \, {r}_{0,\delta,\alpha,k}(\t)}}\,.
\] %\eeq
%% is $(q_0+q_1,\delta,\alpha,k)$-URA for $\f(q_0+q_1,\delta,\alpha;\t)$ 
%and
%$$
% \bar{r}(q_0+q_1,\delta,\alpha,k;\t)= {{r}(0,\delta,\alpha,k;\t)\over {1+(q_0+q_1)\, {r}(0,\delta,\alpha,k;\t)}}\,.
%$$
\end{remark}

%\red{And what about if we combine the  Definition~\ref{pdef0ura} with Remark~\ref{prem0ura}. }
  
\begin{definition}\label{pdef1ura}
Let $q =q_0 + q_1$, $q_0, q_1 >0$, and ${r}_{q_0,\delta,\alpha,k}(\t)$ be 
$(q_0,\delta,\alpha,k)$-BURA of  $\f(q_0,\delta,\alpha;\t)$.
A  rational function %of uniform approximation %(1-URA) 
$\bar{\bar{r}}_{q,\delta,\alpha,k}(\t) \in \mathcal R_k$  is an uniform approximation 
of $ \f(q,\delta,\alpha;\t)$ on $[\delta,1]$,  %for $q=q_0+q_1$ 
called  $(q,\delta,\alpha,k)$-1-URA  approximation,  and its error $\bar{\bar{E}}_{q,\delta,\alpha,k}$, 
are defined as:
\beq\label{1-URA}
%\begin{aligned}
\bar{\bar{r}}_{q,\delta,\alpha,k}(\t) %= \bar{\bar{r}}(q_0+q_1,\delta,\alpha,k;\t) & 
:= {{r}_{q_0,\delta,\alpha,k}(\t)\over {1+q_1 \, r_{q_0,\delta,\alpha,k}(\t)}}
%\end{aligned}
\eeq
and
\[
\begin{aligned}
\bar{\bar{E}}_{q,\delta,\alpha,k}  = \left \| \f(q,\delta,\alpha;\t)   - \bar{\bar{r}}_{q,\delta,\alpha,k}(\t)  \right \|_{L^\infty[\delta,1)} 
 = \sup_{\t\in [\delta,1]} \left| \bar{\bar{\varepsilon}}(q,\delta,\alpha,k;\t) \right| .
\end{aligned}
\] 
\end{definition}

Finally, we present this rational function as a sum of partial fractions
\beq\label{frac-1-URA}
\bar{\bar{r}}_{q,\delta,\alpha,k} (\t)=\bar{\bar c}_0+\sum_{i=1}^k \frac{\bar{\bar c}_i}{\t-\bar{\bar d}_i}
\eeq

We remark that the $(q,\delta,\alpha,k)$-1-URA  approximation gives a possibility to use a
previously computed   $(q_0,\delta,\alpha,k)$-BURA of  $\f(q_0,\delta,\alpha;\t)$ for a fixed $q_0$ 
in order to find an acceptable approximation for $q=q_0 + q_1$ with $q_0, q_1 >0$.

%\endinput

%\input sect1ra-p.tex

%\pagebreak

%    \part{Tables }\label{part:tabl}   %containing various data regarding BURA} 
\chapter{Tables }\label{chapt:tabl}  %containing various data regarding BURA} 

\section{Description of the data provided by our numerical experiments}\label{ss:num-ex}

%Clearly, the BURA is the most accurate approximation, but we need to establish theoretically 
%and experimentally which of 0-URA or 1-URA  is better.

\bigskip

%\subsection{Numerical experiments -- description}\label{ss:num-ex-dscr}

We provide all data for the uniform rational approximation of the function $\f(q,\delta,\alpha;\t)$
for various values of the 
% Our numerical experiments are performed for various values of the 
parameters $q,\delta,\alpha,k$.
 The Tables and the  corresponding files are named according the following encoding:

 \begin{itemize}
 \item[(a)] $ q    \in \{ 0,\,1,\,100,\,200,\,400\}$, coded  by $qQQQ,\ QQQ\in \{000,001,100,200,400\}$, total 5 parameters;
 
 \item[(b)] $\delta \in \{ 0.0,\, 10^{-6},\, 10^{-7},\,10^{-8}\}$, notation  $dD,\ D\in\{0,6,7,8\}$ -- 4 parameters;

 \item[(c)] $\alpha \in \{ 0.250,\, 0.500,\, 0.750\}$, notation   $aAA,\ AA\in\{25,50,75\}$ -- 3 parameters;

 \item[(d)] $ k  \in \{ 3,\,4,\,5,\,6,\,7,\,8 \}$, notation $kK,\ k\in\{3,\ldots,8\}$ -- 6 parameters;

 \item[(e)] $q_0$ and $q_1$ go by   
 pairs $(q_0,q_1)= (0,200)$,     $(q_0,q_1)=(100,100)$,   $(q_0,q_1)=(0,400)$,    $(q_0,q_1)=(200,200) $,
 %$(q_0,q_1)\in \{ (0,200),\, (100,100),\, (0,400),\, (200,200) \}$, 
 total of 4 cases,    
 coded as $\{qq02,\, qq11,\, qq04,\, qq22\}$.

\end{itemize}

 \medskip

 The  computational data is presented in a number of Tables that contain:

 \begin{enumerate}
  \item[(a)] the errors of  $(q,\delta,\alpha,k)$-BURA; 
  \item[(b)] the extreme points of ${\varepsilon}(q,\delta,\alpha,k;\t)$ defined by \eqref{BURA-er-eps};
  \item[(c)] the poles of ${r}_{q,\delta,\alpha,k}(\t)$ defined by \eqref{BURA}; 
  \item[(d)] the coefficients of decomposition $\bigl\{{c}_j\bigr\}_{j=0}^{k}$ of ${r}_{q,\delta,\alpha,k}(\t)$ as a sum of partial fractions. 
  \item[(e)] the poles of $\bar{r}_{q,\delta,\alpha,k}(\t)$ defined by \eqref{0-URA}; 
  \item[(f)]  the coefficients of the decomposition $\bigl\{\bar{c}_j\bigr\}_{j=0}^{k}$ of $\bar{r}_{q,\delta,\alpha,k}(\t)$ as a sum of partial fractions. 
  \item[(g)] the poles of $\bar{\bar{r}}_{q,\delta,\alpha,k}(\t)$ defined by \eqref{1-URA}; 
  \item[(h)] the coefficients of the decomposition $\bigl\{\bar{\bar{c}}_j\bigr\}_{j=0}^{k}$ of 
  $\bar{\bar{r}}_{q,\delta,\alpha,k}(\t)$ as a sum of partial fractions. 
 \end{enumerate}

\noindent 	
Short description of the type-tables:
\begin{description}
\item[(a-b)]  These tables correspond to characterization of the BURA.
\item[(c-d)]  These tables correspond to decomposition of BURA as a sum of partial fractions.
\item[(e-f)]  These tables correspond to decomposition of 0-URA as a sum of partial fractions.
\item[(g-h)]  These tables correspond to decomposition of 1-URA as a sum of partial fractions.
\end{description}

\bigskip

\begin{table}
\caption{ Table of Tables and Data-files} \label{p1ltabl}
\begin{tabular}{||c|c|c|c|c|c||}
\hline \hline 
 Table \& File  &Number of& Rows    & Cols   &Folder   &Number  \\
 qQQQdDaAAkK.txt& Tables  & ( q,k ) & (d, a) & name    &of Files\\ \hline
 0-head-tabl.txt&  ( 1 )  & ( 5*6 ) & (1+12) &BURA-tabl& 360    \\ \hline\hline
 qQQQdDaAAkK.txt&  (q,k)  & (     ) & (d, a) &         &        \\ \hline
 BqQQQkK.txt    &  (5*6)  & (2*K+2) & (2+12) &BURA-tabl& 360    \\ \hline
 CqQQQkK.txt    &  (5*6)  & (  K  ) & (2+12) &BURA-dcmp& 360    \\  
 DqQQQkK.txt    &  (5*6)  & ( K+1 ) & (2+12) &BURA-dcmp& 360    \\ \hline
 EqQQQkK.txt    &  (4*6)  & (  K  ) & (2+12) &0URA-dcmp& 288    \\
 FqQQQkK.txt    &  (4*6)  & ( K+1 ) & (2+12) &0URA-dcmp& 288    \\ \hline
 GqqQQkK.txt    &  (2*6)  & (  K  ) & (2+12) &1URA-dcmp& 144    \\
 HqqQQkK.txt    &  (2*6)  & ( K+1 ) & (2+12) &1URA-dcmp& 144    \\
\hline \hline
\end{tabular}
\end{table}

Total number of Tables is $(1+90+48+24)=163$.

\noindent
More descriptions for the tables and files:
\begin{description}
\item[(a-b)]  The data for these tables in files of types (a),(b) correspond to characterization of the BURA.
  Folder with values of BURA and extreme points has $360$ files, and it is named \texttt{BURA-tabl}.
 Only for 5 cases the program did not finished with a solution.

\item[(c-d)]  The files with normalized coefficients $A$ and $B$, poles $d_j$ (named $U0(j)$)
 and coefficients $c_j$ (named $E(j)$) are in the
 folder  named \texttt{BURA-dcmp}. One sub-folder more was present (\texttt{BURA-dcmp/add/})
with more details about poles and coefficients -- its Imaginary parts. Number of files is also $360$
 and have names \texttt{qQQQdDaAAkK}.

\item[(e-f)]  The files have the same names as in the previous item \textbf{(c-d)} 
and the same structure but folder is other.
The folder is named \texttt{0URA-dcmp/}. One sub-folder more is given (\texttt{0URA-dcmp/add/})
with more details about poles and coefficients -- its Imaginary parts.
 Number of files is $288$, because the cases
 $(0,\delta, \lambda)$-URA coincide with $(0,\delta, \lambda)$-BURA and files with $QQQ=000$ are not present.

\item[(g-h)]  The folder with files is \texttt{1URA-dcmp/} and sub-folder \texttt{1URA-dcmp/add/}.
	The names of files are \texttt{qqQQdDaAAkK} and as the cases \texttt{qqQQ = qq02, qq04} are not present
	because coincide with the cases \texttt{qQQQdDaAAkK}, \texttt{qQQQ = q200, q400} from \texttt{0URA-dcmp},
	and number of files is $144$ -- for \texttt{qqQQ = qq11, qq22} only. 
\end{description}

% portrait tables (only one is doubled and rotated)

%\newpage

%\pagebreak

\section{Tables type \textbf{(a-b)} for BURA-errors and BURA-extreme points}

%\include{sect1ra-pm0}

%\subsubsection{Tables type \textbf{(a-b)} for BURA-errors and BURA-extreme points}

%  \input table-landscape.tex

It is clear from Table \ref{tabl:A0BURAp} that for  fixed parameters $\alpha$, $q$, and $k$, 
the error is increasing, when $\delta \to 0$. However, the  differences are not that pronounced, so for practical purposes one can
use for all $\delta$ the approximations for $\delta=0$. The significance of using $\delta >0$ 
is in the performance of Remez algorithm for computing BURA. 

One should realize that BURA-based  methods involve 
Remez method of finding the best uniform rational approximation
by solving the highly non-linear min-max problem \eqref{bura}. 
It is well known that Remez algorithm is very sensitive to the precision of the computer arithmetic, 
cf. \cite{PGMASA1987,varga1992some,Driscoll2014}.  Various techniques for stabilization
of the method have been used, mostly by using Tchebyshev orthogonal polynomials, cf.
\cite{PGMASA1987,Driscoll2014}. It seems that to achieve high accuracy one needs to use high arithmetic precision.
For example, in \cite{varga1992some}
the errors of best uniform rational approximation of $t^\alpha$ for six values of $\alpha  \in [0,1]$ are reported for
degree $k \le 30$ by using computer arithmetic with 200 significant digits.
In short, for $\delta >0$ the Remez algorithm has substantially better stability and is significantly more reliable.

We also note that in Table \ref{tabl:A0BURAp} there are  5 sets of parameters (all of 
them for $\alpha=0.25$)  for which Remez algorithm dies not 
provide the needed information.  In these cases the convergence in the iterative process for finding the extremal points of 
BURA for these parameters either does not converge or fail to produce equal absolute values at the extremal points with the
desired accuracy.

~ ~ ~

\begin{tiny}
\begin{table}
\setlength{\tabcolsep}{1.6pt}
\renewcommand{\arraystretch}{1.02}
\caption{The error $E_{q, \delta,\alpha, k}$ of $(q,\delta,\alpha,k)$-BURA approximation, 
i.e. the best uniform rational approximation of  
$\f(q,\delta,\alpha;\t)=\t^{\alpha}/(1+q\,\t^{\alpha})$  with functions from $\mathcal R_k$ for 
 $t \in [\delta, 1]$,  $\alpha =0.25, 0.5, 0.75$, $k=3, \dots, 8$, and $q=0,1,100,200,400$
}
\label{tabl:A0BURAp}
%\rotatebox{90}{
\begin{tabular}{||c||c|c|c||c|c|c||c|c|c||c|c|c||}
\hline \hline
d ($\delta$) & $0.00$ & $0.00$ & $0.00$ &$10^{-6}$&$10^{-6}$&$10^{-6}$&$10^{-7}$&$10^{-7}$&$10^{-7}$&$10^{-8}$&$10^{-8}$&$10^{-8}$\\  \hline
q,k  & $\alpha =0.25$   & $\alpha =0.50$   & $\alpha =0.75$   & $\alpha =0.25$   & $\alpha =0.50$   & $\alpha =0.75$   & $\alpha =0.25$   & $\alpha =0.50$   & $\alpha =0.75 $  & $\alpha =0.25 $  & $\alpha =0.50$   & $\alpha =0.75 $  \\ \hline
000,\,3&1.235E-2&2.282E-3&4.041E-4&7.364E-3&2.107E-3&3.993E-4&9.217E-3&2.223E-3&4.032E-4&1.048E-2&2.263E-3&4.039E-4\\
000,\,4&5.566E-3&7.366E-4&9.954E-5&2.270E-3&6.066E-4&9.591E-5&3.263E-3&6.880E-4&9.878E-5&4.082E-3&7.202E-4&9.939E-5\\
000,\,5&2.735E-3&2.690E-4&2.868E-5&6.972E-4&1.823E-4&2.607E-5&1.159E-3&2.309E-4&2.806E-5&1.619E-3&2.550E-4&2.855E-5\\
000,\,6&1.431E-3&1.075E-4&9.252E-6&2.134E-4&5.534E-5&7.548E-6&4.104E-4&8.015E-5&8.773E-6&6.424E-4&9.598E-5&9.146E-6\\
000,\,7&7.865E-4&4.604E-5&3.257E-6&6.522E-5&1.683E-5&2.245E-6&1.450E-4&2.811E-5&2.907E-6&2.545E-4&3.722E-5&3.168E-6\\
000,\,8&4.495E-4&2.085E-5&1.229E-6&1.990E-5&5.121E-6&6.751E-7&5.119E-5&9.881E-6&9.930E-7&1.007E-4&1.460E-5&1.159E-6\\
\hline \hline
\end{tabular}
%}
\end{table}
\end{tiny}

~ ~ ~

\begin{tiny}
\begin{table}
\setlength{\tabcolsep}{2.0pt}
\renewcommand{\arraystretch}{1.02}
\caption{The extreme points of ${\varepsilon}(q,\delta,\alpha,k;\t)$
 of $(q,\delta,\alpha,k)$-BURA approximation, 
i.e. the best uniform rational approximation of  
$\f(q,\delta,\alpha;\t)=\t^{\alpha}/(1+q\,\t^{\alpha})$  with functions from $\mathcal R_k$ on 
 $ [\delta, 1]$,   $\alpha =0.25, 0.5, 0.75$, ${\bf k=3}$, and $q=0$} %,1,100,200,400$}
 \label{tabl:Bq000dDaAAk3p}
%\rotatebox{90}{
\begin{tabular}{||c||c|c|c||c|c|c||c|c|c||c|c|c||}
\hline \hline
{\bf k=3} /$\delta$& $0.00$ & $0.00$ & $0.00$ &$10^{-6}$&$10^{-6}$&$10^{-6}$&$10^{-7}$&$10^{-7}$&$10^{-7}$&$10^{-8}$&$10^{-8}$&$10^{-8}$\\
\ j  & $\alpha =0.25$   & 0.50   & 0.75   & 0.25   & 0.50   & 0.75   & 0.25   & 0.50   & 0.75   & 0.25   & 0.50   & 0.75   \\ \hline
1&0.000E0&0.000E0&0.000E0&1.000E-6&1.000E-6&1.000E-6&1.000E-7&1.000E-7&1.000E-7&1.000E-8&1.000E-8&1.000E-8\\
2&1.342E-6&1.037E-4&6.496E-4&1.645E-5&1.354E-4&6.702E-4&6.389E-6&1.134E-4&6.535E-4&3.403E-6&1.067E-4&6.503E-4\\
3&4.731E-5&1.612E-3&6.907E-3&2.440E-4&1.858E-3&7.007E-3&1.306E-4&1.690E-3&6.926E-3&8.650E-5&1.637E-3&6.911E-3\\
4&1.011E-3&1.263E-2&3.708E-2&2.682E-3&1.374E-2&3.739E-2&1.819E-3&1.299E-2&3.714E-2&1.423E-3&1.274E-2&3.709E-2\\
5&1.125E-2&6.819E-2&1.423E-1&2.100E-2&7.163E-2&1.429E-1&1.635E-2&6.931E-2&1.424E-1&1.397E-2&6.855E-2&1.423E-1\\
6&9.395E-2&2.688E-1&4.016E-1&1.311E-1&2.754E-1&4.024E-1&1.145E-1&2.710E-1&4.017E-1&1.053E-1&2.695E-1&4.016E-1\\
7&4.887E-1&7.013E-1&7.896E-1&5.472E-1&7.061E-1&7.900E-1&5.229E-1&7.029E-1&7.896E-1&5.082E-1&7.018E-1&7.896E-1\\
8&1.000E0&1.000E0&1.000E0&1.000E0&1.000E0&1.000E0&1.000E0&1.000E0&1.000E0&1.000E0&1.000E0&1.000E0\\ \hline  \hline 
\end{tabular}
%}
\end{table}
\end{tiny}

~ ~ ~

\begin{tiny}
\begin{table}
\setlength{\tabcolsep}{2.0pt}
\renewcommand{\arraystretch}{1.02}
\caption{The extreme points of ${\varepsilon}(q,\delta,\alpha,k;\t)$
 of $(q,\delta,\alpha,k)$-BURA approximation, 
i.e. the best uniform rational approximation of  
$\f(q,\delta,\alpha;\t)=\t^{\alpha}/(1+q\,\t^{\alpha})$  with functions from $\mathcal R_k$ on 
 $ [\delta, 1]$,   $\alpha =0.25, 0.5, 0.75$, ${\bf k=4}$, and $q=0$}
 %$(q=000,\delta,\alpha,k=4)$-BURA }
 \label{tabl:Bq000dDaAAk4p}
%\rotatebox{90}{
\begin{tabular}{||c||c|c|c||c|c|c||c|c|c||c|c|c||}
\hline \hline
{\bf k=4 } /$\delta$& $0.00$ & $0.00$ & $0.00$ &$10^{-6}$&$10^{-6}$&$10^{-6}$&$10^{-7}$&$10^{-7}$&$10^{-7}$&$10^{-8}$&$10^{-8}$&$10^{-8}$\\
\ j & $\alpha=0.25$   & $\alpha= 0.50$   & $\alpha= 0.75$   & $\alpha= 0.25$   & $\alpha=0.50$   & $\alpha=0.75$   & $\alpha= 0.25$   & $\alpha= 0.50$   & $\alpha=0.75$   & $\alpha= 0.25$   & $\alpha= 0.50$   & $\alpha= 0.75 $  \\ \hline
1&0.000E0&0.000E0&0.000E0&1.000E-6&1.000E-6&1.000E-6&1.000E-7&1.000E-7&1.000E-7&1.000E-8&1.000E-8&1.000E-8\\
2&5.542E-8&1.080E-5&1.004E-4&5.079E-6&2.285E-5&1.129E-4&1.213E-6&1.427E-5&1.029E-4&4.020E-7&1.186E-5&1.009E-4\\
3&1.954E-6&1.682E-4&1.073E-3&4.485E-5&2.581E-4&1.137E-3&1.621E-5&1.967E-4&1.086E-3&7.517E-6&1.772E-4&1.076E-3\\
4&4.176E-5&1.328E-3&5.879E-3&3.420E-4&1.756E-3&6.102E-3&1.655E-4&1.469E-3&5.925E-3&9.807E-5&1.373E-3&5.888E-3\\
5&4.680E-4&7.419E-3&2.406E-2&2.134E-3&8.969E-3&2.466E-2&1.251E-3&7.944E-3&2.418E-2&8.568E-4&7.590E-3&2.408E-2\\
6&4.099E-3&3.305E-2&8.003E-2&1.167E-2&3.752E-2&8.132E-2&8.011E-3&3.459E-2&8.030E-2&6.165E-3&3.355E-2&8.008E-2\\
7&2.769E-2&1.216E-1&2.213E-1&5.551E-2&1.316E-1&2.234E-1&4.315E-2&1.251E-1&2.217E-1&3.625E-2&1.228E-1&2.214E-1\\
8&1.525E-1&3.585E-1&4.932E-1&2.239E-1&3.729E-1&4.954E-1&1.947E-1&3.636E-1&4.937E-1&1.769E-1&3.602E-1&4.933E-1\\
9&5.753E-1&7.597E-1&8.327E-1&6.515E-1&7.681E-1&8.336E-1&6.231E-1&7.627E-1&8.329E-1&6.040E-1&7.607E-1&8.327E-1\\
10&1.000E0&1.000E0&1.000E0&1.000E0&1.000E0&1.000E0&1.000E0&1.000E0&1.000E0&1.000E0&1.000E0&1.000E0\\ \hline  \hline 
\end{tabular}
%}
\end{table}
\end{tiny}

~ ~ ~

\begin{tiny}
\begin{table}
\setlength{\tabcolsep}{2.0pt}
\renewcommand{\arraystretch}{1.02}
\caption{The extreme points of ${\varepsilon}(q,\delta,\alpha,k;\t)$
 of $(q,\delta,\alpha,k)$-BURA approximation, 
i.e. the best uniform rational approximation of  
$\f(q,\delta,\alpha;\t)=\t^{\alpha}/(1+q\,\t^{\alpha})$  with functions from $\mathcal R_k$ on 
 $ [\delta, 1]$,   $\alpha =0.25, 0.5, 0.75$, ${\bf k=5}$, and $q=0$}
 %$(q=000,\delta,\alpha,k=5)$-BURA }
\label{tabl:Bq000dDaAAk5p}
%\rotatebox{90}{
\begin{tabular}{||c||c|c|c||c|c|c||c|c|c||c|c|c||}
\hline \hline
{\bf k=5} /$\delta$& $0.00$ & $0.00$ & $0.00$ &$10^{-6}$&$10^{-6}$&$10^{-6}$&$10^{-7}$&$10^{-7}$&$10^{-7}$&$10^{-8}$&$10^{-8}$&$10^{-8}$\\
\ j & $\alpha=0.25$   & $\alpha=0.50$   & $\alpha=0.75$   & $\alpha=0.25$   & $\alpha=0.50$   & $\alpha=0.75$   & $\alpha=0.25$   & $\alpha=0.50$   & $\alpha=0.75$   & $\alpha=0.25$   & $\alpha=0.50$   & $\alpha= 0.75 $  \\ \hline
1&0.000E0&0.000E0&0.000E0&1.000E-6&1.000E-6&1.000E-6&1.000E-7&1.000E-7&1.000E-7&1.000E-8&1.000E-8&1.000E-8\\
2&3.229E-9&1.440E-6&1.910E-5&2.756E-6&7.075E-6&2.676E-5&4.753E-7&2.875E-6&2.069E-5&1.049E-7&1.857E-6&1.942E-5\\
3&1.139E-7&2.242E-5&2.044E-4&1.525E-5&5.885E-5&2.437E-4&4.027E-6&3.373E-5&2.129E-4&1.327E-6&2.601E-5&2.061E-4\\
4&2.434E-6&1.772E-4&1.124E-3&8.372E-5&3.454E-4&1.265E-3&2.984E-5&2.338E-4&1.155E-3&1.302E-5&1.958E-4&1.130E-3\\
5&2.728E-5&9.942E-4&4.650E-3&4.129E-4&1.614E-3&5.057E-3&1.827E-4&1.213E-3&4.741E-3&9.609E-5&1.068E-3&4.669E-3\\
6&2.395E-4&4.496E-3&1.597E-2&1.864E-3&6.420E-3&1.697E-2&9.882E-4&5.200E-3&1.620E-2&6.054E-4&4.735E-3&1.602E-2\\
7&1.640E-3&1.739E-2&4.777E-2&7.730E-3&2.254E-2&4.988E-2&4.749E-3&1.933E-2&4.825E-2&3.275E-3&1.806E-2&4.787E-2\\
8&9.687E-3&5.904E-2&1.261E-1&2.972E-2&7.073E-2&1.298E-1&2.079E-2&6.352E-2&1.270E-1&1.587E-2&6.059E-2&1.263E-1\\
9&4.899E-2&1.748E-1&2.903E-1&1.043E-1&1.959E-1&2.954E-1&8.189E-2&1.830E-1&2.915E-1&6.826E-2&1.776E-1&2.905E-1\\
10&2.089E-1&4.301E-1&5.608E-1&3.172E-1&4.555E-1&5.655E-1&2.776E-1&4.402E-1&5.619E-1&2.510E-1&4.336E-1&5.610E-1\\
11&6.373E-1&7.990E-1&8.611E-1&7.260E-1&8.118E-1&8.629E-1&6.969E-1&8.042E-1&8.615E-1&6.754E-1&8.009E-1&8.612E-1\\
12&1.000E0&1.000E0&1.000E0&1.000E0&1.000E0&1.000E0&1.000E0&1.000E0&1.000E0&1.000E0&1.000E0&1.000E0\\ \hline  \hline 
\end{tabular}
%}
\end{table}
\end{tiny}

~ ~ ~

\begin{tiny}
\begin{table}
\setlength{\tabcolsep}{2.0pt}
\renewcommand{\arraystretch}{1.02}
\caption{The extreme points of ${\varepsilon}(q,\delta,\alpha,k;\t)$
 of $(q,\delta,\alpha,k)$-BURA approximation, 
i.e. the best uniform rational approximation of  
$\f(q,\delta,\alpha;\t)=\t^{\alpha}/(1+q\,\t^{\alpha})$  with functions from $\mathcal R_k$ on 
 $ [\delta, 1]$,   $\alpha =0.25, 0.5, 0.75$, ${\bf k=6}$, and $q=0$}
 %$(q=000,\delta,\alpha,k=6)$-BURA }
\label{tabl:Bq000dDaAAk6p}
%\rotatebox{90}{
\begin{tabular}{||c||c|c|c||c|c|c||c|c|c||c|c|c||}
\hline \hline
{\bf k=6} /$\delta$& $0.00$ & $0.00$ & $0.00$ &$10^{-6}$&$10^{-6}$&$10^{-6}$&$10^{-7}$&$10^{-7}$&$10^{-7}$&$10^{-8}$&$10^{-8}$&$10^{-8}$\\
\ j /$\alpha$ & 0.25   & 0.50   & 0.75   & 0.25   & 0.50   & 0.75   & 0.25   & 0.50   & 0.75   & 0.25   & 0.50   & 0.75   \\ \hline
1&0.000E0&0.000E0&0.000E0&1.000E-6&1.000E-6&1.000E-6&1.000E-7&1.000E-7&1.000E-7&1.000E-8&1.000E-8&1.000E-8\\
2&2.422E-10&2.299E-7&4.228E-6&1.971E-6&3.488E-6&9.101E-6&2.804E-7&9.260E-7&5.245E-6&4.709E-8&4.131E-7&4.435E-6\\
3&8.541E-9&3.581E-6&4.524E-5&7.510E-6&2.039E-5&6.892E-5&1.577E-6&8.448E-6&5.075E-5&3.966E-7&5.095E-6&4.640E-5\\
4&1.826E-7&2.830E-5&2.489E-4&3.133E-5&9.958E-5&3.337E-4&8.760E-6&5.191E-5&2.694E-4&2.926E-6&3.617E-5&2.532E-4\\
5&2.046E-6&1.589E-4&1.032E-3&1.252E-4&4.124E-4&1.281E-3&4.359E-5&2.499E-4&1.094E-3&1.787E-5&1.905E-4&1.046E-3\\
6&1.797E-5&7.200E-4&3.569E-3&4.734E-4&1.505E-3&4.205E-3&1.984E-4&1.018E-3&3.729E-3&9.646E-5&8.267E-4&3.604E-3\\
7&1.232E-4&2.806E-3&1.086E-2&1.692E-3&4.975E-3&1.231E-2&8.304E-4&3.666E-3&1.123E-2&4.636E-4&3.120E-3&1.094E-2\\
8&7.307E-4&9.740E-3&2.987E-2&5.754E-3&1.513E-2&3.283E-2&3.247E-3&1.195E-2&3.063E-2&2.044E-3&1.056E-2&3.003E-2\\
9&3.775E-3&3.067E-2&7.502E-2&1.861E-2&4.267E-2&8.039E-2&1.189E-2&3.574E-2&7.641E-2&8.300E-3&3.258E-2&7.532E-2\\
10&1.759E-2&8.786E-2&1.715E-1&5.707E-2&1.110E-1&1.798E-1&4.088E-2&9.790E-2&1.736E-1&3.132E-2&9.169E-2&1.719E-1\\
11&7.316E-2&2.247E-1&3.495E-1&1.621E-1&2.607E-1&3.598E-1&1.294E-1&2.406E-1&3.522E-1&1.081E-1&2.308E-1&3.501E-1\\
12&2.609E-1&4.879E-1&6.126E-1&4.033E-1&5.256E-1&6.210E-1&3.569E-1&5.050E-1&6.148E-1&3.234E-1&4.946E-1&6.130E-1\\
13&6.836E-1&8.273E-1&8.812E-1&7.801E-1&8.442E-1&8.844E-1&7.523E-1&8.351E-1&8.821E-1&7.303E-1&8.303E-1&8.814E-1\\
14&1.000E0&1.000E0&1.000E0&1.000E0&1.000E0&1.000E0&1.000E0&1.000E0&1.000E0&1.000E0&1.000E0&1.000E0\\ \hline  \hline 
\end{tabular}
%}
\end{table}
\end{tiny}

~ ~ ~

\begin{tiny}
\begin{table}
\setlength{\tabcolsep}{2.0pt}
\renewcommand{\arraystretch}{1.02}
\caption{The extreme points of ${\varepsilon}(q,\delta,\alpha,k;\t)$
 of $(q,\delta,\alpha,k)$-BURA approximation, 
i.e. the best uniform rational approximation of  
$\f(q,\delta,\alpha;\t)=\t^{\alpha}/(1+q\,\t^{\alpha})$  with functions from $\mathcal R_k$ on 
 $ [\delta, 1]$,   $\alpha =0.25, 0.5, 0.75$, ${\bf k=7}$, and $q=0$}
 %$(q=000,\delta,\alpha,k=7)$-BURA }
 \label{tabl:Bq000dDaAAk7p}
%\rotatebox{90}{
\begin{tabular}{||c||c|c|c||c|c|c||c|c|c||c|c|c||}
\hline \hline
{\bf k=7} /$\delta$& $0.00$ & $0.00$ & $0.00$ &$10^{-6}$&$10^{-6}$&$10^{-6}$&$10^{-7}$&$10^{-7}$&$10^{-7}$&$10^{-8}$&$10^{-8}$&$10^{-8}$\\
\ j /$\alpha$ & 0.25   & 0.50   & 0.75   & 0.25   & 0.50   & 0.75   & 0.25   & 0.50   & 0.75   & 0.25   & 0.50   & 0.75   \\ \hline
1&0.000E0&0.000E0&0.000E0&1.000E-6&1.000E-6&1.000E-6&1.000E-7&1.000E-7&1.000E-7&1.000E-8&1.000E-8&1.000E-8\\
2&2.209E-11&4.219E-8&1.051E-6&1.618E-6&2.301E-6&4.350E-6&2.055E-7&4.423E-7&1.716E-6&2.914E-8&1.330E-7&1.189E-6\\
3&7.789E-10&6.570E-7&1.124E-5&4.642E-6&9.653E-6&2.558E-5&8.260E-7&2.977E-6&1.473E-5&1.699E-7&1.342E-6&1.202E-5\\
4&1.665E-8&5.192E-6&6.187E-5&1.545E-5&3.892E-5&1.116E-4&3.587E-6&1.568E-5&7.489E-5&9.675E-7&8.661E-6&6.485E-5\\
5&1.866E-7&2.915E-5&2.568E-4&5.160E-5&1.413E-4&4.016E-4&1.480E-5&6.830E-5&2.963E-4&4.905E-6&4.302E-5&2.660E-4\\
6&1.639E-6&1.322E-4&8.891E-4&1.673E-4&4.669E-4&1.262E-3&5.749E-5&2.590E-4&9.939E-4&2.269E-5&1.793E-4&9.136E-4\\
\hline  \hline 
\end{tabular}
%}
\end{table}
\end{tiny}

~ ~ ~

\begin{tiny}
\begin{table}
\setlength{\tabcolsep}{2.0pt}
\renewcommand{\arraystretch}{1.02}
\caption{The extreme points of ${\varepsilon}(q,\delta,\alpha,k;\t)$
 of $(q,\delta,\alpha,k)$-BURA approximation, 
i.e. the best uniform rational approximation of  
$\f(q,\delta,\alpha;\t)=\t^{\alpha}/(1+q\,\t^{\alpha})$  with functions from $\mathcal R_k$ on 
 $ [\delta, 1]$,   $\alpha =0.25, 0.5, 0.75$, ${\bf k=8}$, and $q=0$}
 %$(q=000,\delta,\alpha,k=8)$-BURA }
 \label{tabl:Bq000dDaAAk8p}
%\rotatebox{90}{
\begin{tabular}{||c||c|c|c||c|c|c||c|c|c||c|c|c||}
\hline \hline
{\bf k=8} /$\delta$& $0.00$ & $0.00$ & $0.00$ &$10^{-6}$&$10^{-6}$&$10^{-6}$&$10^{-7}$&$10^{-7}$&$10^{-7}$&$10^{-8}$&$10^{-8}$&$10^{-8}$\\
\ j /$\alpha$ & 0.25   & 0.50   & 0.75   & 0.25   & 0.50   & 0.75   & 0.25   & 0.50   & 0.75   & 0.25   & 0.50   & 0.75   \\ \hline
1&0.000E0&0.000E0&0.000E0&1.000E-6&1.000E-6&1.000E-6&1.000E-7&1.000E-7&1.000E-7&1.000E-8&1.000E-8&1.000E-8\\
2&2.357E-12&8.655E-9&2.865E-7&1.429E-6&1.792E-6&2.706E-6&1.694E-7&2.780E-7&7.363E-7&2.163E-8&6.023E-8&3.800E-7\\
3&8.310E-11&1.348E-7&3.066E-6&3.315E-6&5.682E-6&1.207E-5&5.222E-7&1.378E-6&5.264E-6&9.249E-8&4.698E-7&3.584E-6\\
4&1.776E-9&1.065E-6&1.687E-5&9.172E-6&1.910E-5&4.625E-5&1.844E-6&6.137E-6&2.492E-5&4.208E-7&2.669E-6&1.886E-5\\
5&1.991E-8&5.981E-6&7.002E-5&2.630E-5&6.090E-5&1.534E-4&6.458E-6&2.387E-5&9.440E-5&1.799E-6&1.222E-5&7.622E-5\\
6&1.748E-7&2.712E-5&2.426E-4&7.474E-5&1.820E-4&4.549E-4&2.181E-5&8.316E-5&3.075E-4&7.208E-6&4.808E-5&2.594E-4\\
\hline  \hline 
\end{tabular}
%}
\end{table}
\end{tiny}

~ ~ ~

\begin{tiny}
\begin{table}
\setlength{\tabcolsep}{2.0pt}
\renewcommand{\arraystretch}{1.02}
\caption{The extreme points of ${\varepsilon}(q,\delta,\alpha,k;\t)$
 of $(q,\delta,\alpha,k)$-BURA approximation, 
i.e. the best uniform rational approximation of  
$\f(q,\delta,\alpha;\t)=\t^{\alpha}/(1+q\,\t^{\alpha})$  with functions from $\mathcal R_k$ on 
 $ [\delta, 1]$,   $\alpha =0.25, 0.5, 0.75$, ${\bf k=3}$, and $q=1$}
 %$(q=001,\delta,\alpha,k=3)$-BURA }
 \label{tabl:Bq001dDaAAk3p}
%\rotatebox{90}{
\begin{tabular}{||c||c|c|c||c|c|c||c|c|c||c|c|c||}
\hline \hline
{\bf k=3} /$\delta$& $0.00$ & $0.00$ & $0.00$ &$10^{-6}$&$10^{-6}$&$10^{-6}$&$10^{-7}$&$10^{-7}$&$10^{-7}$&$10^{-8}$&$10^{-8}$&$10^{-8}$\\
\ j /$\alpha$ & 0.25   & 0.50   & 0.75   & 0.25   & 0.50   & 0.75   & 0.25   & 0.50   & 0.75   & 0.25   & 0.50   & 0.75   \\ \hline
1&0.000E0&0.000E0&0.000E0&1.000E-6&1.000E-6&1.000E-6&1.000E-7&1.000E-7&1.000E-7&1.000E-8&1.000E-8&1.000E-8\\
2&3.498E-7&5.726E-5&4.471E-4&9.698E-6&8.172E-5&4.657E-4&3.021E-6&6.467E-5&4.507E-4&1.317E-6&5.957E-5&4.478E-4\\
3&1.288E-5&8.925E-4&4.723E-3&1.229E-4&1.083E-3&4.813E-3&5.501E-5&9.531E-4&4.740E-3&3.130E-5&9.117E-4&4.726E-3\\
4&3.074E-4&7.073E-3&2.520E-2&1.280E-3&7.966E-3&2.549E-2&7.497E-4&7.363E-3&2.526E-2&5.237E-4&7.165E-3&2.521E-2\\
5&3.893E-3&3.986E-2&9.862E-2&1.037E-2&4.293E-2&9.929E-2&7.153E-3&4.087E-2&9.875E-2&5.591E-3&4.018E-2&9.864E-2\\
6&4.179E-2&1.777E-1&3.040E-1&7.497E-2&1.853E-1&3.051E-1&5.992E-2&1.802E-1&3.042E-1&5.172E-2&1.785E-1&3.040E-1\\
7&3.259E-1&5.902E-1&7.089E-1&4.194E-1&5.990E-1&7.097E-1&3.815E-1&5.931E-1&7.090E-1&3.580E-1&5.912E-1&7.089E-1\\
8&1.000E0&1.000E0&1.000E0&1.000E0&1.000E0&1.000E0&1.000E0&1.000E0&1.000E0&1.000E0&1.000E0&1.000E0\\ \hline  \hline 
\end{tabular}
%}
\end{table}
\end{tiny}

~ ~ ~

\begin{tiny}
\begin{table}
\setlength{\tabcolsep}{2.0pt}
\renewcommand{\arraystretch}{1.02}
\caption{The extreme points of ${\varepsilon}(q,\delta,\alpha,k;\t)$
 of $(q,\delta,\alpha,k)$-BURA approximation, 
i.e. the best uniform rational approximation of  
$\f(q,\delta,\alpha;\t)=\t^{\alpha}/(1+q\,\t^{\alpha})$  with functions from $\mathcal R_k$ on 
 $ [\delta, 1]$,   $\alpha =0.25, 0.5, 0.75$, ${\bf k=4}$, and $q=1$}
 %$(q=001,\delta,\alpha,k=4)$-BURA }
 \label{tabl:Bq001dDaAAk4p}
%\rotatebox{90}{
\begin{tabular}{||c||c|c|c||c|c|c||c|c|c||c|c|c||}
\hline \hline
{\bf k=4} /$\delta$& $0.00$ & $0.00$ & $0.00$ &$10^{-6}$&$10^{-6}$&$10^{-6}$&$10^{-7}$&$10^{-7}$&$10^{-7}$&$10^{-8}$&$10^{-8}$&$10^{-8}$\\
\ j /$\alpha$ & 0.25   & 0.50   & 0.75   & 0.25   & 0.50   & 0.75   & 0.25   & 0.50   & 0.75   & 0.25   & 0.50   & 0.75   \\ \hline
1&0.000E0&0.000E0&0.000E0&1.000E-6&1.000E-6&1.000E-6&1.000E-7&1.000E-7&1.000E-7&1.000E-8&1.000E-8&1.000E-8\\
2&1.642E-8&6.397E-6&7.245E-5&3.852E-6&1.617E-5&8.378E-5&7.865E-7&9.141E-6&7.473E-5&2.161E-7&7.229E-6&7.288E-5\\
3&5.905E-7&9.964E-5&7.724E-4&2.867E-5&1.707E-4&8.306E-4&8.925E-6&1.222E-4&7.843E-4&3.538E-6&1.068E-4&7.747E-4\\
4&1.327E-5&7.878E-4&4.211E-3&1.987E-4&1.126E-3&4.414E-3&8.338E-5&9.004E-4&4.253E-3&4.310E-5&8.243E-4&4.219E-3\\
5&1.572E-4&4.430E-3&1.715E-2&1.192E-3&5.677E-3&1.770E-2&6.131E-4&4.857E-3&1.726E-2&3.732E-4&4.569E-3&1.717E-2\\
6&1.518E-3&2.012E-2&5.725E-2&6.554E-3&2.391E-2&5.848E-2&4.006E-3&2.145E-2&5.751E-2&2.794E-3&2.056E-2&5.730E-2\\
7&1.168E-2&7.818E-2&1.637E-1&3.299E-2&8.770E-2&1.659E-1&2.321E-2&8.156E-2&1.641E-1&1.796E-2&7.930E-2&1.638E-1\\
8&8.061E-2&2.607E-1&3.992E-1&1.523E-1&2.784E-1&4.021E-1&1.227E-1&2.671E-1&3.998E-1&1.049E-1&2.629E-1&3.993E-1\\
9&4.322E-1&6.732E-1&7.723E-1&5.553E-1&6.881E-1&7.740E-1&5.114E-1&6.787E-1&7.726E-1&4.808E-1&6.751E-1&7.723E-1\\
10&1.000E0&1.000E0&1.000E0&1.000E0&1.000E0&1.000E0&1.000E0&1.000E0&1.000E0&1.000E0&1.000E0&1.000E0\\ \hline  \hline 
\end{tabular}
%}
\end{table}
\end{tiny}

~ ~ ~

\begin{tiny}
\begin{table}
\setlength{\tabcolsep}{2.0pt}
\renewcommand{\arraystretch}{1.02}
\caption{The extreme points of ${\varepsilon}(q,\delta,\alpha,k;\t)$
 of $(q,\delta,\alpha,k)$-BURA approximation, 
i.e. the best uniform rational approximation of  
$\f(q,\delta,\alpha;\t)=\t^{\alpha}/(1+q\,\t^{\alpha})$  with functions from $\mathcal R_k$ on 
 $ [\delta, 1]$,   $\alpha =0.25, 0.5, 0.75$, ${\bf k=5}$, and $q=1$}
 %$(q=001,\delta,\alpha,k=5)$-BURA }
 \label{tabl:Bq001dDaAAk5p}
%\rotatebox{90}{
\begin{tabular}{||c||c|c|c||c|c|c||c|c|c||c|c|c||}
\hline \hline
{\bf k=5} /$\delta$& $0.00$ & $0.00$ & $0.00$ &$10^{-6}$&$10^{-6}$&$10^{-6}$&$10^{-7}$&$10^{-7}$&$10^{-7}$&$10^{-8}$&$10^{-8}$&$10^{-8}$\\
\ j /$\alpha$ & 0.25   & 0.50   & 0.75   & 0.25   & 0.50   & 0.75   & 0.25   & 0.50   & 0.75   & 0.25   & 0.50   & 0.75   \\ \hline
1&0.000E0&0.000E0&0.000E0&1.000E-6&1.000E-6&1.000E-6&1.000E-7&1.000E-7&1.000E-7&1.000E-8&1.000E-8&1.000E-8\\
2&1.060E-9&8.967E-7&1.424E-5&2.372E-6&5.717E-6&2.124E-5&3.725E-7&2.076E-6&1.569E-5&7.213E-8&1.233E-6&1.452E-5\\
3&3.775E-8&1.397E-5&1.522E-4&1.135E-5&4.360E-5&1.879E-4&2.686E-6&2.305E-5&1.600E-4&7.789E-7&1.684E-5&1.537E-4\\
4&8.273E-7&1.104E-4&8.348E-4&5.672E-5&2.450E-4&9.630E-4&1.800E-5&1.557E-4&8.637E-4&6.937E-6&1.253E-4&8.407E-4\\
5&9.531E-6&6.199E-4&3.443E-3&2.641E-4&1.114E-3&3.812E-3&1.042E-4&7.959E-4&3.528E-3&4.879E-5&6.794E-4&3.461E-3\\
6&8.766E-5&2.813E-3&1.178E-2&1.156E-3&4.360E-3&1.269E-2&5.481E-4&3.386E-3&1.199E-2&3.018E-4&3.010E-3&1.182E-2\\
7&6.370E-4&1.099E-2&3.519E-2&4.766E-3&1.524E-2&3.714E-2&2.636E-3&1.261E-2&3.564E-2&1.651E-3&1.155E-2&3.528E-2\\
8&4.120E-3&3.828E-2&9.395E-2&1.878E-2&4.855E-2&9.759E-2&1.194E-2&4.229E-2&9.480E-2&8.371E-3&3.969E-2&9.413E-2\\
9&2.367E-2&1.204E-1&2.253E-1&7.023E-2&1.415E-1&2.310E-1&5.072E-2&1.288E-1&2.267E-1&3.931E-2&1.234E-1&2.256E-1\\
10&1.242E-1&3.335E-1&4.744E-1&2.410E-1&3.654E-1&4.808E-1&1.980E-1&3.465E-1&4.759E-1&1.695E-1&3.382E-1&4.747E-1\\
11&5.139E-1&7.299E-1&8.140E-1&6.551E-1&7.515E-1&8.172E-1&6.114E-1&7.389E-1&8.148E-1&5.779E-1&7.332E-1&8.141E-1\\
12&1.000E0&1.000E0&1.000E0&1.000E0&1.000E0&1.000E0&1.000E0&1.000E0&1.000E0&1.000E0&1.000E0&1.000E0\\ \hline  \hline 
\end{tabular}
%}
\end{table}
\end{tiny}

~ ~ ~

\begin{tiny}
\begin{table}
\setlength{\tabcolsep}{2.0pt}
\renewcommand{\arraystretch}{1.02}
\caption{The extreme points of ${\varepsilon}(q,\delta,\alpha,k;\t)$
 of $(q,\delta,\alpha,k)$-BURA approximation, 
i.e. the best uniform rational approximation of  
$\f(q,\delta,\alpha;\t)=\t^{\alpha}/(1+q\,\t^{\alpha})$  with functions from $\mathcal R_k$ on 
 $ [\delta, 1]$,   $\alpha =0.25, 0.5, 0.75$, ${\bf k=6}$, and $q=1$}
 %$(q=001,\delta,\alpha,k=6)$-BURA }
\label{tabl:Bq001dDaAAk6p}
%\rotatebox{90}{
\begin{tabular}{||c||c|c|c||c|c|c||c|c|c||c|c|c||}
\hline \hline
{\bf k=6} /$\delta$& $0.00$ & $0.00$ & $0.00$ &$10^{-6}$&$10^{-6}$&$10^{-6}$&$10^{-7}$&$10^{-7}$&$10^{-7}$&$10^{-8}$&$10^{-8}$&$10^{-8}$\\
\ j /$\alpha$ & 0.25   & 0.50   & 0.75   & 0.25   & 0.50   & 0.75   & 0.25   & 0.50   & 0.75   & 0.25   & 0.50   & 0.75   \\ \hline
1&0.000E0&0.000E0&0.000E0&1.000E-6&1.000E-6&1.000E-6&1.000E-7&1.000E-7&1.000E-7&1.000E-8&1.000E-8&1.000E-8\\
2&8.636E-11&1.487E-7&3.227E-6&1.807E-6&3.078E-6&7.729E-6&2.433E-7&7.452E-7&4.164E-6&3.764E-8&3.006E-7&3.419E-6\\
3&3.061E-9&2.316E-6&3.452E-5&6.142E-6&1.650E-5&5.602E-5&1.192E-6&6.319E-6&3.957E-5&2.726E-7&3.550E-6&3.559E-5\\
4&6.631E-8&1.830E-5&1.897E-4&2.360E-5&7.678E-5&2.663E-4&6.037E-6&3.745E-5&2.086E-4&1.823E-6&2.469E-5&1.938E-4\\
5&7.541E-7&1.028E-4&7.859E-4&8.907E-5&3.078E-4&1.010E-3&2.823E-5&1.762E-4&8.426E-4&1.047E-5&1.285E-4&7.982E-4\\
6&6.784E-6&4.661E-4&2.710E-3&3.235E-4&1.098E-3&3.283E-3&1.232E-4&7.070E-4&2.857E-3&5.427E-5&5.532E-4&2.742E-3\\
7&4.795E-5&1.820E-3&8.221E-3&1.128E-3&3.569E-3&9.526E-3&5.033E-4&2.519E-3&8.561E-3&2.557E-4&2.079E-3&8.295E-3\\
8&2.978E-4&6.352E-3&2.255E-2&3.804E-3&1.076E-2&2.524E-2&1.956E-3&8.178E-3&2.326E-2&1.127E-3&7.041E-3&2.270E-2\\
9&1.634E-3&2.027E-2&5.679E-2&1.242E-2&3.041E-2&6.180E-2&7.265E-3&2.460E-2&5.812E-2&4.670E-3&2.193E-2&5.708E-2\\
10&8.305E-3&5.987E-2&1.320E-1&3.936E-2&8.086E-2&1.403E-1&2.599E-2&6.909E-2&1.343E-1&1.849E-2&6.346E-2&1.325E-1\\
11&3.909E-2&1.632E-1&2.816E-1&1.193E-1&2.003E-1&2.930E-1&8.867E-2&1.799E-1&2.846E-1&6.945E-2&1.698E-1&2.822E-1\\
12&1.690E-1&3.959E-1&5.345E-1&3.299E-1&4.438E-1&5.459E-1&2.770E-1&4.181E-1&5.376E-1&2.394E-1&4.048E-1&5.352E-1\\
13&5.771E-1&7.707E-1&8.433E-1&7.277E-1&7.982E-1&8.486E-1&6.869E-1&7.838E-1&8.447E-1&6.536E-1&7.760E-1&8.436E-1\\
14&1.000E0&1.000E0&1.000E0&1.000E0&1.000E0&1.000E0&1.000E0&1.000E0&1.000E0&1.000E0&1.000E0&1.000E0\\ \hline \hline
\end{tabular}
%}
\end{table}
\end{tiny}

~ ~ ~

\begin{tiny}
\begin{table}
\setlength{\tabcolsep}{2.0pt}
\renewcommand{\arraystretch}{1.02}
\caption{The extreme points of ${\varepsilon}(q,\delta,\alpha,k;\t)$
 of $(q,\delta,\alpha,k)$-BURA approximation, 
i.e. the best uniform rational approximation of  
$\f(q,\delta,\alpha;\t)=\t^{\alpha}/(1+q\,\t^{\alpha})$  with functions from $\mathcal R_k$ on 
 $ [\delta, 1]$,   $\alpha =0.25, 0.5, 0.75$, ${\bf k=7}$, and $q=1$}
 %$(q=001,\delta,\alpha,k=7)$-BURA }
\label{tabl:Bq001dDaAAk7p}
%\rotatebox{90}{
\begin{tabular}{||c||c|c|c||c|c|c||c|c|c||c|c|c||}
\hline \hline
{\bf k=7} /$\delta$& $0.00$ & $0.00$ & $0.00$ &$10^{-6}$&$10^{-6}$&$10^{-6}$&$10^{-7}$&$10^{-7}$&$10^{-7}$&$10^{-8}$&$10^{-8}$&$10^{-8}$\\
\ j /$\alpha$ & 0.25   & 0.50   & 0.75   & 0.25   & 0.50   & 0.75   & 0.25   & 0.50   & 0.75   & 0.25   & 0.50   & 0.75   \\ \hline
1&0.000E0&0.000E0&0.000E0&1.000E-6&1.000E-6&1.000E-6&1.000E-7&1.000E-7&1.000E-7&1.000E-8&1.000E-8&1.000E-8\\
2&8.427E-12&2.812E-8&8.172E-7&1.533E-6&2.136E-6&3.909E-6&1.884E-7&3.856E-7&1.433E-6&2.539E-8&1.060E-7&9.458E-7\\
3&2.980E-10&4.380E-7&8.743E-6&4.034E-6&8.309E-6&2.183E-5&6.776E-7&2.401E-6&1.194E-5&1.299E-7&1.006E-6&9.465E-6\\
4&6.417E-9&3.461E-6&4.809E-5&1.252E-5&3.197E-5&9.298E-5&2.715E-6&1.213E-5&6.000E-5&6.764E-7&6.307E-6&5.086E-5\\
5&7.252E-8&1.943E-5&1.995E-4&3.969E-5&1.123E-4&3.297E-4&1.056E-5&5.141E-5&2.356E-4&3.222E-6&3.076E-5&2.080E-4\\
6&6.454E-7&8.812E-5&6.899E-4&1.236E-4&3.620E-4&1.024E-3&3.923E-5&1.910E-4&7.857E-4&1.423E-5&1.265E-4&7.128E-4\\
\hline  \hline 
\end{tabular}
%}
\end{table}
\end{tiny}

~ ~ ~

\begin{tiny}
\begin{table}
\setlength{\tabcolsep}{2.0pt}
\renewcommand{\arraystretch}{1.02}
\caption{The extreme points of ${\varepsilon}(q,\delta,\alpha,k;\t)$
 of $(q,\delta,\alpha,k)$-BURA approximation, 
i.e. the best uniform rational approximation of  
$\f(q,\delta,\alpha;\t)=\t^{\alpha}/(1+q\,\t^{\alpha})$  with functions from $\mathcal R_k$ on 
 $ [\delta, 1]$,   $\alpha =0.25, 0.5, 0.75$, ${\bf k=8}$, and $q=1$}
 %$(q=001,\delta,\alpha,k=8)$-BURA }
\label{tabl:Bq001dDaAAk8p}
%\rotatebox{90}{
\begin{tabular}{||c||c|c|c||c|c|c||c|c|c||c|c|c||}
\hline \hline
{\bf k=8} /$\delta$& $0.00$ & $0.00$ & $0.00$ &$10^{-6}$&$10^{-6}$&$10^{-6}$&$10^{-7}$&$10^{-7}$&$10^{-7}$&$10^{-8}$&$10^{-8}$&$10^{-8}$\\
\ j /$\alpha$ & 0.25   & 0.50   & 0.75   & 0.25   & 0.50   & 0.75   & 0.25   & 0.50   & 0.75   & 0.25   & 0.50   & 0.75   \\ \hline
1&0.000E0&0.000E0&0.000E0&1.000E-6&1.000E-6&1.000E-6&1.000E-7&1.000E-7&1.000E-7&1.000E-8&1.000E-8&1.000E-8\\
2&9.514E-13&5.913E-9&2.263E-7&1.379E-6&1.711E-6&2.529E-6&1.601E-7&2.551E-7&6.464E-7&1.980E-8&5.166E-8&3.132E-7\\
3&3.360E-11&9.209E-8&2.421E-6&2.998E-6&5.106E-6&1.072E-5&4.525E-7&1.177E-6&4.437E-6&7.603E-8&3.761E-7&2.901E-6\\
4&7.213E-10&7.277E-7&1.332E-5&7.825E-6&1.645E-5&4.001E-5&1.491E-6&5.022E-6&2.065E-5&3.198E-7&2.063E-6&1.516E-5\\
5&8.123E-9&4.086E-6&5.528E-5&2.142E-5&5.082E-5&1.304E-4&4.952E-6&1.897E-5&7.744E-5&1.290E-6&9.234E-6&6.101E-5\\
6&7.188E-8&1.853E-5&1.914E-4&5.866E-5&1.482E-4&3.819E-4&1.603E-5&6.467E-5&2.503E-4&4.935E-6&3.571E-5&2.069E-4\\
\hline  \hline 
\end{tabular}
%}
\end{table}
\end{tiny}

~ ~ ~

\begin{tiny}
\begin{table}
\setlength{\tabcolsep}{2.0pt}
\renewcommand{\arraystretch}{1.02}
\caption{The extreme points of ${\varepsilon}(q,\delta,\alpha,k;\t)$
 of $(q,\delta,\alpha,k)$-BURA approximation, 
i.e. the best uniform rational approximation of  
$\f(q,\delta,\alpha;\t)=\t^{\alpha}/(1+q\,\t^{\alpha})$  with functions from $\mathcal R_k$ on 
 $ [\delta, 1]$,   $\alpha =0.25, 0.5, 0.75$, ${\bf k=3}$, and $q=100$}
 %$(q=100,\delta,\alpha,k=3)$-BURA }
\label{tabl:Bq100dDaAAk3p}
%\rotatebox{90}{
\begin{tabular}{||c||c|c|c||c|c|c||c|c|c||c|c|c||}
\hline \hline
{\bf k=3}/$\delta$& $0.00$ & $0.00$ & $0.00$ &$10^{-6}$&$10^{-6}$&$10^{-6}$&$10^{-7}$&$10^{-7}$&$10^{-7}$&$10^{-8}$&$10^{-8}$&$10^{-8}$\\
\ j /$\alpha$ & 0.25   & 0.50   & 0.75   & 0.25   & 0.50   & 0.75   & 0.25   & 0.50   & 0.75   & 0.25   & 0.50   & 0.75   \\ \hline
1&0.000E0&0.000E0&0.000E0&1.000E-6&1.000E-6&1.000E-6&1.000E-7&1.000E-7&1.000E-7&1.000E-8&1.000E-8&1.000E-8\\
2&1.136E-12&2.045E-7&1.336E-5&2.013E-6&3.411E-6&1.999E-5&2.314E-7&8.732E-7&1.473E-5&2.793E-8&3.782E-7&1.363E-5\\
3&4.845E-11&3.280E-6&1.381E-4&9.434E-6&2.082E-5&1.703E-4&1.307E-6&8.162E-6&1.451E-4&1.933E-7&4.776E-6&1.395E-4\\
4&1.779E-9&2.874E-5&7.446E-4&6.248E-5&1.191E-4&8.627E-4&1.028E-5&5.675E-5&7.710E-4&1.789E-6&3.783E-5&7.499E-4\\
5&3.721E-8&2.112E-4&3.403E-3&4.963E-4&6.988E-4&3.824E-3&9.605E-5&3.706E-4&3.498E-3&1.909E-5&2.643E-4&3.422E-3\\
6&1.281E-6&1.721E-3&1.657E-2&5.355E-3&4.942E-3&1.833E-2&1.324E-3&2.818E-3&1.696E-2&3.195E-4&2.094E-3&1.665E-2\\
7&4.953E-5&2.233E-2&1.197E-1&7.502E-2&5.447E-2&1.296E-1&2.527E-2&3.408E-2&1.220E-1&7.590E-3&2.644E-2&1.201E-1\\
8&1.000E0&1.000E0&1.000E0&1.000E0&1.000E0&1.000E0&1.000E0&1.000E0&1.000E0&1.000E0&1.000E0&1.000E0\\ \hline \hline
\end{tabular}
%}
\end{table}
\end{tiny}

~ ~ ~

\begin{tiny}
\begin{table}
\setlength{\tabcolsep}{2.0pt}
\renewcommand{\arraystretch}{1.02}
\caption{The extreme points of ${\varepsilon}(q,\delta,\alpha,k;\t)$
 of $(q,\delta,\alpha,k)$-BURA approximation, 
i.e. the best uniform rational approximation of  
$\f(q,\delta,\alpha;\t)=\t^{\alpha}/(1+q\,\t^{\alpha})$  with functions from $\mathcal R_k$ on 
 $ [\delta, 1]$,   $\alpha =0.25, 0.5, 0.75$, ${\bf k=4}$, and $q=100$}
 %$(q=100,\delta,\alpha,k=4)$-BURA }
 \label{tabl:Bq100dDaAAk4p}
%\rotatebox{90}{
\begin{tabular}{||c||c|c|c||c|c|c||c|c|c||c|c|c||}
\hline \hline
{\bf k=4} /$\delta$& $0.00$ & $0.00$ & $0.00$ &$10^{-6}$&$10^{-6}$&$10^{-6}$&$10^{-7}$&$10^{-7}$&$10^{-7}$&$10^{-8}$&$10^{-8}$&$10^{-8}$\\
\ j /$\alpha$ & 0.25   & 0.50   & 0.75   & 0.25   & 0.50   & 0.75   & 0.25   & 0.50   & 0.75   & 0.25   & 0.50   & 0.75   \\ \hline
1&0.000E0&0.000E0&0.000E0&1.000E-6&1.000E-6&1.000E-6&1.000E-7&1.000E-7&1.000E-7&1.000E-8&1.000E-8&1.000E-8\\
2&1.093E-13&3.667E-8&3.029E-6&1.629E-6&2.260E-6&7.384E-6&1.811E-7&4.220E-7&3.935E-6&2.080E-8&1.227E-7&3.215E-6\\
3&4.275E-12&5.752E-7&3.180E-5&5.177E-6&9.531E-6&5.191E-5&6.855E-7&2.811E-6&3.655E-5&9.535E-8&1.226E-6&3.281E-5\\
4&1.207E-10&4.676E-6&1.700E-4&2.203E-5&4.080E-5&2.397E-4&3.424E-6&1.534E-5&1.872E-4&5.567E-7&8.128E-6&1.737E-4\\
5&1.858E-9&2.848E-5&6.970E-4&1.071E-4&1.707E-4&9.056E-4&1.932E-5&7.445E-5&7.495E-4&3.579E-6&4.420E-5&7.084E-4\\
6&2.981E-8&1.559E-4&2.541E-3&5.870E-4&7.374E-4&3.154E-3&1.253E-4&3.548E-4&2.697E-3&2.679E-5&2.262E-4&2.575E-3\\
7&4.421E-7&8.860E-4&9.434E-3&3.609E-3&3.530E-3&1.138E-2&9.320E-4&1.836E-3&9.933E-3&2.320E-4&1.230E-3&9.543E-3\\
8&1.096E-5&6.101E-3&4.033E-2&2.655E-2&2.047E-2&4.745E-2&9.013E-3&1.160E-2&4.218E-2&2.826E-3&8.158E-3&4.073E-2\\
9&3.599E-4&6.466E-2&2.318E-1&2.199E-1&1.620E-1&2.589E-1&1.089E-1&1.070E-1&2.391E-1&4.626E-2&8.150E-2&2.334E-1\\
10&1.000E0&1.000E0&1.000E0&1.000E0&1.000E0&1.000E0&1.000E0&1.000E0&1.000E0&1.000E0&1.000E0&1.000E0\\ \hline \hline
\end{tabular}
%}
\end{table}
\end{tiny}

~ ~ ~

\begin{tiny}
\begin{table}
\setlength{\tabcolsep}{2.0pt}
\renewcommand{\arraystretch}{1.02}
\caption{The extreme points of ${\varepsilon}(q,\delta,\alpha,k;\t)$
 of $(q,\delta,\alpha,k)$-BURA approximation, 
i.e. the best uniform rational approximation of  
$\f(q,\delta,\alpha;\t)=\t^{\alpha}/(1+q\,\t^{\alpha})$  with functions from $\mathcal R_k$ on 
 $ [\delta, 1]$,   $\alpha =0.25, 0.5, 0.75$, ${\bf k=5}$, and $q=100$}
 %$(q=100,\delta,\alpha,k=5)$-BURA }
 \label{tabl:Bq100dDaAAk5p}
%\rotatebox{90}{
\begin{tabular}{||c||c|c|c||c|c|c||c|c|c||c|c|c||}
\hline \hline
{\bf k=5} /$\delta$& $0.00$ & $0.00$ & $0.00$ &$10^{-6}$&$10^{-6}$&$10^{-6}$&$10^{-7}$&$10^{-7}$&$10^{-7}$&$10^{-8}$&$10^{-8}$&$10^{-8}$\\
\ j /$\alpha$ & 0.25   & 0.50   & 0.75   & 0.25   & 0.50   & 0.75   & 0.25   & 0.50   & 0.75   & 0.25   & 0.50   & 0.75   \\ \hline
1&0.000E0&0.000E0&0.000E0&1.000E-6&1.000E-6&1.000E-6&1.000E-7&1.000E-7&1.000E-7&1.000E-8&1.000E-8&1.000E-8\\
2&1.357E-14&7.520E-9&7.680E-7&1.432E-6&1.770E-6&3.790E-6&1.556E-7&2.695E-7&1.368E-6&1.730E-8&5.683E-8&8.935E-7\\
3&5.084E-13&1.173E-7&8.147E-6&3.505E-6&5.584E-6&2.065E-5&4.454E-7&1.311E-6&1.121E-5&5.898E-8&4.340E-7&8.841E-6\\
4&1.273E-11&9.346E-7&4.407E-5&1.118E-5&1.922E-5&8.592E-5&1.653E-6&5.873E-6&5.521E-5&2.544E-7&2.468E-6&4.667E-5\\
5&1.706E-10&5.381E-6&1.791E-4&4.021E-5&6.536E-5&2.983E-4&6.829E-6&2.377E-5&2.122E-4&1.195E-6&1.163E-5&1.869E-4\\
6&2.066E-9&2.600E-5&6.127E-4&1.560E-4&2.215E-4&9.241E-4&3.053E-5&9.108E-5&7.015E-4&6.058E-6&4.936E-5&6.339E-4\\
7&2.138E-8&1.164E-4&1.928E-3&6.452E-4&7.715E-4&2.727E-3&1.468E-4&3.478E-4&2.159E-3&3.313E-5&2.026E-4&1.984E-3\\
8&2.481E-7&5.254E-4&6.013E-3&2.878E-3&2.862E-3&8.151E-3&7.803E-4&1.397E-3&6.639E-3&2.050E-4&8.601E-4&6.164E-3\\
9&3.066E-6&2.607E-3&2.007E-2&1.389E-2&1.177E-2&2.624E-2&4.613E-3&6.244E-3&2.190E-2&1.442E-3&4.052E-3&2.052E-2\\
10&6.204E-5&1.576E-2&7.685E-2&7.355E-2&5.606E-2&9.593E-2&3.211E-2&3.323E-2&8.265E-2&1.284E-2&2.304E-2&7.827E-2\\
11&1.763E-3&1.343E-1&3.513E-1&3.850E-1&3.107E-1&3.995E-1&2.461E-1&2.249E-1&3.667E-1&1.386E-1&1.760E-1&3.551E-1\\
12&1.000E0&1.000E0&1.000E0&1.000E0&1.000E0&1.000E0&1.000E0&1.000E0&1.000E0&1.000E0&1.000E0&1.000E0\\ \hline \hline
\end{tabular}
%}
\end{table}
\end{tiny}

~ ~ ~

\begin{tiny}
\begin{table}
\setlength{\tabcolsep}{2.0pt}
\renewcommand{\arraystretch}{1.02}
\caption{The extreme points of ${\varepsilon}(q,\delta,\alpha,k;\t)$
 of $(q,\delta,\alpha,k)$-BURA approximation, 
i.e. the best uniform rational approximation of  
$\f(q,\delta,\alpha;\t)=\t^{\alpha}/(1+q\,\t^{\alpha})$  with functions from $\mathcal R_k$ on 
 $ [\delta, 1]$,   $\alpha =0.25, 0.5, 0.75$, ${\bf k=6}$, and $q=100$}
 %$(q=100,\delta,\alpha,k=6)$-BURA }
 \label{tabl:Bq100dDaAAk6p}
%\rotatebox{90}{
\begin{tabular}{||c||c|c|c||c|c|c||c|c|c||c|c|c||}
\hline \hline
{\bf k=6} /$\delta$& $0.00$ & $0.00$ & $0.00$ &$10^{-6}$&$10^{-6}$&$10^{-6}$&$10^{-7}$&$10^{-7}$&$10^{-7}$&$10^{-8}$&$10^{-8}$&$10^{-8}$\\
\ j /$\alpha$ & 0.25   & 0.50   & 0.75   & 0.25   & 0.50   & 0.75   & 0.25   & 0.50   & 0.75   & 0.25   & 0.50   & 0.75   \\ \hline
1&0.000E0&0.000E0&0.000E0&1.000E-6&1.000E-6&1.000E-6&1.000E-7&1.000E-7&1.000E-7&1.000E-8&1.000E-8&1.000E-8\\
2&1.984E-15&1.700E-9&2.129E-7&1.316E-6&1.519E-6&2.479E-6&1.407E-7&2.032E-7&6.245E-7&1.530E-8&3.414E-8&2.979E-7\\
3&7.261E-14&2.649E-8&2.269E-6&2.681E-6&3.825E-6&1.030E-5&3.287E-7&7.514E-7&4.221E-6&4.177E-8&1.970E-7&2.735E-6\\
4&1.708E-12&2.098E-7&1.238E-5&6.953E-6&1.101E-5&3.775E-5&9.818E-7&2.822E-6&1.938E-5&1.440E-7&9.601E-7&1.415E-5\\
5&2.129E-11&1.187E-6&5.069E-5&2.025E-5&3.206E-5&1.208E-4&3.259E-6&9.905E-6&7.145E-5&5.409E-7&4.016E-6&5.609E-5\\
6&2.254E-10&5.499E-6&1.726E-4&6.285E-5&9.274E-5&3.479E-4&1.148E-5&3.290E-5&2.269E-4&2.141E-6&1.513E-5&1.870E-4\\
7&1.962E-9&2.256E-5&5.219E-4&2.040E-4&2.692E-4&9.355E-4&4.241E-5&1.062E-4&6.539E-4&8.876E-6&5.365E-5&5.573E-4\\
8&1.690E-8&8.709E-5&1.473E-3&6.910E-4&7.985E-4&2.433E-3&1.651E-4&3.435E-4&1.785E-3&3.911E-5&1.864E-4&1.558E-3\\
9&1.419E-7&3.343E-4&4.083E-3&2.446E-3&2.470E-3&6.356E-3&6.806E-4&1.149E-3&4.832E-3&1.843E-4&6.608E-4&4.287E-3\\
10&1.374E-6&1.348E-3&1.166E-2&9.106E-3&8.129E-3&1.729E-2&3.017E-3&4.096E-3&1.354E-2&9.555E-4&2.490E-3&1.217E-2\\
11&1.499E-5&6.049E-3&3.592E-2&3.560E-2&2.905E-2&5.059E-2&1.444E-2&1.610E-2&4.094E-2&5.487E-3&1.042E-2&3.731E-2\\
12&2.575E-4&3.236E-2&1.236E-1&1.439E-1&1.139E-1&1.615E-1&7.569E-2&7.206E-2&1.370E-1&3.669E-2&5.075E-2&1.273E-1\\
13&6.286E-3&2.210E-1&4.584E-1&5.251E-1&4.558E-1&5.243E-1&3.906E-1&3.585E-1&4.834E-1&2.657E-1&2.929E-1&4.656E-1\\
14&1.000E0&1.000E0&1.000E0&1.000E0&1.000E0&1.000E0&1.000E0&1.000E0&1.000E0&1.000E0&1.000E0&1.000E0\\  \hline \hline
\end{tabular}
%}
\end{table}
\end{tiny}

~ ~ ~

\begin{tiny}
\begin{table}
\setlength{\tabcolsep}{2.0pt}
\renewcommand{\arraystretch}{1.02}
\caption{The extreme points of ${\varepsilon}(q,\delta,\alpha,k;\t)$
 of $(q,\delta,\alpha,k)$-BURA approximation, 
i.e. the best uniform rational approximation of  
$\f(q,\delta,\alpha;\t)=\t^{\alpha}/(1+q\,\t^{\alpha})$  with functions from $\mathcal R_k$ on 
 $ [\delta, 1]$,   $\alpha =0.25, 0.5, 0.75$, ${\bf k=7}$, and $q=100$}
 %$(q=100,\delta,\alpha,k=7)$-BURA }
 \label{tabl:Bq100dDaAAk7p}
%\rotatebox{90}{
\begin{tabular}{||c||c|c|c||c|c|c||c|c|c||c|c|c||}
\hline \hline
{\bf k=7} /$\delta$& $0.00$ & $0.00$ & $0.00$ &$10^{-6}$&$10^{-6}$&$10^{-6}$&$10^{-7}$&$10^{-7}$&$10^{-7}$&$10^{-8}$&$10^{-8}$&$10^{-8}$\\
\ j /$\alpha$ & 0.25   & 0.50   & 0.75   & 0.25   & 0.50   & 0.75   & 0.25   & 0.50   & 0.75   & 0.25   & 0.50   & 0.75   \\ \hline
1&0.000E0&0.000E0&0.000E0&1.000E-6&1.000E-6&1.000E-6&1.000E-7&1.000E-7&1.000E-7&1.000E-8&1.000E-8&1.000E-8\\
2&3.253E-16&4.158E-10&6.353E-8&1.243E-6&1.373E-6&1.896E-6&1.312E-7&1.692E-7&3.618E-7&1.403E-8&2.440E-8&1.221E-7\\
3&1.174E-14&6.476E-9&6.786E-7&2.213E-6&2.908E-6&6.137E-6&2.634E-7&4.981E-7&1.936E-6&3.232E-8&1.082E-7&9.858E-7\\
4&2.664E-13&5.122E-8&3.719E-6&4.909E-6&7.207E-6&1.963E-5&6.659E-7&1.595E-6&8.059E-6&9.354E-8&4.523E-7&4.877E-6\\
5&3.191E-12&2.883E-7&1.533E-5&1.221E-5&1.849E-5&5.751E-5&1.872E-6&4.951E-6&2.803E-5&2.963E-7&1.692E-6&1.889E-5\\
6&3.141E-11&1.316E-6&5.249E-5&3.224E-5&4.746E-5&1.550E-4&5.559E-6&1.471E-5&8.551E-5&9.815E-7&5.787E-6&6.206E-5\\
\hline  \hline 
\end{tabular}
%}
\end{table}
\end{tiny}

~ ~ ~

\begin{tiny}
\begin{table}
\setlength{\tabcolsep}{2.0pt}
\renewcommand{\arraystretch}{1.02}
\caption{The extreme points of ${\varepsilon}(q,\delta,\alpha,k;\t)$
 of $(q,\delta,\alpha,k)$-BURA approximation, 
i.e. the best uniform rational approximation of  
$\f(q,\delta,\alpha;\t)=\t^{\alpha}/(1+q\,\t^{\alpha})$  with functions from $\mathcal R_k$ on 
 $ [\delta, 1]$,   $\alpha =0.25, 0.5, 0.75$, ${\bf k=8}$, and $q=100$}
 %$(q=100,\delta,\alpha,k=8)$-BURA }
 \label{tabl:Bq100dDaAAk8p}
%\rotatebox{90}{
\begin{tabular}{||c||c|c|c||c|c|c||c|c|c||c|c|c||}
\hline \hline
{\bf k=8} /$\delta$& $0.00$ & $0.00$ & $0.00$ &$10^{-6}$&$10^{-6}$&$10^{-6}$&$10^{-7}$&$10^{-7}$&$10^{-7}$&$10^{-8}$&$10^{-8}$&$10^{-8}$\\
\ j /$\alpha$ & 0.25   & 0.50   & 0.75   & 0.25   & 0.50   & 0.75   & 0.25   & 0.50   & 0.75   & 0.25   & 0.50   & 0.75   \\ \hline
1&&0.000E0&0.000E0&1.000E-6&1.000E-6&1.000E-6&1.000E-7&1.000E-7&1.000E-7&1.000E-8&1.000E-8&1.000E-8\\
2&&1.087E-10&2.018E-8&1.192E-6&1.282E-6&1.596E-6&1.247E-7&1.496E-7&2.507E-7&1.318E-8&1.949E-8&6.184E-8\\
3&&1.692E-9&2.158E-7&1.921E-6&2.374E-6&4.180E-6&2.231E-7&3.661E-7&1.054E-6&2.660E-8&6.855E-8&4.179E-7\\
4&&1.338E-8&1.185E-6&3.773E-6&5.191E-6&1.169E-5&4.941E-7&1.018E-6&3.912E-6&6.684E-8&2.474E-7&1.930E-6\\
5&&7.516E-8&4.902E-6&8.284E-6&1.196E-5&3.128E-5&1.218E-6&2.830E-6&1.267E-5&1.847E-7&8.323E-7&7.181E-6\\
6&&3.416E-7&1.687E-5&1.932E-5&2.785E-5&7.882E-5&3.170E-6&7.646E-6&3.677E-5&5.334E-7&2.607E-6&2.301E-5\\
\hline  \hline 
\end{tabular}
%}
\end{table}
\end{tiny}

~ ~ ~

\begin{tiny}
\begin{table}
\setlength{\tabcolsep}{2.0pt}
\renewcommand{\arraystretch}{1.02}
\caption{The extreme points of ${\varepsilon}(q,\delta,\alpha,k;\t)$
 of $(q,\delta,\alpha,k)$-BURA approximation, 
i.e. the best uniform rational approximation of  
$\f(q,\delta,\alpha;\t)=\t^{\alpha}/(1+q\,\t^{\alpha})$  with functions from $\mathcal R_k$ on 
 $ [\delta, 1]$,   $\alpha =0.25, 0.5, 0.75$, ${\bf k=3}$, and $q=200$}
 %$(q=200,\delta,\alpha,k=3)$-BURA }
 \label{tabl:Bq200dDaAAk3p}
%\rotatebox{90}{
\begin{tabular}{||c||c|c|c||c|c|c||c|c|c||c|c|c||}
\hline \hline
{\bf k=3} /$\delta$& $0.00$ & $0.00$ & $0.00$ &$10^{-6}$&$10^{-6}$&$10^{-6}$&$10^{-7}$&$10^{-7}$&$10^{-7}$&$10^{-8}$&$10^{-8}$&$10^{-8}$\\
\ j /$\alpha$ & 0.25   & 0.50   & 0.75   & 0.25   & 0.50   & 0.75   & 0.25   & 0.50   & 0.75   & 0.25   & 0.50   & 0.75   \\ \hline
1&0.000E0&0.000E0&0.000E0&1.000E-6&1.000E-6&1.000E-6&1.000E-7&1.000E-7&1.000E-7&1.000E-8&1.000E-8&1.000E-8\\
2&7.673E-14&5.630E-8&5.903E-6&1.922E-6&2.490E-6&1.106E-5&2.115E-7&4.920E-7&6.978E-6&2.376E-8&1.572E-7&6.120E-6\\
3&3.286E-12&9.051E-7&6.095E-5&8.531E-6&1.212E-5&8.523E-5&1.086E-6&3.697E-6&6.648E-5&1.413E-7&1.720E-6&6.209E-5\\
4&1.222E-10&7.992E-6&3.299E-4&5.497E-5&6.304E-5&4.196E-4&8.108E-6&2.357E-5&3.509E-4&1.199E-6&1.292E-5&3.342E-4\\
5&2.594E-9&5.980E-5&1.534E-3&4.338E-4&3.601E-4&1.866E-3&7.449E-5&1.501E-4&1.612E-3&1.236E-5&8.935E-5&1.550E-3\\
6&9.325E-8&5.073E-4&7.786E-3&4.723E-3&2.603E-3&9.290E-3&1.033E-3&1.166E-3&8.143E-3&2.060E-4&7.277E-4&7.861E-3\\
7&3.798E-6&7.250E-3&6.395E-2&6.800E-2&3.132E-2&7.437E-2&2.038E-2&1.541E-2&6.646E-2&5.061E-3&1.007E-2&6.448E-2\\
8&1.000E0&1.000E0&1.000E0&1.000E0&1.000E0&1.000E0&1.000E0&1.000E0&1.000E0&1.000E0&1.000E0&1.000E0\\ \hline  \hline 
\end{tabular}
%}
\end{table}
\end{tiny}

~ ~ ~

\begin{tiny}
\begin{table}
\setlength{\tabcolsep}{2.0pt}
\renewcommand{\arraystretch}{1.02}
\caption{The extreme points of ${\varepsilon}(q,\delta,\alpha,k;\t)$
 of $(q,\delta,\alpha,k)$-BURA approximation, 
i.e. the best uniform rational approximation of  
$\f(q,\delta,\alpha;\t)=\t^{\alpha}/(1+q\,\t^{\alpha})$  with functions from $\mathcal R_k$ on 
 $ [\delta, 1]$,   $\alpha =0.25, 0.5, 0.75$, ${\bf k=4}$, and $q=200$}
 %$(q=200,\delta,\alpha,k=4)$-BURA }
 \label{tabl:Bq200dDaAAk4p}
%\rotatebox{90}{
\begin{tabular}{||c||c|c|c||c|c|c||c|c|c||c|c|c||}
\hline \hline
{\bf k=4 } /$\delta$& $0.00$ & $0.00$ & $0.00$ &$10^{-6}$&$10^{-6}$&$10^{-6}$&$10^{-7}$&$10^{-7}$&$10^{-7}$&$10^{-8}$&$10^{-8}$&$10^{-8}$\\
\ j /$\alpha$ & 0.25   & 0.50   & 0.75   & 0.25   & 0.50   & 0.75   & 0.25   & 0.50   & 0.75   & 0.25   & 0.50   & 0.75   \\ \hline
1&0.000E0&0.000E0&0.000E0&1.000E-6&1.000E-6&1.000E-6&1.000E-7&1.000E-7&1.000E-7&1.000E-8&1.000E-8&1.000E-8\\
2&7.638E-15&1.067E-8&1.408E-6&1.583E-6&1.869E-6&4.909E-6&1.712E-7&2.934E-7&2.123E-6&1.875E-8&6.576E-8&1.557E-6\\
3&2.997E-13&1.675E-7&1.476E-5&4.829E-6&6.463E-6&2.991E-5&6.006E-7&1.551E-6&1.843E-5&7.581E-8&5.396E-7&1.557E-5\\
4&8.537E-12&1.367E-6&7.883E-5&2.004E-5&2.489E-5&1.305E-4&2.866E-6&7.575E-6&9.211E-5&4.102E-7&3.277E-6&8.181E-5\\
5&1.327E-10&8.402E-6&3.247E-4&9.640E-5&9.927E-5&4.808E-4&1.584E-5&3.485E-5&3.659E-4&2.539E-6&1.712E-5&3.340E-4\\
6&2.176E-9&4.693E-5&1.202E-3&5.279E-4&4.250E-4&1.678E-3&1.022E-4&1.640E-4&1.329E-3&1.870E-5&8.699E-5&1.231E-3\\
7&3.319E-8&2.769E-4&4.627E-3&3.269E-3&2.081E-3&6.240E-3&7.655E-4&8.666E-4&5.062E-3&1.628E-4&4.846E-4&4.726E-3\\
8&8.824E-7&2.046E-3&2.121E-2&2.446E-2&1.277E-2&2.779E-2&7.595E-3&5.816E-3&2.301E-2&2.047E-3&3.421E-3&2.163E-2\\
9&3.160E-5&2.539E-2&1.455E-1&2.088E-1&1.146E-1&1.783E-1&9.618E-2&6.171E-2&1.548E-1&3.562E-2&3.961E-2&1.476E-1\\
10&1.000E0&1.000E0&1.000E0&1.000E0&1.000E0&1.000E0&1.000E0&1.000E0&1.000E0&1.000E0&1.000E0&1.000E0\\ \hline \hline
\end{tabular}
%}
\end{table}
\end{tiny}

~ ~ ~

\begin{tiny}
\begin{table}[h!t]
\setlength{\tabcolsep}{2.0pt}
\renewcommand{\arraystretch}{1.02}
\caption{The extreme points of ${\varepsilon}(q,\delta,\alpha,k;\t)$
 of $(q,\delta,\alpha,k)$-BURA approximation, 
i.e. the best uniform rational approximation of  
$\f(q,\delta,\alpha;\t)=\t^{\alpha}/(1+q\,\t^{\alpha})$  with functions from $\mathcal R_k$ on 
 $ [\delta, 1]$,   $\alpha =0.25, 0.5, 0.75$, ${\bf k=5}$, and $q=200$}
 %$(q=200,\delta,\alpha,k=5)$-BURA }
 \label{tabl:Bq200dDaAAk5p}
%\rotatebox{90}{
\begin{tabular}{||c||c|c|c||c|c|c||c|c|c||c|c|c||}
\hline \hline
{\bf k=5} /$\delta$& $0.00$ & $0.00$ & $0.00$ &$10^{-6}$&$10^{-6}$&$10^{-6}$&$10^{-7}$&$10^{-7}$&$10^{-7}$&$10^{-8}$&$10^{-8}$&$10^{-8}$\\
\ j /$\alpha$ & 0.25   & 0.50   & 0.75   & 0.25   & 0.50   & 0.75   & 0.25   & 0.50   & 0.75   & 0.25   & 0.50   & 0.75   \\ \hline
1&0.000E0&0.000E0&0.000E0&1.000E-6&1.000E-6&1.000E-6&1.000E-7&1.000E-7&1.000E-7&1.000E-8&1.000E-8&1.000E-8\\
2&9.914E-16&2.321E-9&3.740E-7&1.406E-6&1.571E-6&2.915E-6&1.499E-7&2.141E-7&8.567E-7&1.614E-8&3.740E-8&4.747E-7\\
3&3.723E-14&3.622E-8&3.963E-6&3.335E-6&4.217E-6&1.358E-5&4.043E-7&8.434E-7&6.325E-6&4.970E-8&2.300E-7&4.515E-6\\
4&9.388E-13&2.890E-7&2.140E-5&1.040E-5&1.308E-5&5.246E-5&1.440E-6&3.358E-6&2.988E-5&2.004E-7&1.174E-6&2.346E-5\\
5&1.267E-11&1.671E-6&8.689E-5&3.698E-5&4.200E-5&1.744E-4&5.824E-6&1.270E-5&1.121E-4&9.060E-7&5.198E-6&9.318E-5\\
6&1.558E-10&8.159E-6&2.985E-4&1.428E-4&1.388E-4&5.293E-4&2.578E-5&4.705E-5&3.670E-4&4.505E-6&2.140E-5&3.158E-4\\
7&1.643E-9&3.726E-5&9.524E-4&5.914E-4&4.836E-4&1.563E-3&1.239E-4&1.788E-4&1.137E-3&2.452E-5&8.751E-5&9.991E-4\\
8&1.974E-8&1.742E-4&3.058E-3&2.653E-3&1.833E-3&4.780E-3&6.642E-4&7.327E-4&3.584E-3&1.530E-4&3.795E-4&3.192E-3\\
9&2.554E-7&9.167E-4&1.076E-2&1.293E-2&7.858E-3&1.614E-2&3.991E-3&3.430E-3&1.243E-2&1.100E-3&1.877E-3&1.119E-2\\
10&5.733E-6&6.121E-3&4.509E-2&6.968E-2&4.009E-2&6.403E-2&2.859E-2&1.980E-2&5.114E-2&1.023E-2&1.165E-2&4.666E-2\\
11&1.856E-4&6.433E-2&2.511E-1&3.745E-1&2.533E-1&3.151E-1&2.294E-1&1.579E-1&2.729E-1&1.185E-1&1.069E-1&2.569E-1\\
12&1.000E0&1.000E0&1.000E0&1.000E0&1.000E0&1.000E0&1.000E0&1.000E0&1.000E0&1.000E0&1.000E0&1.000E0\\ \hline  \hline 
\end{tabular}
%}
\end{table}
\end{tiny}

~ ~ ~

\begin{tiny}
\begin{table}
\setlength{\tabcolsep}{2.0pt}
\renewcommand{\arraystretch}{1.02}
\caption{The extreme points of ${\varepsilon}(q,\delta,\alpha,k;\t)$
 of $(q,\delta,\alpha,k)$-BURA approximation, 
i.e. the best uniform rational approximation of  
$\f(q,\delta,\alpha;\t)=\t^{\alpha}/(1+q\,\t^{\alpha})$  with functions from $\mathcal R_k$ on 
 $ [\delta, 1]$,   $\alpha =0.25, 0.5, 0.75$, ${\bf k=6}$, and $q=200$}
 %$(q=200,\delta,\alpha,k=6)$-BURA }
 \label{tabl:Bq200dDaAAk6p}
%\rotatebox{90}{
\begin{tabular}{||c||c|c|c||c|c|c||c|c|c||c|c|c||}
\hline \hline
{\bf k=6} /$\delta$& $0.00$ & $0.00$ & $0.00$ &$10^{-6}$&$10^{-6}$&$10^{-6}$&$10^{-7}$&$10^{-7}$&$10^{-7}$&$10^{-8}$&$10^{-8}$&$10^{-8}$\\
\ j /$\alpha$ & 0.25   & 0.50   & 0.75   & 0.25   & 0.50   & 0.75   & 0.25   & 0.50   & 0.75   & 0.25   & 0.50   & 0.75   \\ \hline
1&0.000E0&0.000E0&0.000E0&1.000E-6&1.000E-6&1.000E-6&1.000E-7&1.000E-7&1.000E-7&1.000E-8&1.000E-8&1.000E-8\\
2&1.531E-16&5.552E-10&1.079E-7&1.300E-6&1.405E-6&2.101E-6&1.371E-7&1.751E-7&4.489E-7&1.457E-8&2.590E-8&1.767E-7\\
3&5.613E-15&8.652E-9&1.150E-6&2.585E-6&3.114E-6&7.534E-6&3.057E-7&5.413E-7&2.665E-6&3.666E-8&1.213E-7&1.518E-6\\
4&1.327E-13&6.857E-8&6.264E-6&6.573E-6&8.149E-6&2.549E-5&8.809E-7&1.813E-6&1.157E-5&1.190E-7&5.269E-7&7.649E-6\\
5&1.665E-12&3.887E-7&2.561E-5&1.892E-5&2.233E-5&7.742E-5&2.862E-6&5.913E-6&4.117E-5&4.311E-7&2.049E-6&2.984E-5\\
6&1.783E-11&1.809E-6&8.712E-5&5.838E-5&6.244E-5&2.157E-4&9.969E-6&1.877E-5&1.278E-4&1.671E-6&7.376E-6&9.847E-5\\
7&1.574E-10&7.499E-6&2.644E-4&1.892E-4&1.789E-4&5.700E-4&3.667E-5&5.934E-5&3.643E-4&6.863E-6&2.556E-5&2.927E-4\\
8&1.388E-9&2.951E-5&7.548E-4&6.419E-4&5.331E-4&1.482E-3&1.430E-4&1.919E-4&9.974E-4&3.024E-5&8.864E-5&8.242E-4\\
9&1.202E-8&1.169E-4&2.142E-3&2.283E-3&1.680E-3&3.941E-3&5.936E-4&6.535E-4&2.751E-3&1.437E-4&3.201E-4&2.317E-3\\
10&1.226E-7&4.955E-4&6.365E-3&8.560E-3&5.709E-3&1.111E-2&2.664E-3&2.418E-3&8.000E-3&7.587E-4&1.254E-3&6.839E-3\\
11&1.429E-6&2.398E-3&2.091E-2&3.382E-2&2.138E-2&3.436E-2&1.298E-2&1.006E-2&2.566E-2&4.476E-3&5.587E-3&2.230E-2\\
12&2.828E-5&1.448E-2&7.967E-2&1.388E-1&8.991E-2&1.196E-1&6.985E-2&4.918E-2&9.442E-2&3.125E-2&3.010E-2&8.407E-2\\
13&8.278E-4&1.258E-1&3.597E-1&5.167E-1&4.039E-1&4.512E-1&3.750E-1&2.873E-1&3.967E-1&2.418E-1&2.105E-1&3.712E-1\\
14&1.000E0&1.000E0&1.000E0&1.000E0&1.000E0&1.000E0&1.000E0&1.000E0&1.000E0&1.000E0&1.000E0&1.000E0\\ \hline  \hline 
\end{tabular}
%}
\end{table}
\end{tiny}

~ ~ ~

\begin{tiny}
\begin{table}
\setlength{\tabcolsep}{2.0pt}
\renewcommand{\arraystretch}{1.02}
\caption{The extreme points of ${\varepsilon}(q,\delta,\alpha,k;\t)$
 of $(q,\delta,\alpha,k)$-BURA approximation, 
i.e. the best uniform rational approximation of  
$\f(q,\delta,\alpha;\t)=\t^{\alpha}/(1+q\,\t^{\alpha})$  with functions from $\mathcal R_k$ on 
 $ [\delta, 1]$,   $\alpha =0.25, 0.5, 0.75$, ${\bf k=7}$, and $q=200$}
 %$(q=200,\delta,\alpha,k=7)$-BURA }
 \label{tabl:Bq200dDaAAk7p}
%\rotatebox{90}{
\begin{tabular}{||c||c|c|c||c|c|c||c|c|c||c|c|c||}
\hline \hline
{\bf k=7} /$\delta$& $0.00$ & $0.00$ & $0.00$ &$10^{-6}$&$10^{-6}$&$10^{-6}$&$10^{-7}$&$10^{-7}$&$10^{-7}$&$10^{-8}$&$10^{-8}$&$10^{-8}$\\
\ j /$\alpha$ & 0.25   & 0.50   & 0.75   & 0.25   & 0.50   & 0.75   & 0.25   & 0.50   & 0.75   & 0.25   & 0.50   & 0.75   \\ \hline
1&&0.000E0&0.000E0&1.000E-6&1.000E-6&1.000E-6&1.000E-7&1.000E-7&1.000E-7&1.000E-8&1.000E-8&1.000E-8\\
2&&1.429E-10&3.334E-8&1.232E-6&1.302E-6&1.706E-6&1.287E-7&1.531E-7&2.895E-7&1.354E-8&2.029E-8&8.152E-8\\
3&&2.226E-9&3.560E-7&2.153E-6&2.496E-6&4.868E-6&2.492E-7&3.894E-7&1.351E-6&2.922E-8&7.475E-8&5.986E-7\\
4&&1.761E-8&1.950E-6&4.696E-6&5.685E-6&1.439E-5&6.104E-7&1.122E-6&5.268E-6&8.018E-8&2.789E-7&2.851E-6\\
5&&9.917E-8&8.024E-6&1.154E-5&1.374E-5&3.992E-5&1.681E-6&3.236E-6&1.756E-5&2.452E-7&9.666E-7&1.078E-5\\
6&&4.538E-7&2.744E-5&3.029E-5&3.395E-5&1.035E-4&4.931E-6&9.145E-6&5.196E-5&7.952E-7&3.134E-6&3.483E-5\\
\hline  \hline 
\end{tabular}
%}
\end{table}
\end{tiny}

~ ~ ~

\begin{tiny}
\begin{table}
\setlength{\tabcolsep}{2.0pt}
\renewcommand{\arraystretch}{1.02}
\caption{The extreme points of ${\varepsilon}(q,\delta,\alpha,k;\t)$
 of $(q,\delta,\alpha,k)$-BURA approximation, 
i.e. the best uniform rational approximation of  
$\f(q,\delta,\alpha;\t)=\t^{\alpha}/(1+q\,\t^{\alpha})$  with functions from $\mathcal R_k$ on 
 $ [\delta, 1]$,   $\alpha =0.25, 0.5, 0.75$, ${\bf k=8}$, and $q=200$}
 %$(q=200,\delta,\alpha,k=8)$-BURA }
 \label{tabl:Bq200dDaAAk8p}
%\rotatebox{90}{
\begin{tabular}{||c||c|c|c||c|c|c||c|c|c||c|c|c||}
\hline \hline
{\bf k=8} /$\delta$& $0.00$ & $0.00$ & $0.00$ &$10^{-6}$&$10^{-6}$&$10^{-6}$&$10^{-7}$&$10^{-7}$&$10^{-7}$&$10^{-8}$&$10^{-8}$&$10^{-8}$\\
\ j /$\alpha$ & 0.25   & 0.50   & 0.75   & 0.25   & 0.50   & 0.75   & 0.25   & 0.50   & 0.75   & 0.25   & 0.50   & 0.75   \\ \hline
1&&0.000E0&0.000E0&1.000E-6&1.000E-6&1.000E-6&1.000E-7&1.000E-7&1.000E-7&1.000E-8&1.000E-8&1.000E-8\\
2&&3.908E-11&1.091E-8&1.184E-6&1.234E-6&1.488E-6&1.229E-7&1.395E-7&2.160E-7&1.283E-8&1.718E-8&4.602E-8\\
3&&6.087E-10&1.166E-7&1.881E-6&2.115E-6&3.519E-6&2.136E-7&3.036E-7&7.989E-7&2.456E-8&5.160E-8&2.773E-7\\
4&&4.812E-9&6.404E-7&3.641E-6&4.298E-6&9.141E-6&4.604E-7&7.681E-7&2.765E-6&5.891E-8&1.674E-7&1.219E-6\\
5&&2.704E-8&2.647E-6&7.909E-6&9.366E-6&2.319E-5&1.113E-6&1.991E-6&8.547E-6&1.575E-7&5.220E-7&4.399E-6\\
6&&1.230E-7&9.100E-6&1.832E-5&2.096E-5&5.612E-5&2.861E-6&5.112E-6&2.398E-5&4.452E-7&1.546E-6&1.378E-5\\
\hline  \hline 
\end{tabular}
%}
\end{table}
\end{tiny}

~ ~ ~

\begin{tiny}
\begin{table}
\setlength{\tabcolsep}{2.0pt}
\renewcommand{\arraystretch}{1.02}
\caption{The extreme points of ${\varepsilon}(q,\delta,\alpha,k;\t)$
 of $(q,\delta,\alpha,k)$-BURA approximation, 
i.e. the best uniform rational approximation of  
$\f(q,\delta,\alpha;\t)=\t^{\alpha}/(1+q\,\t^{\alpha})$  with functions from $\mathcal R_k$ on 
 $ [\delta, 1]$,   $\alpha =0.25, 0.5, 0.75$, ${\bf k=3}$, and $q=400$}
 %$(q=400,\delta,\alpha,k=3)$-BURA }
 \label{tabl:Bq400dDaAAk3p}
%\rotatebox{90}{
\begin{tabular}{||c||c|c|c||c|c|c||c|c|c||c|c|c||}
\hline \hline
{\bf k=3}/$\delta$& $0.00$ & $0.00$ & $0.00$ &$10^{-6}$&$10^{-6}$&$10^{-6}$&$10^{-7}$&$10^{-7}$&$10^{-7}$&$10^{-8}$&$10^{-8}$&$10^{-8}$\\
\ j /$\alpha$ & 0.25   & 0.50   & 0.75   & 0.25   & 0.50   & 0.75   & 0.25   & 0.50   & 0.75   & 0.25   & 0.50   & 0.75   \\ \hline
1&0.000E0&0.000E0&0.000E0&1.000E-6&1.000E-6&1.000E-6&1.000E-7&1.000E-7&1.000E-7&1.000E-8&1.000E-8&1.000E-8\\
2&4.989E-15&1.484E-8&2.509E-6&1.875E-6&1.997E-6&6.499E-6&2.014E-7&3.192E-7&3.334E-6&2.171E-8&7.570E-8&2.679E-6\\
3&2.141E-13&2.390E-7&2.589E-5&8.079E-6&7.922E-6&4.368E-5&9.784E-7&1.884E-6&3.008E-5&1.174E-7&6.794E-7&2.679E-5\\
4&8.012E-12&2.119E-6&1.405E-4&5.128E-5&3.763E-5&2.060E-4&7.090E-6&1.085E-5&1.565E-4&9.439E-7&4.713E-6&1.440E-4\\
5&1.714E-10&1.603E-5&6.613E-4&4.031E-4&2.089E-4&9.101E-4&6.447E-5&6.660E-5&7.228E-4&9.520E-6&3.170E-5&6.746E-4\\
6&6.305E-9&1.391E-4&3.454E-3&4.410E-3&1.527E-3&4.643E-3&8.962E-4&5.205E-4&3.749E-3&1.579E-4&2.601E-4&3.518E-3\\
7&2.638E-7&2.105E-3&3.106E-2&6.442E-2&1.942E-2&4.062E-2&1.799E-2&7.252E-3&3.347E-2&3.945E-3&3.791E-3&3.158E-2\\
8&1.000E0&1.000E0&1.000E0&1.000E0&1.000E0&1.000E0&1.000E0&1.000E0&1.000E0&1.000E0&1.000E0&1.000E0\\ \hline \hline
\end{tabular}
%}
\end{table}
\end{tiny}

~ ~ ~

\begin{tiny}
\begin{table}
\setlength{\tabcolsep}{2.0pt}
\renewcommand{\arraystretch}{1.02}
\caption{The extreme points of ${\varepsilon}(q,\delta,\alpha,k;\t)$
 of $(q,\delta,\alpha,k)$-BURA approximation, 
i.e. the best uniform rational approximation of  
$\f(q,\delta,\alpha;\t)=\t^{\alpha}/(1+q\,\t^{\alpha})$  with functions from $\mathcal R_k$ on 
 $ [\delta, 1]$,   $\alpha =0.25, 0.5, 0.75$, ${\bf k=4}$, and $q=400$}
 %$(q=400,\delta,\alpha,k=4)$-BURA }
 \label{tabl:Bq400dDaAAk4p}
%\rotatebox{90}{
\begin{tabular}{||c||c|c|c||c|c|c||c|c|c||c|c|c||}
\hline \hline
{\bf k=4} /$\delta$& $0.00$ & $0.00$ & $0.00$ &$10^{-6}$&$10^{-6}$&$10^{-6}$&$10^{-7}$&$10^{-7}$&$10^{-7}$&$10^{-8}$&$10^{-8}$&$10^{-8}$\\
\ j /$\alpha$ & 0.25   & 0.50   & 0.75   & 0.25   & 0.50   & 0.75   & 0.25   & 0.50   & 0.75   & 0.25   & 0.50   & 0.75   \\ \hline
 1&0.000E0&0.000E0&0.000E0&1.000E-6&1.000E-6&1.000E-6&1.000E-7&1.000E-7&1.000E-7&1.000E-8&1.000E-8&1.000E-8\\
 2&5.056E-16&2.919E-9&6.241E-7&1.560E-6&1.630E-6&3.451E-6&1.659E-7&2.231E-7&1.179E-6&1.768E-8&4.000E-8&7.402E-7\\
 3&1.987E-14&4.586E-8&6.535E-6&4.650E-6&4.751E-6&1.773E-5&5.574E-7&9.366E-7&9.280E-6&6.616E-8&2.597E-7&7.162E-6\\
 4&5.685E-13&3.750E-7&3.487E-5&1.903E-5&1.663E-5&7.202E-5&2.591E-6&4.067E-6&4.476E-5&3.415E-7&1.414E-6&3.720E-5\\
 5&8.888E-12&2.320E-6&1.442E-4&9.104E-5&6.357E-5&2.569E-4&1.414E-5&1.767E-5&1.752E-4&2.065E-6&6.990E-6&1.516E-4\\
 6&1.474E-10&1.312E-5&5.404E-4&4.983E-4&2.699E-4&8.930E-4&9.091E-5&8.159E-5&6.388E-4&1.505E-5&3.486E-5&5.640E-4\\
 7&2.282E-9&7.934E-5&2.139E-3&3.097E-3&1.343E-3&3.397E-3&6.841E-4&4.362E-4&2.493E-3&1.312E-4&1.961E-4&2.224E-3\\
 8&6.301E-8&6.139E-4&1.035E-2&2.339E-2&8.575E-3&1.594E-2&6.883E-3&3.049E-3&1.195E-2&1.680E-3&1.441E-3&1.073E-2\\
 9&2.363E-6&8.453E-3&8.187E-2&2.030E-1&8.414E-2&1.168E-1&8.945E-2&3.572E-2&9.232E-2&3.025E-2&1.844E-2&8.443E-2\\
 10&1.000E0&1.000E0&1.000E0&1.000E0&1.000E0&1.000E0&1.000E0&1.000E0&1.000E0&1.000E0&1.000E0&1.000E0\\ \hline \hline
\end{tabular}%
%}
\end{table}
\end{tiny}
%
%~ ~ ~
%
%
\begin{tiny}
\begin{table}[h!t]
\setlength{\tabcolsep}{2.0pt}
\renewcommand{\arraystretch}{1.02}
\caption{The extreme points of ${\varepsilon}(q,\delta,\alpha,k;\t)$
 of $(q,\delta,\alpha,k)$-BURA approximation, 
i.e. the best uniform rational approximation of  
$\f(q,\delta,\alpha;\t)=\t^{\alpha}/(1+q\,\t^{\alpha})$  with functions from $\mathcal R_k$ on 
 $ [\delta, 1]$,   $\alpha =0.25, 0.5, 0.75$, ${\bf k=8}$, and $q=400$}
 %$(q=400,\delta,\alpha,k=8)$-BURA }
 \label{tabl:Bq400dDaAAk8p}
%\rotatebox{90}{
%
%}
\end{table}
\end{tiny}

%~ ~ ~\\[8ex]

%\endinput

\break

%\newpage

\section{Tables type \textbf{(c)} for BURA-poles} 
Now we provide the data about the poles of the BURA element $r_{q,\delta,\alpha,k}(\t)$:
\begin{tiny}
\begin{table}[t!]
\setlength{\tabcolsep}{2.0pt}
\renewcommand{\arraystretch}{1.02}
\caption{The poles of $r_{q,\delta,\alpha,k}(\t)$, \eqref{BURA},  of $(q,\delta,\alpha,k)$-BURA approximation, 
i.e. the best uniform rational approximation of  
$\f(q,\delta,\alpha;\t)=\t^{\alpha}/(1+q\,\t^{\alpha})$  with functions from $\mathcal R_k$ on 
 $ [\delta, 1]$,   $\alpha =0.25, 0.5, 0.75$, ${\bf k=3}$, and $q=0$}
 %$(q=000,\delta,\alpha,k=3)$-BURA-poles}
 \label{tabl:Cq000dDaAAk3p}
%\rotatebox{90}{
\begin{tabular}{||c||c|c|c||c|c|c||c|c|c||c|c|c||}
\hline \hline
{\bf k=3} /$\delta$& $0.00$ & $0.00$ & $0.00$ &$10^{-6}$&$10^{-6}$&$10^{-6}$&$10^{-7}$&$10^{-7}$&$10^{-7}$&$10^{-8}$&$10^{-8}$&$10^{-8}$\\
\ j /$\alpha$ & 0.25   & 0.50   & 0.75   & 0.25   & 0.50   & 0.75   & 0.25   & 0.50   & 0.75   & 0.25   & 0.50   & 0.75   \\ \hline
  1&-1.36E-5&-8.81E-4&-5.45E-3&-9.29E-5&-1.04E-3&-5.53E-3&-4.51E-5&-9.32E-4&-5.46E-3&-2.78E-5&-8.97E-4&-5.45E-3\\
  2&-4.78E-3&-4.93E-2&-1.49E-1&-1.01E-2&-5.24E-2&-1.50E-1&-7.49E-3&-5.03E-2&-1.50E-1&-6.21E-3&-4.97E-2&-1.50E-1\\
  3&-4.14E-1&-1.77E0&-6.09E0&-5.48E-1&-1.82E0&-6.11E0&-4.88E-1&-1.79E0&-6.10E0&-4.55E-1&-1.78E0&-6.10E0\\
\hline \hline
\end{tabular}
%}
\end{table}
\end{tiny}

~ ~ ~

\begin{tiny}
\begin{table}
\setlength{\tabcolsep}{2.0pt}
\renewcommand{\arraystretch}{1.02}
\caption{The poles of  $r_{q,\delta,\alpha,k}(\t)$, %$r(q,\delta,\alpha,k;\t)$, 
\eqref{BURA},  of $(q,\delta,\alpha,k)$-BURA approximation, 
i.e. the best uniform rational approximation of  
$\f(q,\delta,\alpha;\t)=\t^{\alpha}/(1+q\,\t^{\alpha})$  with functions from $\mathcal R_k$ on 
 $ [\delta, 1]$,   $\alpha =0.25, 0.5, 0.75$, ${\bf k=4}$, and $q=0$}
 %$(q=000,\delta,\alpha,k=4)$-BURA-poles}
 \label{tabl:Cq000dDaAAk4p}
%\rotatebox{90}{
\begin{tabular}{||c||c|c|c||c|c|c||c|c|c||c|c|c||}
\hline \hline
{\bf k=4} /$\delta$& $0.00$ & $0.00$ & $0.00$ &$10^{-6}$&$10^{-6}$&$10^{-6}$&$10^{-7}$&$10^{-7}$&$10^{-7}$&$10^{-8}$&$10^{-8}$&$10^{-8}$\\
\ j /$\alpha$ & 0.25   & 0.50   & 0.75   & 0.25   & 0.50   & 0.75   & 0.25   & 0.50   & 0.75   & 0.25   & 0.50   & 0.75   \\ \hline
  1&-5.61E-7&-9.17E-5&-8.37E-4&-1.95E-5&-1.49E-4&-8.91E-4&-6.35E-6&-1.10E-4&-8.48E-4&-2.67E-6&-9.75E-5&-8.39E-4\\
  2&-1.97E-4&-4.98E-3&-2.11E-2&-1.10E-3&-6.15E-3&-2.17E-2&-6.01E-4&-5.38E-3&-2.13E-2&-3.92E-4&-5.11E-3&-2.12E-2\\
  3&-1.41E-2&-1.01E-1&-2.68E-1&-3.27E-2&-1.12E-1&-2.72E-1&-2.41E-2&-1.05E-1&-2.69E-1&-1.95E-2&-1.03E-1&-2.68E-1\\
  4&-6.27E-1&-2.49E0&-8.40E0&-9.10E-1&-2.62E0&-8.47E0&-7.90E-1&-2.53E0&-8.42E0&-7.20E-1&-2.50E0&-8.41E0\\
\hline \hline
\end{tabular}
%}
\end{table}
\end{tiny}

~ ~ ~
\begin{tiny}
\begin{table}
\setlength{\tabcolsep}{2.0pt}
\renewcommand{\arraystretch}{1.02}
\caption{The poles of  $r_{q,\delta,\alpha,k}(\t)$, %${r}(q,\delta,\alpha,k;\t)$, 
\eqref{BURA},  of $(q,\delta,\alpha,k)$-BURA approximation, 
i.e. the best uniform rational approximation of  
$\f(q,\delta,\alpha;\t)=\t^{\alpha}/(1+q\,\t^{\alpha})$  with functions from $\mathcal R_k$ on 
 $ [\delta, 1]$,   $\alpha =0.25, 0.5, 0.75$, ${\bf k=5}$, and $q=0$}
 %$(q=000,\delta,\alpha,k=5)$-BURA-poles}
 \label{tabl:Cq000dDaAAk5p}
%\rotatebox{90}{
\begin{tabular}{||c||c|c|c||c|c|c||c|c|c||c|c|c||}
\hline \hline
{\bf k=5} /$\delta$& $0.00$ & $0.00$ & $0.00$ &$10^{-6}$&$10^{-6}$&$10^{-6}$&$10^{-7}$&$10^{-7}$&$10^{-7}$&$10^{-8}$&$10^{-8}$&$10^{-8}$\\
\ j /$\alpha$ & 0.25   & 0.50   & 0.75   & 0.25   & 0.50   & 0.75   & 0.25   & 0.50   & 0.75   & 0.25   & 0.50   & 0.75   \\ \hline
  1&-3.27E-8&-1.22E-5&-1.59E-4&-7.21E-6&-3.56E-5&-1.92E-4&-1.77E-6&-1.94E-5&-1.66E-4&-5.28E-7&-1.45E-5&-1.60E-4\\
  2&-1.15E-5&-6.62E-4&-3.97E-3&-2.29E-4&-1.12E-3&-4.33E-3&-9.42E-5&-8.22E-4&-4.05E-3&-4.65E-5&-7.15E-4&-3.98E-3\\
  3&-8.15E-4&-1.28E-2&-4.47E-2&-4.60E-3&-1.70E-2&-4.69E-2&-2.68E-3&-1.44E-2&-4.52E-2&-1.77E-3&-1.33E-2&-4.48E-2\\
  4&-2.81E-2&-1.63E-1&-3.97E-1&-7.13E-2&-1.90E-1&-4.08E-1&-5.27E-2&-1.73E-1&-4.00E-1&-4.21E-2&-1.66E-1&-3.98E-1\\
  5&-8.47E-1&-3.21E0&-1.08E1&-1.35E0&-3.52E0&-1.09E1&-1.15E0&-3.33E0&-1.08E1&-1.03E0&-3.25E0&-1.08E1\\
\hline \hline
\end{tabular}
%}
\end{table}
\end{tiny}

~ ~ ~

\begin{tiny}
\begin{table}
\setlength{\tabcolsep}{2.0pt}
\renewcommand{\arraystretch}{1.02}
\caption{The poles of  $r_{q,\delta,\alpha,k}(\t)$,  %${r}(q,\delta,\alpha,k;\t)$, 
\eqref{BURA},  of $(q,\delta,\alpha,k)$-BURA approximation, 
i.e. the best uniform rational approximation of  
$\f(q,\delta,\alpha;\t)=\t^{\alpha}/(1+q\,\t^{\alpha})$  with functions from $\mathcal R_k$ on 
 $ [\delta, 1]$,   $\alpha =0.25, 0.5, 0.75$, ${\bf k=6}$, and $q=0$}
 %$(q=000,\delta,\alpha,k=6)$-BURA-poles}
 \label{tabl:Cq000dDaAAk6p}
%\rotatebox{90}{
\begin{tabular}{||c||c|c|c||c|c|c||c|c|c||c|c|c||}
\hline \hline
{\bf k=6} /$\delta$& $0.00$ & $0.00$ & $0.00$ &$10^{-6}$&$10^{-6}$&$10^{-6}$&$10^{-7}$&$10^{-7}$&$10^{-7}$&$10^{-8}$&$10^{-8}$&$10^{-8}$\\
\ j /$\alpha$ & 0.25   & 0.50   & 0.75   & 0.25   & 0.50   & 0.75   & 0.25   & 0.50   & 0.75   & 0.25   & 0.50   & 0.75   \\ \hline
  1&-2.45E-9&-1.95E-6&-3.52E-5&-3.62E-6&-1.28E-5&-5.49E-5&-7.44E-7&-5.05E-6&-3.98E-5&-1.74E-7&-2.90E-6&-3.62E-5\\
  2&-8.61E-7&-1.06E-4&-8.75E-4&-7.42E-5&-2.93E-4&-1.10E-3&-2.41E-5&-1.72E-4&-9.30E-4&-9.22E-6&-1.28E-4&-8.87E-4\\
  3&-6.11E-5&-2.03E-3&-9.65E-3&-1.06E-3&-3.74E-3&-1.10E-2&-4.90E-4&-2.70E-3&-1.00E-2&-2.60E-4&-2.27E-3&-9.73E-3\\
  4&-2.07E-3&-2.40E-2&-7.40E-2&-1.22E-2&-3.47E-2&-8.01E-2&-7.43E-3&-2.84E-2&-7.56E-2&-4.98E-3&-2.56E-2&-7.44E-2\\
  5&-4.58E-2&-2.30E-1&-5.32E-1&-1.26E-1&-2.85E-1&-5.58E-1&-9.40E-2&-2.53E-1&-5.39E-1&-7.46E-2&-2.39E-1&-5.34E-1\\
  6&-1.07E0&-3.95E0&-1.31E1&-1.86E0&-4.52E0&-1.36E1&-1.57E0&-4.20E0&-1.32E1&-1.38E0&-4.04E0&-1.32E1\\
\hline \hline
\end{tabular}
%}
\end{table}
\end{tiny}

~ ~ ~

\begin{tiny}
\begin{table}
\setlength{\tabcolsep}{2.0pt}
\renewcommand{\arraystretch}{1.02}
\caption{The poles of  $r_{q,\delta,\alpha,k}(\t)$,  %${r}(q,\delta,\alpha,k;\t)$, 
\eqref{BURA},  of $(q,\delta,\alpha,k)$-BURA approximation, 
i.e. the best uniform rational approximation of  
$\f(q,\delta,\alpha;\t)=\t^{\alpha}/(1+q\,\t^{\alpha})$  with functions from $\mathcal R_k$ on 
 $ [\delta, 1]$,   $\alpha =0.25, 0.5, 0.75$, ${\bf k=7}$, and $q=0$}
 %$(q=000,\delta,\alpha,k=7)$-BURA-poles}
 \label{tabl:Cq000dDaAAk7p}
%\rotatebox{90}{
\begin{tabular}{||c||c|c|c||c|c|c||c|c|c||c|c|c||}
\hline \hline
{\bf k=7} /$\delta$& $0.00$ & $0.00$ & $0.00$ &$10^{-6}$&$10^{-6}$&$10^{-6}$&$10^{-7}$&$10^{-7}$&$10^{-7}$&$10^{-8}$&$10^{-8}$&$10^{-8}$\\
\ j /$\alpha$ & 0.25   & 0.50   & 0.75   & 0.25   & 0.50   & 0.75   & 0.25   & 0.50   & 0.75   & 0.25   & 0.50   & 0.75   \\ \hline
  1&-2.24E-10&-3.58E-7&-8.75E-6&-2.17E-6&-6.05E-6&-2.05E-5&-3.98E-7&-1.84E-6&-1.17E-5&-7.97E-8&-7.90E-7&-9.40E-6\\
  2&-7.85E-8&-1.94E-5&-2.17E-4&-3.21E-5&-1.03E-4&-3.45E-4&-8.69E-6&-4.79E-5&-2.52E-4&-2.70E-6&-2.94E-5&-2.25E-4\\
  3&-5.57E-6&-3.72E-4&-2.39E-3&-3.42E-4&-1.09E-3&-3.18E-3&-1.30E-4&-6.56E-4&-2.61E-3&-5.65E-5&-4.80E-4&-2.44E-3\\
  4&-1.89E-4&-4.34E-3&-1.77E-2&-3.10E-3&-8.76E-3&-2.13E-2&-1.56E-3&-6.26E-3&-1.88E-2&-8.72E-4&-5.11E-3&-1.80E-2\\
  5&-4.08E-3&-3.80E-2&-1.08E-1&-2.51E-2&-5.94E-2&-1.21E-1&-1.58E-2&-4.79E-2&-1.12E-1&-1.08E-2&-4.21E-2&-1.09E-1\\
  6&-6.66E-2&-3.01E-1&-6.71E-1&-1.97E-1&-3.96E-1&-7.24E-1&-1.48E-1&-3.46E-1&-6.87E-1&-1.17E-1&-3.20E-1&-6.75E-1\\
  7&-1.30E00&-4.69E00&-1.55E01&-2.45E00&-5.64E00&-1.64E01&-2.05E00&-5.15E00&-1.58E01&-1.78E00&-4.88E00&-1.56E01\\
\hline \hline
\end{tabular}
%}
\end{table}
\end{tiny}

~ ~ ~

\begin{tiny}
\begin{table}
\setlength{\tabcolsep}{2.0pt}
\renewcommand{\arraystretch}{1.02}
\caption{The poles of $r_{q,\delta,\alpha,k}(\t)$,  %${r}(q,\delta,\alpha,k;\t)$, 
\eqref{BURA},  of $(q,\delta,\alpha,k)$-BURA approximation, 
i.e. the best uniform rational approximation of  
$\f(q,\delta,\alpha;\t)=\t^{\alpha}/(1+q\,\t^{\alpha})$  with functions from $\mathcal R_k$ on 
 $ [\delta, 1]$,   $\alpha =0.25, 0.5, 0.75$, ${\bf k=8}$, and $q=0$}
 %$(q=000,\delta,\alpha,k=8)$-BURA-poles}
 \label{tabl:Cq000dDaAAk8p}
%\rotatebox{90}{
\begin{tabular}{||c||c|c|c||c|c|c||c|c|c||c|c|c||}
\hline \hline
{\bf k=8} /$\delta$& $0.00$ & $0.00$ & $0.00$ &$10^{-6}$&$10^{-6}$&$10^{-6}$&$10^{-7}$&$10^{-7}$&$10^{-7}$&$10^{-8}$&$10^{-8}$&$10^{-8}$\\
\ j /$\alpha$ & 0.25   & 0.50   & 0.75   & 0.25   & 0.50   & 0.75   & 0.25   & 0.50   & 0.75   & 0.25   & 0.50   & 0.75   \\ \hline
  1&-2.39E-11&-7.35E-8&-2.38E-6&-1.45E-6&-3.45E-6&-9.54E-6&-2.47E-7&-8.67E-7&-4.20E-6&-4.45E-8&-2.86E-7&-2.82E-6\\
  2&-8.37E-9&-3.98E-6&-5.93E-5&-1.69E-5&-4.52E-5&-1.33E-4&-3.97E-6&-1.71E-5&-8.07E-5&-1.04E-6&-8.49E-6&-6.47E-5\\
  3&-5.95E-7&-7.62E-5&-6.50E-4&-1.41E-4&-3.97E-4&-1.10E-3&-4.56E-5&-2.00E-4&-7.91E-4&-1.66E-5&-1.24E-4&-6.87E-4\\
  4&-2.01E-5&-8.89E-4&-4.79E-3&-1.04E-3&-2.79E-3&-6.84E-3&-4.41E-4&-1.72E-3&-5.47E-3&-2.08E-4&-1.24E-3&-4.97E-3\\
  5&-4.34E-4&-7.65E-3&-2.79E-2&-6.99E-3&-1.67E-2&-3.55E-2&-3.73E-3&-1.19E-2&-3.05E-2&-2.17E-3&-9.55E-3&-2.86E-2\\
  6&-6.88E-3&-5.44E-2&-1.44E-1&-4.40E-2&-9.13E-2&-1.69E-1&-2.85E-2&-7.28E-2&-1.53E-1&-1.97E-2&-6.28E-2&-1.47E-1\\
  7&-8.99E-2&-3.75E-1&-8.13E-1&-2.83E-1&-5.23E-1&-9.06E-1&-2.13E-1&-4.51E-1&-8.46E-1&-1.69E-1&-4.10E-1&-8.22E-1\\
  8&-1.53E00&-5.43E00&-1.79E01&-3.11E00&-6.88E00&-1.95E01&-2.58E00&-6.18E00&-1.85E01&-2.23E00&-5.78E00&-1.81E01\\
\hline \hline
\end{tabular}
%}
\end{table}
\end{tiny}

~ ~ ~

\begin{tiny}
\begin{table}
\setlength{\tabcolsep}{2.0pt}
\renewcommand{\arraystretch}{1.02}
\caption{The poles of $r_{q,\delta,\alpha,k}(\t)$,  %${r}(q,\delta,\alpha,k;\t)$, 
\eqref{BURA},  of $(q,\delta,\alpha,k)$-BURA approximation, 
i.e. the best uniform rational approximation of  
$\f(q,\delta,\alpha;\t)=\t^{\alpha}/(1+q\,\t^{\alpha})$  with functions from $\mathcal R_k$ on 
 $ [\delta, 1]$,   $\alpha =0.25, 0.5, 0.75$, ${\bf k=3}$, and $q=1$}
 %$(q=001,\delta,\alpha,k=3)$-BURA-poles}
 \label{tabl:Cq001dDaAAk3p}
%\rotatebox{90}{
\begin{tabular}{||c||c|c|c||c|c|c||c|c|c||c|c|c||}
\hline \hline
{\bf k=3} /$\delta$& $0.00$ & $0.00$ & $0.00$ &$10^{-6}$&$10^{-6}$&$10^{-6}$&$10^{-7}$&$10^{-7}$&$10^{-7}$&$10^{-8}$&$10^{-8}$&$10^{-8}$\\
\ j /$\alpha$ & 0.25   & 0.50   & 0.75   & 0.25   & 0.50   & 0.75   & 0.25   & 0.50   & 0.75   & 0.25   & 0.50   & 0.75   \\ \hline
  1&-3.53E-6&-4.86E-4&-3.78E-3&-4.67E-5&-6.08E-4&-3.87E-3&-1.89E-5&-5.25E-4&-3.80E-3&-9.90E-6&-4.98E-4&-3.79E-3\\
  2&-1.43E-3&-2.61E-2&-9.83E-2&-4.41E-3&-2.84E-2&-9.90E-2&-2.88E-3&-2.69E-2&-9.84E-2&-2.17E-3&-2.64E-2&-9.83E-2\\
  3&-1.49E-1&-5.52E-1&-9.17E-1&-2.37E-1&-5.71E-1&-9.20E-1&-1.97E-1&-5.58E-1&-9.18E-1&-1.76E-1&-5.54E-1&-9.17E-1\\
\hline \hline
\end{tabular}
%}
\end{table}
\end{tiny}

~ ~ ~

\begin{tiny}
\begin{table}
\setlength{\tabcolsep}{2.0pt}
\renewcommand{\arraystretch}{1.02}
\caption{The poles of $r_{q,\delta,\alpha,k}(\t)$,  %${r}(q,\delta,\alpha,k;\t)$,  
\eqref{BURA},  of $(q,\delta,\alpha,k)$-BURA approximation, 
i.e. the best uniform rational approximation of  
$\f(q,\delta,\alpha;\t)=\t^{\alpha}/(1+q\,\t^{\alpha})$  with functions from $\mathcal R_k$ on 
 $ [\delta, 1]$,   $\alpha =0.25, 0.5, 0.75$, ${\bf k=4}$, and $q=1$}
 %$(q=001,\delta,\alpha,k=4)$-BURA-poles}
 \label{tabl:Cq001dDaAAk4p}
%\rotatebox{90}{
\begin{tabular}{||c||c|c|c||c|c|c||c|c|c||c|c|c||}
\hline \hline
{\bf k=4} /$\delta$& $0.00$ & $0.00$ & $0.00$ &$10^{-6}$&$10^{-6}$&$10^{-6}$&$10^{-7}$&$10^{-7}$&$10^{-7}$&$10^{-8}$&$10^{-8}$&$10^{-8}$\\
\ j /$\alpha$ & 0.25   & 0.50   & 0.75   & 0.25   & 0.50   & 0.75   & 0.25   & 0.50   & 0.75   & 0.25   & 0.50   & 0.75   \\ \hline
  1&-1.66E-7&-5.43E-5&-6.06E-4&-1.26E-5&-9.98E-5&-6.56E-4&-3.57E-6&-6.86E-5&-6.16E-4&-1.28E-6&-5.88E-5&-6.08E-4\\
  2&-6.21E-5&-2.94E-3&-1.53E-2&-5.87E-4&-3.85E-3&-1.59E-2&-2.82E-4&-3.25E-3&-1.54E-2&-1.63E-4&-3.04E-3&-1.53E-2\\
  3&-5.06E-3&-5.62E-2&-1.75E-1&-1.68E-2&-6.42E-2&-1.78E-1&-1.11E-2&-5.90E-2&-1.76E-1&-8.30E-3&-5.72E-2&-1.75E-1\\
  4&-2.51E-1&-7.73E-1&-1.13E00&-4.39E-1&-8.23E-1&-1.13E00&-3.61E-1&-7.91E-1&-1.13E00&-3.14E-1&-7.79E-1&-1.13E00\\
\hline \hline
\end{tabular}
%}
\end{table}
\end{tiny}

~ ~ ~

\begin{tiny}
\begin{table}
\setlength{\tabcolsep}{2.0pt}
\renewcommand{\arraystretch}{1.02}
\caption{The poles of $r_{q,\delta,\alpha,k}(\t)$,  %${ r}(q,\delta,\alpha,k;\t)$,  
\eqref{BURA},   of $(q,\delta,\alpha,k)$-BURA approximation, 
i.e. the best uniform rational approximation of  
$\f(q,\delta,\alpha;\t)=\t^{\alpha}/(1+q\,\t^{\alpha})$  with functions from $\mathcal R_k$ on 
 $ [\delta, 1]$,   $\alpha =0.25, 0.5, 0.75$, ${\bf k=5}$, and $q=1$}
 %$(q=001,\delta,\alpha,k=5)$-BURA-poles}
 \label{tabl:Cq001dDaAAk5p}
%\rotatebox{90}{
\begin{tabular}{||c||c|c|c||c|c|c||c|c|c||c|c|c||}
\hline \hline
{\bf k=5} /$\delta$& $0.00$ & $0.00$ & $0.00$ &$10^{-6}$&$10^{-6}$&$10^{-6}$&$10^{-7}$&$10^{-7}$&$10^{-7}$&$10^{-8}$&$10^{-8}$&$10^{-8}$\\
\ j /$\alpha$ & 0.25   & 0.50   & 0.75   & 0.25   & 0.50   & 0.75   & 0.25   & 0.50   & 0.75   & 0.25   & 0.50   & 0.75   \\ \hline
  1&-1.07E-8&-7.62E-6&-1.19E-4&-5.32E-6&-2.67E-5&-1.49E-4&-1.20E-6&-1.34E-5&-1.25E-4&-3.19E-7&-9.42E-6&-1.20E-4\\
  2&-3.88E-6&-4.12E-4&-2.97E-3&-1.44E-4&-7.74E-4&-3.30E-3&-5.31E-5&-5.39E-4&-3.04E-3&-2.33E-5&-4.55E-4&-2.98E-3\\
  3&-2.94E-4&-7.88E-3&-3.32E-2&-2.67E-3&-1.12E-2&-3.52E-2&-1.40E-3&-9.14E-3&-3.37E-2&-8.39E-4&-8.32E-3&-3.33E-2\\
  4&-1.14E-2&-9.28E-2&-2.55E-1&-4.08E-2&-1.12E-1&-2.62E-1&-2.77E-2&-1.00E-1&-2.57E-1&-2.06E-2&-9.55E-2&-2.55E-1\\
  5&-3.65E-1&-9.93E-1&-1.33E00&-6.96E-1&-1.10E00&-1.35E00&-5.67E-1&-1.04E00&-1.33E00&-4.87E-1&-1.01E00&-1.33E00\\
\hline \hline
\end{tabular}
%}
\end{table}
\end{tiny}

~ ~ ~

\begin{tiny}
\begin{table}
\setlength{\tabcolsep}{2.0pt}
\renewcommand{\arraystretch}{1.02}
\caption{The poles of $r_{q,\delta,\alpha,k}(\t)$,  %${ r}(q,\delta,\alpha,k;\t)$,  
\eqref{BURA},   of $(q,\delta,\alpha,k)$-BURA approximation, 
i.e. the best uniform rational approximation of  
$\f(q,\delta,\alpha;\t)=\t^{\alpha}/(1+q\,\t^{\alpha})$  with functions from $\mathcal R_k$ on 
 $ [\delta, 1]$,   $\alpha =0.25, 0.5, 0.75$, ${\bf k=6}$, and $q=1$}
 %$(q=001,\delta,\alpha,k=6)$-BURA-poles}
 \label{tabl:Cq001dDaAAk6p}
%\rotatebox{90}{
\begin{tabular}{||c||c|c|c||c|c|c||c|c|c||c|c|c||}
\hline \hline
{\bf k=6} /$\delta$& $0.00$ & $0.00$ & $0.00$ &$10^{-6}$&$10^{-6}$&$10^{-6}$&$10^{-7}$&$10^{-7}$&$10^{-7}$&$10^{-8}$&$10^{-8}$&$10^{-8}$\\
\ j /$\alpha$ & 0.25   & 0.50   & 0.75   & 0.25   & 0.50   & 0.75   & 0.25   & 0.50   & 0.75   & 0.25   & 0.50   & 0.75   \\ \hline
  1&-8.74E-10&-1.26E-6&-2.69E-5&-2.89E-6&-1.04E-5&-4.47E-5&-5.63E-7&-3.81E-6&-3.11E-5&-1.22E-7&-2.04E-6&-2.78E-5\\
  2&-3.12E-7&-6.83E-5&-6.69E-4&-5.25E-5&-2.20E-4&-8.69E-4&-1.56E-5&-1.22E-4&-7.20E-4&-5.42E-6&-8.68E-5&-6.80E-4\\
  3&-2.29E-5&-1.31E-3&-7.39E-3&-6.84E-4&-2.67E-3&-8.63E-3&-2.89E-4&-1.85E-3&-7.71E-3&-1.40E-4&-1.51E-3&-7.46E-3\\
  4&-8.29E-4&-1.52E-2&-5.59E-2&-7.59E-3&-2.36E-2&-6.13E-2&-4.23E-3&-1.88E-2&-5.73E-2&-2.61E-3&-1.66E-2&-5.62E-2\\
  5&-2.04E-2&-1.34E-1&-3.33E-1&-7.79E-2&-1.72E-1&-3.50E-1&-5.40E-2&-1.51E-1&-3.38E-1&-4.02E-2&-1.40E-1&-3.34E-1\\
  6&-4.85E-1&-1.21E00&-1.53E00&-1.00E00&-1.41E00&-1.57E00&-8.13E-1&-1.30E00&-1.54E00&-6.91E-1&-1.25E00&-1.53E00\\
\hline \hline
\end{tabular}
%}
\end{table}
\end{tiny}

~ ~ ~

\begin{tiny}
\begin{table}
\setlength{\tabcolsep}{2.0pt}
\renewcommand{\arraystretch}{1.02}
\caption{The poles of $r_{q,\delta,\alpha,k}(\t)$,  %${ r}(q,\delta,\alpha,k;\t)$,  
\eqref{BURA},   of $(q,\delta,\alpha,k)$-BURA approximation, 
i.e. the best uniform rational approximation of  
$\f(q,\delta,\alpha;\t)=\t^{\alpha}/(1+q\,\t^{\alpha})$  with functions from $\mathcal R_k$ on 
 $ [\delta, 1]$,   $\alpha =0.25, 0.5, 0.75$, ${\bf k=7}$, and $q=1$}
 %$(q=001,\delta,\alpha,k=7)$-BURA-poles}
 \label{tabl:Cq001dDaAAk7p}
%\rotatebox{90}{
\begin{tabular}{||c||c|c|c||c|c|c||c|c|c||c|c|c||}
\hline \hline
{\bf k=7} /$\delta$& $0.00$ & $0.00$ & $0.00$ &$10^{-6}$&$10^{-6}$&$10^{-6}$&$10^{-7}$&$10^{-7}$&$10^{-7}$&$10^{-8}$&$10^{-8}$&$10^{-8}$\\
\ j /$\alpha$ & 0.25   & 0.50   & 0.75   & 0.25   & 0.50   & 0.75   & 0.25   & 0.50   & 0.75   & 0.25   & 0.50   & 0.75   \\ \hline
  1&-8.53E-11&-2.39E-7&-6.80E-6&-1.82E-6&-5.18E-6&-1.75E-5&-3.22E-7&-1.50E-6&-9.46E-6&-6.14E-8&-5.98E-7&-7.41E-6\\
  2&-3.02E-8&-1.29E-5&-1.69E-4&-2.46E-5&-8.22E-5&-2.84E-4&-6.22E-6&-3.63E-5&-2.01E-4&-1.79E-6&-2.11E-5&-1.77E-4\\
  3&-2.19E-6&-2.48E-4&-1.86E-3&-2.41E-4&-8.29E-4&-2.57E-3&-8.50E-5&-4.78E-4&-2.07E-3&-3.40E-5&-3.36E-4&-1.91E-3\\
  4&-7.68E-5&-2.88E-3&-1.38E-2&-2.08E-3&-6.42E-3&-1.71E-2&-9.65E-4&-4.43E-3&-1.48E-2&-4.99E-4&-3.51E-3&-1.41E-2\\
  5&-1.77E-3&-2.48E-2&-8.19E-2&-1.65E-2&-4.14E-2&-9.35E-2&-9.62E-3&-3.25E-2&-8.55E-2&-6.11E-3&-2.81E-2&-8.28E-2\\
  6&-3.18E-2&-1.77E-1&-4.09E-1&-1.28E-1&-2.43E-1&-4.40E-1&-9.06E-2&-2.09E-1&-4.19E-1&-6.80E-2&-1.91E-1&-4.12E-1\\
  7&-6.10E-1&-1.43E00&-1.73E00&-1.36E00&-1.75E00&-1.82E00&-1.10E00&-1.59E00&-1.76E00&-9.27E-1&-1.50E00&-1.74E00\\
\hline \hline
\end{tabular}
%}
\end{table}
\end{tiny}

~ ~ ~

\begin{tiny}
\begin{table}
\setlength{\tabcolsep}{2.0pt}
\renewcommand{\arraystretch}{1.02}
\caption{The poles of $r_{q,\delta,\alpha,k}(\t)$,  %${ r}(q,\delta,\alpha,k;\t)$,  
\eqref{BURA},   of $(q,\delta,\alpha,k)$-BURA approximation, 
i.e. the best uniform rational approximation of  
$\f(q,\delta,\alpha;\t)=\t^{\alpha}/(1+q\,\t^{\alpha})$  with functions from $\mathcal R_k$ on 
 $ [\delta, 1]$,   $\alpha =0.25, 0.5, 0.75$, ${\bf k=8}$, and $q=1$}
 %$(q=001,\delta,\alpha,k=8)$-BURA-poles}
 \label{tabl:Cq001dDaAAk8p}
%\rotatebox{90}{
\begin{tabular}{||c||c|c|c||c|c|c||c|c|c||c|c|c||}
\hline \hline
{\bf k=8} /$\delta$& $0.00$ & $0.00$ & $0.00$ &$10^{-6}$&$10^{-6}$&$10^{-6}$&$10^{-7}$&$10^{-7}$&$10^{-7}$&$10^{-8}$&$10^{-8}$&$10^{-8}$\\
\ j /$\alpha$ & 0.25   & 0.50   & 0.75   & 0.25   & 0.50   & 0.75   & 0.25   & 0.50   & 0.75   & 0.25   & 0.50   & 0.75   \\ \hline
  1&-9.64E-12&-5.02E-8&-1.88E-6&-1.25E-6&-3.06E-6&-8.44E-6&-2.09E-7&-7.40E-7&-3.54E-6&-3.64E-8&-2.31E-7&-2.28E-6\\
  2&-3.40E-9&-2.72E-6&-4.68E-5&-1.36E-5&-3.79E-5&-1.13E-4&-3.05E-6&-1.37E-5&-6.63E-5&-7.56E-7&-6.44E-6&-5.18E-5\\
  3&-2.44E-7&-5.21E-5&-5.14E-4&-1.07E-4&-3.18E-4&-9.17E-4&-3.23E-5&-1.53E-4&-6.42E-4&-1.09E-5&-9.15E-5&-5.48E-4\\
  4&-8.43E-6&-6.07E-4&-3.79E-3&-7.45E-4&-2.16E-3&-5.63E-3&-2.94E-4&-1.28E-3&-4.41E-3&-1.29E-4&-8.90E-4&-3.96E-3\\
  5&-1.89E-4&-5.20E-3&-2.21E-2&-4.85E-3&-1.25E-2&-2.89E-2&-2.41E-3&-8.65E-3&-2.44E-2&-1.31E-3&-6.75E-3&-2.27E-2\\
  6&-3.18E-3&-3.61E-2&-1.10E-1&-3.01E-2&-6.48E-2&-1.31E-1&-1.83E-2&-5.05E-2&-1.18E-1&-1.19E-2&-4.28E-2&-1.12E-1\\
  7&-4.54E-2&-2.23E-1&-4.83E-1&-1.92E-1&-3.24E-1&-5.33E-1&-1.37E-1&-2.75E-1&-5.01E-1&-1.04E-1&-2.48E-1&-4.88E-1\\
  8&-7.38E-1&-1.65E00&-1.93E00&-1.76E00&-2.11E00&-2.08E00&-1.42E00&-1.89E00&-1.99E00&-1.19E00&-1.77E00&-1.95E00\\
\hline \hline
\end{tabular}
%}
\end{table}
\end{tiny}

~ ~ ~

\begin{tiny}
\begin{table}
\setlength{\tabcolsep}{2.0pt}
\renewcommand{\arraystretch}{1.02}
\caption{The poles of $r_{q,\delta,\alpha,k}(\t)$,  %${ r}(q,\delta,\alpha,k;\t)$,  
\eqref{BURA},   of $(q,\delta,\alpha,k)$-BURA approximation, 
i.e. the best uniform rational approximation of  
$\f(q,\delta,\alpha;\t)=\t^{\alpha}/(1+q\,\t^{\alpha})$  with functions from $\mathcal R_k$ on 
 $ [\delta, 1]$,   $\alpha =0.25, 0.5, 0.75$, ${\bf k=3}$, and $q=100$}
 %$(q=100,\delta,\alpha,k=3)$-BURA-poles}
 \label{tabl:Cq100dDaAAk3p}
%\rotatebox{90}{
\begin{tabular}{||c||c|c|c||c|c|c||c|c|c||c|c|c||}
\hline \hline
{\bf k=3} /$\delta$& $0.00$ & $0.00$ & $0.00$ &$10^{-6}$&$10^{-6}$&$10^{-6}$&$10^{-7}$&$10^{-7}$&$10^{-7}$&$10^{-8}$&$10^{-8}$&$10^{-8}$\\
\ j /$\alpha$ & 0.25   & 0.50   & 0.75   & 0.25   & 0.50   & 0.75   & 0.25   & 0.50   & 0.75   & 0.25   & 0.50   & 0.75   \\ \hline
  1&-1.13E-11&-1.69E-6&-1.17E-4&-2.60E-6&-1.11E-5&-1.46E-4&-3.72E-7&-4.44E-6&-1.23E-4&-5.70E-8&-2.54E-6&-1.18E-4\\
  2&-8.24E-9&-8.15E-5&-1.81E-3&-1.36E-4&-2.44E-4&-1.95E-3&-2.44E-5&-1.36E-4&-1.84E-3&-4.62E-6&-1.00E-4&-1.81E-3\\
  3&-5.38E-6&-3.29E-3&-1.92E-2&-1.41E-2&-8.79E-3&-2.12E-2&-3.85E-3&-5.19E-3&-1.96E-2&-1.01E-3&-3.94E-3&-1.93E-2\\
\hline \hline
\end{tabular}
%}
\end{table}
\end{tiny}

~ ~ ~

\begin{tiny}
\begin{table}
\setlength{\tabcolsep}{2.0pt}
\renewcommand{\arraystretch}{1.02}
\caption{The poles of $r_{q,\delta,\alpha,k}(\t)$,  %${ r}(q,\delta,\alpha,k;\t)$, 
 \eqref{BURA},   of $(q,\delta,\alpha,k)$-BURA approximation, 
i.e. the best uniform rational approximation of  
$\f(q,\delta,\alpha;\t)=\t^{\alpha}/(1+q\,\t^{\alpha})$  with functions from $\mathcal R_k$ on 
 $ [\delta, 1]$,   $\alpha =0.25, 0.5, 0.75$, ${\bf k=4}$, and $q=100$}
 %$(q=100,\delta,\alpha,k=4)$-BURA-poles}
 \label{tabl:Cq100dDaAAk4p}
%\rotatebox{90}{
\begin{tabular}{||c||c|c|c||c|c|c||c|c|c||c|c|c||}
\hline \hline
{\bf k=4} /$\delta$& $0.00$ & $0.00$ & $0.00$ &$10^{-6}$&$10^{-6}$&$10^{-6}$&$10^{-7}$&$10^{-7}$&$10^{-7}$&$10^{-8}$&$10^{-8}$&$10^{-8}$\\
\ j /$\alpha$ & 0.25   & 0.50   & 0.75   & 0.25   & 0.50   & 0.75   & 0.25   & 0.50   & 0.75   & 0.25   & 0.50   & 0.75   \\ \hline
  1&-1.09E-12&-3.09E-7&-2.59E-5&-1.51E-6&-5.49E-6&-4.36E-5&-2.14E-7&-1.68E-6&-3.01E-5&-3.18E-8&-7.08E-7&-2.68E-5\\
  2&-5.53E-10&-1.57E-5&-5.91E-4&-3.97E-5&-8.48E-5&-7.44E-4&-6.76E-6&-3.92E-5&-6.31E-4&-1.21E-6&-2.40E-5&-6.00E-4\\
  3&-1.01E-7&-3.05E-4&-3.66E-3&-1.16E-3&-1.18E-3&-4.25E-3&-2.68E-4&-6.15E-4&-3.81E-3&-6.14E-5&-4.17E-4&-3.69E-3\\
  4&-4.08E-5&-1.07E-2&-4.71E-2&-6.07E-2&-3.40E-2&-5.57E-2&-2.27E-2&-1.98E-2&-4.93E-2&-7.80E-3&-1.41E-2&-4.76E-2\\
\hline \hline
\end{tabular}
%}
\end{table}
\end{tiny}

~ ~ ~

\begin{tiny}
\begin{table}
\setlength{\tabcolsep}{2.0pt}
\renewcommand{\arraystretch}{1.02}
\caption{The poles of ${ r}(q,\delta,\alpha,k;\t)$,  
\eqref{BURA},   of $r_{q,\delta,\alpha,k}(\t)$,  %$(q,\delta,\alpha,k)$-BURA approximation, 
i.e. the best uniform rational approximation of  
$\f(q,\delta,\alpha;\t)=\t^{\alpha}/(1+q\,\t^{\alpha})$  with functions from $\mathcal R_k$ on 
 $ [\delta, 1]$,   $\alpha =0.25, 0.5, 0.75$, ${\bf k=5}$, and $q=100$}
 %$(q=100,\delta,\alpha,k=5)$-BURA-poles}
 \label{tabl:Cq100dDaAAk5p}
%\rotatebox{90}{
\begin{tabular}{||c||c|c|c||c|c|c||c|c|c||c|c|c||}
\hline \hline
{\bf k=5} /$\delta$& $0.00$ & $0.00$ & $0.00$ &$10^{-6}$&$10^{-6}$&$10^{-6}$&$10^{-7}$&$10^{-7}$&$10^{-7}$&$10^{-8}$&$10^{-8}$&$10^{-8}$\\
\ j /$\alpha$ & 0.25   & 0.50   & 0.75   & 0.25   & 0.50   & 0.75   & 0.25   & 0.50   & 0.75   & 0.25   & 0.50   & 0.75   \\ \hline
  1&-1.36E-13&-6.38E-8&-6.47E-6&-1.00E-6&-3.21E-6&-1.70E-5&-1.41E-7&-8.09E-7&-9.09E-6&-2.05E-8&-2.63E-7&-7.06E-6\\
  2&-5.88E-11&-3.37E-6&-1.61E-4&-1.75E-5&-3.90E-5&-2.71E-4&-2.86E-6&-1.49E-5&-1.92E-4&-4.87E-7&-7.36E-6&-1.69E-4\\
  3&-6.77E-9&-5.97E-5&-1.36E-3&-2.68E-4&-3.54E-4&-1.76E-3&-5.62E-5&-1.66E-4&-1.48E-3&-1.19E-5&-1.00E-4&-1.39E-3\\
  4&-7.18E-7&-8.70E-4&-7.06E-3&-5.21E-3&-4.17E-3&-9.24E-3&-1.52E-3&-2.12E-3&-7.69E-3&-4.28E-4&-1.36E-3&-7.21E-3\\
  5&-2.13E-4&-2.64E-2&-9.23E-2&-1.56E-1&-8.97E-2&-1.17E-1&-7.23E-2&-5.40E-2&-9.97E-2&-3.13E-2&-3.80E-2&-9.41E-2\\
\hline \hline
\end{tabular}
%}
\end{table}
\end{tiny}

~ ~ ~

\begin{tiny}
\begin{table}
\setlength{\tabcolsep}{2.0pt}
\renewcommand{\arraystretch}{1.02}
\caption{The poles of $r_{q,\delta,\alpha,k}(\t)$,  %${ r}(q,\delta,\alpha,k;\t)$,  
\eqref{BURA},   of $(q,\delta,\alpha,k)$-BURA approximation, 
i.e. the best uniform rational approximation of  
$\f(q,\delta,\alpha;\t)=\t^{\alpha}/(1+q\,\t^{\alpha})$  with functions from $\mathcal R_k$ on 
 $ [\delta, 1]$,   $\alpha =0.25, 0.5, 0.75$, ${\bf k=6}$, and $q=100$}
 %$(q=100,\delta,\alpha,k=6)$-BURA-poles}
 \label{tabl:Cq100dDaAAk6p}
%\rotatebox{90}{
\begin{tabular}{||c||c|c|c||c|c|c||c|c|c||c|c|c||}
\hline \hline
{\bf k=6} /$\delta$& $0.00$ & $0.00$ & $0.00$ &$10^{-6}$&$10^{-6}$&$10^{-6}$&$10^{-7}$&$10^{-7}$&$10^{-7}$&$10^{-8}$&$10^{-8}$&$10^{-8}$\\
\ j /$\alpha$ & 0.25   & 0.50   & 0.75   & 0.25   & 0.50   & 0.75   & 0.25   & 0.50   & 0.75   & 0.25   & 0.50   & 0.75   \\ \hline
  1&-2.00E-14&-1.44E-8&-1.78E-6&-7.16E-7&-2.10E-6&-8.23E-6&-1.01E-7&-4.61E-7&-3.41E-6&-1.45E-8&-1.23E-7&-2.17E-6\\
  2&-7.94E-12&-7.75E-7&-4.47E-5&-9.60E-6&-2.11E-5&-1.10E-4&-1.52E-6&-6.78E-6&-6.39E-5&-2.49E-7&-2.74E-6&-4.96E-5\\
  3&-7.40E-10&-1.42E-5&-4.68E-4&-9.83E-5&-1.50E-4&-7.97E-4&-1.92E-5&-6.25E-5&-5.78E-4&-3.83E-6&-3.27E-5&-4.98E-4\\
  4&-4.60E-8&-1.60E-4&-2.35E-3&-1.12E-3&-1.13E-3&-3.30E-3&-2.84E-4&-5.24E-4&-2.67E-3&-7.14E-5&-3.06E-4&-2.44E-3\\
  5&-3.61E-6&-2.06E-3&-1.29E-2&-1.54E-2&-1.13E-2&-1.90E-2&-5.45E-3&-5.86E-3&-1.49E-2&-1.84E-3&-3.65E-3&-1.35E-2\\
  6&-8.23E-4&-5.27E-2&-1.54E-1&-3.01E-1&-1.83E-1&-2.08E-1&-1.60E-1&-1.15E-1&-1.73E-1&-8.17E-2&-8.14E-2&-1.59E-1\\
\hline \hline
\end{tabular}
%}
\end{table}
\end{tiny}

~ ~ ~

\begin{tiny}
\begin{table}
\setlength{\tabcolsep}{2.0pt}
\renewcommand{\arraystretch}{1.02}
\caption{The poles of $r_{q,\delta,\alpha,k}(\t)$,  %${ r}(q,\delta,\alpha,k;\t)$,  
\eqref{BURA},   of $(q,\delta,\alpha,k)$-BURA approximation, 
i.e. the best uniform rational approximation of  
$\f(q,\delta,\alpha;\t)=\t^{\alpha}/(1+q\,\t^{\alpha})$  with functions from $\mathcal R_k$ on 
 $ [\delta, 1]$,   $\alpha =0.25, 0.5, 0.75$, ${\bf k=7}$, and $q=100$}
 %$(q=100,\delta,\alpha,k=7)$-BURA-poles}
\label{tabl:Cq100dDaAAk7p}
%\rotatebox{90}{
\begin{tabular}{||c||c|c|c||c|c|c||c|c|c||c|c|c||}
\hline \hline
{\bf k=7} /$\delta$& $0.00$ & $0.00$ & $0.00$ &$10^{-6}$&$10^{-6}$&$10^{-6}$&$10^{-7}$&$10^{-7}$&$10^{-7}$&$10^{-8}$&$10^{-8}$&$10^{-8}$\\
\ j /$\alpha$ & 0.25   & 0.50   & 0.75   & 0.25   & 0.50   & 0.75   & 0.25   & 0.50   & 0.75   & 0.25   & 0.50   & 0.75   \\ \hline
  1&-3.29E-15&-3.53E-9&-5.30E-7&-5.41E-7&-1.48E-6&-4.68E-6&-7.58E-8&-2.95E-7&-1.56E-6&-1.08E-8&-6.78E-8&-7.86E-7\\
  2&-1.24E-12&-1.91E-7&-1.32E-5&-6.04E-6&-1.28E-5&-5.16E-5&-9.28E-7&-3.55E-6&-2.48E-5&-1.47E-7&-1.20E-6&-1.65E-5\\
  3&-1.04E-10&-3.59E-6&-1.45E-4&-4.69E-5&-7.71E-5&-3.61E-4&-8.65E-6&-2.83E-5&-2.19E-4&-1.64E-6&-1.27E-5&-1.67E-4\\
  4&-4.98E-9&-3.95E-5&-9.42E-4&-3.76E-4&-4.43E-4&-1.62E-3&-8.67E-5&-1.89E-4&-1.20E-3&-2.01E-5&-1.00E-4&-1.02E-3\\
  5&-2.20E-7&-3.56E-4&-3.73E-3&-3.34E-3&-2.94E-3&-6.03E-3&-1.00E-3&-1.38E-3&-4.56E-3&-2.94E-4&-7.87E-4&-3.99E-3\\
  6&-1.40E-5&-4.18E-3&-2.17E-2&-3.45E-2&-2.48E-2&-3.51E-2&-1.42E-2&-1.34E-2&-2.67E-2&-5.60E-3&-8.32E-3&-2.33E-2\\
  7&-2.48E-3&-9.00E-2&-2.30E-1&-4.93E-1&-3.13E-1&-3.29E-1&-2.87E-1&-2.05E-1&-2.69E-1&-1.63E-1&-1.47E-1&-2.42E-1\\
\hline \hline
\end{tabular}
%}
\end{table}
\end{tiny}

~ ~ ~

\begin{tiny}
\begin{table}
\setlength{\tabcolsep}{2.0pt}
\renewcommand{\arraystretch}{1.02}
\caption{The poles of $r_{q,\delta,\alpha,k}(\t)$,  %${ r}(q,\delta,\alpha,k;\t)$,  
\eqref{BURA},   of $(q,\delta,\alpha,k)$-BURA approximation, 
i.e. the best uniform rational approximation of  
$\f(q,\delta,\alpha;\t)=\t^{\alpha}/(1+q\,\t^{\alpha})$  with functions from $\mathcal R_k$ on 
 $ [\delta, 1]$,   $\alpha =0.25, 0.5, 0.75$, ${\bf k=8}$, and $q=100$}
 %$(q=100,\delta,\alpha,k=8)$-BURA-poles}
 \label{tabl:Cq100dDaAAk8p}
%\rotatebox{90}{
\begin{tabular}{||c||c|c|c||c|c|c||c|c|c||c|c|c||}
\hline \hline
{\bf k=8} /$\delta$& $0.00$ & $0.00$ & $0.00$ &$10^{-6}$&$10^{-6}$&$10^{-6}$&$10^{-7}$&$10^{-7}$&$10^{-7}$&$10^{-8}$&$10^{-8}$&$10^{-8}$\\
\ j /$\alpha$ & 0.25   & 0.50   & 0.75   & 0.25   & 0.50   & 0.75   & 0.25   & 0.50   & 0.75   & 0.25   & 0.50   & 0.75   \\ \hline
  1&-7.52E-16&-9.23E-10&-1.68E-7&-4.25E-7&-1.10E-6&-2.98E-6&-5.94E-8&-2.05E-7&-8.32E-7&-8.39E-9&-4.22E-8&-3.35E-7\\
  2& 3.96E-14&-4.99E-8&-4.19E-6&-4.15E-6&-8.40E-6&-2.78E-5&-6.24E-7&-2.08E-6&-1.11E-5&-9.63E-8&-6.07E-7&-6.22E-6\\
  3& 1.25E-11&-9.50E-7&-4.62E-5&-2.64E-5&-4.46E-5&-1.75E-4&-4.65E-6&-1.45E-5&-9.03E-5&-8.40E-7&-5.63E-6&-6.02E-5\\
  4&-3.72E-9&-1.07E-5&-3.33E-4&-1.65E-4&-2.15E-4&-8.29E-4&-3.53E-5&-8.39E-5&-5.27E-4&-7.64E-6&-3.99E-5&-4.00E-4\\
  5& 3.72E-9&-8.78E-5&-1.54E-3&-1.09E-3&-1.10E-3&-2.78E-3&-2.92E-4&-4.76E-4&-2.05E-3&-7.74E-5&-2.51E-4&-1.72E-3\\
  6& 3.35E-6&-7.02E-4&-5.75E-3&-7.88E-3&-6.51E-3&-1.07E-2&-2.72E-3&-3.13E-3&-7.69E-3&-9.14E-4&-1.78E-3&-6.41E-3\\
  7& 2.46E-5&-7.55E-3&-3.36E-2&-6.42E-2&-4.66E-2&-5.90E-2&-2.97E-2&-2.61E-2&-4.40E-2&-1.32E-2&-1.65E-2&-3.73E-2\\
  8& 1.03E-2&-1.38E-1&-3.18E-1&-7.30E-1&-4.80E-1&-4.79E-1&-4.50E-1&-3.24E-1&-3.87E-1&-2.75E-1&-2.37E-1&-3.43E-1\\
\hline \hline
\end{tabular}
%}
\end{table}
\end{tiny}

~ ~ ~

\begin{tiny}
\begin{table}
\setlength{\tabcolsep}{2.0pt}
\renewcommand{\arraystretch}{1.02}
\caption{The poles of $r_{q,\delta,\alpha,k}(\t)$,  %${ r}(q,\delta,\alpha,k;\t)$,  
\eqref{BURA},   of $(q,\delta,\alpha,k)$-BURA approximation, 
i.e. the best uniform rational approximation of  
$\f(q,\delta,\alpha;\t)=\t^{\alpha}/(1+q\,\t^{\alpha})$  with functions from $\mathcal R_k$ on 
 $ [\delta, 1]$,   $\alpha =0.25, 0.5, 0.75$, ${\bf k=3}$, and $q=200$}
 %$(q=200,\delta,\alpha,k=3)$-BURA-poles}
 \label{tabl:Cq200dDaAAk3p}
%\rotatebox{90}{
\begin{tabular}{||c||c|c|c||c|c|c||c|c|c||c|c|c||}
\hline \hline
{\bf k=3} /$\delta$& $0.00$ & $0.00$ & $0.00$ &$10^{-6}$&$10^{-6}$&$10^{-6}$&$10^{-7}$&$10^{-7}$&$10^{-7}$&$10^{-8}$&$10^{-8}$&$10^{-8}$\\
\ j /$\alpha$ & 0.25   & 0.50   & 0.75   & 0.25   & 0.50   & 0.75   & 0.25   & 0.50   & 0.75   & 0.25   & 0.50   & 0.75   \\ \hline
  1&-7.62E-13&-4.65E-7&-5.18E-5&-2.21E-6&-5.95E-6&-7.36E-5&-2.87E-7&-1.99E-6&-5.68E-5&-3.87E-8&-9.21E-7&-5.28E-5\\
  2&-5.67E-10&-2.24E-5&-7.69E-4&-1.17E-4&-1.16E-4&-8.72E-4&-1.85E-5&-5.15E-5&-7.94E-4&-2.92E-6&-3.22E-5&-7.74E-4\\
  3&-3.98E-7&-9.78E-4&-9.03E-3&-1.24E-2&-4.63E-3&-1.08E-2&-3.00E-3&-2.15E-3&-9.44E-3&-6.48E-4&-1.37E-3&-9.12E-3\\
\hline \hline
\end{tabular}
%}
\end{table}
\end{tiny}

~ ~ ~

\begin{tiny}
\begin{table}
\setlength{\tabcolsep}{2.0pt}
\renewcommand{\arraystretch}{1.02}
\caption{The poles of $r_{q,\delta,\alpha,k}(\t)$,  %${ r}(q,\delta,\alpha,k;\t)$,  
\eqref{BURA},    %$(q,\delta,\alpha,k)$-BURA approximation, 
i.e. the best uniform rational approximation of  
$\f(q,\delta,\alpha;\t)=\t^{\alpha}/(1+q\,\t^{\alpha})$  with functions from $\mathcal R_k$ on 
 $ [\delta, 1]$,   $\alpha =0.25, 0.5, 0.75$, ${\bf k=4}$, and $q=200$}
 %$(q=200,\delta,\alpha,k=4)$-BURA-poles}
 \label{tabl:Cq200dDaAAk4p}
%\rotatebox{90}{
\begin{tabular}{||c||c|c|c||c|c|c||c|c|c||c|c|c||}
\hline \hline
{\bf k=4} /$\delta$& $0.00$ & $0.00$ & $0.00$ &$10^{-6}$&$10^{-6}$&$10^{-6}$&$10^{-7}$&$10^{-7}$&$10^{-7}$&$10^{-8}$&$10^{-8}$&$10^{-8}$\\
\ j /$\alpha$ & 0.25   & 0.50   & 0.75   & 0.25   & 0.50   & 0.75   & 0.25   & 0.50   & 0.75   & 0.25   & 0.50   & 0.75   \\ \hline
  1&-7.64E-14&-9.00E-8&-1.21E-5&-1.31E-6&-3.41E-6&-2.54E-5&-1.72E-7&-9.12E-7&-1.53E-5&-2.30E-8&-3.17E-7&-1.28E-5\\
  2&-3.91E-11&-4.53E-6&-2.69E-4&-3.51E-5&-4.49E-5&-3.74E-4&-5.42E-6&-1.72E-5&-2.98E-4&-8.33E-7&-8.90E-6&-2.76E-4\\
  3&-7.43E-9&-9.08E-5&-1.65E-3&-1.04E-3&-6.54E-4&-2.11E-3&-2.16E-4&-2.71E-4&-1.78E-3&-4.20E-5&-1.54E-4&-1.68E-3\\
  4&-3.37E-6&-3.68E-3&-2.49E-2&-5.61E-2&-2.14E-2&-3.28E-2&-1.92E-2&-1.00E-2&-2.71E-2&-5.71E-3&-6.03E-3&-2.54E-2\\
\hline \hline
\end{tabular}
%}
\end{table}
\end{tiny}

~ ~ ~

\begin{tiny}
\begin{table}
\setlength{\tabcolsep}{2.0pt}
\renewcommand{\arraystretch}{1.02}
\caption{The poles of $r_{q,\delta,\alpha,k}(\t)$,  %${ r}(q,\delta,\alpha,k;\t)$, 
 \eqref{BURA},   of $(q,\delta,\alpha,k)$-BURA approximation, 
i.e. the best uniform rational approximation of  
$\f(q,\delta,\alpha;\t)=\t^{\alpha}/(1+q\,\t^{\alpha})$  with functions from $\mathcal R_k$ on 
 $ [\delta, 1]$,   $\alpha =0.25, 0.5, 0.75$, ${\bf k=5}$, and $q=200$}
 %$(q=200,\delta,\alpha,k=5)$-BURA-poles}
 \label{tabl:Cq200dDaAAk5p}
%\rotatebox{90}{
\begin{tabular}{||c||c|c|c||c|c|c||c|c|c||c|c|c||}
\hline \hline
{\bf k=5} /$\delta$& $0.00$ & $0.00$ & $0.00$ &$10^{-6}$&$10^{-6}$&$10^{-6}$&$10^{-7}$&$10^{-7}$&$10^{-7}$&$10^{-8}$&$10^{-8}$&$10^{-8}$\\
\ j /$\alpha$ & 0.25   & 0.50   & 0.75   & 0.25   & 0.50   & 0.75   & 0.25   & 0.50   & 0.75   & 0.25   & 0.50   & 0.75   \\ \hline
  1&-9.97E-15&-1.97E-8&-3.15E-6&-8.77E-7&-2.20E-6&-1.12E-5&-1.16E-7&-5.07E-7&-5.18E-6&-1.55E-8&-1.41E-7&-3.63E-6\\
  2&-4.33E-12&-1.04E-6&-7.83E-5&-1.58E-5&-2.28E-5&-1.56E-4&-2.38E-6&-7.56E-6&-1.01E-4&-3.58E-7&-3.21E-6&-8.41E-5\\
  3&-5.11E-10&-1.83E-5&-6.27E-4&-2.43E-4&-2.06E-4&-8.91E-4&-4.67E-5&-7.88E-5&-7.11E-4&-8.63E-6&-4.03E-5&-6.48E-4\\
  4&-5.77E-8&-2.87E-4&-3.50E-3&-4.79E-3&-2.64E-3&-5.30E-3&-1.29E-3&-1.10E-3&-4.04E-3&-3.18E-4&-5.89E-4&-3.64E-3\\
  5&-2.05E-5&-1.05E-2&-5.40E-2&-1.48E-1&-6.45E-2&-7.77E-2&-6.47E-2&-3.27E-2&-6.15E-2&-2.53E-2&-1.96E-2&-5.60E-2\\
\hline \hline
\end{tabular}
%}
\end{table}
\end{tiny}

~ ~ ~

\begin{tiny}
\begin{table}
\setlength{\tabcolsep}{2.0pt}
\renewcommand{\arraystretch}{1.02}
\caption{The poles of $r_{q,\delta,\alpha,k}(\t)$,  %${ r}(q,\delta,\alpha,k;\t)$, 
 \eqref{BURA},   of $(q,\delta,\alpha,k)$-BURA approximation, 
i.e. the best uniform rational approximation of  
$\f(q,\delta,\alpha;\t)=\t^{\alpha}/(1+q\,\t^{\alpha})$  with functions from $\mathcal R_k$ on 
 $ [\delta, 1]$,   $\alpha =0.25, 0.5, 0.75$, ${\bf k=6}$, and $q=200$}
 %$(q=200,\delta,\alpha,k=6)$-BURA-poles}
 \label{tabl:Cq200dDaAAk6p}
%\rotatebox{90}{
\begin{tabular}{||c||c|c|c||c|c|c||c|c|c||c|c|c||}
\hline \hline
{\bf k=6} /$\delta$& $0.00$ & $0.00$ & $0.00$ &$10^{-6}$&$10^{-6}$&$10^{-6}$&$10^{-7}$&$10^{-7}$&$10^{-7}$&$10^{-8}$&$10^{-8}$&$10^{-8}$\\
\ j /$\alpha$ & 0.25   & 0.50   & 0.75   & 0.25   & 0.50   & 0.75   & 0.25   & 0.50   & 0.75   & 0.25   & 0.50   & 0.75   \\ \hline
  1&-1.54E-15&-4.71E-9&-9.03E-7&-6.32E-7&-1.54E-6&-5.94E-6&-8.37E-8&-3.19E-7&-2.16E-6&-1.12E-8&-7.57E-8&-1.21E-6\\
  2&-6.16E-13&-2.53E-7&-2.27E-5&-8.77E-6&-1.35E-5&-7.07E-5&-1.29E-6&-3.90E-6&-3.71E-5&-1.91E-7&-1.38E-6&-2.66E-5\\
  3&-5.84E-11&-4.59E-6&-2.32E-4&-9.01E-5&-9.19E-5&-4.49E-4&-1.63E-5&-3.23E-5&-3.09E-4&-2.90E-6&-1.47E-5&-2.55E-4\\
  4&-3.79E-9&-5.29E-5&-1.12E-3&-1.03E-3&-7.30E-4&-1.83E-3&-2.44E-4&-2.79E-4&-1.36E-3&-5.46E-5&-1.38E-4&-1.19E-3\\
  5&-3.27E-7&-7.57E-4&-7.00E-3&-1.44E-2&-7.94E-3&-1.22E-2&-4.81E-3&-3.46E-3&-8.77E-3&-1.47E-3&-1.84E-3&-7.51E-3\\
  6&-9.57E-5&-2.42E-2&-9.77E-2&-2.89E-1&-1.43E-1&-1.51E-1&-1.48E-1&-7.87E-2&-1.17E-1&-7.01E-2&-4.89E-2&-1.03E-1\\
\hline \hline
\end{tabular}
%}
\end{table}
\end{tiny}

~ ~ ~

\begin{tiny}
\begin{table}
\setlength{\tabcolsep}{2.0pt}
\renewcommand{\arraystretch}{1.02}
\caption{The poles of $r_{q,\delta,\alpha,k}(\t)$,  %${ r}(q,\delta,\alpha,k;\t)$, 
 \eqref{BURA},   of $(q,\delta,\alpha,k)$-BURA approximation, 
i.e. the best uniform rational approximation of  
$\f(q,\delta,\alpha;\t)=\t^{\alpha}/(1+q\,\t^{\alpha})$  with functions from $\mathcal R_k$ on 
 $ [\delta, 1]$,   $\alpha =0.25, 0.5, 0.75$, ${\bf k=7}$, and $q=200$}
 %$(q=200,\delta,\alpha,k=7)$-BURA-poles}
 \label{tabl:Cq200dDaAAk7p}
%\rotatebox{90}{
\begin{tabular}{||c||c|c|c||c|c|c||c|c|c||c|c|c||}
\hline \hline
{\bf k=7} /$\delta$& $0.00$ & $0.00$ & $0.00$ &$10^{-6}$&$10^{-6}$&$10^{-6}$&$10^{-7}$&$10^{-7}$&$10^{-7}$&$10^{-8}$&$10^{-8}$&$10^{-8}$\\
\ j /$\alpha$ & 0.25   & 0.50   & 0.75   & 0.25   & 0.50   & 0.75   & 0.25   & 0.50   & 0.75   & 0.25   & 0.50   & 0.75   \\ \hline
  1&-9.85E-15&-1.21E-9&-2.78E-7&-4.79E-7&-1.14E-6&-3.61E-6&-6.37E-8&-2.18E-7&-1.08E-6&-8.54E-9&-4.61E-8&-4.80E-7\\
  2&-4.28E-12&-6.55E-8&-6.96E-6&-5.57E-6&-8.79E-6&-3.61E-5&-8.07E-7&-2.26E-6&-1.57E-5&-1.17E-7&-6.85E-7&-9.46E-6\\
  3&-5.07E-10&-1.23E-6&-7.61E-5&-4.34E-5&-4.96E-5&-2.28E-4&-7.50E-6&-1.59E-5&-1.30E-4&-1.28E-6&-6.38E-6&-9.29E-5\\
  4& 1.65E-9&-1.35E-5&-4.68E-4&-3.49E-4&-2.92E-4&-9.06E-4&-7.56E-5&-1.04E-4&-6.39E-4&-1.57E-5&-4.78E-5&-5.26E-4\\
  5&-5.79E-8&-1.28E-4&-1.89E-3&-3.14E-3&-2.07E-3&-3.74E-3&-8.89E-4&-8.18E-4&-2.56E-3&-2.36E-4&-3.99E-4&-2.10E-3\\
  6&-2.05E-5&-1.71E-3&-1.27E-2&-3.28E-2&-1.88E-2&-2.45E-2&-1.29E-2&-8.74E-3&-1.72E-2&-4.68E-3&-4.73E-3&-1.41E-2\\
  7& 1.51E-3&-4.63E-2&-1.55E-1&-4.78E-1&-2.60E-1&-2.54E-1&-2.70E-1&-1.53E-1&-1.95E-1&-1.45E-1&-9.90E-2&-1.69E-1\\
\hline \hline
\end{tabular}
%}
\end{table}
\end{tiny}

~ ~ ~

\begin{tiny}
\begin{table}
\setlength{\tabcolsep}{2.0pt}
\renewcommand{\arraystretch}{1.02}
\caption{The poles of $r_{q,\delta,\alpha,k}(\t)$,  %${ r}(q,\delta,\alpha,k;\t)$, 
 \eqref{BURA},   of $(q,\delta,\alpha,k)$-BURA approximation, 
i.e. the best uniform rational approximation of  
$\f(q,\delta,\alpha;\t)=\t^{\alpha}/(1+q\,\t^{\alpha})$  with functions from $\mathcal R_k$ on 
 $ [\delta, 1]$,   $\alpha =0.25, 0.5, 0.75$, ${\bf k=8}$, and $q=200$}
 %$(q=200,\delta,\alpha,k=8)$-BURA-poles}
 \label{tabl:Cq200dDaAAk8p}
%\rotatebox{90}{
\begin{tabular}{||c||c|c|c||c|c|c||c|c|c||c|c|c||}
\hline \hline
{\bf k=8} /$\delta$& $0.00$ & $0.00$ & $0.00$ &$10^{-6}$&$10^{-6}$&$10^{-6}$&$10^{-7}$&$10^{-7}$&$10^{-7}$&$10^{-8}$&$10^{-8}$&$10^{-8}$\\
\ j /$\alpha$ & 0.25   & 0.50   & 0.75   & 0.25   & 0.50   & 0.75   & 0.25   & 0.50   & 0.75   & 0.25   & 0.50   & 0.75   \\ \hline
  1& 6.67E-15&-3.32E-10&-9.09E-8&-3.77E-7&-8.79E-7&-2.41E-6&-5.03E-8&-1.59E-7&-6.21E-7&-6.75E-9&-3.08E-8&-2.23E-7\\
  2& 2.42E-13&-1.79E-8&-2.27E-6&-3.86E-6&-6.13E-6&-2.07E-5&-5.50E-7&-1.43E-6&-7.57E-6&-7.86E-8&-3.81E-7&-3.83E-6\\
  3& 2.39E-9&-3.41E-7&-2.50E-5&-2.46E-5&-3.02E-5&-1.21E-4&-4.09E-6&-8.86E-6&-5.78E-5&-6.76E-7&-3.12E-6&-3.57E-5\\
  4& 1.00E-8&-3.82E-6&-1.78E-4&-1.54E-4&-1.45E-4&-5.07E-4&-3.11E-5&-4.87E-5&-3.12E-4&-6.13E-6&-2.04E-5&-2.26E-4\\
  5& 4.96E-8&-3.16E-5&-7.72E-4&-1.03E-3&-7.80E-4&-1.66E-3&-2.60E-4&-2.86E-4&-1.12E-3&-6.29E-5&-1.30E-4&-8.99E-4\\
  6& 5.10E-7&-2.74E-4&-3.13E-3&-7.49E-3&-4.88E-3&-7.28E-3&-2.47E-3&-2.03E-3&-4.75E-3&-7.62E-4&-9.99E-4&-3.70E-3\\
  7& 3.81E-5&-3.40E-3&-2.08E-2&-6.18E-2&-3.73E-2&-4.38E-2&-2.75E-2&-1.85E-2&-3.03E-2&-1.15E-2&-1.03E-2&-2.43E-2\\
  8& 9.99E-3&-7.77E-2&-2.25E-1&-7.11E-1&-4.13E-1&-3.87E-1&-4.29E-1&-2.57E-1&-2.97E-1&-2.52E-1&-1.72E-1&-2.52E-1\\
\hline \hline
\end{tabular}
%}
\end{table}
\end{tiny}

~ ~ ~

\begin{tiny}
\begin{table}
\setlength{\tabcolsep}{2.0pt}
\renewcommand{\arraystretch}{1.02}
\caption{The poles of $r_{q,\delta,\alpha,k}(\t)$,  %${ r}(q,\delta,\alpha,k;\t)$, 
 \eqref{BURA},   of $(q,\delta,\alpha,k)$-BURA approximation, 
i.e. the best uniform rational approximation of  
$\f(q,\delta,\alpha;\t)=\t^{\alpha}/(1+q\,\t^{\alpha})$  with functions from $\mathcal R_k$ on 
 $ [\delta, 1]$,   $\alpha =0.25, 0.5, 0.75$, ${\bf k=3}$, and $q=400$}
 %$(q=400,\delta,\alpha,k=3)$-BURA-poles}
 \label{tabl:Cq400dDaAAk3p}
%\rotatebox{90}{
\begin{tabular}{||c||c|c|c||c|c|c||c|c|c||c|c|c||}
\hline \hline
{\bf k=3} /$\delta$& $0.00$ & $0.00$ & $0.00$ &$10^{-6}$&$10^{-6}$&$10^{-6}$&$10^{-7}$&$10^{-7}$&$10^{-7}$&$10^{-8}$&$10^{-8}$&$10^{-8}$\\
\ j /$\alpha$ & 0.25   & 0.50   & 0.75   & 0.25   & 0.50   & 0.75   & 0.25   & 0.50   & 0.75   & 0.25   & 0.50   & 0.75   \\ \hline
  1&-4.95E-14&-1.22E-7&-2.20E-5&-2.02E-6&-3.37E-6&-3.77E-5&-2.46E-7&-9.68E-7&-2.58E-5&-3.02E-8&-3.65E-7&-2.28E-5\\
  2&-3.72E-11&-5.90E-6&-3.19E-4&-1.08E-4&-6.26E-5&-3.92E-4&-1.58E-5&-2.14E-5&-3.38E-4&-2.22E-6&-1.09E-5&-3.23E-4\\
  3&-2.71E-8&-2.70E-4&-4.02E-3&-1.16E-2&-2.72E-3&-5.41E-3&-2.60E-3&-9.55E-4&-4.37E-3&-4.94E-4&-4.89E-4&-4.10E-3\\
\hline \hline
\end{tabular}
%}
\end{table}
\end{tiny}

~ ~ ~

\begin{tiny}
\begin{table}
\setlength{\tabcolsep}{2.0pt}
\renewcommand{\arraystretch}{1.02}
\caption{The poles of $r_{q,\delta,\alpha,k}(\t)$,  %${ r}(q,\delta,\alpha,k;\t)$, 
 \eqref{BURA},   of $(q,\delta,\alpha,k)$-BURA approximation, 
i.e. the best uniform rational approximation of  
$\f(q,\delta,\alpha;\t)=\t^{\alpha}/(1+q\,\t^{\alpha})$  with functions from $\mathcal R_k$ on 
 $ [\delta, 1]$,   $\alpha =0.25, 0.5, 0.75$, ${\bf k=4}$, and $q=400$}
 %$(q=400,\delta,\alpha,k=4)$-BURA-poles}
 \label{tabl:Cq400dDaAAk4p}
%\rotatebox{90}{
\begin{tabular}{||c||c|c|c||c|c|c||c|c|c||c|c|c||}
\hline \hline
{\bf k=4} /$\delta$& $0.00$ & $0.00$ & $0.00$ &$10^{-6}$&$10^{-6}$&$10^{-6}$&$10^{-7}$&$10^{-7}$&$10^{-7}$&$10^{-8}$&$10^{-8}$&$10^{-8}$\\
\ j /$\alpha$ & 0.25   & 0.50   & 0.75   & 0.25   & 0.50   & 0.75   & 0.25   & 0.50   & 0.75   & 0.25   & 0.50   & 0.75   \\ \hline
  1&-5.06E-15&-2.46E-8&-5.36E-6&-1.21E-6&-2.15E-6&-1.51E-5&-1.50E-7&-5.23E-7&-7.80E-6&-1.85E-8&-1.53E-7&-5.92E-6\\
  2&-2.61E-12&-1.24E-6&-1.17E-4&-3.29E-5&-2.62E-5&-1.86E-4&-4.78E-6&-8.07E-6&-1.38E-4&-6.66E-7&-3.47E-6&-1.22E-4\\
  3&-5.04E-10&-2.52E-5&-7.21E-4&-9.75E-4&-4.05E-4&-1.07E-3&-1.91E-4&-1.30E-4&-8.17E-4&-3.34E-5&-5.94E-5&-7.44E-4\\
  4&-2.44E-7&-1.12E-3&-1.23E-2&-5.37E-2&-1.45E-2&-1.90E-2&-1.75E-2&-5.33E-3&-1.42E-2&-4.71E-3&-2.57E-3&-1.27E-2\\
\hline \hline
\end{tabular}
%}
\end{table}
\end{tiny}

~ ~ ~

\begin{tiny}
\begin{table}
\setlength{\tabcolsep}{2.0pt}
\renewcommand{\arraystretch}{1.02}
\caption{The poles of $r_{q,\delta,\alpha,k}(\t)$,  %${ r}(q,\delta,\alpha,k;\t)$, 
 \eqref{BURA},   of $(q,\delta,\alpha,k)$-BURA approximation, 
i.e. the best uniform rational approximation of  
$\f(q,\delta,\alpha;\t)=\t^{\alpha}/(1+q\,\t^{\alpha})$  with functions from $\mathcal R_k$ on 
 $ [\delta, 1]$,   $\alpha =0.25, 0.5, 0.75$, ${\bf k=5}$, and $q=400$}
% $(q=400,\delta,\alpha,k=5)$-BURA-poles}
 \label{tabl:Cq400dDaAAk5p}
%\rotatebox{90}{
\begin{tabular}{||c||c|c|c||c|c|c||c|c|c||c|c|c||}
\hline \hline
{\bf k=5} /$\delta$& $0.00$ & $0.00$ & $0.00$ &$10^{-6}$&$10^{-6}$&$10^{-6}$&$10^{-7}$&$10^{-7}$&$10^{-7}$&$10^{-8}$&$10^{-8}$&$10^{-8}$\\
\ j /$\alpha$ & 0.25   & 0.50   & 0.75   & 0.25   & 0.50   & 0.75   & 0.25   & 0.50   & 0.75   & 0.25   & 0.50   & 0.75   \\ \hline
  1&-6.76E-16&-5.63E-9&-1.46E-6&-8.11E-7&-1.50E-6&-7.48E-6&-1.02E-7&-3.26E-7&-2.98E-6&-1.27E-8&-8.01E-8&-1.83E-6\\
  2&-2.95E-13&-2.96E-7&-3.61E-5&-1.49E-5&-1.43E-5&-8.95E-5&-2.14E-6&-4.00E-6&-5.26E-5&-2.97E-7&-1.45E-6&-4.05E-5\\
  3&-3.53E-11&-5.24E-6&-2.77E-4&-2.30E-4&-1.33E-4&-4.54E-4&-4.20E-5&-4.03E-5&-3.35E-4&-7.10E-6&-1.69E-5&-2.93E-4\\
  4&-4.12E-9&-8.64E-5&-1.66E-3&-4.58E-3&-1.83E-3&-3.08E-3&-1.18E-3&-6.04E-4&-2.10E-3&-2.66E-4&-2.63E-4&-1.78E-3\\
  5&-1.62E-6&-3.66E-3&-2.94E-2&-1.43E-1&-4.89E-2&-5.10E-2&-6.07E-2&-2.03E-2&-3.65E-2&-2.21E-2&-1.00E-2&-3.14E-2\\
\hline \hline
\end{tabular}
%}
\end{table}
\end{tiny}

~ ~ ~

\begin{tiny}
\begin{table}
\setlength{\tabcolsep}{2.0pt}
\renewcommand{\arraystretch}{1.02}
\caption{The poles of $r_{q,\delta,\alpha,k}(\t)$,  %${ r}(q,\delta,\alpha,k;\t)$, 
 \eqref{BURA},   of $(q,\delta,\alpha,k)$-BURA approximation, 
i.e. the best uniform rational approximation of  
$\f(q,\delta,\alpha;\t)=\t^{\alpha}/(1+q\,\t^{\alpha})$  with functions from $\mathcal R_k$ on 
 $ [\delta, 1]$,   $\alpha =0.25, 0.5, 0.75$, ${\bf k=6}$, and $q=400$}
 %$(q=400,\delta,\alpha,k=6)$-BURA-poles}
 \label{tabl:Cq400dDaAAk6p}
%\rotatebox{90}{
\begin{tabular}{||c||c|c|c||c|c|c||c|c|c||c|c|c||}
\hline \hline
{\bf k=6} /$\delta$& $0.00$ & $0.00$ & $0.00$ &$10^{-6}$&$10^{-6}$&$10^{-6}$&$10^{-7}$&$10^{-7}$&$10^{-7}$&$10^{-8}$&$10^{-8}$&$10^{-8}$\\
\ j /$\alpha$ & 0.25   & 0.50   & 0.75   & 0.25   & 0.50   & 0.75   & 0.25   & 0.50   & 0.75   & 0.25   & 0.50   & 0.75   \\ \hline
  1&-1.08E-16&-1.42E-9&-4.34E-7&-5.86E-7&-1.12E-6&-4.33E-6&-7.42E-8&-2.23E-7&-1.39E-6&-9.35E-9&-4.82E-8&-6.73E-7\\
  2&-4.33E-14&-7.59E-8&-1.09E-5&-8.35E-6&-8.95E-6&-4.53E-5&-1.18E-6&-2.29E-6&-2.14E-5&-1.63E-7&-7.17E-7&-1.39E-5\\
  3&-4.14E-12&-1.37E-6&-1.09E-4&-8.60E-5&-6.14E-5&-2.47E-4&-1.49E-5&-1.76E-5&-1.60E-4&-2.45E-6&-6.77E-6&-1.25E-4\\
  4&-2.75E-10&-1.61E-5&-5.15E-4&-9.90E-4&-5.15E-4&-1.04E-3&-2.24E-4&-1.60E-4&-6.91E-4&-4.64E-5&-6.48E-5&-5.68E-4\\
  5&-2.51E-8&-2.51E-4&-3.60E-3&-1.40E-2&-5.95E-3&-7.87E-3&-4.48E-3&-2.14E-3&-5.08E-3&-1.28E-3&-9.42E-4&-4.06E-3\\
  6&-8.52E-6&-9.61E-3&-5.81E-2&-2.83E-1&-1.17E-1&-1.09E-1&-1.41E-1&-5.52E-2&-7.70E-2&-6.39E-2&-2.91E-2&-6.41E-2\\
\hline \hline
\end{tabular}
%}
\end{table}
\end{tiny}

~ ~ ~

\begin{tiny}
\begin{table}
\setlength{\tabcolsep}{2.0pt}
\renewcommand{\arraystretch}{0.95}
\caption{The poles of $r_{q,\delta,\alpha,k}(\t)$,  %${ r}(q,\delta,\alpha,k;\t)$,  
\eqref{BURA},   of $(q,\delta,\alpha,k)$-BURA approximation, 
i.e. the best uniform rational approximation of  
$\f(q,\delta,\alpha;\t)=\t^{\alpha}/(1+q\,\t^{\alpha})$  with functions from $\mathcal R_k$ on 
 $ [\delta, 1]$,   $\alpha =0.25, 0.5, 0.75$, ${\bf k=7}$, and $q=400$}
 %$(q=400,\delta,\alpha,k=7)$-BURA-poles}
 \label{tabl:Cq400dDaAAk7p}
%\rotatebox{90}{
\begin{tabular}{||c||c|c|c||c|c|c||c|c|c||c|c|c||}
\hline \hline
{\bf k=7} /$\delta$& $0.00$ & $0.00$ & $0.00$ &$10^{-6}$&$10^{-6}$&$10^{-6}$&$10^{-7}$&$10^{-7}$&$10^{-7}$&$10^{-8}$&$10^{-8}$&$10^{-8}$\\
\ j /$\alpha$ & 0.25   & 0.50   & 0.75   & 0.25   & 0.50   & 0.75   & 0.25   & 0.50   & 0.75   & 0.25   & 0.50   & 0.75   \\ \hline
  1& 9.05E-14&-3.83E-10&-1.38E-7&-4.45E-7&-8.62E-7&-2.79E-6&-5.67E-8&-1.62E-7&-7.58E-7&-7.20E-9&-3.20E-8&-2.93E-7\\
  2&-1.37E-13&-2.06E-8&-3.47E-6&-5.32E-6&-6.13E-6&-2.52E-5&-7.43E-7&-1.45E-6&-9.86E-6&-1.01E-7&-3.97E-7&-5.34E-6\\
  3& 3.25E-13&-3.85E-7&-3.77E-5&-4.16E-5&-3.42E-5&-1.40E-4&-6.92E-6&-9.29E-6&-7.50E-5&-1.10E-6&-3.24E-6&-5.00E-5\\
  4& 4.90E-10&-4.22E-6&-2.21E-4&-3.36E-4&-2.10E-4&-5.11E-4&-6.99E-5&-6.17E-5&-3.32E-4&-1.36E-5&-2.36E-5&-2.60E-4\\
  5& 3.57E-9&-4.21E-5&-9.30E-4&-3.03E-3&-1.56E-3&-2.39E-3&-8.29E-4&-5.15E-4&-1.45E-3&-2.07E-4&-2.10E-4&-1.10E-3\\
  6& 7.90E-8&-6.29E-4&-7.01E-3&-3.20E-2&-1.50E-2&-1.72E-2&-1.22E-2&-5.95E-3&-1.09E-2&-4.21E-3&-2.73E-3&-8.34E-3\\
  7& 2.77E-4&-2.10E-2&-9.93E-2&-4.70E-1&-2.23E-1&-1.96E-1&-2.61E-1&-1.17E-1&-1.39E-1&-1.36E-1&-6.63E-2&-1.14E-1\\
\hline \hline
\end{tabular}
%}
\end{table}
\end{tiny}

~ ~ ~

\begin{tiny}
\begin{table}[h!t]
\setlength{\tabcolsep}{1.0pt}
\renewcommand{\arraystretch}{1.02}
\caption{The poles of $r_{q,\delta,\alpha,k}(\t)$,  %${ r}(q,\delta,\alpha,k;\t)$,  
\eqref{BURA},   of $(q,\delta,\alpha,k)$-BURA approximation, 
i.e. the best uniform rational approximation of  
$\f(q,\delta,\alpha;\t)=\t^{\alpha}/(1+q\,\t^{\alpha})$  with functions from $\mathcal R_k$ on 
 $ [\delta, 1]$,   $\alpha =0.25, 0.5, 0.75$, ${\bf k=8}$, and $q=400$}
 %$(q=400,\delta,\alpha,k=8)$-BURA-poles}
 \label{tabl:Cq400dDaAAk8p}
%\rotatebox{90}{
\begin{tabular}{||c||c|c|c||c|c|c||c|c|c||c|c|c||}
\hline \hline
{\bf k=8} /$\delta$& $0.00$ & $0.00$ & $0.00$ &$10^{-6}$&$10^{-6}$&$10^{-6}$&$10^{-7}$&$10^{-7}$&$10^{-7}$&$10^{-8}$&$10^{-8}$&$10^{-8}$\\
\ j /$\alpha$ & 0.25   & 0.50   & 0.75   & 0.25   & 0.50   & 0.75   & 0.25   & 0.50   & 0.75   & 0.25   & 0.50   & 0.75   \\ \hline
  1& 1.55E-13&-1.10E-10&-4.66E-8&-3.51E-7&-6.88E-7&-1.94E-6&-4.49E-8&-1.23E-7&-4.65E-7&-5.73E-9&-2.27E-8&-1.49E-7\\
  2& 2.42E-11&-5.92E-9&-1.16E-6&-3.70E-6&-4.47E-6&-1.53E-5&-5.10E-7&-9.81E-7&-5.13E-6&-6.91E-8&-2.42E-7&-2.34E-6\\
  3&-1.05E-10&-1.12E-7&-1.28E-5&-2.37E-5&-2.14E-5&-8.14E-5&-3.81E-6&-5.51E-6&-3.64E-5&-5.92E-7&-1.74E-6&-2.07E-5\\
  4& 3.87E-10&-1.25E-6&-8.93E-5&-1.49E-4&-1.06E-4&-3.00E-4&-2.90E-5&-2.97E-5&-1.77E-4&-5.37E-6&-1.07E-5&-1.23E-4\\
  5&-3.93E-9&-1.05E-5&-3.71E-4&-9.94E-4&-5.94E-4&-1.03E-3&-2.44E-4&-1.83E-4&-6.18E-4&-5.56E-5&-7.02E-5&-4.62E-4\\
  6&-1.59E-7&-9.78E-5&-1.65E-3&-7.29E-3&-3.88E-3&-5.05E-3&-2.33E-3&-1.38E-3&-2.95E-3&-6.84E-4&-5.79E-4&-2.13E-3\\
  7&-1.55E-5&-1.37E-3&-1.23E-2&-6.06E-2&-3.12E-2&-3.26E-2&-2.64E-2&-1.36E-2&-2.07E-2&-1.06E-2&-6.57E-3&-1.55E-2\\
  8&-3.37E-3&-3.93E-2&-1.52E-1&-7.01E-1&-3.66E-1&-3.13E-1&-4.17E-1&-2.08E-1&-2.24E-1&-2.39E-1&-1.25E-1&-1.81E-1\\
\hline \hline
\end{tabular}
%}
\end{table}
\end{tiny}

~ ~ ~ 

%\vspace{2cm}
%\endinput

%\newpage

\section{Tables type \textbf{(d)} for the coefficients of partial fractions representation of BURA}
%
%\include{sect1ra-pm2}
%\subsubsection{Tables type \textbf{(d)} for BURA-decomposition coefficients}
\begin{tiny}
%\vspace{-2cm}
\begin{table}[h!t]
\setlength{\tabcolsep}{1.0pt}
\renewcommand{\arraystretch}{1.02}
\caption{The coefficients $c_j$, $ j=0, \dots, k$ 
of the partial fraction representation defined by \eqref{trg}
of $(q,\delta,\alpha,k)$-BURA approximation for  $\delta=0, 10^{-6}, 10^{-7}, 10^{-8}$, $\alpha= 0.25, 0.50, 0.75$,
${\bf k=3}$,  and $q=0$}
  \label{tabl:Dq000dDaAAk3p}
%\rotatebox{90}{
\begin{tabular}{||c||c|c|c||c|c|c||c|c|c||c|c|c||}
\hline \hline
{\bf k=3} /$\delta$& $0.00$ & $0.00$ & $0.00$ &$10^{-6}$&$10^{-6}$&$10^{-6}$&$10^{-7}$&$10^{-7}$&$10^{-7}$&$10^{-8}$&$10^{-8}$&$10^{-8}$\\
\ j /$\alpha$ & 0.25   & 0.50   & 0.75   & 0.25   & 0.50   & 0.75   & 0.25   & 0.50   & 0.75   & 0.25   & 0.50   & 0.75   \\ \hline
  0& 1.23E00& 2.28E00& 6.15E00& 1.29E00& 2.31E00& 6.17E00& 1.27E00& 2.29E00& 6.16E00& 1.25E00& 2.28E00& 6.15E00\\
  1&-1.40E-6&-4.18E-5&-9.88E-5&-1.12E-5&-5.15E-5&-1.01E-4&-5.10E-6&-4.48E-5&-9.92E-5&-3.02E-6&-4.28E-5&-9.88E-5\\
  2&-1.37E-3&-1.20E-2&-2.39E-2&-2.99E-3&-1.29E-2&-2.41E-2&-2.19E-3&-1.23E-2&-2.39E-2&-1.80E-3&-1.21E-2&-2.39E-2\\
  3&-3.44E-1&-3.52E00&-3.64E01&-4.60E-1&-3.65E00&-3.66E01&-4.07E-1&-3.56E00&-3.65E01&-3.79E-1&-3.54E00&-3.64E01\\
\hline \hline \end{tabular}
%}
\end{table}
\end{tiny}

~ ~ ~

\begin{tiny}
\begin{table}
\setlength{\tabcolsep}{1.0pt}
\renewcommand{\arraystretch}{1.02}
\caption{The coefficients $c_j$, $ j=0, \dots, k$   %$c_0, c_1, \dots, c_k$ 
of the partial fraction representation defined by \eqref{trg}
of $(q,\delta,\alpha,k)$-BURA approximation for  $\delta=0, 10^{-6}, 10^{-7}, 10^{-8}$, $\alpha= 0.25, 0.50, 0.75$,
${\bf k=4}$,  and $q=0$}
 % $(q=000,\delta,\alpha,k=4)$-BURA-decomposition}
 \label{tabl:Dq000dDaAAk4p}
%\rotatebox{90}{
%
%}
\end{table}
\end{tiny}

~ ~ ~

\begin{tiny}
\begin{table}
\setlength{\tabcolsep}{1.0pt}
\renewcommand{\arraystretch}{1.02}
\caption{The coefficients $c_j$, $ j=0, \dots, k$   %$c_0, c_1, \dots, c_k$ 
of the partial fraction representation defined by \eqref{trg}
of $(q,\delta,\alpha,k)$-BURA approximation for  $\delta=0, 10^{-6}, 10^{-7}, 10^{-8}$, $\alpha= 0.25, 0.50, 0.75$,
${\bf k=5}$,  and $q=0$}
 %$(q=000,\delta,\alpha,k=5)$-BURA-decomposition}\
\label{tabl:Dq000dDaAAk5p}
%\rotatebox{90}{
%
%}
\end{table}
\end{tiny}

~ ~ ~

\begin{tiny}
\begin{table}
\setlength{\tabcolsep}{1.0pt}
\renewcommand{\arraystretch}{1.02}
\caption{The coefficients $c_j$, $ j=0, \dots, k$   % $c_0, c_1, \dots, c_k$ 
of the partial fraction representation defined by \eqref{trg}
of $(q,\delta,\alpha,k)$-BURA approximation for  $\delta=0, 10^{-6}, 10^{-7}, 10^{-8}$, $\alpha= 0.25, 0.50, 0.75$,
${\bf k=6}$,  and $q=0$}
 %$(q=000,\delta,\alpha,k=6)$-BURA-decomposition}
 \label{tabl:Dq000dDaAAk6p}
%\rotatebox{90}{
%
%}
\end{table}
\end{tiny}

~ ~ ~

\begin{tiny}
\begin{table}
\setlength{\tabcolsep}{1.0pt}
\renewcommand{\arraystretch}{1.02}
\caption{The coefficients $c_j$, $ j=0, \dots, k$   %$c_0, c_1, \dots, c_k$ 
of the partial fraction representation defined by \eqref{trg}
of $(q,\delta,\alpha,k)$-BURA approximation for  $\delta=0, 10^{-6}, 10^{-7}, 10^{-8}$, $\alpha= 0.25, 0.50, 0.75$,
${\bf k=7}$,  and $q=0$}
 %$(q=000,\delta,\alpha,k=7)$-BURA-decomposition}
 \label{tabl:Dq000dDaAAk7p}
%\rotatebox{90}{
%
%}
\end{table}
\end{tiny}

~ ~ ~

\begin{tiny}
\begin{table}
\setlength{\tabcolsep}{1.0pt}
\renewcommand{\arraystretch}{1.02}
\caption{The coefficients $c_j$, $ j=0, \dots, k$   %$c_0, c_1, \dots, c_k$ 
of the partial fraction representation defined by \eqref{trg}
of $(q,\delta,\alpha,k)$-BURA approximation for  $\delta=0, 10^{-6}, 10^{-7}, 10^{-8}$, $\alpha= 0.25, 0.50, 0.75$,
${\bf k=8}$,  and $q=0$}
 %$(q=000,\delta,\alpha,k=8)$-BURA-decomposition}
 \label{tabl:Dq000dDaAAk8p}
%\rotatebox{90}{
%
%}
\end{table}
\end{tiny}

~ ~ ~

\begin{tiny}
\begin{table}
\setlength{\tabcolsep}{1.0pt}
\renewcommand{\arraystretch}{1.02}
\caption{The coefficients $c_j$, $ j=0, \dots, k$ 
of the partial fraction representation defined by \eqref{trg}
of $(q,\delta,\alpha,k)$-BURA approximation for  $\delta=0, 10^{-6}, 10^{-7}, 10^{-8}$, $\alpha= 0.25, 0.50, 0.75$,
${\bf k=3}$,  and $q=1$}
% $(q=001,\delta,\alpha,k=3)$-BURA-decomposition}
\label{tabl:Dq001dDaAAk3p}
%\rotatebox{90}{
%
%}
\end{table}
\end{tiny}

~ ~ ~

\begin{tiny}
\begin{table}
\setlength{\tabcolsep}{1.0pt}
\renewcommand{\arraystretch}{1.02}
\caption{The coefficients $c_j$, $ j=0, \dots, k$ 
of the partial fraction representation defined by \eqref{trg}
of $(q,\delta,\alpha,k)$-BURA approximation for  $\delta=0, 10^{-6}, 10^{-7}, 10^{-8}$, $\alpha= 0.25, 0.50, 0.75$,
${\bf k=4}$,  and $q=1$}
 %$(q=001,\delta,\alpha,k=4)$-BURA-decomposition}
 \label{tabl:Dq001dDaAAk4p}
%\rotatebox{90}{
%
%}
\end{table}
\end{tiny}

~ ~ ~

\begin{tiny}
\begin{table}
\setlength{\tabcolsep}{1.0pt}
\renewcommand{\arraystretch}{1.02}
\caption{The coefficients $c_j$, $ j=0, \dots, k$ 
of the partial fraction representation defined by \eqref{trg}
of $(q,\delta,\alpha,k)$-BURA approximation for  $\delta=0, 10^{-6}, 10^{-7}, 10^{-8}$, $\alpha= 0.25, 0.50, 0.75$,
${\bf k=5}$,  and $q=1$}
 %$(q=001,\delta,\alpha,k=5)$-BURA-decomposition}
 \label{tabl:Dq001dDaAAk5p}
%\rotatebox{90}{
%
%}
\end{table}
\end{tiny}

~ ~ ~

\begin{tiny}
\begin{table}
\setlength{\tabcolsep}{1.0pt}
\renewcommand{\arraystretch}{1.02}
\caption{The coefficients $c_j$, $ j=0, \dots, k$ 
of the partial fraction representation defined by \eqref{trg}
of $(q,\delta,\alpha,k)$-BURA approximation for  $\delta=0, 10^{-6}, 10^{-7}, 10^{-8}$, $\alpha= 0.25, 0.50, 0.75$,
${\bf k=6}$,  and $q=1$}
 %$(q=001,\delta,\alpha,k=6)$-BURA-decomposition}
 \label{tabl:Dq001dDaAAk6p}
%\rotatebox{90}{
%
%}
\end{table}
\end{tiny}

~ ~ ~ 

%\endinput

%\newpage

%\bigskip

%\subsubsection{Tables type \textbf{(e)} for 0-URA-poles}

%\include{sect1ra-pm3}

\section{Tables type \textbf{(e)} for 0-URA-poles}

\begin{tiny}
\begin{table}[h!t]
\setlength{\tabcolsep}{1.0pt}
\renewcommand{\arraystretch}{1.02}
\caption{The poles of ${\bar r}_{q,\delta,\alpha,k}(\t)$,  % 
 \eqref{0-URA},  of $(q,\delta,\alpha,k)$-0-URA approximation, 
i.e. the best uniform rational approximation of  
$\f(q,\delta,\alpha;\t)=\t^{\alpha}/(1+q\,\t^{\alpha})$ 
with functions from $\mathcal R_k$   on 
 $ [\delta, 1]$,   $\alpha =0.25, 0.5, 0.75$, ${\bf k=3}$, and $q=1$}
 %$(q=001,\delta,\alpha,k=3)$-0URA-poles}
 \label{tabl:Eq001dDaAAk3p}
%\rotatebox{90}{
\begin{tabular}{||c||c|c|c||c|c|c||c|c|c||c|c|c||}
\hline \hline
{\bf k=3} /$\delta$& $0.00$ & $0.00$ & $0.00$ &$10^{-6}$&$10^{-6}$&$10^{-6}$&$10^{-7}$&$10^{-7}$&$10^{-7}$&$10^{-8}$&$10^{-8}$&$10^{-8}$\\
\ j /$\alpha$ & 0.25   & 0.50   & 0.75   & 0.25   & 0.50   & 0.75   & 0.25   & 0.50   & 0.75   & 0.25   & 0.50   & 0.75   \\ \hline
  1&-1.23E-5&-8.41E-4&-5.35E-3&-8.32E-5&-9.92E-4&-5.43E-3&-4.06E-5&-8.89E-4&-5.36E-3&-2.51E-5&-8.57E-4&-5.35E-3\\
  2&-3.80E-3&-3.98E-2&-1.27E-1&-8.05E-3&-4.21E-2&-1.27E-1&-5.96E-3&-4.05E-2&-1.27E-1&-4.93E-3&-4.00E-2&-1.27E-1\\
  3&-2.60E-1&-7.04E-1&-1.02E00&-3.48E-1&-7.21E-1&-1.03E00&-3.09E-1&-7.10E-1&-1.02E00&-2.87E-1&-7.06E-1&-1.02E00\\
\hline \hline \end{tabular}
%}
\end{table}
\end{tiny}

~ ~ ~

\begin{tiny}
\begin{table}
\setlength{\tabcolsep}{1.0pt}
\renewcommand{\arraystretch}{1.02}
\caption{The poles of ${\bar r}_{q,\delta,\alpha,k}(\t)$,  % ${\bar r}(q,\delta,\alpha,k;\t)$,  
 \eqref{0-URA},  of $(q,\delta,\alpha,k)$-0-URA approximation, 
i.e. the best uniform rational approximation of  $\f(q,\delta,\alpha;\t)=\t^{\alpha}/(1+q\,\t^{\alpha})$ 
with functions from $\mathcal R_k$   on 
 $ [\delta, 1]$,   $\alpha =0.25, 0.5, 0.75$, ${\bf k=4}$, and $q=1$}
 %$(q=001,\delta,\alpha,k=4)$-0URA-poles}
 \label{tabl:Eq001dDaAAk4p}
%\rotatebox{90}{
\begin{tabular}{||c||c|c|c||c|c|c||c|c|c||c|c|c||}
\hline \hline
{\bf k=4} /$\delta$& $0.00$ & $0.00$ & $0.00$ &$10^{-6}$&$10^{-6}$&$10^{-6}$&$10^{-7}$&$10^{-7}$&$10^{-7}$&$10^{-8}$&$10^{-8}$&$10^{-8}$\\
\ j /$\alpha$ & 0.25   & 0.50   & 0.75   & 0.25   & 0.50   & 0.75   & 0.25   & 0.50   & 0.75   & 0.25   & 0.50   & 0.75   \\ \hline
  1&-5.37E-7&-9.04E-5&-8.33E-4&-1.83E-5&-1.47E-4&-8.87E-4&-6.00E-6&-1.08E-4&-8.44E-4&-2.53E-6&-9.60E-5&-8.35E-4\\
  2&-1.75E-4&-4.64E-3&-2.04E-2&-9.63E-4&-5.70E-3&-2.10E-2&-5.31E-4&-4.99E-3&-2.05E-2&-3.47E-4&-4.75E-3&-2.04E-2\\
  3&-1.12E-2&-7.92E-2&-2.15E-1&-2.61E-2&-8.72E-2&-2.17E-1&-1.92E-2&-8.20E-2&-2.15E-1&-1.55E-2&-8.01E-2&-2.15E-1\\
  4&-4.00E-1&-9.55E-1&-1.25E00&-5.87E-1&-1.00E00&-1.25E00&-5.08E-1&-9.71E-1&-1.25E00&-4.61E-1&-9.61E-1&-1.25E00\\
\hline \hline \end{tabular}
%}
\end{table}
\end{tiny}

~ ~ ~

\begin{tiny}
\begin{table}
\setlength{\tabcolsep}{1.0pt}
\renewcommand{\arraystretch}{1.02}
\caption{The poles of ${\bar r}_{q,\delta,\alpha,k}(\t)$,  % ${\bar r}(q,\delta,\alpha,k;\t)$, 
 \eqref{0-URA},  of $(q,\delta,\alpha,k)$-0-URA approximation, 
i.e. the rational approximation of  
$\f(q,\delta,\alpha;\t)=\t^{\alpha}/(1+q\,\t^{\alpha})$ 
with functions from $\mathcal R_k$  on 
 $ [\delta, 1]$,   $\alpha =0.25, 0.5, 0.75$, ${\bf k=5}$, and $q=1$}
 %$(q=001,\delta,\alpha,k=5)$-0URA-poles}
 \label{tabl:Eq001dDaAAk5p}
%\rotatebox{90}{
\begin{tabular}{||c||c|c|c||c|c|c||c|c|c||c|c|c||}
\hline \hline
{\bf k=5} /$\delta$& $0.00$ & $0.00$ & $0.00$ &$10^{-6}$&$10^{-6}$&$10^{-6}$&$10^{-7}$&$10^{-7}$&$10^{-7}$&$10^{-8}$&$10^{-8}$&$10^{-8}$\\
\ j /$\alpha$ & 0.25   & 0.50   & 0.75   & 0.25   & 0.50   & 0.75   & 0.25   & 0.50   & 0.75   & 0.25   & 0.50   & 0.75   \\ \hline
  1&-3.20E-8&-1.22E-5&-1.59E-4&-6.90E-6&-3.53E-5&-1.92E-4&-1.71E-6&-1.93E-5&-1.66E-4&-5.12E-7&-1.44E-5&-1.60E-4\\
  2&-1.08E-5&-6.45E-4&-3.93E-3&-2.12E-4&-1.08E-3&-4.29E-3&-8.75E-5&-7.98E-4&-4.01E-3&-4.34E-5&-6.96E-4&-3.95E-3\\
  3&-7.18E-4&-1.17E-2&-4.25E-2&-4.02E-3&-1.54E-2&-4.46E-2&-2.34E-3&-1.31E-2&-4.30E-2&-1.55E-3&-1.22E-2&-4.26E-2\\
  4&-2.24E-2&-1.24E-1&-3.02E-1&-5.75E-2&-1.44E-1&-3.09E-1&-4.23E-2&-1.32E-1&-3.04E-1&-3.37E-2&-1.27E-1&-3.02E-1\\
  5&-5.46E-1&-1.20E00&-1.47E00&-8.77E-1&-1.30E00&-1.49E00&-7.46E-1&-1.24E00&-1.47E00&-6.65E-1&-1.22E00&-1.47E00\\
\hline \hline \end{tabular}
%}
\end{table}
\end{tiny}

~ ~ ~

\begin{tiny}
\begin{table}
\setlength{\tabcolsep}{1.0pt}
\renewcommand{\arraystretch}{1.02}
\caption{The poles of ${\bar r}_{q,\delta,\alpha,k}(\t)$,  % ${\bar r}(q,\delta,\alpha,k;\t)$,  
\eqref{0-URA},  of $(q,\delta,\alpha,k)$-0-URA approximation, 
i.e. the rational approximation of  
$\f(q,\delta,\alpha;\t)=\t^{\alpha}/(1+q\,\t^{\alpha})$ 
with functions from $\mathcal R_k$   on 
 $ [\delta, 1]$,   $\alpha =0.25, 0.5, 0.75$, ${\bf k=6}$, and $q=1$}
 %$(q=001,\delta,\alpha,k=6)$-0URA-poles}
 \label{tabl:Eq001dDaAAk6p}
%\rotatebox{90}{
\begin{tabular}{||c||c|c|c||c|c|c||c|c|c||c|c|c||}
\hline \hline
{\bf k=6} /$\delta$& $0.00$ & $0.00$ & $0.00$ &$10^{-6}$&$10^{-6}$&$10^{-6}$&$10^{-7}$&$10^{-7}$&$10^{-7}$&$10^{-8}$&$10^{-8}$&$10^{-8}$\\
\ j /$\alpha$ & 0.25   & 0.50   & 0.75   & 0.25   & 0.50   & 0.75   & 0.25   & 0.50   & 0.75   & 0.25   & 0.50   & 0.75   \\ \hline
  1&-2.42E-9&-1.95E-6&-3.52E-5&-3.50E-6&-1.27E-5&-5.48E-5&-7.25E-7&-5.03E-6&-3.98E-5&-1.71E-7&-2.90E-6&-3.61E-5\\
  2&-8.33E-7&-1.05E-4&-8.73E-4&-7.03E-5&-2.88E-4&-1.09E-3&-2.30E-5&-1.70E-4&-9.27E-4&-8.84E-6&-1.27E-4&-8.84E-4\\
  3&-5.70E-5&-1.95E-3&-9.51E-3&-9.66E-4&-3.58E-3&-1.08E-2&-4.50E-4&-2.59E-3&-9.84E-3&-2.40E-4&-2.19E-3&-9.58E-3\\
  4&-1.82E-3&-2.15E-2&-6.94E-2&-1.06E-2&-3.08E-2&-7.48E-2&-6.47E-3&-2.54E-2&-7.08E-2&-4.35E-3&-2.30E-2&-6.97E-2\\
  5&-3.68E-2&-1.73E-1&-3.86E-1&-1.03E-1&-2.12E-1&-4.01E-1&-7.61E-2&-1.90E-1&-3.90E-1&-6.02E-2&-1.79E-1&-3.87E-1\\
  6&-6.94E-1&-1.45E00&-1.68E00&-1.22E00&-1.63E00&-1.73E00&-1.02E00&-1.53E00&-1.70E00&-8.99E-1&-1.48E00&-1.69E00\\
\hline \hline \end{tabular}
%}
\end{table}
\end{tiny}

~ ~ ~

\begin{tiny}
\begin{table}
\setlength{\tabcolsep}{1.0pt}
\renewcommand{\arraystretch}{1.02}
\caption{The poles of ${\bar r}_{q,\delta,\alpha,k}(\t)$,  % ${\bar r}(q,\delta,\alpha,k;\t)$,  
\eqref{0-URA},   of $(q,\delta,\alpha,k)$-0-URA approximation, 
i.e. the rational approximation of  
$\f(q,\delta,\alpha;\t)=\t^{\alpha}/(1+q\,\t^{\alpha})$ 
with functions from $\mathcal R_k$  %of the type \eqref{0-URA} on 
on  $ [\delta, 1]$,   $\alpha =0.25, 0.5, 0.75$, ${\bf k=7}$, and $q=1$}
 %$(q=001,\delta,\alpha,k=7)$-0URA-poles}
 \label{tabl:Eq001dDaAAk7p}
%\rotatebox{90}{
\begin{tabular}{||c||c|c|c||c|c|c||c|c|c||c|c|c||}
\hline \hline
{\bf k=7} /$\delta$& $0.00$ & $0.00$ & $0.00$ &$10^{-6}$&$10^{-6}$&$10^{-6}$&$10^{-7}$&$10^{-7}$&$10^{-7}$&$10^{-8}$&$10^{-8}$&$10^{-8}$\\
\ j /$\alpha$ & 0.25   & 0.50   & 0.75   & 0.25   & 0.50   & 0.75   & 0.25   & 0.50   & 0.75   & 0.25   & 0.50   & 0.75   \\ \hline
  1&-2.22E-10&-3.58E-7&-8.74E-6&-2.11E-6&-6.03E-6&-2.05E-5&-3.90E-7&-1.84E-6&-1.16E-5&-7.85E-8&-7.89E-7&-9.40E-6\\
  2&-7.71E-8&-1.93E-5&-2.17E-4&-3.09E-5&-1.02E-4&-3.44E-4&-8.41E-6&-4.76E-5&-2.52E-4&-2.62E-6&-2.92E-5&-2.25E-4\\
  3&-5.35E-6&-3.66E-4&-2.37E-3&-3.21E-4&-1.06E-3&-3.16E-3&-1.23E-4&-6.43E-4&-2.60E-3&-5.36E-5&-4.72E-4&-2.43E-3\\
  4&-1.75E-4&-4.14E-3&-1.74E-2&-2.82E-3&-8.27E-3&-2.08E-2&-1.42E-3&-5.95E-3&-1.84E-2&-7.99E-4&-4.87E-3&-1.76E-2\\
  5&-3.56E-3&-3.37E-2&-9.93E-2&-2.19E-2&-5.20E-2&-1.11E-1&-1.38E-2&-4.22E-2&-1.03E-1&-9.39E-3&-3.72E-2&-1.00E-1\\
  6&-5.37E-2&-2.23E-1&-4.66E-1&-1.61E-1&-2.90E-1&-4.95E-1&-1.20E-1&-2.55E-1&-4.75E-1&-9.51E-2&-2.37E-1&-4.68E-1\\
  7&-8.45E-1&-1.69E00&-1.90E00&-1.60E00&-1.99E00&-1.98E00&-1.34E00&-1.83E00&-1.93E00&-1.16E00&-1.75E00&-1.91E00\\
\hline \hline \end{tabular}
%}
\end{table}
\end{tiny}

~ ~ ~

\begin{tiny}
\begin{table}
\setlength{\tabcolsep}{1.0pt}
\renewcommand{\arraystretch}{1.02}
\caption{The poles of ${\bar r}_{q,\delta,\alpha,k}(\t)$,  % ${\bar r}(q,\delta,\alpha,k;\t)$,  
\eqref{0-URA},   of $(q,\delta,\alpha,k)$-0-URA approximation, 
i.e. the rational approximation of  
$\f(q,\delta,\alpha;\t)=\t^{\alpha}/(1+q\,\t^{\alpha})$ 
with functions from $\mathcal R_k$   on 
 $ [\delta, 1]$,   $\alpha =0.25, 0.5, 0.75$, ${\bf k=8}$, and $q=1$}
 %$(q=001,\delta,\alpha,k=8)$-0URA-poles}
 \label{tabl:Eq001dDaAAk8p}
%\rotatebox{90}{
\begin{tabular}{||c||c|c|c||c|c|c||c|c|c||c|c|c||}
\hline \hline
{\bf k=8} /$\delta$& $0.00$ & $0.00$ & $0.00$ &$10^{-6}$&$10^{-6}$&$10^{-6}$&$10^{-7}$&$10^{-7}$&$10^{-7}$&$10^{-8}$&$10^{-8}$&$10^{-8}$\\
\ j /$\alpha$ & 0.25   & 0.50   & 0.75   & 0.25   & 0.50   & 0.75   & 0.25   & 0.50   & 0.75   & 0.25   & 0.50   & 0.75   \\ \hline
  1&-2.38E-11&-7.35E-8&-2.38E-6&-1.41E-6&-3.44E-6&-9.54E-6&-2.43E-7&-8.66E-7&-4.20E-6&-4.40E-8&-2.86E-7&-2.82E-6\\
  2&-8.29E-9&-3.97E-6&-5.92E-5&-1.64E-5&-4.50E-5&-1.32E-4&-3.88E-6&-1.71E-5&-8.06E-5&-1.02E-6&-8.47E-6&-6.47E-5\\
  3&-5.81E-7&-7.56E-5&-6.48E-4&-1.35E-4&-3.92E-4&-1.10E-3&-4.38E-5&-1.98E-4&-7.89E-4&-1.60E-5&-1.23E-4&-6.85E-4\\
  4&-1.92E-5&-8.70E-4&-4.76E-3&-9.72E-4&-2.70E-3&-6.77E-3&-4.14E-4&-1.67E-3&-5.42E-3&-1.96E-4&-1.21E-3&-4.93E-3\\
  5&-3.99E-4&-7.24E-3&-2.72E-2&-6.34E-3&-1.56E-2&-3.45E-2&-3.39E-3&-1.12E-2&-2.97E-2&-1.98E-3&-9.00E-3&-2.79E-2\\
  6&-6.00E-3&-4.78E-2&-1.31E-1&-3.84E-2&-7.88E-2&-1.53E-1&-2.48E-2&-6.33E-2&-1.39E-1&-1.72E-2&-5.49E-2&-1.33E-1\\
  7&-7.27E-2&-2.75E-1&-5.43E-1&-2.32E-1&-3.78E-1&-5.90E-1&-1.75E-1&-3.28E-1&-5.60E-1&-1.38E-1&-3.00E-1&-5.47E-1\\
  8&-9.97E-1&-1.92E00&-2.12E00&-2.03E00&-2.38E00&-2.27E00&-1.68E00&-2.16E00&-2.17E00&-1.45E00&-2.04E00&-2.14E00\\
\hline \hline \end{tabular}
%}
\end{table}
\end{tiny}

~ ~ ~

\begin{tiny}
\begin{table}
\setlength{\tabcolsep}{1.0pt}
\renewcommand{\arraystretch}{1.02}
\caption{The poles of ${\bar r}_{q,\delta,\alpha,k}(\t)$,  % ${\bar r}(q,\delta,\alpha,k;\t)$,  
\eqref{0-URA},   of $(q,\delta,\alpha,k)$-0-URA approximation, 
i.e. the rational approximation of  
$\f(q,\delta,\alpha;\t)=\t^{\alpha}/(1+q\,\t^{\alpha})$ 
with functions from $\mathcal R_k$   on 
 $ [\delta, 1]$,   $\alpha =0.25, 0.5, 0.75$, ${\bf k=3}$, and $q=100$}
 %$(q=100,\delta,\alpha,k=3)$-0URA-poles}
 \label{tabl:Eq100dDaAAk3p}
%\rotatebox{90}{
\begin{tabular}{||c||c|c|c||c|c|c||c|c|c||c|c|c||}
\hline \hline
{\bf k=3} /$\delta$& $0.00$ & $0.00$ & $0.00$ &$10^{-6}$&$10^{-6}$&$10^{-6}$&$10^{-7}$&$10^{-7}$&$10^{-7}$&$10^{-8}$&$10^{-8}$&$10^{-8}$\\
\ j /$\alpha$ & 0.25   & 0.50   & 0.75   & 0.25   & 0.50   & 0.75   & 0.25   & 0.50   & 0.75   & 0.25   & 0.50   & 0.75   \\ \hline
  1&-2.41E-6&-1.69E-4&-1.55E-3&-2.60E-5&-2.02E-4&-1.57E-3&-1.09E-5&-1.79E-4&-1.55E-3&-6.03E-6&-1.72E-4&-1.55E-3\\
  2&-1.44E-3&-9.07E-3&-1.70E-2&-3.64E-3&-9.91E-3&-1.71E-2&-2.52E-3&-9.34E-3&-1.70E-2&-2.00E-3&-9.16E-3&-1.70E-2\\
  3&-1.39E-1&-2.70E-1&-3.19E-1&-1.99E-1&-2.80E-1&-3.20E-1&-1.72E-1&-2.73E-1&-3.19E-1&-1.57E-1&-2.71E-1&-3.19E-1\\
\hline \hline \end{tabular}
%}
\end{table}
\end{tiny}

~ ~ ~

\begin{tiny}
\begin{table}
\setlength{\tabcolsep}{1.0pt}
\renewcommand{\arraystretch}{1.02}
\caption{The poles of ${\bar r}_{q,\delta,\alpha,k}(\t)$,  % ${\bar r}(q,\delta,\alpha,k;\t)$, 
 \eqref{0-URA},   of $(q,\delta,\alpha,k)$-0-URA approximation, 
i.e. the rational approximation of  
$\f(q,\delta,\alpha;\t)=\t^{\alpha}/(1+q\,\t^{\alpha})$ 
with functions from $\mathcal R_k$   on 
 $ [\delta, 1]$,   $\alpha =0.25, 0.5, 0.75$, ${\bf k=4}$, and $q=100$}
 %$(q=100,\delta,\alpha,k=4)$-0URA-poles}
 \label{tabl:Eq100dDaAAk4p}
%\rotatebox{90}{
\begin{tabular}{||c||c|c|c||c|c|c||c|c|c||c|c|c||}
\hline \hline
{\bf k=4} /$\delta$& $0.00$ & $0.00$ & $0.00$ &$10^{-6}$&$10^{-6}$&$10^{-6}$&$10^{-7}$&$10^{-7}$&$10^{-7}$&$10^{-8}$&$10^{-8}$&$10^{-8}$\\
\ j /$\alpha$ & 0.25   & 0.50   & 0.75   & 0.25   & 0.50   & 0.75   & 0.25   & 0.50   & 0.75   & 0.25   & 0.50   & 0.75   \\ \hline
  1&-1.40E-7&-3.64E-5&-5.45E-4&-7.16E-6&-5.61E-5&-5.73E-4&-2.11E-6&-4.27E-5&-5.51E-4&-8.11E-7&-3.84E-5&-5.46E-4\\
  2&-6.34E-5&-1.17E-3&-4.22E-3&-4.56E-4&-1.50E-3&-4.33E-3&-2.30E-4&-1.28E-3&-4.24E-3&-1.41E-4&-1.20E-3&-4.23E-3\\
  3&-5.44E-3&-2.54E-2&-4.39E-2&-1.48E-2&-2.92E-2&-4.48E-2&-1.03E-2&-2.67E-2&-4.40E-2&-8.03E-3&-2.59E-2&-4.39E-2\\
  4&-2.36E-1&-4.32E-1&-5.15E-1&-3.71E-1&-4.63E-1&-5.21E-1&-3.13E-1&-4.43E-1&-5.16E-1&-2.79E-1&-4.36E-1&-5.15E-1\\
\hline \hline \end{tabular}
%}
\end{table}
\end{tiny}

~ ~ ~

\begin{tiny}
\begin{table}
\setlength{\tabcolsep}{1.0pt}
\renewcommand{\arraystretch}{1.02}
\caption{The poles of ${\bar r}_{q,\delta,\alpha,k}(\t)$,  % ${\bar r}(q,\delta,\alpha,k;\t)$, 
 \eqref{0-URA},   of $(q,\delta,\alpha,k)$-0-URA approximation, 
i.e. the rational approximation of  
$\f(q,\delta,\alpha;\t)=\t^{\alpha}/(1+q\,\t^{\alpha})$ 
with functions from $\mathcal R_k$  of the type \eqref{0-URA} on 
 $ [\delta, 1]$,   $\alpha =0.25, 0.5, 0.75$, ${\bf k=5}$, and $q=100$}
 %$(q=100,\delta,\alpha,k=5)$-0URA-poles}
 \label{tabl:Eq100dDaAAk5p}
%\rotatebox{90}{
\begin{tabular}{||c||c|c|c||c|c|c||c|c|c||c|c|c||}
\hline \hline
{\bf k=5} /$\delta$& $0.00$ & $0.00$ & $0.00$ &$10^{-6}$&$10^{-6}$&$10^{-6}$&$10^{-7}$&$10^{-7}$&$10^{-7}$&$10^{-8}$&$10^{-8}$&$10^{-8}$\\
\ j /$\alpha$ & 0.25   & 0.50   & 0.75   & 0.25   & 0.50   & 0.75   & 0.25   & 0.50   & 0.75   & 0.25   & 0.50   & 0.75   \\ \hline
  1&-1.17E-8&-7.81E-6&-1.40E-4&-3.15E-6&-2.02E-5&-1.67E-4&-7.36E-7&-1.18E-5&-1.46E-4&-2.10E-7&-9.09E-6&-1.41E-4\\
  2&-4.09E-6&-2.23E-4&-1.68E-3&-1.09E-4&-3.77E-4&-1.79E-3&-4.13E-5&-2.76E-4&-1.70E-3&-1.91E-5&-2.41E-4&-1.69E-3\\
  3&-3.26E-4&-3.67E-3&-9.61E-3&-2.31E-3&-5.22E-3&-1.02E-2&-1.26E-3&-4.23E-3&-9.74E-3&-7.86E-4&-3.86E-3&-9.64E-3\\
  4&-1.24E-2&-4.85E-2&-8.01E-2&-3.69E-2&-5.96E-2&-8.34E-2&-2.60E-2&-5.27E-2&-8.09E-2&-2.00E-2&-5.00E-2&-8.03E-2\\
  5&-3.40E-1&-6.03E-1&-7.21E-1&-5.88E-1&-6.76E-1&-7.38E-1&-4.89E-1&-6.31E-1&-7.25E-1&-4.28E-1&-6.13E-1&-7.22E-1\\
\hline \hline \end{tabular}
%}
\end{table}
\end{tiny}

~ ~ ~

\begin{tiny}
\begin{table}
\setlength{\tabcolsep}{1.0pt}
\renewcommand{\arraystretch}{1.02}
\caption{The poles of ${\bar r}_{q,\delta,\alpha,k}(\t)$,  % ${\bar r}(q,\delta,\alpha,k;\t)$,  
\eqref{0-URA},   of $(q,\delta,\alpha,k)$-0-URA approximation, 
i.e. the rational approximation of  
$\f(q,\delta,\alpha;\t)=\t^{\alpha}/(1+q\,\t^{\alpha})$ 
with functions from $\mathcal R_k$  of the type \eqref{0-URA} on 
 $ [\delta, 1]$,   $\alpha =0.25, 0.5, 0.75$, ${\bf k=6}$, and $q=100$}
 %$(q=100,\delta,\alpha,k=6)$-0URA-poles}
 \label{tabl:Eq100dDaAAk6p}
%\rotatebox{90}{
\begin{tabular}{||c||c|c|c||c|c|c||c|c|c||c|c|c||}
\hline \hline
{\bf k=6} /$\delta$& $0.00$ & $0.00$ & $0.00$ &$10^{-6}$&$10^{-6}$&$10^{-6}$&$10^{-7}$&$10^{-7}$&$10^{-7}$&$10^{-8}$&$10^{-8}$&$10^{-8}$\\
\ j /$\alpha$ & 0.25   & 0.50   & 0.75   & 0.25   & 0.50   & 0.75   & 0.25   & 0.50   & 0.75   & 0.25   & 0.50   & 0.75   \\ \hline
  1&-1.20E-9&-1.59E-6&-3.38E-5&-1.77E-6&-9.04E-6&-5.21E-5&-3.59E-7&-3.88E-6&-3.81E-5&-8.39E-8&-2.32E-6&-3.47E-5\\
  2&-3.53E-7&-5.32E-5&-6.32E-4&-3.94E-5&-1.34E-4&-7.57E-4&-1.20E-5&-8.24E-5&-6.64E-4&-4.31E-6&-6.34E-5&-6.39E-4\\
  3&-2.57E-5&-7.00E-4&-3.21E-3&-5.77E-4&-1.36E-3&-3.58E-3&-2.50E-4&-9.53E-4&-3.31E-3&-1.26E-4&-7.91E-4&-3.23E-3\\
  4&-9.41E-4&-7.95E-3&-1.83E-2&-6.86E-3&-1.25E-2&-2.03E-2&-3.94E-3&-9.79E-3&-1.89E-2&-2.52E-3&-8.63E-3&-1.85E-2\\
  5&-2.21E-2&-7.66E-2&-1.23E-1&-7.12E-2&-1.01E-1&-1.31E-1&-5.08E-2&-8.70E-2&-1.25E-1&-3.89E-2&-8.05E-2&-1.23E-1\\
  6&-4.50E-1&-7.80E-1&-9.33E-1&-8.48E-1&-9.19E-1&-9.73E-1&-6.99E-1&-8.40E-1&-9.43E-1&-6.04E-1&-8.03E-1&-9.35E-1\\
\hline \hline \end{tabular}
%}
\end{table}
\end{tiny}

~ ~ ~

\begin{tiny}
\begin{table}
\setlength{\tabcolsep}{1.0pt}
\renewcommand{\arraystretch}{1.02}
\caption{The poles of ${\bar r}_{q,\delta,\alpha,k}(\t)$,  % ${\bar r}(q,\delta,\alpha,k;\t)$, 
 \eqref{0-URA},   of $(q,\delta,\alpha,k)$-0-URA approximation, 
i.e. the rational approximation of  
$\f(q,\delta,\alpha;\t)=\t^{\alpha}/(1+q\,\t^{\alpha})$ 
with functions from $\mathcal R_k$   on 
 $ [\delta, 1]$,   $\alpha =0.25, 0.5, 0.75$, ${\bf k=7}$, and $q=100$}
 % $(q=100,\delta,\alpha,k=7)$-0URA-poles}
 \label{tabl:Eq100dDaAAk7p}
%\rotatebox{90}{
\begin{tabular}{||c||c|c|c||c|c|c||c|c|c||c|c|c||}
\hline \hline
{\bf k=7} /$\delta$& $0.00$ & $0.00$ & $0.00$ &$10^{-6}$&$10^{-6}$&$10^{-6}$&$10^{-7}$&$10^{-7}$&$10^{-7}$&$10^{-8}$&$10^{-8}$&$10^{-8}$\\
\ j /$\alpha$ & 0.25   & 0.50   & 0.75   & 0.25   & 0.50   & 0.75   & 0.25   & 0.50   & 0.75   & 0.25   & 0.50   & 0.75   \\ \hline
  1&-1.39E-10&-3.27E-7&-8.62E-6&-1.14E-6&-4.82E-6&-2.00E-5&-2.12E-7&-1.58E-6&-1.15E-5&-4.35E-8&-7.04E-7&-9.26E-6\\
  2&-3.80E-8&-1.33E-5&-1.94E-4&-1.86E-5&-5.98E-5&-2.97E-4&-4.77E-6&-3.01E-5&-2.23E-4&-1.42E-6&-1.94E-5&-2.01E-4\\
  3&-2.53E-6&-1.68E-4&-1.37E-3&-2.01E-4&-4.79E-4&-1.68E-3&-7.19E-5&-2.90E-4&-1.46E-3&-2.96E-5&-2.14E-4&-1.39E-3\\
  4&-8.85E-5&-1.61E-3&-5.57E-3&-1.85E-3&-3.54E-3&-6.73E-3&-8.80E-4&-2.42E-3&-5.91E-3&-4.69E-4&-1.93E-3&-5.65E-3\\
  5&-2.02E-3&-1.40E-2&-3.00E-2&-1.52E-2&-2.41E-2&-3.51E-2&-9.12E-3&-1.85E-2&-3.15E-2&-5.96E-3&-1.58E-2&-3.04E-2\\
  6&-3.41E-2&-1.08E-1&-1.70E-1&-1.18E-1&-1.53E-1&-1.88E-1&-8.51E-2&-1.29E-1&-1.75E-1&-6.53E-2&-1.17E-1&-1.71E-1\\
  7&-5.64E-1&-9.60E-1&-1.15E00&-1.15E00&-1.19E00&-1.23E00&-9.41E-1&-1.07E00&-1.17E00&-8.06E-1&-1.01E00&-1.15E00\\
\hline \hline \end{tabular}
%}
\end{table}
\end{tiny}

~ ~ ~

\begin{tiny}
\begin{table}
\setlength{\tabcolsep}{1.0pt}
\renewcommand{\arraystretch}{1.02}
\caption{The poles of ${\bar r}_{q,\delta,\alpha,k}(\t)$,  % ${\bar r}(q,\delta,\alpha,k;\t)$, 
 \eqref{0-URA},   of $(q,\delta,\alpha,k)$-0-URA approximation, 
i.e. the rational approximation of  
$\f(q,\delta,\alpha;\t)=\t^{\alpha}/(1+q\,\t^{\alpha})$ 
with functions from $\mathcal R_k$   on 
 $ [\delta, 1]$,   $\alpha =0.25, 0.5, 0.75$, ${\bf k=8}$, and $q=100$}
 %$(q=100,\delta,\alpha,k=8)$-0URA-poles}
 \label{tabl:Eq100dDaAAk8p}
%\rotatebox{90}{
\begin{tabular}{||c||c|c|c||c|c|c||c|c|c||c|c|c||}
\hline \hline
{\bf k=8} /$\delta$& $0.00$ & $0.00$ & $0.00$ &$10^{-6}$&$10^{-6}$&$10^{-6}$&$10^{-7}$&$10^{-7}$&$10^{-7}$&$10^{-8}$&$10^{-8}$&$10^{-8}$\\
\ j /$\alpha$ & 0.25   & 0.50   & 0.75   & 0.25   & 0.50   & 0.75   & 0.25   & 0.50   & 0.75   & 0.25   & 0.50   & 0.75   \\ \hline
  1&-1.76E-11&-7.05E-8&-2.37E-6&-7.99E-7&-2.93E-6&-9.43E-6&-1.41E-7&-7.86E-7&-4.17E-6&-2.65E-8&-2.68E-7&-2.80E-6\\
  2&-4.81E-9&-3.29E-6&-5.69E-5&-1.04E-5&-3.09E-5&-1.24E-4&-2.36E-6&-1.28E-5&-7.67E-5&-6.03E-7&-6.68E-6&-6.19E-5\\
  3&-2.98E-7&-4.57E-5&-5.23E-4&-8.81E-5&-2.09E-4&-8.14E-4&-2.70E-5&-1.10E-4&-6.19E-4&-9.39E-6&-7.12E-5&-5.48E-4\\
  4&-9.89E-6&-3.91E-4&-2.27E-3&-6.55E-4&-1.27E-3&-2.99E-3&-2.63E-4&-7.65E-4&-2.51E-3&-1.18E-4&-5.46E-4&-2.34E-3\\
  5&-2.20E-4&-3.04E-3&-9.04E-3&-4.45E-3&-7.41E-3&-1.19E-2&-2.26E-3&-5.04E-3&-9.99E-3&-1.26E-3&-3.91E-3&-9.29E-3\\
  6&-3.62E-3&-2.17E-2&-4.42E-2&-2.81E-2&-4.04E-2&-5.44E-2&-1.75E-2&-3.08E-2&-4.77E-2&-1.17E-2&-2.58E-2&-4.51E-2\\
  7&-4.82E-2&-1.43E-1&-2.20E-1&-1.76E-1&-2.16E-1&-2.54E-1&-1.29E-1&-1.80E-1&-2.32E-1&-9.93E-2&-1.60E-1&-2.23E-1\\
  8&-6.79E-1&-1.14E00&-1.37E00&-1.49E00&-1.50E00&-1.51E00&-1.21E00&-1.33E00&-1.42E00&-1.03E00&-1.23E00&-1.38E00\\
\hline \hline \end{tabular}
%}
\end{table}
\end{tiny}

~ ~ ~

\begin{tiny}
\begin{table}
\setlength{\tabcolsep}{1.0pt}
\renewcommand{\arraystretch}{1.02}
\caption{The poles of ${\bar r}_{q,\delta,\alpha,k}(\t)$,  % ${\bar r}(q,\delta,\alpha,k;\t)$, 
 \eqref{0-URA},   of $(q,\delta,\alpha,k)$-0-URA approximation, 
i.e. the rational approximation of  
$\f(q,\delta,\alpha;\t)=\t^{\alpha}/(1+q\,\t^{\alpha})$ 
with functions from $\mathcal R_k$  of the type \eqref{0-URA} on 
 $ [\delta, 1]$,   $\alpha =0.25, 0.5, 0.75$, ${\bf k=3}$, and $q=200$}
 %$(q=200,\delta,\alpha,k=3)$-0URA-poles}
 \label{tabl:Eq200dDaAAk3p}
%\rotatebox{90}{
\begin{tabular}{||c||c|c|c||c|c|c||c|c|c||c|c|c||}
\hline \hline
{\bf k=3} /$\delta$& $0.00$ & $0.00$ & $0.00$ &$10^{-6}$&$10^{-6}$&$10^{-6}$&$10^{-7}$&$10^{-7}$&$10^{-7}$&$10^{-8}$&$10^{-8}$&$10^{-8}$\\
\ j /$\alpha$ & 0.25   & 0.50   & 0.75   & 0.25   & 0.50   & 0.75   & 0.25   & 0.50   & 0.75   & 0.25   & 0.50   & 0.75   \\ \hline
  1&-1.95E-6&-1.08E-4&-8.95E-4&-2.40E-5&-1.34E-4&-9.05E-4&-9.77E-6&-1.16E-4&-8.97E-4&-5.23E-6&-1.11E-4&-8.95E-4\\
  2&-1.40E-3&-8.52E-3&-1.54E-2&-3.57E-3&-9.35E-3&-1.56E-2&-2.46E-3&-8.79E-3&-1.54E-2&-1.95E-3&-8.61E-3&-1.54E-2\\
  3&-1.38E-1&-2.67E-1&-3.16E-1&-1.98E-1&-2.77E-1&-3.18E-1&-1.71E-1&-2.70E-1&-3.16E-1&-1.56E-1&-2.68E-1&-3.16E-1\\
\hline \hline \end{tabular}
%}
\end{table}
\end{tiny}

~ ~ ~

\begin{tiny}
\begin{table}
\setlength{\tabcolsep}{1.0pt}
\renewcommand{\arraystretch}{1.02}
\caption{The poles of ${\bar r}_{q,\delta,\alpha,k}(\t)$,  % ${\bar r}(q,\delta,\alpha,k;\t)$, 
 \eqref{0-URA},   of $(q,\delta,\alpha,k)$-0-URA approximation, 
i.e. the rational approximation of  
$\f(q,\delta,\alpha;\t)=\t^{\alpha}/(1+q\,\t^{\alpha})$ 
with functions from $\mathcal R_k$  on 
 $ [\delta, 1]$,   $\alpha =0.25, 0.5, 0.75$, ${\bf k=4}$, and $q=200$}
 %$(q=200,\delta,\alpha,k=4)$-0URA-poles}
 \label{tabl:Eq200dDaAAk4p}
%\rotatebox{90}{
\begin{tabular}{||c||c|c|c||c|c|c||c|c|c||c|c|c||}
\hline \hline
{\bf k=4} /$\delta$& $0.00$ & $0.00$ & $0.00$ &$10^{-6}$&$10^{-6}$&$10^{-6}$&$10^{-7}$&$10^{-7}$&$10^{-7}$&$10^{-8}$&$10^{-8}$&$10^{-8}$\\
\ j /$\alpha$ & 0.25   & 0.50   & 0.75   & 0.25   & 0.50   & 0.75   & 0.25   & 0.50   & 0.75   & 0.25   & 0.50   & 0.75   \\ \hline
  1&-1.03E-7&-2.36E-5&-3.87E-4&-6.57E-6&-3.74E-5&-4.04E-4&-1.86E-6&-2.80E-5&-3.91E-4&-6.84E-7&-2.50E-5&-3.88E-4\\
  2&-5.98E-5&-1.00E-3&-3.12E-3&-4.44E-4&-1.32E-3&-3.23E-3&-2.22E-4&-1.11E-3&-3.14E-3&-1.35E-4&-1.03E-3&-3.13E-3\\
  3&-5.36E-3&-2.47E-2&-4.24E-2&-1.47E-2&-2.85E-2&-4.33E-2&-1.02E-2&-2.60E-2&-4.26E-2&-7.93E-3&-2.51E-2&-4.25E-2\\
  4&-2.34E-1&-4.29E-1&-5.12E-1&-3.69E-1&-4.60E-1&-5.18E-1&-3.11E-1&-4.40E-1&-5.13E-1&-2.78E-1&-4.32E-1&-5.13E-1\\
\hline \hline \end{tabular}
%}
\end{table}
\end{tiny}

~ ~ ~

\begin{tiny}
\begin{table}
\setlength{\tabcolsep}{1.0pt}
\renewcommand{\arraystretch}{1.02}
\caption{The poles of ${\bar r}_{q,\delta,\alpha,k}(\t)$,  % ${\bar r}(q,\delta,\alpha,k;\t)$, 
 \eqref{0-URA},   of $(q,\delta,\alpha,k)$-0-URA approximation, 
i.e. the rational approximation of  
$\f(q,\delta,\alpha;\t)=\t^{\alpha}/(1+q\,\t^{\alpha})$ 
with functions from $\mathcal R_k$   on 
 $ [\delta, 1]$,   $\alpha =0.25, 0.5, 0.75$, ${\bf k=5}$, and $q=200$}
 %$(q=200,\delta,\alpha,k=5)$-0URA-poles}
 \label{tabl:Eq200dDaAAk5p}
%\rotatebox{90}{
\begin{tabular}{||c||c|c|c||c|c|c||c|c|c||c|c|c||}
\hline \hline
{\bf k=5} /$\delta$& $0.00$ & $0.00$ & $0.00$ &$10^{-6}$&$10^{-6}$&$10^{-6}$&$10^{-7}$&$10^{-7}$&$10^{-7}$&$10^{-8}$&$10^{-8}$&$10^{-8}$\\
\ j /$\alpha$ & 0.25   & 0.50   & 0.75   & 0.25   & 0.50   & 0.75   & 0.25   & 0.50   & 0.75   & 0.25   & 0.50   & 0.75   \\ \hline
  1&-8.24E-9&-5.76E-6&-1.23E-4&-2.88E-6&-1.46E-5&-1.45E-4&-6.49E-7&-8.62E-6&-1.28E-4&-1.77E-7&-6.68E-6&-1.24E-4\\
  2&-3.70E-6&-1.71E-4&-1.10E-3&-1.06E-4&-3.06E-4&-1.17E-3&-3.95E-5&-2.17E-4&-1.12E-3&-1.80E-5&-1.86E-4&-1.10E-3\\
  3&-3.16E-4&-3.41E-3&-8.55E-3&-2.27E-3&-4.93E-3&-9.14E-3&-1.23E-3&-3.96E-3&-8.68E-3&-7.69E-4&-3.60E-3&-8.57E-3\\
  4&-1.23E-2&-4.76E-2&-7.88E-2&-3.67E-2&-5.87E-2&-8.20E-2&-2.58E-2&-5.18E-2&-7.95E-2&-1.98E-2&-4.91E-2&-7.89E-2\\
  5&-3.39E-1&-6.00E-1&-7.18E-1&-5.85E-1&-6.72E-1&-7.35E-1&-4.87E-1&-6.28E-1&-7.22E-1&-4.26E-1&-6.09E-1&-7.19E-1\\
\hline \hline \end{tabular}
%}
\end{table}
\end{tiny}

~ ~ ~

\begin{tiny}
\begin{table}
\setlength{\tabcolsep}{1.0pt}
\renewcommand{\arraystretch}{1.02}
\caption{The poles of ${\bar r}_{q,\delta,\alpha,k}(\t)$,  % ${\bar r}(q,\delta,\alpha,k;\t)$, 
 \eqref{0-URA},   of $(q,\delta,\alpha,k)$-0-URA approximation, 
i.e. the rational approximation of  
$\f(q,\delta,\alpha;\t)=\t^{\alpha}/(1+q\,\t^{\alpha})$ 
with functions from $\mathcal R_k$  on 
 $ [\delta, 1]$,   $\alpha =0.25, 0.5, 0.75$, ${\bf k=6}$, and $q=200$}
 % $(q=200,\delta,\alpha,k=6)$-0URA-poles}
 \label{tabl:Eq200dDaAAk6p}
%\rotatebox{90}{
\begin{tabular}{||c||c|c|c||c|c|c||c|c|c||c|c|c||}
\hline \hline
{\bf k=6} /$\delta$& $0.00$ & $0.00$ & $0.00$ &$10^{-6}$&$10^{-6}$&$10^{-6}$&$10^{-7}$&$10^{-7}$&$10^{-7}$&$10^{-8}$&$10^{-8}$&$10^{-8}$\\
\ j /$\alpha$ & 0.25   & 0.50   & 0.75   & 0.25   & 0.50   & 0.75   & 0.25   & 0.50   & 0.75   & 0.25   & 0.50   & 0.75   \\ \hline
  1&-8.56E-10&-1.34E-6&-3.24E-5&-1.62E-6&-7.10E-6&-4.94E-5&-3.17E-7&-3.16E-6&-3.64E-5&-7.10E-8&-1.93E-6&-3.33E-5\\
  2&-3.04E-7&-3.93E-5&-4.74E-4&-3.81E-5&-1.05E-4&-5.55E-4&-1.13E-5&-6.23E-5&-4.95E-4&-4.01E-6&-4.72E-5&-4.79E-4\\
  3&-2.44E-5&-6.05E-4&-2.40E-3&-5.67E-4&-1.23E-3&-2.75E-3&-2.44E-4&-8.43E-4&-2.49E-3&-1.22E-4&-6.90E-4&-2.42E-3\\
  4&-9.22E-4&-7.61E-3&-1.73E-2&-6.80E-3&-1.21E-2&-1.94E-2&-3.89E-3&-9.44E-3&-1.79E-2&-2.49E-3&-8.29E-3&-1.75E-2\\
  5&-2.19E-2&-7.56E-2&-1.21E-1&-7.09E-2&-1.00E-1&-1.30E-1&-5.05E-2&-8.60E-2&-1.24E-1&-3.87E-2&-7.95E-2&-1.22E-1\\
  6&-4.48E-1&-7.76E-1&-9.30E-1&-8.45E-1&-9.15E-1&-9.71E-1&-6.96E-1&-8.36E-1&-9.41E-1&-6.02E-1&-7.99E-1&-9.32E-1\\
\hline \hline \end{tabular}
%}
\end{table}
\end{tiny}

~ ~ ~

\begin{tiny}
\begin{table}
\setlength{\tabcolsep}{1.0pt}
\renewcommand{\arraystretch}{1.02}
\caption{The poles of ${\bar r}_{q,\delta,\alpha,k}(\t)$,  % ${\bar r}(q,\delta,\alpha,k;\t)$, 
 \eqref{0-URA},   of $(q,\delta,\alpha,k)$-0-URA approximation, 
i.e. the rational approximation of  
$\f(q,\delta,\alpha;\t)=\t^{\alpha}/(1+q\,\t^{\alpha})$ 
with functions from $\mathcal R_k$   on 
 $ [\delta, 1]$,   $\alpha =0.25, 0.5, 0.75$, ${\bf k=7}$, and $q=200$}
 %$(q=200,\delta,\alpha,k=7)$-0URA-poles}
 \label{tabl:Eq200dDaAAk7p}
%\rotatebox{90}{
\begin{tabular}{||c||c|c|c||c|c|c||c|c|c||c|c|c||}
\hline \hline
{\bf k=7} /$\delta$& $0.00$ & $0.00$ & $0.00$ &$10^{-6}$&$10^{-6}$&$10^{-6}$&$10^{-7}$&$10^{-7}$&$10^{-7}$&$10^{-8}$&$10^{-8}$&$10^{-8}$\\
\ j /$\alpha$ & 0.25   & 0.50   & 0.75   & 0.25   & 0.50   & 0.75   & 0.25   & 0.50   & 0.75   & 0.25   & 0.50   & 0.75   \\ \hline
  1&-1.05E-10&-3.00E-7&-8.50E-6&-1.04E-6&-4.03E-6&-1.95E-5&-1.88E-7&-1.39E-6&-1.13E-5&-3.72E-8&-6.34E-7&-9.12E-6\\
  2&-3.15E-8&-1.04E-5&-1.73E-4&-1.79E-5&-4.70E-5&-2.55E-4&-4.51E-6&-2.33E-5&-1.96E-4&-1.31E-6&-1.50E-5&-1.78E-4\\
  3&-2.33E-6&-1.34E-4&-9.56E-4&-1.97E-4&-4.17E-4&-1.19E-3&-6.99E-5&-2.42E-4&-1.02E-3&-2.85E-5&-1.75E-4&-9.72E-4\\
  4&-8.55E-5&-1.47E-3&-4.74E-3&-1.83E-3&-3.36E-3&-5.91E-3&-8.67E-4&-2.26E-3&-5.08E-3&-4.60E-4&-1.78E-3&-4.82E-3\\
  5&-1.99E-3&-1.36E-2&-2.91E-2&-1.51E-2&-2.36E-2&-3.41E-2&-9.05E-3&-1.81E-2&-3.06E-2&-5.91E-3&-1.54E-2&-2.94E-2\\
  6&-3.39E-2&-1.07E-1&-1.68E-1&-1.17E-1&-1.52E-1&-1.87E-1&-8.47E-2&-1.28E-1&-1.74E-1&-6.50E-2&-1.16E-1&-1.70E-1\\
  7&-5.61E-1&-9.56E-1&-1.15E00&-1.14E00&-1.19E00&-1.23E00&-9.38E-1&-1.07E00&-1.17E00&-8.04E-1&-1.00E00&-1.15E00\\
\hline \hline \end{tabular}
%}
\end{table}
\end{tiny}

~ ~ ~

\begin{tiny}
\begin{table}
\setlength{\tabcolsep}{1.0pt}
\renewcommand{\arraystretch}{1.02}
\caption{The poles of ${\bar r}_{q,\delta,\alpha,k}(\t)$,  % ${\bar r}(q,\delta,\alpha,k;\t)$, 
 \eqref{0-URA},   of $(q,\delta,\alpha,k)$-0-URA approximation, 
i.e. the rational approximation of  
$\f(q,\delta,\alpha;\t)=\t^{\alpha}/(1+q\,\t^{\alpha})$ 
with functions from $\mathcal R_k$  of the type \eqref{0-URA} on 
 $ [\delta, 1]$,   $\alpha =0.25, 0.5, 0.75$, ${\bf k=8}$, and $q=200$}
 %$(q=200,\delta,\alpha,k=8)$-0URA-poles}
 \label{tabl:Eq200dDaAAk8p}
%\rotatebox{90}{
\begin{tabular}{||c||c|c|c||c|c|c||c|c|c||c|c|c||}
\hline \hline
{\bf k=8} /$\delta$& $0.00$ & $0.00$ & $0.00$ &$10^{-6}$&$10^{-6}$&$10^{-6}$&$10^{-7}$&$10^{-7}$&$10^{-7}$&$10^{-8}$&$10^{-8}$&$10^{-8}$\\
\ j /$\alpha$ & 0.25   & 0.50   & 0.75   & 0.25   & 0.50   & 0.75   & 0.25   & 0.50   & 0.75   & 0.25   & 0.50   & 0.75   \\ \hline
  1&-1.41E-11&-6.76E-8&-2.36E-6&-7.31E-7&-2.55E-6&-9.32E-6&-1.25E-7&-7.20E-7&-4.14E-6&-2.28E-8&-2.52E-7&-2.78E-6\\
  2&-3.90E-9&-2.81E-6&-5.45E-5&-1.00E-5&-2.50E-5&-1.16E-4&-2.23E-6&-1.04E-5&-7.29E-5&-5.56E-7&-5.57E-6&-5.92E-5\\
  3&-2.65E-7&-3.57E-5&-4.21E-4&-8.62E-5&-1.77E-4&-6.22E-4&-2.62E-5&-8.93E-5&-4.88E-4&-8.98E-6&-5.65E-5&-4.39E-4\\
  4&-9.36E-6&-3.35E-4&-1.67E-3&-6.47E-4&-1.18E-3&-2.33E-3&-2.58E-4&-6.87E-4&-1.88E-3&-1.16E-4&-4.80E-4&-1.73E-3\\
  5&-2.14E-4&-2.87E-3&-8.25E-3&-4.41E-3&-7.18E-3&-1.11E-2&-2.23E-3&-4.83E-3&-9.20E-3&-1.24E-3&-3.72E-3&-8.50E-3\\
  6&-3.57E-3&-2.12E-2&-4.33E-2&-2.80E-2&-3.99E-2&-5.35E-2&-1.74E-2&-3.03E-2&-4.68E-2&-1.16E-2&-2.53E-2&-4.42E-2\\
  7&-4.79E-2&-1.42E-1&-2.19E-1&-1.76E-1&-2.15E-1&-2.52E-1&-1.28E-1&-1.79E-1&-2.30E-1&-9.89E-2&-1.59E-1&-2.22E-1\\
  8&-6.77E-1&-1.14E00&-1.36E00&-1.48E00&-1.49E00&-1.51E00&-1.21E00&-1.32E00&-1.41E00&-1.03E00&-1.22E00&-1.38E00\\
\hline \hline \end{tabular}
%}
\end{table}
\end{tiny}

~ ~ ~

\begin{tiny}
\begin{table}
\setlength{\tabcolsep}{1.0pt}
\renewcommand{\arraystretch}{1.02}
\caption{The poles of ${\bar r}_{q,\delta,\alpha,k}(\t)$,  % ${\bar r}(q,\delta,\alpha,k;\t)$, 
 \eqref{0-URA},   of $(q,\delta,\alpha,k)$-0-URA approximation, 
i.e. the rational approximation of  
$\f(q,\delta,\alpha;\t)=\t^{\alpha}/(1+q\,\t^{\alpha})$ 
with functions from $\mathcal R_k$  of the type \eqref{0-URA} on 
 $ [\delta, 1]$,   $\alpha =0.25, 0.5, 0.75$, ${\bf k=3}$, and $q=400$}
 % $(q=400,\delta,\alpha,k=3)$-0URA-poles}
\label{tabl:Eq400dDaAAk3p}
%\rotatebox{90}{
\begin{tabular}{||c||c|c|c||c|c|c||c|c|c||c|c|c||}
\hline \hline
{\bf k=3} /$\delta$& $0.00$ & $0.00$ & $0.00$ &$10^{-6}$&$10^{-6}$&$10^{-6}$&$10^{-7}$&$10^{-7}$&$10^{-7}$&$10^{-8}$&$10^{-8}$&$10^{-8}$\\
\ j /$\alpha$ & 0.25   & 0.50   & 0.75   & 0.25   & 0.50   & 0.75   & 0.25   & 0.50   & 0.75   & 0.25   & 0.50   & 0.75   \\ \hline
  1&-1.70E-6&-7.38E-5&-5.07E-4&-2.29E-5&-9.56E-5&-5.14E-4&-9.16E-6&-8.06E-5&-5.09E-4&-4.80E-6&-7.59E-5&-5.07E-4\\
  2&-1.38E-3&-8.24E-3&-1.47E-2&-3.53E-3&-9.06E-3&-1.48E-2&-2.43E-3&-8.51E-3&-1.47E-2&-1.92E-3&-8.33E-3&-1.47E-2\\
  3&-1.38E-1&-2.65E-1&-3.15E-1&-1.97E-1&-2.76E-1&-3.16E-1&-1.70E-1&-2.69E-1&-3.15E-1&-1.56E-1&-2.67E-1&-3.15E-1\\
\hline \hline \end{tabular}
%}
\end{table}
\end{tiny}

~ ~ ~

\begin{tiny}
\begin{table}
\setlength{\tabcolsep}{1.0pt}
\renewcommand{\arraystretch}{1.02}
\caption{The poles of ${\bar r}_{q,\delta,\alpha,k}(\t)$,  % ${\bar r}(q,\delta,\alpha,k;\t)$, 
 \eqref{0-URA},   of $(q,\delta,\alpha,k)$-0-URA approximation, 
i.e. the rational approximation of  
$\f(q,\delta,\alpha;\t)=\t^{\alpha}/(1+q\,\t^{\alpha})$ 
with functions from $\mathcal R_k$  of the type \eqref{0-URA} on 
 $ [\delta, 1]$,   $\alpha =0.25, 0.5, 0.75$, ${\bf k=4}$, and $q=400$}
 %$(q=400,\delta,\alpha,k=4)$-0URA-poles}
 \label{tabl:Eq400dDaAAk4p}
%\rotatebox{90}{
\begin{tabular}{||c||c|c|c||c|c|c||c|c|c||c|c|c||}
\hline \hline
{\bf k=4} /$\delta$& $0.00$ & $0.00$ & $0.00$ &$10^{-6}$&$10^{-6}$&$10^{-6}$&$10^{-7}$&$10^{-7}$&$10^{-7}$&$10^{-8}$&$10^{-8}$&$10^{-8}$\\
\ j /$\alpha$ & 0.25   & 0.50   & 0.75   & 0.25   & 0.50   & 0.75   & 0.25   & 0.50   & 0.75   & 0.25   & 0.50   & 0.75   \\ \hline
  1&-8.26E-8&-1.49E-5&-2.42E-4&-6.25E-6&-2.52E-5&-2.51E-4&-1.73E-6&-1.81E-5&-2.44E-4&-6.14E-7&-1.59E-5&-2.42E-4\\
  2&-5.80E-5&-9.14E-4&-2.60E-3&-4.38E-4&-1.22E-3&-2.70E-3&-2.18E-4&-1.01E-3&-2.62E-3&-1.32E-4&-9.46E-4&-2.60E-3\\
  3&-5.31E-3&-2.43E-2&-4.17E-2&-1.46E-2&-2.81E-2&-4.27E-2&-1.01E-2&-2.56E-2&-4.19E-2&-7.87E-3&-2.48E-2&-4.18E-2\\
  4&-2.34E-1&-4.27E-1&-5.11E-1&-3.68E-1&-4.58E-1&-5.17E-1&-3.10E-1&-4.38E-1&-5.12E-1&-2.77E-1&-4.31E-1&-5.11E-1\\
\hline \hline \end{tabular}
%}
\end{table}
\end{tiny}

~ ~ ~

\begin{tiny}
\begin{table}
\setlength{\tabcolsep}{1.0pt}
\renewcommand{\arraystretch}{1.02}
\caption{The poles of ${\bar r}_{q,\delta,\alpha,k}(\t)$,  % ${\bar r}(q,\delta,\alpha,k;\t)$, 
\eqref{0-URA},   of $(q,\delta,\alpha,k)$-0-URA approximation, 
i.e. the rational approximation of  
$\f(q,\delta,\alpha;\t)=\t^{\alpha}/(1+q\,\t^{\alpha})$ 
with functions from $\mathcal R_k$  of the type \eqref{0-URA} on 
 $ [\delta, 1]$,   $\alpha =0.25, 0.5, 0.75$, ${\bf k=5}$, and $q=400$}
 %$(q=400,\delta,\alpha,k=5)$-0URA-poles}
 \label{tabl:Eq400dDaAAk5p}
%\rotatebox{90}{
\begin{tabular}{||c||c|c|c||c|c|c||c|c|c||c|c|c||}
\hline \hline
{\bf k=5} /$\delta$& $0.00$ & $0.00$ & $0.00$ &$10^{-6}$&$10^{-6}$&$10^{-6}$&$10^{-7}$&$10^{-7}$&$10^{-7}$&$10^{-8}$&$10^{-8}$&$10^{-8}$\\
\ j /$\alpha$ & 0.25   & 0.50   & 0.75   & 0.25   & 0.50   & 0.75   & 0.25   & 0.50   & 0.75   & 0.25   & 0.50   & 0.75   \\ \hline
  1&-6.07E-9&-3.85E-6&-9.78E-5&-2.74E-6&-1.01E-5&-1.13E-4&-5.99E-7&-5.82E-6&-1.01E-4&-1.57E-7&-4.48E-6&-9.85E-5\\
  2&-3.49E-6&-1.42E-4&-7.46E-4&-1.04E-4&-2.68E-4&-8.13E-4&-3.85E-5&-1.84E-4&-7.61E-4&-1.74E-5&-1.56E-4&-7.49E-4\\
  3&-3.11E-4&-3.28E-3&-8.06E-3&-2.26E-3&-4.78E-3&-8.66E-3&-1.22E-3&-3.83E-3&-8.20E-3&-7.61E-4&-3.47E-3&-8.09E-3\\
  4&-1.22E-2&-4.72E-2&-7.81E-2&-3.66E-2&-5.82E-2&-8.13E-2&-2.57E-2&-5.14E-2&-7.88E-2&-1.97E-2&-4.87E-2&-7.82E-2\\
  5&-3.38E-1&-5.98E-1&-7.17E-1&-5.84E-1&-6.70E-1&-7.34E-1&-4.86E-1&-6.26E-1&-7.21E-1&-4.25E-1&-6.08E-1&-7.18E-1\\
\hline \hline \end{tabular}
%}
\end{table}
\end{tiny}

~ ~ ~

\begin{tiny}
\begin{table}
\setlength{\tabcolsep}{1.0pt}
\renewcommand{\arraystretch}{1.02}
\caption{The poles of ${\bar r}_{q,\delta,\alpha,k}(\t)$,  % ${\bar r}(q,\delta,\alpha,k;\t)$, 
 \eqref{0-URA},   of $(q,\delta,\alpha,k)$-0-URA approximation, 
i.e. the rational approximation of  
$\f(q,\delta,\alpha;\t)=\t^{\alpha}/(1+q\,\t^{\alpha})$ 
with functions from $\mathcal R_k$  of the type \eqref{0-URA} on 
 $ [\delta, 1]$,   $\alpha =0.25, 0.5, 0.75$, ${\bf k=6}$, and $q=400$}
 % $(q=400,\delta,\alpha,k=6)$-0URA-poles}
 \label{tabl:Eq400dDaAAk6p}
%\rotatebox{90}{
\begin{tabular}{||c||c|c|c||c|c|c||c|c|c||c|c|c||}
\hline \hline
{\bf k=6} /$\delta$& $0.00$ & $0.00$ & $0.00$ &$10^{-6}$&$10^{-6}$&$10^{-6}$&$10^{-7}$&$10^{-7}$&$10^{-7}$&$10^{-8}$&$10^{-8}$&$10^{-8}$\\
\ j /$\alpha$ & 0.25   & 0.50   & 0.75   & 0.25   & 0.50   & 0.75   & 0.25   & 0.50   & 0.75   & 0.25   & 0.50   & 0.75   \\ \hline
  1&-6.05E-10&-1.03E-6&-2.98E-5&-1.53E-6&-5.18E-6&-4.45E-5&-2.93E-7&-2.34E-6&-3.34E-5&-6.31E-8&-1.45E-6&-3.06E-5\\
  2&-2.76E-7&-2.96E-5&-3.19E-4&-3.73E-5&-8.74E-5&-3.71E-4&-1.10E-5&-4.92E-5&-3.32E-4&-3.85E-6&-3.62E-5&-3.21E-4\\
  3&-2.37E-5&-5.54E-4&-2.01E-3&-5.62E-4&-1.17E-3&-2.37E-3&-2.41E-4&-7.86E-4&-2.10E-3&-1.20E-4&-6.36E-4&-2.03E-3\\
  4&-9.12E-4&-7.45E-3&-1.69E-2&-6.77E-3&-1.19E-2&-1.89E-2&-3.87E-3&-9.26E-3&-1.74E-2&-2.47E-3&-8.12E-3&-1.70E-2\\
  5&-2.18E-2&-7.51E-2&-1.21E-1&-7.07E-2&-9.95E-2&-1.29E-1&-5.03E-2&-8.55E-2&-1.23E-1&-3.85E-2&-7.90E-2&-1.21E-1\\
  6&-4.47E-1&-7.74E-1&-9.29E-1&-8.43E-1&-9.13E-1&-9.69E-1&-6.95E-1&-8.34E-1&-9.39E-1&-6.01E-1&-7.97E-1&-9.31E-1\\
\hline \hline \end{tabular}
%}
\end{table}
\end{tiny}

~ ~ ~

\begin{tiny}
\begin{table}
\setlength{\tabcolsep}{1.0pt}
\renewcommand{\arraystretch}{0.95}
\caption{The poles of ${\bar r}_{q,\delta,\alpha,k}(\t)$,  % ${\bar r}(q,\delta,\alpha,k;\t)$,  
\eqref{0-URA},   of $(q,\delta,\alpha,k)$-0-URA approximation, 
i.e. the rational approximation of  
$\f(q,\delta,\alpha;\t)=\t^{\alpha}/(1+q\,\t^{\alpha})$ 
with functions from $\mathcal R_k$  of the type \eqref{0-URA} on 
 $ [\delta, 1]$,   $\alpha =0.25, 0.5, 0.75$, ${\bf k=7}$, and $q=400$}
 %$(q=400,\delta,\alpha,k=7)$-0URA-poles}
 \label{tabl:Eq400dDaAAk7p}
%\rotatebox{90}{
\begin{tabular}{||c||c|c|c||c|c|c||c|c|c||c|c|c||}
\hline \hline
{\bf k=7} /$\delta$& $0.00$ & $0.00$ & $0.00$ &$10^{-6}$&$10^{-6}$&$10^{-6}$&$10^{-7}$&$10^{-7}$&$10^{-7}$&$10^{-8}$&$10^{-8}$&$10^{-8}$\\
\ j /$\alpha$ & 0.25   & 0.50   & 0.75   & 0.25   & 0.50   & 0.75   & 0.25   & 0.50   & 0.75   & 0.25   & 0.50   & 0.75   \\ \hline
  1&-7.44E-11&-2.58E-7&-8.25E-6&-9.87E-7&-3.10E-6&-1.86E-5&-1.73E-7&-1.12E-6&-1.09E-5&-3.31E-8&-5.30E-7&-8.85E-6\\
  2&-2.73E-8&-7.72E-6&-1.39E-4&-1.75E-5&-3.81E-5&-1.94E-4&-4.37E-6&-1.80E-5&-1.55E-4&-1.25E-6&-1.13E-5&-1.42E-4\\
  3&-2.21E-6&-1.14E-4&-6.80E-4&-1.95E-4&-3.83E-4&-8.89E-4&-6.88E-5&-2.16E-4&-7.38E-4&-2.79E-5&-1.52E-4&-6.94E-4\\
  4&-8.39E-5&-1.40E-3&-4.37E-3&-1.82E-3&-3.27E-3&-5.55E-3&-8.60E-4&-2.18E-3&-4.71E-3&-4.55E-4&-1.71E-3&-4.45E-3\\
  5&-1.97E-3&-1.34E-2&-2.86E-2&-1.51E-2&-2.34E-2&-3.37E-2&-9.01E-3&-1.79E-2&-3.01E-2&-5.88E-3&-1.52E-2&-2.90E-2\\
  6&-3.38E-2&-1.07E-1&-1.68E-1&-1.17E-1&-1.52E-1&-1.86E-1&-8.45E-2&-1.28E-1&-1.73E-1&-6.48E-2&-1.16E-1&-1.69E-1\\
  7&-5.60E-1&-9.54E-1&-1.14E00&-1.14E00&-1.19E00&-1.23E00&-9.37E-1&-1.07E00&-1.17E00&-8.02E-1&-1.00E00&-1.15E00\\
\hline \hline \end{tabular}
%}
\end{table}
\end{tiny}

~ ~ ~

\begin{tiny}
\begin{table}
\setlength{\tabcolsep}{1.0pt}
\renewcommand{\arraystretch}{1.02}
\caption{The poles of ${\bar r}_{q,\delta,\alpha,k}(\t)$,  % ${\bar r}(q,\delta,\alpha,k;\t)$, 
 \eqref{0-URA},   of $(q,\delta,\alpha,k)$-0-URA approximation, 
i.e. the rational approximation of  
$\f(q,\delta,\alpha;\t)=\t^{\alpha}/(1+q\,\t^{\alpha})$ 
with functions from $\mathcal R_k$  of the type \eqref{0-URA} on 
 $ [\delta, 1]$,   $\alpha =0.25, 0.5, 0.75$, ${\bf k=8}$, and $q=400$}
 %$(q=400,\delta,\alpha,k=8)$-0URA-poles}
 \label{tabl:Eq400dDaAAk8p}
%\rotatebox{90}{
\begin{tabular}{||c||c|c|c||c|c|c||c|c|c||c|c|c||}
\hline \hline
{\bf k=8} /$\delta$& $0.00$ & $0.00$ & $0.00$ &$10^{-6}$&$10^{-6}$&$10^{-6}$&$10^{-7}$&$10^{-7}$&$10^{-7}$&$10^{-8}$&$10^{-8}$&$10^{-8}$\\
\ j /$\alpha$ & 0.25   & 0.50   & 0.75   & 0.25   & 0.50   & 0.75   & 0.25   & 0.50   & 0.75   & 0.25   & 0.50   & 0.75   \\ \hline
  1&-1.05E-11&-6.26E-8&-2.33E-6&-6.92E-7&-2.06E-6&-9.09E-6&-1.15E-7&-6.17E-7&-4.08E-6&-2.03E-8&-2.25E-7&-2.75E-6\\
  2&-3.25E-9&-2.22E-6&-4.99E-5&-9.83E-6&-2.02E-5&-1.00E-4&-2.16E-6&-8.17E-6&-6.56E-5&-5.29E-7&-4.34E-6&-5.40E-5\\
  3&-2.45E-7&-2.83E-5&-3.03E-4&-8.53E-5&-1.59E-4&-4.43E-4&-2.57E-5&-7.62E-5&-3.48E-4&-8.76E-6&-4.65E-5&-3.15E-4\\
  4&-9.07E-6&-3.05E-4&-1.36E-3&-6.42E-4&-1.13E-3&-2.03E-3&-2.56E-4&-6.46E-4&-1.57E-3&-1.14E-4&-4.44E-4&-1.42E-3\\
  5&-2.11E-4&-2.78E-3&-7.90E-3&-4.39E-3&-7.06E-3&-1.08E-2&-2.22E-3&-4.73E-3&-8.86E-3&-1.24E-3&-3.62E-3&-8.15E-3\\
  6&-3.55E-3&-2.10E-2&-4.28E-2&-2.79E-2&-3.96E-2&-5.31E-2&-1.73E-2&-3.00E-2&-4.63E-2&-1.15E-2&-2.51E-2&-4.38E-2\\
  7&-4.78E-2&-1.41E-1&-2.18E-1&-1.75E-1&-2.14E-1&-2.52E-1&-1.28E-1&-1.78E-1&-2.30E-1&-9.87E-2&-1.58E-1&-2.21E-1\\
  8&-6.76E-1&-1.14E00&-1.36E00&-1.48E00&-1.49E00&-1.51E00&-1.21E00&-1.32E00&-1.41E00&-1.03E00&-1.22E00&-1.38E00\\
\hline \hline \end{tabular}
%}
\end{table}
\end{tiny}

~ ~ ~

%\endinput

%\newpage

\section{Tables type \textbf{(f)} for 0-URA-decomposition coefficients}

%\include{sect1ra-pm4}

%\subsubsection{Tables type \textbf{(f)} for 0URA-decomposition coefficients}

\begin{tiny}
\begin{table}[b!]
\setlength{\tabcolsep}{1.0pt}
\renewcommand{\arraystretch}{1.02}
\caption{The coefficients $ \bar c_j,  j=0, \dots, k$ 
of the partial fraction representation \eqref{trg-bar}
of $(q,\delta,\alpha,k)$-0-URA approximation 
\eqref{0-URA} for  $\delta=0, 10^{-6}, 10^{-7}, 10^{-8}$, $\alpha= 0.25, 0.50, 0.75$,
${\bf k=3}$,  and $q=1$}
 %$(q=001,\delta,\alpha,k=3)$-0-URA-decomposition}
\label{tabl:Fq001dDaAAk3p}
%\rotatebox{90}{

%}
\end{table}
\end{tiny}

~ ~ ~

\begin{tiny}
\begin{table}
\setlength{\tabcolsep}{1.0pt}
\renewcommand{\arraystretch}{1.02}
\caption{The coefficients $  \bar c_j,  j=0, \dots, k$ 
of the partial fraction representation \eqref{trg-bar}
of $(q,\delta,\alpha,k)$-0-URA approximation 
\eqref{0-URA} for  $\delta=0, 10^{-6}, 10^{-7}, 10^{-8}$, $\alpha= 0.25, 0.50, 0.75$,
${\bf k=4}$,  and $q=1$}
 %$(q=001,\delta,\alpha,k=4)$-0-URA-decomposition}
\label{tabl:Fq001dDaAAk4p}
%\rotatebox{90}{
%
%}
\end{table}
\end{tiny}

~ ~ ~

\begin{tiny}
\begin{table}
\setlength{\tabcolsep}{1.0pt}
\renewcommand{\arraystretch}{1.02}
\caption{The coefficients $  \bar c_j,  j=0, \dots, k$ 
of the partial fraction representation \eqref{trg-bar}
of $(q,\delta,\alpha,k)$-0-URA approximation 
\eqref{0-URA} for  $\delta=0, 10^{-6}, 10^{-7}, 10^{-8}$, $\alpha= 0.25, 0.50, 0.75$,
${\bf k=5}$,  and $q=1$}
 %$(q=001,\delta,\alpha,k=5)$-0URA-decomposition}
\label{tabl:Fq001dDaAAk5p}
%\rotatebox{90}{
%
%}
\end{table}
\end{tiny}

~ ~ ~

\begin{tiny}
\begin{table}
\setlength{\tabcolsep}{1.0pt}
\renewcommand{\arraystretch}{1.02}
\caption{The coefficients $  \bar c_j,  j=0, \dots, k$ 
of the partial fraction representation \eqref{trg-bar}
of $(q,\delta,\alpha,k)$-0-URA approximation 
\eqref{0-URA} for  $\delta=0, 10^{-6}, 10^{-7}, 10^{-8}$, $\alpha= 0.25, 0.50, 0.75$,
${\bf k=6}$,  and $q=1$}
 % $(q=001,\delta,\alpha,k=6)$-0URA-decomposition}
\label{tabl:Fq001dDaAAk6p}
%\rotatebox{90}{
%
%}
\end{table}
\end{tiny}

~ ~ ~

\begin{tiny}
\begin{table}
\setlength{\tabcolsep}{1.0pt}
\renewcommand{\arraystretch}{1.02}
\caption{The coefficients $  \bar c_j,  j=0, \dots, k$ 
of the partial fraction representation \eqref{trg-bar}
of $(q,\delta,\alpha,k)$-0-URA approximation 
\eqref{0-URA} for  $\delta=0, 10^{-6}, 10^{-7}, 10^{-8}$, $\alpha= 0.25, 0.50, 0.75$,
${\bf k=7}$,  and $q=1$}
 %$(q=001,\delta,\alpha,k=7)$-0URA-decomposition}
\label{tabl:Fq001dDaAAk7p}
%\rotatebox{90}{
%
%}
\end{table}
\end{tiny}

~ ~ ~

\begin{tiny}
\begin{table}
\setlength{\tabcolsep}{1.0pt}
\renewcommand{\arraystretch}{1.02}
\caption{The coefficients $  \bar c_j,  j=0, \dots, k$ 
of the partial fraction representation \eqref{trg-bar}
of $(q,\delta,\alpha,k)$-0-URA approximation 
\eqref{0-URA} for  $\delta=0, 10^{-6}, 10^{-7}, 10^{-8}$, $\alpha= 0.25, 0.50, 0.75$,
${\bf k=8}$,  and $q=1$}
 %$(q=001,\delta,\alpha,k=8)$-0URA-decomposition}
\label{tabl:Fq001dDaAAk8p}
%\rotatebox{90}{
%
%}
\end{table}
\end{tiny}

~ ~ ~

\begin{tiny}
\begin{table}
\setlength{\tabcolsep}{1.0pt}
\renewcommand{\arraystretch}{1.02}
\caption{The coefficients $ \bar c_j,  j=0, \dots, k$ 
of the partial fraction representation \eqref{trg-bar}
of $(q,\delta,\alpha,k)$-0-URA approximation 
\eqref{0-URA} for  $\delta=0, 10^{-6}, 10^{-7}, 10^{-8}$, $\alpha= 0.25, 0.50, 0.75$,
${\bf k=6}$,  and $q=200$}
 % $(q=200,\delta,\alpha,k=6)$-0URA-decomposition}
\label{tabl:Fq200dDaAAk6p}
%\rotatebox{90}{
%
%}
\end{table}
\end{tiny}

~ ~ ~

\begin{tiny}
\begin{table}
\setlength{\tabcolsep}{1.0pt}
\renewcommand{\arraystretch}{1.02}
\caption{The coefficients $  \bar c_j,  j=0, \dots, k$ 
of the partial fraction representation \eqref{trg-bar}
of $(q,\delta,\alpha,k)$-0-URA approximation 
\eqref{0-URA} for  $\delta=0, 10^{-6}, 10^{-7}, 10^{-8}$, $\alpha= 0.25, 0.50, 0.75$,
${\bf k=3}$,  and $q=400$}
  % $(q=400,\delta,\alpha,k=3)$-0URA-decomposition}
\label{tabl:Fq400dDaAAk3p}
%\rotatebox{90}{
%
%}
\end{table}
\end{tiny}

~ ~ ~

%\endinput

%\newpage
%
%\subsubsection{Tables type \textbf{(g)} for 1-URA-poles}

%\include{sect1ra-pm5}

\section{Tables type \textbf{(g)} for 1-URA-poles}
%\vspace{-3mm}
 
\begin{tiny}
\begin{table}[b!]
\setlength{\tabcolsep}{1.0pt}
\renewcommand{\arraystretch}{1.02}
\caption{The poles of $\bar{\bar r}_{q,\delta,\alpha,k}(\t)$,  %
\eqref{1-URA},   of $(q,\delta,\alpha,k)$-1-URA approximation, 
i.e. the rational approximation of  
$\f(q,\delta,\alpha;\t)=\t^{\alpha}/(1+q\,\t^{\alpha})$ 
with functions from $\mathcal R_k$   on 
 $ [\delta, 1]$,   $\alpha =0.25, 0.5, 0.75$, ${\bf k=3}$, and $q_0=q_1=100$}
   % $(q_0=q_1=100,\delta,\alpha,k=3)$-1URA-poles}
\label{tabl:Gqq11dDaAAk3p}
%\rotatebox{90}{
\begin{tabular}{||c||c|c|c||c|c|c||c|c|c||c|c|c||}
\hline \hline
{\bf k=3} /$\delta$& $0.00$ & $0.00$ & $0.00$ &$10^{-6}$&$10^{-6}$&$10^{-6}$&$10^{-7}$&$10^{-7}$&$10^{-7}$&$10^{-8}$&$10^{-8}$&$10^{-8}$\\
\ j /$\alpha$ & 0.25   & 0.50   & 0.75   & 0.25   & 0.50   & 0.75   & 0.25   & 0.50   & 0.75   & 0.25   & 0.50   & 0.75   \\ \hline
  0& 4.88E-3& 4.95E-3& 4.98E-3& 4.97E-3& 4.97E-3& 4.98E-3& 4.97E-3& 4.96E-3& 4.98E-3& 4.96E-3& 4.95E-3& 4.98E-3\\
  1&-1.73E-14&-2.35E-9&-1.18E-7&-1.21E-9&-1.92E-8&-1.64E-7&-2.51E-10&-7.07E-9&-1.28E-7&-5.17E-11&-3.76E-9&-1.20E-7\\
  2&-1.15E-11&-1.41E-7&-3.74E-6&-2.68E-8&-3.10E-7&-3.93E-6&-7.54E-9&-2.08E-7&-3.79E-6&-2.12E-9&-1.65E-7&-3.75E-6\\
  3&-3.78E-9&-2.02E-6&-8.77E-6&-1.19E-6&-3.24E-6&-8.85E-6&-4.97E-7&-2.53E-6&-8.79E-6&-1.94E-7&-2.21E-6&-8.78E-6\\
\hline \hline \end{tabular}
%}
\end{table}
\end{tiny}

~ ~ ~

\begin{tiny}
\begin{table}
\setlength{\tabcolsep}{1.0pt}
\renewcommand{\arraystretch}{1.02}
\caption{The poles of $\bar{\bar r}_{q,\delta,\alpha,k}(\t)$,  %$\bar{\bar r}(q,\delta,\alpha,k;\t)$, 
 \eqref{1-URA},   of $(q,\delta,\alpha,k)$-1-URA approximation, 
i.e. the rational approximation of  
$\f(q,\delta,\alpha;\t)=\t^{\alpha}/(1+q\,\t^{\alpha})$ 
with functions from $\mathcal R_k$   on 
 $ [\delta, 1]$,   $\alpha =0.25, 0.5, 0.75$, ${\bf k=4}$, and $q_0=q_1=100$}
  % $(q_0=q_1=100,\delta,\alpha,k=4)$-1URA-poles}
  \label{tabl:Gqq11dDaAAk4p}
%\rotatebox{90}{
\begin{tabular}{||c||c|c|c||c|c|c||c|c|c||c|c|c||}
\hline \hline
{\bf k=4} /$\delta$& $0.00$ & $0.00$ & $0.00$ &$10^{-6}$&$10^{-6}$&$10^{-6}$&$10^{-7}$&$10^{-7}$&$10^{-7}$&$10^{-8}$&$10^{-8}$&$10^{-8}$\\
\ j /$\alpha$ & 0.25   & 0.50   & 0.75   & 0.25   & 0.50   & 0.75   & 0.25   & 0.50   & 0.75   & 0.25   & 0.50   & 0.75   \\ \hline
  0& 4.92E-3& 4.97E-3& 4.99E-3& 4.98E-3& 4.98E-3& 4.99E-3& 4.97E-3& 4.97E-3& 4.99E-3& 4.97E-3& 4.97E-3& 4.99E-3\\
  1&-1.23E-15&-2.31E-10&-8.64E-9&-7.02E-10&-7.79E-9&-1.94E-8&-1.42E-10&-1.88E-9&-1.09E-8&-2.77E-11&-6.48E-10&-9.11E-9\\
  2&-7.73E-13&-2.67E-8&-1.14E-6&-8.06E-9&-1.19E-7&-1.49E-6&-2.16E-9&-6.32E-8&-1.23E-6&-5.73E-10&-4.03E-8&-1.16E-6\\
  3&-9.92E-11&-3.57E-7&-4.87E-6&-1.21E-7&-7.57E-7&-4.94E-6&-4.35E-8&-5.37E-7&-4.89E-6&-1.49E-8&-4.31E-7&-4.88E-6\\
  4&-1.80E-8&-3.55E-6&-9.58E-6&-3.00E-6&-5.82E-6&-9.74E-6&-1.62E-6&-4.65E-6&-9.62E-6&-8.05E-7&-4.02E-6&-9.59E-6\\
\hline \hline \end{tabular}
%}
\end{table}
\end{tiny}

~ ~ ~

\begin{tiny}
\begin{table}
\setlength{\tabcolsep}{1.0pt}
\renewcommand{\arraystretch}{1.02}
\caption{The poles of $\bar{\bar r}_{q,\delta,\alpha,k}(\t)$,  %$\bar{\bar r}(q,\delta,\alpha,k;\t)$, 
 \eqref{1-URA},   of $(q,\delta,\alpha,k)$-1-URA approximation, 
i.e. the rational approximation of  
$\f(q,\delta,\alpha;\t)=\t^{\alpha}/(1+q\,\t^{\alpha})$ 
with functions from $\mathcal R_k$   on 
 $ [\delta, 1]$,   $\alpha =0.25, 0.5, 0.75$, ${\bf k=5}$, and $q_0=q_1=100$}
   %$(q_0=q_1=100,\delta,\alpha,k=5)$-1URA-poles}
   \label{tabl:Gqq11dDaAAk5p}
%\rotatebox{90}{
\begin{tabular}{||c||c|c|c||c|c|c||c|c|c||c|c|c||}
\hline \hline
{\bf k=5} /$\delta$& $0.00$ & $0.00$ & $0.00$ &$10^{-6}$&$10^{-6}$&$10^{-6}$&$10^{-7}$&$10^{-7}$&$10^{-7}$&$10^{-8}$&$10^{-8}$&$10^{-8}$\\
\ j /$\alpha$ & 0.25   & 0.50   & 0.75   & 0.25   & 0.50   & 0.75   & 0.25   & 0.50   & 0.75   & 0.25   & 0.50   & 0.75   \\ \hline
  0& 4.94E-3& 4.98E-3& 4.99E-3& 4.98E-3& 4.98E-3& 4.99E-3& 4.98E-3& 4.98E-3& 4.99E-3& 4.97E-3& 4.98E-3& 4.99E-3\\
  1&-1.07E-16&-2.38E-11&-7.55E-10&-4.69E-10&-3.72E-9&-3.42E-9&-9.24E-11&-6.43E-10&-1.28E-9&-1.74E-11&-1.49E-10&-8.64E-10\\
  2&-7.29E-14&-4.26E-9&-1.52E-7&-3.56E-9&-5.39E-8&-3.38E-7&-9.16E-10&-2.17E-8&-1.99E-7&-2.31E-10&-1.03E-8&-1.63E-7\\
  3&-7.66E-12&-8.99E-8&-2.66E-6&-3.11E-8&-3.02E-7&-3.14E-6&-1.02E-8&-1.90E-7&-2.83E-6&-3.28E-9&-1.34E-7&-2.70E-6\\
  4&-4.91E-10&-6.47E-7&-5.02E-6&-3.31E-7&-1.37E-6&-5.05E-6&-1.46E-7&-1.01E-6&-5.03E-6&-6.06E-8&-8.12E-7&-5.02E-6\\
  5&-6.21E-8&-5.25E-6&-1.02E-5&-5.23E-6&-8.47E-6&-1.04E-5&-3.33E-6&-7.00E-6&-1.03E-5&-1.99E-6&-6.09E-6&-1.02E-5\\
\hline \hline \end{tabular}
%}
\end{table}
\end{tiny}

~ ~ ~

\begin{tiny}
\begin{table}
\setlength{\tabcolsep}{1.0pt}
\renewcommand{\arraystretch}{1.02}
\caption{The poles of $\bar{\bar r}_{q,\delta,\alpha,k}(\t)$,  %$\bar{\bar r}(q,\delta,\alpha,k;\t)$, 
 \eqref{1-URA},   of $(q,\delta,\alpha,k)$-1-URA approximation, 
i.e. the rational approximation of  
$\f(q,\delta,\alpha;\t)=\t^{\alpha}/(1+q\,\t^{\alpha})$ 
with functions from $\mathcal R_k$   on 
 $ [\delta, 1]$,   $\alpha =0.25, 0.5, 0.75$, ${\bf k=6}$, and $q_0=q_1=100$}
   % $(q_0=q_1=100,\delta,\alpha,k=6)$-1URA-poles}
   \label{tabl:Gqq11dDaAAk6p}
%\rotatebox{90}{
\begin{tabular}{||c||c|c|c||c|c|c||c|c|c||c|c|c||}
\hline \hline
{\bf k=6} /$\delta$& $0.00$ & $0.00$ & $0.00$ &$10^{-6}$&$10^{-6}$&$10^{-6}$&$10^{-7}$&$10^{-7}$&$10^{-7}$&$10^{-8}$&$10^{-8}$&$10^{-8}$\\
\ j /$\alpha$ & 0.25   & 0.50   & 0.75   & 0.25   & 0.50   & 0.75   & 0.25   & 0.50   & 0.75   & 0.25   & 0.50   & 0.75   \\ \hline
  0& 4.96E-3& 4.98E-3& 4.99E-3& 4.98E-3& 4.99E-3& 4.99E-3& 4.98E-3& 4.98E-3& 4.99E-3& 4.98E-3& 4.98E-3& 4.99E-3\\
  1&-1.07E-17&-2.67E-12&-7.85E-11&-3.40E-10&-2.01E-9&-8.77E-10&-6.57E-11&-2.71E-10&-2.15E-10&-1.20E-11&-4.57E-11&-1.07E-10\\
  2&-8.01E-15&-5.91E-10&-1.57E-8&-1.95E-9&-2.70E-8&-6.76E-8&-4.83E-10&-8.16E-9&-2.80E-8&-1.16E-10&-2.83E-9&-1.86E-8\\
  3&-8.51E-13&-2.06E-8&-7.47E-7&-1.20E-8&-1.48E-7&-1.45E-6&-3.73E-9&-7.90E-8&-9.88E-7&-1.13E-9&-4.59E-8&-8.13E-7\\
  4&-4.02E-11&-1.85E-7&-3.54E-6&-8.25E-8&-5.51E-7&-3.72E-6&-3.23E-8&-3.74E-7&-3.64E-6&-1.22E-8&-2.77E-7&-3.57E-6\\
  5&-1.75E-9&-9.93E-7&-5.10E-6&-6.58E-7&-2.10E-6&-5.18E-6&-3.41E-7&-1.60E-6&-5.12E-6&-1.66E-7&-1.30E-6&-5.10E-6\\
  6&-1.67E-7&-6.93E-6&-1.07E-5&-7.62E-6&-1.09E-5&-1.10E-5&-5.31E-6&-9.26E-6&-1.08E-5&-3.58E-6&-8.17E-6&-1.07E-5\\
\hline \hline \end{tabular}
%}
\end{table}
\end{tiny}

~ ~ ~

\begin{tiny}
\begin{table}
\setlength{\tabcolsep}{1.0pt}
\renewcommand{\arraystretch}{1.02}
\caption{The poles of $\bar{\bar r}_{q,\delta,\alpha,k}(\t)$,  %$\bar{\bar r}(q,\delta,\alpha,k;\t)$, 
 \eqref{1-URA},   of $(q,\delta,\alpha,k)$-1-URA approximation, 
i.e. the rational approximation of  
$\f(q,\delta,\alpha;\t)=\t^{\alpha}/(1+q\,\t^{\alpha})$ 
with functions from $\mathcal R_k$   on 
 $ [\delta, 1]$,   $\alpha =0.25, 0.5, 0.75$, ${\bf k=7}$, and $q_0=q_1=100$}
   % $(q_0=q_1=100,\delta,\alpha,k=7)$-1URA-poles}
   \label{tabl:Gqq11dDaAAk7p}
%\rotatebox{90}{
\begin{tabular}{||c||c|c|c||c|c|c||c|c|c||c|c|c||}
\hline \hline
{\bf k=7} /$\delta$& $0.00$ & $0.00$ & $0.00$ &$10^{-6}$&$10^{-6}$&$10^{-6}$&$10^{-7}$&$10^{-7}$&$10^{-7}$&$10^{-8}$&$10^{-8}$&$10^{-8}$\\
\ j /$\alpha$ & 0.25   & 0.50   & 0.75   & 0.25   & 0.50   & 0.75   & 0.25   & 0.50   & 0.75   & 0.25   & 0.50   & 0.75   \\ \hline
  0& 4.96E-3& 4.98E-3& 4.99E-3& 4.98E-3& 4.99E-3& 5.00E-3& 4.98E-3& 4.99E-3& 5.00E-3& 4.98E-3& 4.99E-3& 4.99E-3\\
  1&-1.19E-18&-3.29E-13&-9.39E-12&-2.61E-10&-1.19E-9&-3.00E-10&-4.95E-11&-1.34E-10&-5.06E-11&-8.77E-12&-1.77E-11&-1.72E-11\\
  2&-9.61E-16&-8.00E-11&-1.81E-9&-1.22E-9&-1.46E-8&-1.65E-8&-2.92E-10&-3.39E-9&-4.95E-9&-6.74E-11&-8.79E-10&-2.57E-9\\
  3&-1.10E-13&-3.90E-9&-1.07E-7&-5.90E-9&-8.05E-8&-4.55E-7&-1.73E-9&-3.51E-8&-2.08E-7&-4.97E-10&-1.61E-8&-1.35E-7\\
  4&-4.90E-12&-5.41E-8&-1.72E-6&-3.04E-8&-2.77E-7&-2.61E-6&-1.10E-8&-1.71E-7&-2.14E-6&-3.84E-9&-1.12E-7&-1.87E-6\\
  5&-1.47E-10&-3.03E-7&-3.74E-6&-1.70E-7&-8.53E-7&-3.70E-6&-7.68E-8&-6.06E-7&-3.73E-6&-3.31E-8&-4.62E-7&-3.74E-6\\
  6&-4.96E-9&-1.38E-6&-5.21E-6&-1.07E-6&-2.85E-6&-5.35E-6&-6.26E-7&-2.24E-6&-5.27E-6&-3.47E-7&-1.85E-6&-5.23E-6\\
  7&-3.66E-7&-8.48E-6&-1.11E-5&-1.01E-5&-1.31E-5&-1.14E-5&-7.41E-6&-1.13E-5&-1.12E-5&-5.36E-6&-1.01E-5&-1.11E-5\\
\hline \hline \end{tabular}
%}
\end{table}
\end{tiny}

~ ~ ~

\begin{tiny}
\begin{table}
\setlength{\tabcolsep}{1.0pt}
\renewcommand{\arraystretch}{1.02}
\caption{The poles of $\bar{\bar r}_{q,\delta,\alpha,k}(\t)$,  %$\bar{\bar r}(q,\delta,\alpha,k;\t)$, 
 \eqref{1-URA},   of $(q,\delta,\alpha,k)$-1-URA approximation, 
i.e. the rational approximation of  
$\f(q,\delta,\alpha;\t)=\t^{\alpha}/(1+q\,\t^{\alpha})$ 
with functions from $\mathcal R_k$   on 
 $ [\delta, 1]$,   $\alpha =0.25, 0.5, 0.75$, ${\bf k=8}$, and $q_0=q_1=100$}
   % $(q_0=q_1=100,\delta,\alpha,k=8)$-1URA-poles}
\label{tabl:Gqq11dDaAAk8p}
%\rotatebox{90}{
\begin{tabular}{||c||c|c|c||c|c|c||c|c|c||c|c|c||}
\hline \hline
{\bf k=8} /$\delta$& $0.00$ & $0.00$ & $0.00$ &$10^{-6}$&$10^{-6}$&$10^{-6}$&$10^{-7}$&$10^{-7}$&$10^{-7}$&$10^{-8}$&$10^{-8}$&$10^{-8}$\\
\ j /$\alpha$ & 0.25   & 0.50   & 0.75   & 0.25   & 0.50   & 0.75   & 0.25   & 0.50   & 0.75   & 0.25   & 0.50   & 0.75   \\ \hline
  0& 3.49E-3& 4.99E-3& 5.00E-3& 4.98E-3& 4.99E-3& 5.00E-3& 4.98E-3& 4.99E-3& 5.00E-3& 4.98E-3& 4.99E-3& 5.00E-3\\
  1& 7.72E-20&-4.44E-14&-1.26E-12&-2.08E-10&-7.65E-10&-1.27E-10&-3.89E-11&-7.45E-11&-1.57E-11&-6.74E-12&-8.18E-12&-3.67E-12\\
  2& 1.80E-17&-1.12E-11&-2.38E-10&-8.36E-10&-8.45E-9&-5.05E-9&-1.94E-10&-1.56E-9&-1.14E-9&-4.31E-11&-3.16E-10&-4.46E-10\\
  3& 1.09E-20&-6.50E-10&-1.37E-8&-3.36E-9&-4.65E-8&-1.25E-7&-9.45E-10&-1.63E-8&-4.13E-8&-2.59E-10&-5.78E-9&-2.11E-8\\
  4&-5.72E-27&-1.34E-8&-4.03E-7&-1.41E-8&-1.57E-7&-1.37E-6&-4.78E-9&-8.53E-8&-7.86E-7&-1.57E-9&-4.75E-8&-5.30E-7\\
  5& 1.50E-13&-1.02E-7&-2.55E-6&-6.23E-8&-4.34E-7&-3.02E-6&-2.56E-8&-2.88E-7&-2.87E-6&-1.01E-8&-2.02E-7&-2.69E-6\\
  6& 1.27E-25&-4.36E-7&-3.71E-6&-2.94E-7&-1.19E-6&-3.67E-6&-1.49E-7&-8.76E-7&-3.68E-6&-7.20E-8&-6.81E-7&-3.70E-6\\
  7& 3.58E-25&-1.78E-6&-5.34E-6&-1.55E-6&-3.59E-6&-5.52E-6&-9.82E-7&-2.91E-6&-5.42E-6&-5.99E-7&-2.44E-6&-5.37E-6\\
  8& 7.10E-24&-9.87E-6&-1.13E-5&-1.28E-5&-1.51E-5&-1.17E-5&-9.60E-6&-1.32E-5&-1.15E-5&-7.23E-6&-1.19E-5&-1.14E-5\\
\hline \hline \end{tabular}
%}
\end{table}
\end{tiny}

~ ~ ~

\begin{tiny}
\begin{table}
\setlength{\tabcolsep}{1.0pt}
\renewcommand{\arraystretch}{1.02}
\caption{The poles of $\bar{\bar r}_{q,\delta,\alpha,k}(\t)$,  %$\bar{\bar r}(q,\delta,\alpha,k;\t)$, 
 \eqref{1-URA},   of $(q,\delta,\alpha,k)$-1-URA approximation, 
i.e. the rational approximation of  
$\f(q,\delta,\alpha;\t)=\t^{\alpha}/(1+q\,\t^{\alpha})$ 
with functions from $\mathcal R_k$   on 
 $ [\delta, 1]$,   $\alpha =0.25, 0.5, 0.75$, ${\bf k=3}$, and $q_0=q_1=200$}
   %$(q_0=q_1=200,\delta,\alpha,k=3)$-1URA-poles}
   \label{tabl:Gqq22dDaAAk3p}
%\rotatebox{90}{
\begin{tabular}{||c||c|c|c||c|c|c||c|c|c||c|c|c||}
\hline \hline
{\bf k=3} /$\delta$& $0.00$ & $0.00$ & $0.00$ &$10^{-6}$&$10^{-6}$&$10^{-6}$&$10^{-7}$&$10^{-7}$&$10^{-7}$&$10^{-8}$&$10^{-8}$&$10^{-8}$\\
\ j /$\alpha$ & 0.25   & 0.50   & 0.75   & 0.25   & 0.50   & 0.75   & 0.25   & 0.50   & 0.75   & 0.25   & 0.50   & 0.75   \\ \hline
  0& 2.45E-3& 2.48E-3& 2.49E-3& 2.49E-3& 2.49E-3& 2.49E-3& 2.49E-3& 2.49E-3& 2.49E-3& 2.49E-3& 2.48E-3& 2.49E-3\\
  1&-5.88E-16&-3.33E-10&-2.80E-8&-3.03E-10&-5.31E-9&-4.72E-8&-6.15E-11&-1.70E-9&-3.22E-8&-1.22E-11&-7.25E-10&-2.89E-8\\
  2&-3.94E-13&-1.92E-8&-7.86E-7&-6.25E-9&-5.92E-8&-8.46E-7&-1.63E-9&-3.54E-8&-8.02E-7&-4.09E-10&-2.54E-8&-7.90E-7\\
  3&-1.37E-10&-2.87E-7&-1.86E-6&-2.78E-7&-6.13E-7&-1.90E-6&-1.07E-7&-4.24E-7&-1.87E-6&-3.61E-8&-3.40E-7&-1.86E-6\\
\hline \hline \end{tabular}
%}
\end{table}
\end{tiny}

~ ~ ~

\begin{tiny}
\begin{table}
\setlength{\tabcolsep}{1.0pt}
\renewcommand{\arraystretch}{1.02}
\caption{The poles of $\bar{\bar r}_{q,\delta,\alpha,k}(\t)$,  %$\bar{\bar r}(q,\delta,\alpha,k;\t)$, 
 \eqref{1-URA},   of $(q,\delta,\alpha,k)$-1-URA approximation, 
i.e. the rational approximation of  
$\f(q,\delta,\alpha;\t)=\t^{\alpha}/(1+q\,\t^{\alpha})$ 
with functions from $\mathcal R_k$   on 
 $ [\delta, 1]$,   $\alpha =0.25, 0.5, 0.75$, ${\bf k=4}$, and $q_0=q_1=200$}
   % $(q_0=q_1=200,\delta,\alpha,k=4)$-1URA-poles}
   \label{tabl:Gqq22dDaAAk4p}
%\rotatebox{90}{
\begin{tabular}{||c||c|c|c||c|c|c||c|c|c||c|c|c||}
\hline \hline
{\bf k=4} /$\delta$& $0.00$ & $0.00$ & $0.00$ &$10^{-6}$&$10^{-6}$&$10^{-6}$&$10^{-7}$&$10^{-7}$&$10^{-7}$&$10^{-8}$&$10^{-8}$&$10^{-8}$\\
\ j /$\alpha$ & 0.25   & 0.50   & 0.75   & 0.25   & 0.50   & 0.75   & 0.25   & 0.50   & 0.75   & 0.25   & 0.50   & 0.75   \\ \hline
  0& 2.47E-3& 2.49E-3& 2.50E-3& 2.49E-3& 2.49E-3& 2.50E-3& 2.49E-3& 2.49E-3& 2.50E-3& 2.49E-3& 2.49E-3& 2.50E-3\\
  1&-4.37E-17&-3.58E-11&-2.27E-9&-1.82E-10&-2.74E-9&-7.28E-9&-3.67E-11&-6.13E-10&-3.29E-9&-7.15E-12&-1.70E-10&-2.49E-9\\
  2&-2.74E-14&-3.92E-9&-2.70E-7&-1.94E-9&-2.61E-8&-3.85E-7&-4.95E-10&-1.29E-8&-3.03E-7&-1.21E-10&-7.39E-9&-2.78E-7\\
  3&-3.60E-12&-5.10E-8&-1.01E-6&-2.87E-8&-1.49E-7&-1.03E-6&-9.62E-9&-9.51E-8&-1.02E-6&-2.96E-9&-6.97E-8&-1.02E-6\\
  4&-7.23E-10&-5.50E-7&-2.12E-6&-7.20E-7&-1.21E-6&-2.20E-6&-3.69E-7&-8.72E-7&-2.15E-6&-1.65E-7&-6.92E-7&-2.13E-6\\
\hline \hline \end{tabular}
%}
\end{table}
\end{tiny}

~ ~ ~

\begin{tiny}
\begin{table}
\setlength{\tabcolsep}{1.0pt}
\renewcommand{\arraystretch}{1.02}
\caption{The poles of $\bar{\bar r}_{q,\delta,\alpha,k}(\t)$,  %$\bar{\bar r}(q,\delta,\alpha,k;\t)$, 
 \eqref{1-URA},   of $(q,\delta,\alpha,k)$-1-URA approximation, 
i.e. the rational approximation of  
$\f(q,\delta,\alpha;\t)=\t^{\alpha}/(1+q\,\t^{\alpha})$ 
with functions from $\mathcal R_k$   on 
 $ [\delta, 1]$,   $\alpha =0.25, 0.5, 0.75$, ${\bf k=5}$, and $q_0=q_1=200$}
   % $(q_0=q_1=200,\delta,\alpha,k=5)$-1URA-poles}
   \label{tabl:Gqq22dDaAAk5p}
%\rotatebox{90}{
\begin{tabular}{||c||c|c|c||c|c|c||c|c|c||c|c|c||}
\hline \hline
{\bf k=5} /$\delta$& $0.00$ & $0.00$ & $0.00$ &$10^{-6}$&$10^{-6}$&$10^{-6}$&$10^{-7}$&$10^{-7}$&$10^{-7}$&$10^{-8}$&$10^{-8}$&$10^{-8}$\\
\ j /$\alpha$ & 0.25   & 0.50   & 0.75   & 0.25   & 0.50   & 0.75   & 0.25   & 0.50   & 0.75   & 0.25   & 0.50   & 0.75   \\ \hline
  0& 2.48E-3& 2.49E-3& 2.50E-3& 2.49E-3& 2.50E-3& 2.50E-3& 2.49E-3& 2.49E-3& 2.50E-3& 2.49E-3& 2.49E-3& 2.50E-3\\
  1&-4.03E-18&-4.04E-12&-2.15E-10&-1.24E-10&-1.58E-9&-1.59E-9&-2.49E-11&-2.67E-10&-4.65E-10&-4.77E-12&-5.19E-11&-2.67E-10\\
  2&-2.72E-15&-6.92E-10&-4.25E-8&-8.80E-10&-1.39E-8&-1.20E-7&-2.19E-10&-5.53E-9&-6.31E-8&-5.26E-11&-2.39E-9&-4.74E-8\\
  3&-2.87E-13&-1.36E-8&-5.98E-7&-7.48E-9&-6.33E-8&-7.16E-7&-2.35E-9&-3.70E-8&-6.46E-7&-6.91E-10&-2.42E-8&-6.12E-7\\
  4&-1.92E-11&-9.79E-8&-1.06E-6&-7.96E-8&-2.86E-7&-1.09E-6&-3.33E-8&-1.91E-7&-1.07E-6&-1.26E-8&-1.41E-7&-1.06E-6\\
  5&-2.84E-9&-8.91E-7&-2.34E-6&-1.27E-6&-1.88E-6&-2.45E-6&-7.85E-7&-1.44E-6&-2.38E-6&-4.40E-7&-1.17E-6&-2.35E-6\\
\hline \hline \end{tabular}
%}
\end{table}
\end{tiny}

~ ~ ~

\begin{tiny}
\begin{table}
\setlength{\tabcolsep}{1.0pt}
\renewcommand{\arraystretch}{1.02}
\caption{The poles of $\bar{\bar r}_{q,\delta,\alpha,k}(\t)$,  %$\bar{\bar r}(q,\delta,\alpha,k;\t)$, 
 \eqref{1-URA},   of $(q,\delta,\alpha,k)$-1-URA approximation, 
i.e. the rational approximation of  
$\f(q,\delta,\alpha;\t)=\t^{\alpha}/(1+q\,\t^{\alpha})$ 
with functions from $\mathcal R_k$   on 
 $ [\delta, 1]$,   $\alpha =0.25, 0.5, 0.75$, ${\bf k=6}$, and $q_0=q_1=200$}
   % $(q_0=q_1=200,\delta,\alpha,k=6)$-1URA-poles}
   \label{tabl:Gqq22dDaAAk6p}
%\rotatebox{90}{
\begin{tabular}{||c||c|c|c||c|c|c||c|c|c||c|c|c||}
\hline \hline
{\bf k=6} /$\delta$& $0.00$ & $0.00$ & $0.00$ &$10^{-6}$&$10^{-6}$&$10^{-6}$&$10^{-7}$&$10^{-7}$&$10^{-7}$&$10^{-8}$&$10^{-8}$&$10^{-8}$\\
\ j /$\alpha$ & 0.25   & 0.50   & 0.75   & 0.25   & 0.50   & 0.75   & 0.25   & 0.50   & 0.75   & 0.25   & 0.50   & 0.75   \\ \hline
  0& 2.48E-3& 2.49E-3& 2.50E-3& 2.50E-3& 2.50E-3& 2.50E-3& 2.49E-3& 2.50E-3& 2.50E-3& 2.49E-3& 2.50E-3& 2.50E-3\\
  1&-4.32E-19&-4.96E-13&-2.40E-11&-9.16E-11&-9.85E-10&-4.79E-10&-1.82E-11&-1.34E-10&-9.38E-11&-3.45E-12&-1.98E-11&-3.78E-11\\
  2&-3.19E-16&-1.07E-10&-4.83E-9&-4.91E-10&-8.19E-9&-3.07E-8&-1.20E-10&-2.58E-9&-1.08E-8&-2.80E-11&-8.37E-10&-6.26E-9\\
  3&-3.37E-14&-3.43E-9&-2.03E-7&-2.95E-9&-3.36E-8&-4.16E-7&-8.78E-10&-1.76E-8&-2.86E-7&-2.49E-10&-9.76E-9&-2.28E-7\\
  4&-1.62E-12&-2.89E-8&-7.60E-7&-2.00E-8&-1.17E-7&-7.85E-7&-7.50E-9&-7.37E-8&-7.77E-7&-2.62E-9&-5.10E-8&-7.67E-7\\
  5&-7.59E-11&-1.60E-7&-1.11E-6&-1.60E-7&-4.57E-7&-1.16E-6&-7.97E-8&-3.21E-7&-1.13E-6&-3.63E-8&-2.42E-7&-1.11E-6\\
  6&-8.98E-9&-1.27E-6&-2.51E-6&-1.87E-6&-2.51E-6&-2.62E-6&-1.28E-6&-2.02E-6&-2.56E-6&-8.24E-7&-1.69E-6&-2.52E-6\\
\hline \hline \end{tabular}
%}
\end{table}
\end{tiny}

~ ~ ~

\begin{tiny}
\begin{table}
\setlength{\tabcolsep}{1.0pt}
\renewcommand{\arraystretch}{0.95}
\caption{The poles of $\bar{\bar r}_{q,\delta,\alpha,k}(\t)$,  %$\bar{\bar r}(q,\delta,\alpha,k;\t)$, 
 \eqref{1-URA},   of $(q,\delta,\alpha,k)$-1-URA approximation, 
i.e. the rational approximation of  
$\f(q,\delta,\alpha;\t)=\t^{\alpha}/(1+q\,\t^{\alpha})$ 
with functions from $\mathcal R_k$   on 
 $ [\delta, 1]$,   $\alpha =0.25, 0.5, 0.75$, ${\bf k=7}$, and $q_0=q_1=200$}
   % $(q_0=q_1=200,\delta,\alpha,k=7)$-1URA-poles}
   \label{tabl:Gqq22dDaAAk7p}
%\rotatebox{90}{
\begin{tabular}{||c||c|c|c||c|c|c||c|c|c||c|c|c||}
\hline \hline
{\bf k=7} /$\delta$& $0.00$ & $0.00$ & $0.00$ &$10^{-6}$&$10^{-6}$&$10^{-6}$&$10^{-7}$&$10^{-7}$&$10^{-7}$&$10^{-8}$&$10^{-8}$&$10^{-8}$\\
\ j /$\alpha$ & 0.25   & 0.50   & 0.75   & 0.25   & 0.50   & 0.75   & 0.25   & 0.50   & 0.75   & 0.25   & 0.50   & 0.75   \\ \hline
  0& 2.48E-3& 2.49E-3& 2.50E-3& 2.50E-3& 2.50E-3& 2.50E-3& 2.50E-3& 2.50E-3& 2.50E-3& 2.49E-3& 2.50E-3& 2.50E-3\\
  1&-3.97E-18&-6.61E-14&-3.04E-12&-7.13E-11&-6.54E-10&-1.85E-10&-1.41E-11&-7.54E-11&-2.58E-11&-2.63E-12&-9.04E-12&-7.08E-12\\
  2&-2.69E-15&-1.58E-11&-5.92E-10&-3.12E-10&-5.15E-9&-8.76E-9&-7.44E-11&-1.30E-9&-2.20E-9&-1.70E-11&-3.21E-10&-9.69E-10\\
  3&-2.86E-13&-7.27E-10&-3.45E-8&-1.46E-9&-2.01E-8&-1.81E-7&-4.19E-10&-9.13E-9&-8.03E-8&-1.15E-10&-4.15E-9&-4.76E-8\\
  4& 9.97E-25&-9.08E-9&-4.30E-7&-7.44E-9&-6.11E-8&-6.14E-7&-2.59E-9&-3.60E-8&-5.36E-7&-8.52E-10&-2.27E-8&-4.72E-7\\
  5&-1.93E-11&-4.89E-8&-7.86E-7&-4.15E-8&-1.86E-7&-7.89E-7&-1.81E-8&-1.23E-7&-7.86E-7&-7.33E-9&-8.80E-8&-7.86E-7\\
  6&-2.85E-9&-2.35E-7&-1.16E-6&-2.63E-7&-6.44E-7&-1.24E-6&-1.49E-7&-4.76E-7&-1.20E-6&-7.84E-8&-3.68E-7&-1.18E-6\\
  7& 9.43E-25&-1.65E-6&-2.63E-6&-2.50E-6&-3.07E-6&-2.75E-6&-1.80E-6&-2.56E-6&-2.69E-6&-1.26E-6&-2.20E-6&-2.65E-6\\
\hline \hline \end{tabular}
%}
\end{table}
\end{tiny}

~ ~ ~

\begin{tiny}
\begin{table}
\setlength{\tabcolsep}{1.0pt}
\renewcommand{\arraystretch}{1.00}
\caption{The poles of $\bar{\bar r}_{q,\delta,\alpha,k}(\t)$,  %$\bar{\bar r}(q,\delta,\alpha,k;\t)$,  
\eqref{1-URA},   of $(q,\delta,\alpha,k)$-1-URA approximation, 
i.e. the rational approximation of  
$\f(q,\delta,\alpha;\t)=\t^{\alpha}/(1+q\,\t^{\alpha})$ 
with functions from $\mathcal R_k$   on 
 $ [\delta, 1]$,   $\alpha =0.25, 0.5, 0.75$, ${\bf k=8}$, and $q_0=q_1=200$}
   % $(q_0=q_1=200,\delta,\alpha,k=8)$-1URA-poles}
   \label{tabl:Gqq22dDaAAk8p}
%\rotatebox{90}{
\begin{tabular}{||c||c|c|c||c|c|c||c|c|c||c|c|c||}
\hline \hline
{\bf k=8} /$\delta$& $0.00$ & $0.00$ & $0.00$ &$10^{-6}$&$10^{-6}$&$10^{-6}$&$10^{-7}$&$10^{-7}$&$10^{-7}$&$10^{-8}$&$10^{-8}$&$10^{-8}$\\
\ j /$\alpha$ & 0.25   & 0.50   & 0.75   & 0.25   & 0.50   & 0.75   & 0.25   & 0.50   & 0.75   & 0.25   & 0.50   & 0.75   \\ \hline
  0& 1.66E-3& 2.50E-3& 2.50E-3& 2.50E-3& 2.50E-3& 2.50E-3& 2.50E-3& 2.50E-3& 2.50E-3& 2.50E-3& 2.50E-3& 2.50E-3\\
  1& 1.52E-22&-9.55E-15&-4.29E-13&-5.76E-11&-4.55E-10&-8.51E-11&-1.13E-11&-4.61E-11&-9.06E-12&-2.09E-12&-4.72E-12&-1.74E-12\\
  2&-1.70E-22&-2.40E-12&-8.16E-11&-2.17E-10&-3.40E-9&-3.00E-9&-5.05E-11&-7.02E-10&-5.68E-10&-1.13E-11&-1.37E-10&-1.89E-10\\
  3& 8.43E-27&-1.34E-10&-4.73E-9&-8.44E-10&-1.30E-8&-6.40E-8&-2.33E-10&-5.01E-9&-1.90E-8&-6.20E-11&-1.82E-9&-8.55E-9\\
  4&-1.24E-27&-2.52E-9&-1.28E-7&-3.49E-9&-3.65E-8&-4.15E-7&-1.15E-9&-1.98E-8&-2.64E-7&-3.60E-10&-1.10E-8&-1.78E-7\\
  5& 2.05E-27&-1.73E-8&-5.85E-7&-1.53E-8&-9.62E-8&-6.38E-7&-6.08E-9&-6.04E-8&-6.34E-7&-2.28E-9&-4.05E-8&-6.13E-7\\
  6& 7.47E-27&-7.30E-8&-7.87E-7&-7.20E-8&-2.67E-7&-8.09E-7&-3.55E-8&-1.85E-7&-7.94E-7&-1.63E-8&-1.35E-7&-7.89E-7\\
  7& 3.06E-26&-3.19E-7&-1.22E-6&-3.81E-7&-8.31E-7&-1.30E-6&-2.37E-7&-6.40E-7&-1.26E-6&-1.39E-7&-5.09E-7&-1.24E-6\\
  8& 6.51E-24&-2.01E-6&-2.72E-6&-3.15E-6&-3.59E-6&-2.86E-6&-2.35E-6&-3.06E-6&-2.79E-6&-1.73E-6&-2.67E-6&-2.75E-6\\
\hline \hline \end{tabular}
%}
\end{table}
\end{tiny}

~ ~ ~

%\endinput

%\newpage

\section{Tables type \textbf{(h)} for 1-URA-decomposition coefficients}

%\include{sect1ra-pm6}

%\subsubsection{Tables type \textbf{(h)} for 1URA-decomposition coefficients}

%
%\begin{tiny}
%\begin{table}
%\setlength{\tabcolsep}{1.0pt}
\renewcommand{\arraystretch}{1.02}
%\caption{The poles of ${\bar{\bar r}}(q,\delta,\alpha,k;\t)$,  \eqref{1-URA},   of $(q,\delta,\alpha,k)$-0-URA approximation, 
%i.e. the rational approximation of  
%$\f(q,\delta,\alpha;\t)=\t^{\alpha}/(1+q\,\t^{\alpha})$ 
%with functions from $\mathcal R_k$   on 
% $ [\delta, 1]$,   $\alpha =0.25, 0.5, 0.75$, ${\bf k=3}$, and $q_0=q_1=100$
% % $(q_0=q_1=100,\delta,\alpha,k=3)$-1URA-decomposition}
% \label{tabl:Hqq11dDaAAk3p}
%%\rotatebox{90}{
%\begin{tabular}{||c||c|c|c||c|c|c||c|c|c||c|c|c||}
%\hline \hline
%k3 /$\delta$& $0.00$ & $0.00$ & $0.00$ &$10^{-6}$&$10^{-6}$&$10^{-6}$&$10^{-7}$&$10^{-7}$&$10^{-7}$&$10^{-8}$&$10^{-8}$&$10^{-8}$\\
%\ j /$\alpha$ & 0.25   & 0.50   & 0.75   & 0.25   & 0.50   & 0.75   & 0.25   & 0.50   & 0.75   & 0.25   & 0.50   & 0.75   \\ \hline
%  0& 4.88E-3& 4.95E-3& 4.98E-3& 4.97E-3& 4.97E-3& 4.98E-3& 4.97E-3& 4.96E-3& 4.98E-3& 4.96E-3& 4.95E-3& 4.98E-3\\
%  1&-1.73E-14&-2.35E-9&-1.18E-7&-1.21E-9&-1.92E-8&-1.64E-7&-2.51E-10&-7.07E-9&-1.28E-7&-5.17E-11&-3.76E-9&-1.20E-7\\
%  2&-1.15E-11&-1.41E-7&-3.74E-6&-2.68E-8&-3.10E-7&-3.93E-6&-7.54E-9&-2.08E-7&-3.79E-6&-2.12E-9&-1.65E-7&-3.75E-6\\
%  3&-3.78E-9&-2.02E-6&-8.77E-6&-1.19E-6&-3.24E-6&-8.85E-6&-4.97E-7&-2.53E-6&-8.79E-6&-1.94E-7&-2.21E-6&-8.78E-6\\
%\hline \hline \end{tabular}
%%}
%\end{table}
%\end{tiny}
%

\begin{tiny}
\begin{table}[b!]
\setlength{\tabcolsep}{1.0pt}
\renewcommand{\arraystretch}{1.02}
\caption{The coefficients $\bar{\bar c}_j,  j=0, \dots, k$  of the partial fraction representation \eqref{frac-1-URA}
of $(q,\delta,\alpha,k)$-1-URA approximation 
\eqref{1-URA} for  $\delta=0, 10^{-6}, 10^{-7}, 10^{-8}$, $\alpha= 0.25, 0.50, 0.75$,
${\bf k=4}$,  and $q_0=q_1=100$ } 
% \red{Maybe we need to introduce the poles and the coefficients ????}
%$(q_0=q_1=100,\delta,\alpha,k=4)$-1URA-decomposition}
  \label{tabl:Hqq11dDaAAk4p}
%\rotatebox{90}{
\begin{tabular}{||c||c|c|c||c|c|c||c|c|c||c|c|c||}
\hline \hline
{\bf k=4} /$\delta$& $0.00$ & $0.00$ & $0.00$ &$10^{-6}$&$10^{-6}$&$10^{-6}$&$10^{-7}$&$10^{-7}$&$10^{-7}$&$10^{-8}$&$10^{-8}$&$10^{-8}$\\
\ j /$\alpha$ & 0.25   & 0.50   & 0.75   & 0.25   & 0.50   & 0.75   & 0.25   & 0.50   & 0.75   & 0.25   & 0.50   & 0.75   \\ \hline
  0& 4.92E-3& 4.97E-3& 4.99E-3& 4.98E-3& 4.98E-3& 4.99E-3& 4.97E-3& 4.97E-3& 4.99E-3& 4.97E-3& 4.97E-3& 4.99E-3\\
  1&-1.23E-15&-2.31E-10&-8.64E-9&-7.02E-10&-7.79E-9&-1.94E-8&-1.42E-10&-1.88E-9&-1.09E-8&-2.77E-11&-6.48E-10&-9.11E-9\\
  2&-7.73E-13&-2.67E-8&-1.14E-6&-8.06E-9&-1.19E-7&-1.49E-6&-2.16E-9&-6.32E-8&-1.23E-6&-5.73E-10&-4.03E-8&-1.16E-6\\
  3&-9.92E-11&-3.57E-7&-4.87E-6&-1.21E-7&-7.57E-7&-4.94E-6&-4.35E-8&-5.37E-7&-4.89E-6&-1.49E-8&-4.31E-7&-4.88E-6\\
  4&-1.80E-8&-3.55E-6&-9.58E-6&-3.00E-6&-5.82E-6&-9.74E-6&-1.62E-6&-4.65E-6&-9.62E-6&-8.05E-7&-4.02E-6&-9.59E-6\\
\hline \hline \end{tabular}
%}
\end{table}
\end{tiny}

~ ~ ~

\begin{tiny}
\begin{table}
\setlength{\tabcolsep}{1.0pt}
\renewcommand{\arraystretch}{1.02}
\caption{The coefficients $\bar{\bar c}_j,  j=0, \dots, k$  of the partial fraction representation \eqref{frac-1-URA}
of $(q,\delta,\alpha,k)$-1-URA approximation 
\eqref{1-URA} for  $\delta=0, 10^{-6}, 10^{-7}, 10^{-8}$, $\alpha= 0.25, 0.50, 0.75$,
${\bf k=5}$,  and $q_0=q_1=100$ } 
%$(q_0=q_1=100,\delta,\alpha,k=5)$-1URA-decomposition}
\label{tabl:Hqq11dDaAAk5p}
%\rotatebox{90}{
\begin{tabular}{||c||c|c|c||c|c|c||c|c|c||c|c|c||}
\hline \hline
{\bf k=5} /$\delta$& $0.00$ & $0.00$ & $0.00$ &$10^{-6}$&$10^{-6}$&$10^{-6}$&$10^{-7}$&$10^{-7}$&$10^{-7}$&$10^{-8}$&$10^{-8}$&$10^{-8}$\\
\ j /$\alpha$ & 0.25   & 0.50   & 0.75   & 0.25   & 0.50   & 0.75   & 0.25   & 0.50   & 0.75   & 0.25   & 0.50   & 0.75   \\ \hline
  0& 4.94E-3& 4.98E-3& 4.99E-3& 4.98E-3& 4.98E-3& 4.99E-3& 4.98E-3& 4.98E-3& 4.99E-3& 4.97E-3& 4.98E-3& 4.99E-3\\
  1&-1.07E-16&-2.38E-11&-7.55E-10&-4.69E-10&-3.72E-9&-3.42E-9&-9.24E-11&-6.43E-10&-1.28E-9&-1.74E-11&-1.49E-10&-8.64E-10\\
  2&-7.29E-14&-4.26E-9&-1.52E-7&-3.56E-9&-5.39E-8&-3.38E-7&-9.16E-10&-2.17E-8&-1.99E-7&-2.31E-10&-1.03E-8&-1.63E-7\\
  3&-7.66E-12&-8.99E-8&-2.66E-6&-3.11E-8&-3.02E-7&-3.14E-6&-1.02E-8&-1.90E-7&-2.83E-6&-3.28E-9&-1.34E-7&-2.70E-6\\
  4&-4.91E-10&-6.47E-7&-5.02E-6&-3.31E-7&-1.37E-6&-5.05E-6&-1.46E-7&-1.01E-6&-5.03E-6&-6.06E-8&-8.12E-7&-5.02E-6\\
  5&-6.21E-8&-5.25E-6&-1.02E-5&-5.23E-6&-8.47E-6&-1.04E-5&-3.33E-6&-7.00E-6&-1.03E-5&-1.99E-6&-6.09E-6&-1.02E-5\\
\hline \hline \end{tabular}
%}
\end{table}
\end{tiny}

~ ~ ~

\begin{tiny}
\begin{table}
\setlength{\tabcolsep}{1.0pt}
\renewcommand{\arraystretch}{1.02}
\caption{The coefficients $\bar{\bar c}_j,  j=0, \dots, k$  of the partial fraction representation \eqref{frac-1-URA}
of $(q,\delta,\alpha,k)$-1-URA approximation 
\eqref{1-URA} for  $\delta=0, 10^{-6}, 10^{-7}, 10^{-8}$, $\alpha= 0.25, 0.50, 0.75$,
${\bf k=6}$,  and $q_0=q_1=100$ } 
  % $(q_0=q_1=100,\delta,\alpha,k=6)$-1URA-decomposition}
\label{tabl:Hqq11dDaAAk6p}
%\rotatebox{90}{
\begin{tabular}{||c||c|c|c||c|c|c||c|c|c||c|c|c||}
\hline \hline
{\bf k=6} /$\delta$& $0.00$ & $0.00$ & $0.00$ &$10^{-6}$&$10^{-6}$&$10^{-6}$&$10^{-7}$&$10^{-7}$&$10^{-7}$&$10^{-8}$&$10^{-8}$&$10^{-8}$\\
\ j /$\alpha$ & 0.25   & 0.50   & 0.75   & 0.25   & 0.50   & 0.75   & 0.25   & 0.50   & 0.75   & 0.25   & 0.50   & 0.75   \\ \hline
  0& 4.96E-3& 4.98E-3& 4.99E-3& 4.98E-3& 4.99E-3& 4.99E-3& 4.98E-3& 4.98E-3& 4.99E-3& 4.98E-3& 4.98E-3& 4.99E-3\\
  1&-1.07E-17&-2.67E-12&-7.85E-11&-3.40E-10&-2.01E-9&-8.77E-10&-6.57E-11&-2.71E-10&-2.15E-10&-1.20E-11&-4.57E-11&-1.07E-10\\
  2&-8.01E-15&-5.91E-10&-1.57E-8&-1.95E-9&-2.70E-8&-6.76E-8&-4.83E-10&-8.16E-9&-2.80E-8&-1.16E-10&-2.83E-9&-1.86E-8\\
  3&-8.51E-13&-2.06E-8&-7.47E-7&-1.20E-8&-1.48E-7&-1.45E-6&-3.73E-9&-7.90E-8&-9.88E-7&-1.13E-9&-4.59E-8&-8.13E-7\\
  4&-4.02E-11&-1.85E-7&-3.54E-6&-8.25E-8&-5.51E-7&-3.72E-6&-3.23E-8&-3.74E-7&-3.64E-6&-1.22E-8&-2.77E-7&-3.57E-6\\
  5&-1.75E-9&-9.93E-7&-5.10E-6&-6.58E-7&-2.10E-6&-5.18E-6&-3.41E-7&-1.60E-6&-5.12E-6&-1.66E-7&-1.30E-6&-5.10E-6\\
  6&-1.67E-7&-6.93E-6&-1.07E-5&-7.62E-6&-1.09E-5&-1.10E-5&-5.31E-6&-9.26E-6&-1.08E-5&-3.58E-6&-8.17E-6&-1.07E-5\\
\hline \hline \end{tabular}
%}
\end{table}
\end{tiny}

~ ~ ~

\begin{tiny}
\begin{table}
\setlength{\tabcolsep}{1.0pt}
\renewcommand{\arraystretch}{1.02}
\caption{The coefficients $\bar{\bar c}_j,  j=0, \dots, k$  of the partial fraction representation \eqref{frac-1-URA}
of $(q,\delta,\alpha,k)$-1-URA approximation 
\eqref{1-URA} for  $\delta=0, 10^{-6}, 10^{-7}, 10^{-8}$, $\alpha= 0.25, 0.50, 0.75$,
${\bf k=7}$,  and $q_0=q_1=100$ } 
  % $(q_0=q_1=100,\delta,\alpha,k=7)$-1URA-decomposition}
\label{tabl:Hqq11dDaAAk7p}
%\rotatebox{90}{
\begin{tabular}{||c||c|c|c||c|c|c||c|c|c||c|c|c||}
\hline \hline
{\bf k=7} /$\delta$& $0.00$ & $0.00$ & $0.00$ &$10^{-6}$&$10^{-6}$&$10^{-6}$&$10^{-7}$&$10^{-7}$&$10^{-7}$&$10^{-8}$&$10^{-8}$&$10^{-8}$\\
\ j /$\alpha$ & 0.25   & 0.50   & 0.75   & 0.25   & 0.50   & 0.75   & 0.25   & 0.50   & 0.75   & 0.25   & 0.50   & 0.75   \\ \hline
  0& 4.96E-3& 4.98E-3& 4.99E-3& 4.98E-3& 4.99E-3& 5.00E-3& 4.98E-3& 4.99E-3& 5.00E-3& 4.98E-3& 4.99E-3& 4.99E-3\\
  1&-1.19E-18&-3.29E-13&-9.39E-12&-2.61E-10&-1.19E-9&-3.00E-10&-4.95E-11&-1.34E-10&-5.06E-11&-8.77E-12&-1.77E-11&-1.72E-11\\
  2&-9.61E-16&-8.00E-11&-1.81E-9&-1.22E-9&-1.46E-8&-1.65E-8&-2.92E-10&-3.39E-9&-4.95E-9&-6.74E-11&-8.79E-10&-2.57E-9\\
  3&-1.10E-13&-3.90E-9&-1.07E-7&-5.90E-9&-8.05E-8&-4.55E-7&-1.73E-9&-3.51E-8&-2.08E-7&-4.97E-10&-1.61E-8&-1.35E-7\\
  4&-4.90E-12&-5.41E-8&-1.72E-6&-3.04E-8&-2.77E-7&-2.61E-6&-1.10E-8&-1.71E-7&-2.14E-6&-3.84E-9&-1.12E-7&-1.87E-6\\
  5&-1.47E-10&-3.03E-7&-3.74E-6&-1.70E-7&-8.53E-7&-3.70E-6&-7.68E-8&-6.06E-7&-3.73E-6&-3.31E-8&-4.62E-7&-3.74E-6\\
  6&-4.96E-9&-1.38E-6&-5.21E-6&-1.07E-6&-2.85E-6&-5.35E-6&-6.26E-7&-2.24E-6&-5.27E-6&-3.47E-7&-1.85E-6&-5.23E-6\\
  7&-3.66E-7&-8.48E-6&-1.11E-5&-1.01E-5&-1.31E-5&-1.14E-5&-7.41E-6&-1.13E-5&-1.12E-5&-5.36E-6&-1.01E-5&-1.11E-5\\
\hline \hline \end{tabular}
%}
\end{table}
\end{tiny}

~ ~ ~

\begin{tiny}
\begin{table}
\setlength{\tabcolsep}{1.0pt}
\renewcommand{\arraystretch}{1.02}
\caption{The coefficients $\bar{\bar c}_j,  j=0, \dots, k$  of the partial fraction representation \eqref{frac-1-URA}
of $(q,\delta,\alpha,k)$-1-URA approximation 
\eqref{1-URA} for  $\delta=0, 10^{-6}, 10^{-7}, 10^{-8}$, $\alpha= 0.25, 0.50, 0.75$,
${\bf k=8}$,  and $q_0=q_1=100$ } 
  % $(q_0=q_1=100,\delta,\alpha,k=8)$-1URA-decomposition}
\label{tabl:Hqq11dDaAAk8p}
%\rotatebox{90}{
\begin{tabular}{||c||c|c|c||c|c|c||c|c|c||c|c|c||}
\hline \hline
{\bf k=8} /$\delta$& $0.00$ & $0.00$ & $0.00$ &$10^{-6}$&$10^{-6}$&$10^{-6}$&$10^{-7}$&$10^{-7}$&$10^{-7}$&$10^{-8}$&$10^{-8}$&$10^{-8}$\\
\ j /$\alpha$ & 0.25   & 0.50   & 0.75   & 0.25   & 0.50   & 0.75   & 0.25   & 0.50   & 0.75   & 0.25   & 0.50   & 0.75   \\ \hline
  0& 3.49E-3& 4.99E-3& 5.00E-3& 4.98E-3& 4.99E-3& 5.00E-3& 4.98E-3& 4.99E-3& 5.00E-3& 4.98E-3& 4.99E-3& 5.00E-3\\
  1& 7.72E-20&-4.44E-14&-1.26E-12&-2.08E-10&-7.65E-10&-1.27E-10&-3.89E-11&-7.45E-11&-1.57E-11&-6.74E-12&-8.18E-12&-3.67E-12\\
  2& 1.80E-17&-1.12E-11&-2.38E-10&-8.36E-10&-8.45E-9&-5.05E-9&-1.94E-10&-1.56E-9&-1.14E-9&-4.31E-11&-3.16E-10&-4.46E-10\\
  3& 1.09E-20&-6.50E-10&-1.37E-8&-3.36E-9&-4.65E-8&-1.25E-7&-9.45E-10&-1.63E-8&-4.13E-8&-2.59E-10&-5.78E-9&-2.11E-8\\
  4&-5.72E-27&-1.34E-8&-4.03E-7&-1.41E-8&-1.57E-7&-1.37E-6&-4.78E-9&-8.53E-8&-7.86E-7&-1.57E-9&-4.75E-8&-5.30E-7\\
  5& 1.50E-13&-1.02E-7&-2.55E-6&-6.23E-8&-4.34E-7&-3.02E-6&-2.56E-8&-2.88E-7&-2.87E-6&-1.01E-8&-2.02E-7&-2.69E-6\\
  6& 1.27E-25&-4.36E-7&-3.71E-6&-2.94E-7&-1.19E-6&-3.67E-6&-1.49E-7&-8.76E-7&-3.68E-6&-7.20E-8&-6.81E-7&-3.70E-6\\
  7& 3.58E-25&-1.78E-6&-5.34E-6&-1.55E-6&-3.59E-6&-5.52E-6&-9.82E-7&-2.91E-6&-5.42E-6&-5.99E-7&-2.44E-6&-5.37E-6\\
  8& 7.10E-24&-9.87E-6&-1.13E-5&-1.28E-5&-1.51E-5&-1.17E-5&-9.60E-6&-1.32E-5&-1.15E-5&-7.23E-6&-1.19E-5&-1.14E-5\\
\hline \hline \end{tabular}
%}
\end{table}
\end{tiny}

~ ~ ~

\begin{tiny}
\begin{table}
\setlength{\tabcolsep}{1.0pt}
\renewcommand{\arraystretch}{1.02}
\caption{The coefficients $\bar{\bar c}_j,  j=0, \dots, k$  of the partial fraction representation \eqref{frac-1-URA}
of $(q,\delta,\alpha,k)$-1-URA approximation 
\eqref{1-URA} for  $\delta=0, 10^{-6}, 10^{-7}, 10^{-8}$, $\alpha= 0.25, 0.50, 0.75$,
${\bf k=3}$,  and $q_0=q_1=200$ } 
  % $(q_0=q_1=200,\delta,\alpha,k=3)$-1URA-decomposition}
\label{tabl:Hqq22dDaAAk3p}
%\rotatebox{90}{
\begin{tabular}{||c||c|c|c||c|c|c||c|c|c||c|c|c||}
\hline \hline
{\bf k=3} /$\delta$& $0.00$ & $0.00$ & $0.00$ &$10^{-6}$&$10^{-6}$&$10^{-6}$&$10^{-7}$&$10^{-7}$&$10^{-7}$&$10^{-8}$&$10^{-8}$&$10^{-8}$\\
\ j /$\alpha$ & 0.25   & 0.50   & 0.75   & 0.25   & 0.50   & 0.75   & 0.25   & 0.50   & 0.75   & 0.25   & 0.50   & 0.75   \\ \hline
  0& 2.45E-3& 2.48E-3& 2.49E-3& 2.49E-3& 2.49E-3& 2.49E-3& 2.49E-3& 2.49E-3& 2.49E-3& 2.49E-3& 2.48E-3& 2.49E-3\\
  1&-5.88E-16&-3.33E-10&-2.80E-8&-3.03E-10&-5.31E-9&-4.72E-8&-6.15E-11&-1.70E-9&-3.22E-8&-1.22E-11&-7.25E-10&-2.89E-8\\
  2&-3.94E-13&-1.92E-8&-7.86E-7&-6.25E-9&-5.92E-8&-8.46E-7&-1.63E-9&-3.54E-8&-8.02E-7&-4.09E-10&-2.54E-8&-7.90E-7\\
  3&-1.37E-10&-2.87E-7&-1.86E-6&-2.78E-7&-6.13E-7&-1.90E-6&-1.07E-7&-4.24E-7&-1.87E-6&-3.61E-8&-3.40E-7&-1.86E-6\\
\hline \hline \end{tabular}
%}
\end{table}
\end{tiny}

~ ~ ~

\begin{tiny}
\begin{table}
\setlength{\tabcolsep}{1.0pt}
\renewcommand{\arraystretch}{1.02}
\caption{The coefficients $\bar{\bar c}_j,  j=0, \dots, k$  of the partial fraction representation \eqref{frac-1-URA}
of $(q,\delta,\alpha,k)$-1-URA approximation 
\eqref{1-URA} for  $\delta=0, 10^{-6}, 10^{-7}, 10^{-8}$, $\alpha= 0.25, 0.50, 0.75$,
${\bf k=4}$,  and $q_0=q_1=200$ } 
  % $(q_0=q_1=200,\delta,\alpha,k=4)$-1URA-decomposition}
\label{tabl:Hqq22dDaAAk4p}
%\rotatebox{90}{
\begin{tabular}{||c||c|c|c||c|c|c||c|c|c||c|c|c||}
\hline \hline
{\bf k=4} /$\delta$& $0.00$ & $0.00$ & $0.00$ &$10^{-6}$&$10^{-6}$&$10^{-6}$&$10^{-7}$&$10^{-7}$&$10^{-7}$&$10^{-8}$&$10^{-8}$&$10^{-8}$\\
\ j /$\alpha$ & 0.25   & 0.50   & 0.75   & 0.25   & 0.50   & 0.75   & 0.25   & 0.50   & 0.75   & 0.25   & 0.50   & 0.75   \\ \hline
  0& 2.47E-3& 2.49E-3& 2.50E-3& 2.49E-3& 2.49E-3& 2.50E-3& 2.49E-3& 2.49E-3& 2.50E-3& 2.49E-3& 2.49E-3& 2.50E-3\\
  1&-4.37E-17&-3.58E-11&-2.27E-9&-1.82E-10&-2.74E-9&-7.28E-9&-3.67E-11&-6.13E-10&-3.29E-9&-7.15E-12&-1.70E-10&-2.49E-9\\
  2&-2.74E-14&-3.92E-9&-2.70E-7&-1.94E-9&-2.61E-8&-3.85E-7&-4.95E-10&-1.29E-8&-3.03E-7&-1.21E-10&-7.39E-9&-2.78E-7\\
  3&-3.60E-12&-5.10E-8&-1.01E-6&-2.87E-8&-1.49E-7&-1.03E-6&-9.62E-9&-9.51E-8&-1.02E-6&-2.96E-9&-6.97E-8&-1.02E-6\\
  4&-7.23E-10&-5.50E-7&-2.12E-6&-7.20E-7&-1.21E-6&-2.20E-6&-3.69E-7&-8.72E-7&-2.15E-6&-1.65E-7&-6.92E-7&-2.13E-6\\
\hline \hline \end{tabular}
%}
\end{table}
\end{tiny}

~ ~ ~

\begin{tiny}
\begin{table}
\setlength{\tabcolsep}{1.0pt}
\renewcommand{\arraystretch}{1.02}
\caption{The coefficients $\bar{\bar c}_j,  j=0, \dots, k$  of the partial fraction representation \eqref{frac-1-URA}
of $(q,\delta,\alpha,k)$-1-URA approximation 
\eqref{1-URA} for  $\delta=0, 10^{-6}, 10^{-7}, 10^{-8}$, $\alpha= 0.25, 0.50, 0.75$,
${\bf k=5}$,  and $q_0=q_1=200$ } 
  % $(q_0=q_1=200,\delta,\alpha,k=5)$-1URA-decomposition}
\label{tabl:Hqq22dDaAAk5p}
%\rotatebox{90}{
\begin{tabular}{||c||c|c|c||c|c|c||c|c|c||c|c|c||}
\hline \hline
{\bf k=5} /$\delta$& $0.00$ & $0.00$ & $0.00$ &$10^{-6}$&$10^{-6}$&$10^{-6}$&$10^{-7}$&$10^{-7}$&$10^{-7}$&$10^{-8}$&$10^{-8}$&$10^{-8}$\\
\ j /$\alpha$ & 0.25   & 0.50   & 0.75   & 0.25   & 0.50   & 0.75   & 0.25   & 0.50   & 0.75   & 0.25   & 0.50   & 0.75   \\ \hline
  0& 2.48E-3& 2.49E-3& 2.50E-3& 2.49E-3& 2.50E-3& 2.50E-3& 2.49E-3& 2.49E-3& 2.50E-3& 2.49E-3& 2.49E-3& 2.50E-3\\
  1&-4.03E-18&-4.04E-12&-2.15E-10&-1.24E-10&-1.58E-9&-1.59E-9&-2.49E-11&-2.67E-10&-4.65E-10&-4.77E-12&-5.19E-11&-2.67E-10\\
  2&-2.72E-15&-6.92E-10&-4.25E-8&-8.80E-10&-1.39E-8&-1.20E-7&-2.19E-10&-5.53E-9&-6.31E-8&-5.26E-11&-2.39E-9&-4.74E-8\\
  3&-2.87E-13&-1.36E-8&-5.98E-7&-7.48E-9&-6.33E-8&-7.16E-7&-2.35E-9&-3.70E-8&-6.46E-7&-6.91E-10&-2.42E-8&-6.12E-7\\
  4&-1.92E-11&-9.79E-8&-1.06E-6&-7.96E-8&-2.86E-7&-1.09E-6&-3.33E-8&-1.91E-7&-1.07E-6&-1.26E-8&-1.41E-7&-1.06E-6\\
  5&-2.84E-9&-8.91E-7&-2.34E-6&-1.27E-6&-1.88E-6&-2.45E-6&-7.85E-7&-1.44E-6&-2.38E-6&-4.40E-7&-1.17E-6&-2.35E-6\\
\hline \hline \end{tabular}
%}
\end{table}
\end{tiny}

~ ~ ~

\begin{tiny}
\begin{table}
\setlength{\tabcolsep}{1.0pt}
\renewcommand{\arraystretch}{1.02}
\caption{The coefficients $\bar{\bar c}_j,  j=0, \dots, k$  of the partial fraction representation \eqref{frac-1-URA}
of $(q,\delta,\alpha,k)$-1-URA approximation 
\eqref{1-URA} for  $\delta=0, 10^{-6}, 10^{-7}, 10^{-8}$, $\alpha= 0.25, 0.50, 0.75$,
${\bf k=6}$,  and $q_0=q_1=200$ } 
  % $(q_0=q_1=200,\delta,\alpha,k=6)$-1URA-decomposition}
\label{tabl:Hqq22dDaAAk6p}
%\rotatebox{90}{
\begin{tabular}{||c||c|c|c||c|c|c||c|c|c||c|c|c||}
\hline \hline
{\bf k=6} /$\delta$& $0.00$ & $0.00$ & $0.00$ &$10^{-6}$&$10^{-6}$&$10^{-6}$&$10^{-7}$&$10^{-7}$&$10^{-7}$&$10^{-8}$&$10^{-8}$&$10^{-8}$\\
\ j /$\alpha$ & 0.25   & 0.50   & 0.75   & 0.25   & 0.50   & 0.75   & 0.25   & 0.50   & 0.75   & 0.25   & 0.50   & 0.75   \\ \hline
  0& 2.48E-3& 2.49E-3& 2.50E-3& 2.50E-3& 2.50E-3& 2.50E-3& 2.49E-3& 2.50E-3& 2.50E-3& 2.49E-3& 2.50E-3& 2.50E-3\\
  1&-4.32E-19&-4.96E-13&-2.40E-11&-9.16E-11&-9.85E-10&-4.79E-10&-1.82E-11&-1.34E-10&-9.38E-11&-3.45E-12&-1.98E-11&-3.78E-11\\
  2&-3.19E-16&-1.07E-10&-4.83E-9&-4.91E-10&-8.19E-9&-3.07E-8&-1.20E-10&-2.58E-9&-1.08E-8&-2.80E-11&-8.37E-10&-6.26E-9\\
  3&-3.37E-14&-3.43E-9&-2.03E-7&-2.95E-9&-3.36E-8&-4.16E-7&-8.78E-10&-1.76E-8&-2.86E-7&-2.49E-10&-9.76E-9&-2.28E-7\\
  4&-1.62E-12&-2.89E-8&-7.60E-7&-2.00E-8&-1.17E-7&-7.85E-7&-7.50E-9&-7.37E-8&-7.77E-7&-2.62E-9&-5.10E-8&-7.67E-7\\
  5&-7.59E-11&-1.60E-7&-1.11E-6&-1.60E-7&-4.57E-7&-1.16E-6&-7.97E-8&-3.21E-7&-1.13E-6&-3.63E-8&-2.42E-7&-1.11E-6\\
  6&-8.98E-9&-1.27E-6&-2.51E-6&-1.87E-6&-2.51E-6&-2.62E-6&-1.28E-6&-2.02E-6&-2.56E-6&-8.24E-7&-1.69E-6&-2.52E-6\\
\hline \hline \end{tabular}
%}
\end{table}
\end{tiny}

~ ~ ~

\begin{tiny}
\begin{table}
\setlength{\tabcolsep}{1.0pt}
\renewcommand{\arraystretch}{0.95}
\caption{The coefficients $\bar{\bar c}_j,  j=0, \dots, k$  of the partial fraction representation \eqref{frac-1-URA}
of $(q,\delta,\alpha,k)$-1-URA approximation 
\eqref{1-URA} for  $\delta=0, 10^{-6}, 10^{-7}, 10^{-8}$, $\alpha= 0.25, 0.50, 0.75$,
${\bf k=7}$,  and $q_0=q_1=200$ } 
  % $(q_0=q_1=200,\delta,\alpha,k=7)$-1URA-decomposition}
\label{tabl:Hqq22dDaAAk7p}
%\rotatebox{90}{
\begin{tabular}{||c||c|c|c||c|c|c||c|c|c||c|c|c||}
\hline \hline
{\bf k=7} /$\delta$& $0.00$ & $0.00$ & $0.00$ &$10^{-6}$&$10^{-6}$&$10^{-6}$&$10^{-7}$&$10^{-7}$&$10^{-7}$&$10^{-8}$&$10^{-8}$&$10^{-8}$\\
\ j /$\alpha$ & 0.25   & 0.50   & 0.75   & 0.25   & 0.50   & 0.75   & 0.25   & 0.50   & 0.75   & 0.25   & 0.50   & 0.75   \\ \hline
  0& 2.48E-3& 2.49E-3& 2.50E-3& 2.50E-3& 2.50E-3& 2.50E-3& 2.50E-3& 2.50E-3& 2.50E-3& 2.49E-3& 2.50E-3& 2.50E-3\\
  1&-3.97E-18&-6.61E-14&-3.04E-12&-7.13E-11&-6.54E-10&-1.85E-10&-1.41E-11&-7.54E-11&-2.58E-11&-2.63E-12&-9.04E-12&-7.08E-12\\
  2&-2.69E-15&-1.58E-11&-5.92E-10&-3.12E-10&-5.15E-9&-8.76E-9&-7.44E-11&-1.30E-9&-2.20E-9&-1.70E-11&-3.21E-10&-9.69E-10\\
  3&-2.86E-13&-7.27E-10&-3.45E-8&-1.46E-9&-2.01E-8&-1.81E-7&-4.19E-10&-9.13E-9&-8.03E-8&-1.15E-10&-4.15E-9&-4.76E-8\\
  4& 9.97E-25&-9.08E-9&-4.30E-7&-7.44E-9&-6.11E-8&-6.14E-7&-2.59E-9&-3.60E-8&-5.36E-7&-8.52E-10&-2.27E-8&-4.72E-7\\
  5&-1.93E-11&-4.89E-8&-7.86E-7&-4.15E-8&-1.86E-7&-7.89E-7&-1.81E-8&-1.23E-7&-7.86E-7&-7.33E-9&-8.80E-8&-7.86E-7\\
  6&-2.85E-9&-2.35E-7&-1.16E-6&-2.63E-7&-6.44E-7&-1.24E-6&-1.49E-7&-4.76E-7&-1.20E-6&-7.84E-8&-3.68E-7&-1.18E-6\\
  7& 9.43E-25&-1.65E-6&-2.63E-6&-2.50E-6&-3.07E-6&-2.75E-6&-1.80E-6&-2.56E-6&-2.69E-6&-1.26E-6&-2.20E-6&-2.65E-6\\
\hline \hline \end{tabular}
%}
\end{table}
\end{tiny}

~ ~ ~

\begin{tiny}
\begin{table}
\setlength{\tabcolsep}{1.0pt}
\renewcommand{\arraystretch}{1.02}
\caption{The coefficients $\bar{\bar c}_j,  j=0, \dots, k$  of the partial fraction representation \eqref{frac-1-URA}
of $(q,\delta,\alpha,k)$-1-URA approximation 
\eqref{1-URA} for  $\delta=0, 10^{-6}, 10^{-7}, 10^{-8}$, $\alpha= 0.25, 0.50, 0.75$,
${\bf k=8}$,  and $q_0=q_1=200$ } 
  % $(q_0=q_1=200,\delta,\alpha,k=8)$-1URA-decomposition}
\label{tabl:Hqq22dDaAAk8p}
%\rotatebox{90}{
\begin{tabular}{||c||c|c|c||c|c|c||c|c|c||c|c|c||}
\hline \hline
{\bf k=8} /$\delta$& $0.00$ & $0.00$ & $0.00$ &$10^{-6}$&$10^{-6}$&$10^{-6}$&$10^{-7}$&$10^{-7}$&$10^{-7}$&$10^{-8}$&$10^{-8}$&$10^{-8}$\\
\ j /$\alpha$ & 0.25   & 0.50   & 0.75   & 0.25   & 0.50   & 0.75   & 0.25   & 0.50   & 0.75   & 0.25   & 0.50   & 0.75   \\ \hline
  0& 1.66E-3& 2.50E-3& 2.50E-3& 2.50E-3& 2.50E-3& 2.50E-3& 2.50E-3& 2.50E-3& 2.50E-3& 2.50E-3& 2.50E-3& 2.50E-3\\
  1& 1.52E-22&-9.55E-15&-4.29E-13&-5.76E-11&-4.55E-10&-8.51E-11&-1.13E-11&-4.61E-11&-9.06E-12&-2.09E-12&-4.72E-12&-1.74E-12\\
  2&-1.70E-22&-2.40E-12&-8.16E-11&-2.17E-10&-3.40E-9&-3.00E-9&-5.05E-11&-7.02E-10&-5.68E-10&-1.13E-11&-1.37E-10&-1.89E-10\\
  3& 8.43E-27&-1.34E-10&-4.73E-9&-8.44E-10&-1.30E-8&-6.40E-8&-2.33E-10&-5.01E-9&-1.90E-8&-6.20E-11&-1.82E-9&-8.55E-9\\
  4&-1.24E-27&-2.52E-9&-1.28E-7&-3.49E-9&-3.65E-8&-4.15E-7&-1.15E-9&-1.98E-8&-2.64E-7&-3.60E-10&-1.10E-8&-1.78E-7\\
  5& 2.05E-27&-1.73E-8&-5.85E-7&-1.53E-8&-9.62E-8&-6.38E-7&-6.08E-9&-6.04E-8&-6.34E-7&-2.28E-9&-4.05E-8&-6.13E-7\\
  6& 7.47E-27&-7.30E-8&-7.87E-7&-7.20E-8&-2.67E-7&-8.09E-7&-3.55E-8&-1.85E-7&-7.94E-7&-1.63E-8&-1.35E-7&-7.89E-7\\
  7& 3.06E-26&-3.19E-7&-1.22E-6&-3.81E-7&-8.31E-7&-1.30E-6&-2.37E-7&-6.40E-7&-1.26E-6&-1.39E-7&-5.09E-7&-1.24E-6\\
  8& 6.51E-24&-2.01E-6&-2.72E-6&-3.15E-6&-3.59E-6&-2.86E-6&-2.35E-6&-3.06E-6&-2.79E-6&-1.73E-6&-2.67E-6&-2.75E-6\\
\hline \hline \end{tabular}
%}
\end{table}
\end{tiny}

%~ ~ ~

%\endinput

\section{How to access the BURA data from the repository}

\subsection{Data generated by Remez algorithm}\label{data-soft}

Remez algorithm is implemented according to the results reported in \cite{PGMASA1987}.
It is part of a software developed in IICT, Bulgarian Academy of Sciences, concerting the 
more general problem of the best uniform rational approximation of a given function of $ x \in [0,1]$
by $P_m(x)/Q_k(x)$, where $P_m$ and $Q_k$ are polynomials of degree  $m \ge k$, correspondingly.
In this work we use only the particular case when $m=k$.
To make the writing more concise we use the following notations
$$ 
P_n(x)  \quad \mbox{denoted by} \quad p(n,F;x) = \sum_{j=0}^{n} F(j)*x^j,
$$
where $n$ is the degree of the polynomial and $F(j)$, $j=0,\dots, n$, are the  coefficients. 
Further, we denote  the zeros of  $ Q_k(x)$ by  $  U0(j))$, $j=0, \dots, k$, which in general are 
complex numbers 
$$
U0(j)) = Re\{U0(j))\} + i Im\{ U0(j))\}, \ \ j=1, \dots, k, \ \ i^2=-1.
$$ 
Accordingly, the best uniform rational approximation $P_m(x)/Q_k(x)=  p(m,A;x) / p(k,B;x)  $  in the data files 
will be denoted by and represented as a sum of partial fractions as
\begin{equation}\label{BURA-gen-data}
\begin{aligned}
%r(m,k,A,B, x) & = %P_m(A;x) / Q_k(B;x) = 
p(m,A;x) / p(k,B;x) 
& = p(m-k,C;x) + p(k-1,D;x)/p(k,B;x) \\
& =  \sum_{j=0}^{m-k} C(j) * x^j+ \sum_{j=1}^{k} E(j) / (x-U0(j)).
 \end{aligned}
\end{equation} 
where the coefficients $C(j)$ and $E(j)$ are in general complex numbers. Of course, when the 
poles are real, then the complex parts, since obtained by computational procedure for finding 
roots of a polynomial, will have very small imaginary parts. 

In our case $k=m$ and the above representation becomes
$$
r(m,k,A,B, x) =  C(0) + \sum_{j=1}^{k} E(j) / (x-U0(j)).
$$

\subsection{Examples how to use the  data in the repository}\label{example4use}

Now we shall illustrate on a couple of examples how to solve the system $\wcalAt^\alpha \tiluh = \tilfh$
arising in the approximation of the Example 1 that produces the matrix $ \wcalAt$ defined by \eqref{FD-matrix-1D-k}
by using BURA method \eqref{wh}.  For definiteness and convenience, we assume that $h=10^{-3}$, 
$\min a(x)= 0.25/\pi^2$ and 
$\max a(x)= 1$, so that $ \lambda_1 \ge 1$ and  $ \lambda_1/\lambda_N = h^2/4= 0.25*10^{-6}$. 
We assume also that the data $f$ is such that the step-size $h=10^{-3}$ allows the following estimate
for the semi-discrete error $ \|u - u_h \| \le 10^{-3} \|f\|$.

The task we have is to find an approximation $\tilwh$ of the solution $\tiluh$ with 
relative error $10^{-3}$ by using \eqref{wh} so that (simplified due to $\lambda_1 =1$)
%$$
%\tilwh = \lambda_1^{-\alpha} r_{\alpha,k}(\lambda_1^{-\alpha} \wcalAt) \tilfh =
%\left ( c_0 + \sum_{i=1}^k c_i((\lambda_1^{-\alpha} \wcalAt)^{-1}  -d_i)^{-1} \right ) \lambda_1^{-\alpha} \tilfh.
%$$
%Since $\lambda_1 =1$ we get the following simple representation of the solution:
$$
\tilwh = r_{\alpha,k}(\wcalAt^{-1}) \tilfh = \left ( c_0 + \sum_{i=1}^k c_i( \wcalAt^{-1} -d_i)^{-1} \right ) \tilfh
=  \left (\widetilde c_0+\sum_{i=1}^k {\widetilde c_i}{(\wcalAt -\widetilde d_i)^{-1}}  \right ) \tilfh.
$$
%Note that we got some simplification due to the fact that  $\lambda_1 =1$.
%

The coefficients $c_j$, $j=0, \dots, k$, and the poles $d_j$, $k=1, \dots, k$ will be
obtained from the repository \verb! http://parallel.bas.bg/~pencho/BURA/!.  Moreover, in the repository one can find 
more information about $r_{\alpha,k}(t)$, namely, the extremal points of the error  $t^\alpha -  r_{\alpha,k}(t)$, 
the zeros of the numerator, etc. 
%data needed in order to use the BURA $r_{\alpha,k}(t)$ of $t^\alpha$ and most important the coefficients 
%For this one should have 
%available a solution procedure for solving the systems of the type $(\wcalAt - c \wcalIt)\tiluh = \tilfh$, $c \le 0$. 

\paragraph{\underline{Example 1.} We need to solve $\wcalAt^\alpha \tiluh = \tilfh$ with relative accuracy $10^{-3}$ 
for $\delta=0$ and $\alpha=0.25$.}

\noindent 
We visit  Table \ref{tabl:A0BURAp} and look for the corresponding error for $q=0$.
It appears that error below $10^{-3}$  is achieved for $k=7$. In order to implement our algorithm we need to 
get the poles $d_j$ from the second row of Table \ref{tabl:Cq000dDaAAk7p} (under $\delta=0$ and $\alpha=0.25$), 
and the coefficients $c_j$ from the second row of Table \ref{tabl:Dq000dDaAAk7p}  (under  $\delta=0$ and $\alpha=0.25$)
of  BURA. According to Table \ref{p1ltabl},
the corresponding tables are encoded as $Cq000d0a25k7p$ and $ Dq000d0a25k7p$.
To recover the data, in a browser open  the site of the repository %\begin{verbatim} 
\verb! http://parallel.bas.bg/~pencho/BURA/ !
%\end{verbatim} 
which will get you to the html-repository. 
\begin{footnotesize}

\begin{figure}[h!t]
\begin{verbatim}
Index
mode   size  last-canged   name
dr-x         Jul 18 09:26  BURA-tabl/ - with extremal points
dr-x         Jul 18 10:50  BURA-dcmp/ - with poles an decomposition coefficients for BURA
dr-x         Jul 18 11:22  0URA-dcmp/ - with poles an decomposition coefficients for 0URA
dr-x         Jul 18 11:35  1URA-dcmp/ - with poles an decomposition coefficients for 1URA
-r--   672k  Sep 21 11:51  BURA-data-report.pdf - 60 pages document
-r--   3.2M  Sep 21 11:51  all-in-one.zip  - archive with all files (approx 2000) 
-r--   165k  Sep 21 11:51  all-in-one_long.lst  -  long  list of all folders and sub-folders
-r--    720  Sep 21 11:51  all-in-one_short.lst -  short list of the main folder 
\end{verbatim}
\caption{Index of the file repository}\label{html--data}
\end{figure}
\end{footnotesize}

The coefficients in the partial fraction representation of BURA are in the folder 
\texttt{BURA-dcmp/}. % folder where are the files with the data. 
%Consequently, from \texttt{http://parallel.bas.bg/~pencho/BURA/BURA-dcmp/} you get 
Consequently, from \verb! http://parallel.bas.bg/~pencho/BURA/BURA-dcmp/! you get 
a long list of data files where you can extract the data from.

\begin{footnotesize}
\begin{figure}[h!t]
\begin{verbatim}
Index
mode   size  last-canged   name
dr-x         Aug  7 12:19  ../
dr-x         Aug  7 13:16  add/
-r--    229  Jul 18 07:26  0-data-info.txt
-r--    20k  Aug 18 07:26  0-data-info.pdf
-r--    901  Jul 18 09:28  q000d0a25k3.tab
-r--   1.1k  Jul 18 09:28  q000d0a25k4.tab
-r--   1.3k  Jul 18 09:28  q000d0a25k5.tab
-r--   1.5k  Jul 18 09:28  q000d0a25k6.tab
-r--   1.7k  Jul 18 09:28  q000d0a25k7.tab
-r--   1.9k  Jul 18 09:28  q000d0a25k8.tab
. . . 
\end{verbatim}
\caption{The list of file containing the BURA data}\label{fig:data-files}
\end{figure}
\end{footnotesize}
%The file \verb! pdf!  contains a pdf file with short explanation hoe to interpret the data in the files listed below.

Note that the name of the file corresponds for the data you need. 
In order to understand better 
how to recover the needed coefficient we recommend you read to file \texttt{0-data-info.pdf}. 
Now, if you need BURA for $q=0$, 
$\delta=0$, $\alpha=0.25$, and $k=7$ you open the file \texttt{q000d0a25k7.tab} to get 
the data shown on Figure \ref{data-example-1}. 
\begin{footnotesize}
\begin{figure}[h!t]
%\label{data-exp1}
\begin{verbatim}
........................
    0 C(j), j=0,M-K 
    0,  1.5251198659461835471669134E+0000,
.......................
    7 Re{U0(j)}, j=1,K 
    1, -2.2376872996078341567977828E-0010,
    2, -7.8486161880478847863452285E-0008,
    3, -5.5731465614475310569044313E-0006,
    4, -1.8879914699444446592534226E-0004,
    5, -4.0807806486154275568196059E-0003,
    6, -6.6631003383177437492683783E-0002,
    7, -1.3009123561707709448297558E+0000,
    7 Re{E(j)}, j=1,K 
    1, -1.4685691492826336539266828E-0012,
    2, -1.4250266060744046359349294E-0009,
    3, -2.3302001273674185907801704E-0007,
    4, -1.6270756180817065256760833E-0005,
    5, -6.7439995512918544429211685E-0004,
    6, -2.0780142734679021643365256E-0002,
    7, -1.1636545994360130292202863E+0000.
    \end{verbatim}
    \vspace{-5mm}
    \caption{BURA Data for Example 1}\label{data-example-1}
        \end{figure}
 \end{footnotesize}       
     %
%where  $p(n,F;x) = \sum_{j=0}^{n} F(j)*x^j $ is polynomial of degree $n$.    
 Comparing this explained by formula \eqref{BURA-gen-data}
 with the representation of $r_{\alpha, k}( x)$ by \eqref{trg}  we realize that 
the poles $d_j$ are in  $ Re\{U0(j)\}$, while the 
coefficients $c_j$ are in $ Re\{E(j)\}$, $j=1, \dots,7$ except $c_0= C(0)$. 
% that since $c_0= C(0)$,
% $j=0,\dots, k-m$, 

{\small
%Thus, we have  \\
\begin{figure}\label{coeffs1}
\centering
\begin{subfigure}{0.47 \textwidth}
\centering
\begin{tabular}{ll}
 $c_0$=&\ $C(0)$\ \ \ \ =$~\ ~ 1.5251198659461835$  \\ %3.2524043182367136372344287E+0000$. \\
 $c_1$=&$Re\{E(1)\}$=$-1.468569149282634*10^{-12} $,   \\ %-4.3583180671197569815569583E-0009$, \\
 $c_2$=&$Re\{E(2)\}$=$-1.425026606074404*10^{-9}$,  \\ %-1.1477817769645348761874750E-0006$, \\
 $c_3$=&$Re\{E(3)\}$=$-2.330200127367418*10^{-7}$, \\ %-7.7384665890689124783540290E-0005,$ \\
 $c_4$=&$Re\{E(4)\}$=$-1.627075618081706*10^{-5} $, \\ % -2.7381282249725339508289331E-0003,$ \\
 $c_5$=&$Re\{E(5)\}$=$-6.743999551291854*10^{-4}$, \\ % -8.0144182846238070945938169E-0002,$ \\
 $c_6$=&$Re\{E(6)\}$=$-2.078014273467902*10^{-2}$, \\ %-1.0807726759451478394181312E+0001$.\\
 $c_7$=&$Re\{E(7)\}$=$-1.1636545994360130$,
\end{tabular}
   \caption{Coefficients $c_j$  with 15 significant digits}
 \end{subfigure}
  \hfill
 \begin{subfigure}{0.47 \textwidth} %\label{poles1}
 \centering
 (recall $d_j= Re\{U0(j)\}$, $j=1, \dots,7$) \\
\begin{tabular}{ll}
 $d_1$=&$Re\{U0(1)\}$=$-2.237687299607834*10^{-10}$,  \\ %-1.9530271231113077713630524E-0006$, \\
 $d_2$=&$Re\{U0(2)\}$=$-7.848616188047885*10^{-8}$, \\ %  -1.0566469257506061685150566E-0004,$ \\
 $d_3$=&$Re\{U0(3)\}$=$-5.573146561447531*10^{-6}$, \\ % -2.0274404802332074762073092E-0003,$ \\
 $d_4$=&$Re\{U0(4)\}$=$-1.887991469944444*10^{-4}$, \\ % -2.7381282249725339508289331E-0003,$ \\
 $d_5$=&$Re\{U0(5)\}$=$-4.080780648615427*10^{-3}$, \\ %-8.0144182846238070945938169E-0002,$ \\ 
 $d_6$=&$Re\{U0(6)\}$=$-6.663100338317743*10^{-2}$, \\ %-1.0807726759451478394181312E+0001$.\\
 $d_7$=&$Re\{U0(7)\}$=$-1.300912356170771$. \\ %-5.3692125077979476196430882E-0002$ 
\end{tabular}
  \caption{The poles $d_j$ with 15 significant digits}
  \end{subfigure}
  \caption{The coefficients and poles of BURA for Example~1}
 \end{figure}
  }

Note that the poles of the rational function $r_{\alpha,k}(t)$ are computed by finding the roots of its
denominator. In the general case these are complex numbers with real and complex parts. However, in our case 
the roots are real, so the fact that the complex parts are near zero is another indication that our computations are correct.
The same remark is valid also for the coefficients in the partial fraction representation. 
One can recover the data including the complex parts of the roots and coefficients from the file with the same name abut extension
\verb!.txt!. Thus, instead of file \verb! q000d0a25k7.tab!  we need to access file \verb!q000d0a25k7.txt!  
The full path is 
\verb! http://parallel.bas.bg/~pencho/BURA/BURA-dcmp/add/!.
According Remark~{\ref{lemma:BURA}} and \eqref{r-compute} we obtain for the coefficients
 $\{\widetilde c_j\}_{j=0}^{k}$
% (see  Fig.~\ref{coef-example-1-infty})
 and for the poles 
 $\{\widetilde d_j\}_{j=1}^{k}$, shown on Figure ~\ref{coef-example-1-infty}.%  and ~\ref{pole-example-1-infty}, respectively.
% (see  Fig.~\ref{pole-example-1-infty})

\begin{footnotesize} %tiny}
\begin{figure}[h!t]
\begin{verbatim}
   j,  coeffs  \tilde c_j, j=0,...,k                       j,    poles \tilde d_j, j=1,...,k          
  ---------------------------------------               --------------------------------------
    0,  7.8649908986141400766246111E-0004,
    1,  6.8758756374255199468223446E-0001,               1, -7.6869129212016636558276595E-0001,
    2,  4.6805384769874068750186188E+0000,	               2, -1.5008028533643146339273005E+0001,
    3,  4.0497762636194054740290213E+0001,	               3, -2.4505115224443418015306629E+0002,
    4,  4.5646520861012335445240629E+0002,	               4, -5.2966340998851315319920697E+0003,
    5,  7.5022631473305133482030654E+0003,	               5, -1.7943185038727329359280164E+0005,
    6,  2.3133257355683798168384452E+0005,	               6, -1.2741099526854566986420878E+0007,
    7,  2.9328888647334126024639520E+0007,	               7, -4.4688996544568804592027580E+0009.
\end{verbatim}
   \vspace{-3mm}
     \caption{The coefficients $\widetilde c_j$ (left) and the poles $\widetilde d_j$ (right) of BURA in $[1,\infty )$  for Example 1; recall $k=7$}
     \label{coef-example-1-infty}
        \end{figure}
\end{footnotesize} %tiny}        

%\begin{footnotesize} %tiny}
%\begin{figure}[h!]
%\begin{verbatim}
%    j,    poles \tilde d_j, j=1,...,k  
%  ---------------------------------------
%    1, -7.6869129212016636558276595E-0001,
%    2, -1.5008028533643146339273005E+0001,
%    3, -2.4505115224443418015306629E+0002,
%    4, -5.2966340998851315319920697E+0003,
%    5, -1.7943185038727329359280164E+0005,
%    6, -1.2741099526854566986420878E+0007,
%    7, -4.4688996544568804592027580E+0009,
%\end{verbatim}
%   \vspace{-3mm}
%     \caption{The poles $\widetilde d_j$ of BURA  in $[1,\infty )$  for Example 1; recall $k=7$}
%     \label{pole-example-1-infty}
%        \end{figure}
%\end{footnotesize} %tiny}        
%
        
%\pagebreak     

\paragraph{\underline{Example 2}} 
We want to solve $\wcalAt^\alpha \tiluh = \tilfh$ with relative accuracy $10^{-4}$ for $\alpha=0.50$.
Since $\lambda_N/\lambda_1 =2.5 * 10^{-7} > 10^{-7}$. Thus, we can choose 
any $\delta \le 10^{-7}$. We shall choose $\delta =10^{-8} $, which indeed is smaller than $\lambda_N/\lambda_1$.
%\red{We can use also $\delta = 10^{-7}$ ?????}
%$\delta=3.4*10^{-7} > 10^{-8}$. 
% consider the case (2) assume that after getting values for $ \underline{\lambda}_1$ and $\bar \lambda_N$ we 
%produce $\delta=3.4*10^{-7} > 10^{-6}$. 
Visiting Table \ref{tabl:A0BURAp} and looking for the corresponding error for $q=0$
and $\delta=10^{-8}$
%Then we can go to  Table \ref{tabl:Cq000dDaAAk5p}. 
%for $\delta=10^{-6}$  and 
we see that $k=6$ gives the desired accuracy.
Then again, to implement the algorithm through the form 
 \eqref{trg} we need to get the poles $d_j, j=1,\dots,6$
from Table \ref{tabl:Cq000dDaAAk6p}, 
row 6 (under $\delta=10^{-8}$ and $\alpha=0.50$), 
 and the coefficients $c_j, j=0,1, \dots, 6$, from Table \ref{tabl:Dq000dDaAAk6p},   
 row 6 (under $\delta=10^{-8}$ and $\alpha=0.50$). In order to recover them 
 from the html-repository with 25 significant digits we read  the file \texttt{q000d8a50k6.tab}:
% thus, for implementing \eqref{trg} we need the coefficients $c_j$  and the poles $d_j$.%
\begin{footnotesize} %tiny}
\begin{figure}[h!t]
 \begin{verbatim}
............................................
 0 C(j), j=0,M-K 
    0,  3.2884110965902564982570418E+0000,
.............................................
    6 Re{U0(j)}, j=1,K 
    1, -2.9026635633198611614392944E-0006,
    2, -1.2844649138280910668994002E-0004,
    3, -2.2714606179594160697493911E-0003,
    4, -2.5644587248425397513332244E-0002,
    5, -2.3875648347625806119814403E-0001,
    6, -4.0426255143071180974523916E+0000,
    6 Re{E(j)}, j=1,K 
    1, -7.3047263730915787436076796E-0009,
    2, -1.4856610579547499079286817E-0006,
    3, -8.9729276335619099229017533E-0005,
    4, -2.9805865012026772176372313E-0003,
    5, -8.4097511581798008760373107E-0002,
    6, -1.1182633938371627571773423E+0001.
 \end{verbatim}
 \vspace{-3mm}
     \caption{The BURA Data for Example 2: recall $m=k=6$}
     \label{data-example-2}
        \end{figure}
\end{footnotesize} %tiny}
 
 One can recover the data including the complex parts of the roots and coefficients from the file with the same name abut extension
\verb!.txt!. Thus, instead of file \verb!q000d8a50k6.tab!  we need to access file \verb!q000d8a50k6.txt!  to get (recall $m=k=6$).
As we see, the complex parts of the poles $U0(j)$ and of the coefficients $E(j)$ are almost zero. 
Therefore, we ignore them in table  \verb!q000d8a50k6.tab!  as seen on Figure \ref{data-example-2}.
{\small
%Thus, we have  \\
\begin{figure}[h!t]\label{coeffs2}
\centering
\begin{subfigure}{0.47 \textwidth}
\centering
\begin{tabular}{ll}
 $c_0$=&\ $C(0)$\ \ \ \ =$~\ ~ 3.288411096590256 $ \\ 
 %1.5251198659461835$  \\ %3.2524043182367136372344287E+0000$. \\
 $c_1$=&$Re\{E(1)\}$=$-7.304726373091578*10^{-9}$\\ 
 %1.46856914928263*10^{-12} $,   \\ %-4.3583180671197569815569583E-0009$, \\
 $c_2$=&$Re\{E(2)\}$=$-1.485661057954749*10^{-6} $ \\
 %1.42502660607440*10^{-9}$,  \\ %-1.1477817769645348761874750E-0006$, \\
 $c_3$=&$Re\{E(3)\}$=$-8.972927633561909*10^{-5} $ \\
 %2.330200127367418*10^{-7}$, \\ %-7.7384665890689124783540290E-0005,$ \\
 $c_4$=&$Re\{E(4)\}$=$-2.980586501202677*10^{-3}$ \\ %1.627075618081706*10^{-5} $, \\ % -2.7381282249725339508289331E-0003,$ \\
 $c_5$=&$Re\{E(5)\}$=$-8.409751158179800*10^{-2} $ \\
 %6.743999551291854*10^{-4}$, \\ % -8.0144182846238070945938169E-0002,$ \\
 $c_6$=&$Re\{E(6)\}$=$-1.118263393837162*10 $ \\
 %2.078014273467902*10^{-2}$, \\ %-1.0807726759451478394181312E+0001$.\\
% $c_7$=&$Re\{E(7)\}$=&$-1.1636545994360130$,
\end{tabular}
   \caption{Coefficients $c_j$  with 15 significant digits}
 \end{subfigure}
 ~ %\hfill
 \begin{subfigure}{0.47 \textwidth} %\label{poles1}
 \centering
 (recall $d_j= Re\{U0(j)\}$, $j=1, \dots,7$) \\
\begin{tabular}{ll}
 $d_1$=&$Re\{U0(1)\}$=$-2.902663563319861*10^{-6} $ \\
 %2.23768729960783*10^{-10}$,  \\ %-1.9530271231113077713630524E-0006$, \\
 $d_2$=&$Re\{U0(2)\}$=$-1.284464913828091*10^{-4}$ \\
 %7.84861618804788*10^{-8}$, \\ %  -1.0566469257506061685150566E-0004,$ \\
 $d_3$=&$Re\{U0(3)\}$=$-2.271460617959416*10^{-3} $ \\
 %5.573146561447531*10^{-6}$, \\ % -2.0274404802332074762073092E-0003,$ \\
 $d_4$=&$Re\{U0(4)\}$=$-2.564458724842539*10^{-2} $ \\
 %1.887991469944444*10^{-4}$, \\ % -2.7381282249725339508289331E-0003,$ \\
 $d_5$=&$Re\{U0(5)\}$=$-2.387564834762580*10^{-1} $ \\ 
 %4.080780648615427*10^{-3}$, \\ %-8.0144182846238070945938169E-0002,$ \\ 
 $d_6$=&$Re\{U0(6)\}$=$-4.042625514307118*10 $\\
 %6.663100338317743*10^{-2}$, \\ %-1.0807726759451478394181312E+0001$.\\
% $d_7$=&$Re\{U0(7)\}$=&$-1.300912356170771$. \\ %-5.3692125077979476196430882E-0002$ 
\end{tabular}
  \caption{The poles $d_j$ with 15 significant digits}
  \end{subfigure}
  \caption{The coefficients and the poles for Example~2}
 \end{figure}
  }
\begin{footnotesize} %tiny}
\begin{figure}
\begin{verbatim} 
........................
0 C(j), j=0,M-K 
    0,  3.2884110965902564982570418E+0000,
................................
    6 U0(j), j=1,K 
    1, -2.9026635633198611614392944E-0006, -1.7520806294041483876573147E-0037
    2, -1.2844649138280910668994002E-0004,  1.8242239649606423471743024E-0037
    3, -2.2714606179594160697493911E-0003, -7.2143334828724323886476075E-0039
    4, -2.5644587248425397513332244E-0002, -7.2776229715659618428404602E-0047
    5, -2.3875648347625806119814403E-0001, -7.3384739154186567125843208E-0052
    6, -4.0426255143071180974523916E+0000, -8.8654988085625587825735093E-0067
    6 E(j), j=1,K 
    1, -7.3047263730915787436076796E-0009,  5.7573111695991243362896879E-0037
    2, -1.4856610579547499079286817E-0006, -5.9992385756159082541115791E-0037
    3, -8.9729276335619099229017533E-0005,  2.4155906802799248832468797E-0038
    4, -2.9805865012026772176372313E-0003,  2.2938808341026098990357342E-0041
    5, -8.4097511581798008760373107E-0002,  9.4285923737400470688723456E-0042
    6, -1.1182633938371627571773423E+0001,  4.4663981643768036611005825E-0042 
 \end{verbatim}
   \vspace{-3mm}
     \caption{The BURA Data for Example 2 in complex numbers; recall $m=k=6$}
     \label{data-example-2-full}
        \end{figure}
\end{footnotesize} %tiny}    
%    
%\medskip%\noindent
%As we see, the complex parts of the poles $U0(j)$ and of the coefficients $E(j)$ are almost zero. 
%Therefore, we ignore them in table  \verb!q000d8a50k6.tab!  as seen  on Figure \ref{data-example-2}.
 
According Remark~{\ref{lemma:BURA}} and \eqref{r-compute} we obtain  the coefficients
 $\{\widetilde c_j\}_{j=0}^{k}$
% (see  Fig.~\ref{coef-example-2-infty})   
 and the poles 
 $\{\widetilde d_j\}_{j=1}^{k}$, shown on 
 Figure ~\ref{example-2-infty}. % and ~\ref{pole-example-2-infty}, respectively.  
% (see  Fig.~\ref{pole-example-2-infty})
 %
% % corrected by Raytcho %%%%%%%%%%%%%%%
%\begin{footnotesize} %tiny}
%\begin{figure}[h!]
%\begin{verbatim}
%    j,  coeffs  \tilde c_j, j=0,...,k    j,                       poles \tilde d_j, j=1,...,k  
%   ---------------------------------------                     ---------------------------------------
%    0,  1.8619525184162842197896828E-0004,   
%    1,  6.8425358785323701384064700E-0001,                  1, -2.4736399561644631729004006E-0001,
%    2,  1.4752743712419678037242149E+0000,                  2, -4.1883679364017772326828007E+0000,
%    3,  4.5322129077686442621147161E+0000,                  3, -3.8994583547504784590380748E+0001,
%    4,  1.7390967446300787239267527E+0001,                  4, -4.4024536110969760472793286E+0002,
%    5,  9.0048244054276093870632027E+0001,                  5, -7.7853430579096145080903339E+0003,
%    6,  8.6698293594856505069077514E+0002,                  6, -3.4451116300101647683572511E+0005.
%\end{verbatim}
%   \vspace{-3mm}
%     \caption{The coefficients $\widetilde c_j$ (left) and the poles $\widetilde d_j$ (right) of BURA  
%     in $[1,\infty )$ for Example 2; recall $k=6$}
%     \label{coef-example-2-infty}
%        \end{figure}
%\end{footnotesize} %tiny}        

\begin{footnotesize} %tiny}
\begin{figure}
\centering
\begin{subfigure}{0.45 \textwidth}
\centering
\begin{verbatim}
    j,  coeffs  \tilde c_j, j=0,...,k    j,                     
   ---------------------------------------                   
    0,  1.8619525184162842197896828E-0004,   
    1,  6.8425358785323701384064700E-0001,       
    2,  1.4752743712419678037242149E+0000,     
    3,  4.5322129077686442621147161E+0000,      
    4,  1.7390967446300787239267527E+0001,     
    5,  9.0048244054276093870632027E+0001,   
    6,  8.6698293594856505069077514E+0002.
\end{verbatim}
    \caption{The coefficients $\widetilde c_j$} % of BURA  in $[1,\infty )$ for Example 2,  recall $k=6$}
\end{subfigure}%
 \hfill%
\begin{subfigure}{0.45 \textwidth}
\centering
\begin{verbatim}
    j,    poles \tilde d_j, j=1,...,k  
  ---------------------------------------
    1, -2.4736399561644631729004006E-0001,
    2, -4.1883679364017772326828007E+0000,
    3, -3.8994583547504784590380748E+0001,
    4, -4.4024536110969760472793286E+0002,
    5, -7.7853430579096145080903339E+0003,
    6, -3.4451116300101647683572511E+0005,
\end{verbatim}
\caption{The poles $\widetilde d_j$} % of BURA  in $[1,\infty )$ for Example 2}
     \end{subfigure}
\caption{Coefficients $\widetilde c_j$ and poles $\widetilde d_j$ of BURA  in $[1,\infty )$ for Example 2, $k=6$.}     
     \label{example-2-infty} 
\end{figure}
\end{footnotesize} %tiny}    
%   

%\pagebreak

\section* {Acknowledgments}

\par\noindent
This work has been partly support by the Grant No BG05M2OP001-1.001-0003, financed by the 
Science and Education for Smart Growth Operational Program (2014-2020) and co-financed by 
the EU through the European structural and Investment funds.

%\tableofcontents

\bibliographystyle{abbrv} % {plain}
\bibliography{BURA_refer_exp}  %,references-R} %,NSF-references}

%\tableofcontents

\end{document}